\documentclass{article}
\usepackage{amsfonts, color}
\usepackage{textgreek}
\usepackage[dvipsnames]{xcolor}
\usepackage{geometry}
\usepackage{enumitem}
\usepackage[normalem]{ulem}
\usepackage{relsize}
\usepackage{comment}
\usepackage{lineno}
\usepackage{xstring}
\usepackage{braket}
\usepackage{amstext, amsthm, amssymb, amsmath, mathtools}
\numberwithin{equation}{section}
\usepackage{mathrsfs}
\usepackage{bbm}
\usepackage{graphicx}
\usepackage{float}
\usepackage[hypertexnames=false]{hyperref}
\usepackage[capitalize,nameinlink,noabbrev]{cleveref}
\usepackage[style=alphabetic,isbn=false]{biblatex}
\usepackage[colorinlistoftodos,textsize=small]{todonotes}
\makeatletter\def\@todonotes@font{\rmfamily}\makeatother
\presetkeys{todonotes}{color=yellow!50}{}

\usepackage{caption}
\usepackage{subcaption} % in your preamble
\usepackage{adjustbox}
\usepackage{extpfeil}

\usepackage{longtable}
\usepackage{booktabs}
\usepackage{tikz}
\usepackage{tikz-cd}
\usepackage{bm}
\LoadFontDefinitionFile{U}{stmry}
\DeclareFontShape{U}{stmry}{b}{n}{<->ssub * stmry/m/n}{}
\usepackage{stackengine}
\usepackage{thmtools}
\usepackage{bookmark}
\setstackEOL{\\}
\usetikzlibrary{calc,decorations.pathreplacing}
\newcommand{\FDiff}{%
  \begin{tikzpicture}[baseline=-0.1ex,scale=0.05]
    \draw[line width=0.7pt] (0,0) -- (2,3.464) -- (4,0) -- cycle;
  \end{tikzpicture}%
}
\newcommand{\FDiffDiff}{%
  \begin{tikzpicture}[baseline=-0.1ex,scale=0.05]
    \draw[line width=0.7pt] (0,0) -- (2,3.464) -- (4,0) -- cycle;
    \draw[line width=0.7pt] (1,0) -- (2.5, 2.598);
  \end{tikzpicture}%
}
\usepackage{crossreftools}

\newcommand{\supp}{\operatorname{supp}}
\newcommand{\dist}{\operatorname{dist}}
\newcommand{\diff}{\mathrm{d}}

\DeclareMathOperator{\sn}{sn}

\newcommand{\Pp}{\mathcal{P}}
\newcommand{\R}{\mathbb{R}}

\newcommand{\NN}{{\mathbb N}}
\newcommand{\Q}{\mathbb{Q}}

\newcommand{\ACt}{\ensuremath{AC_{\mathrm{II} } } }
\newcommand{\EpsE}{\ensuremath{\epsilon\text{-}\mathfrak{E}}}
\NewDocumentCommand{\Tang}{ O{\varepsilon} }{
  \mathcal{T}_{#1} X
}

\DeclareMathOperator*{\argmin}{arg\,min}
\newcommand{\limsupTauTS}{\limsup_{\substack{T \searrow s,\, t \searrow s, \\ T \geq t}}}
\newcommand{\AC}{\ensuremath{AC}}

\usepackage{xcolor}

\newif\ifshowcomments
\showcommentstrue

\newcommand{\FundingLogos}{%
  \raisebox{0pt}{\includegraphics[height=1.5cm]{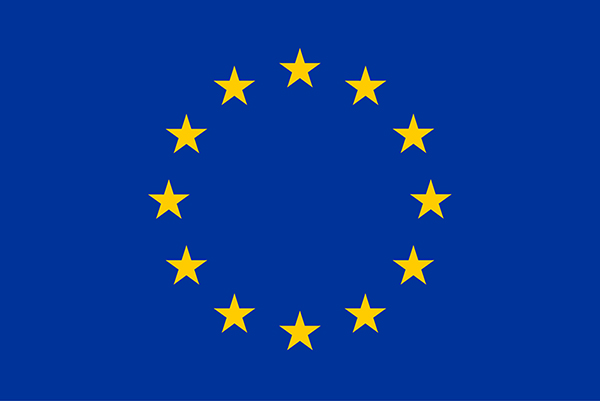}}%
  \hspace{1em}%
  \raisebox{0pt}{\includegraphics[height=1.5cm]{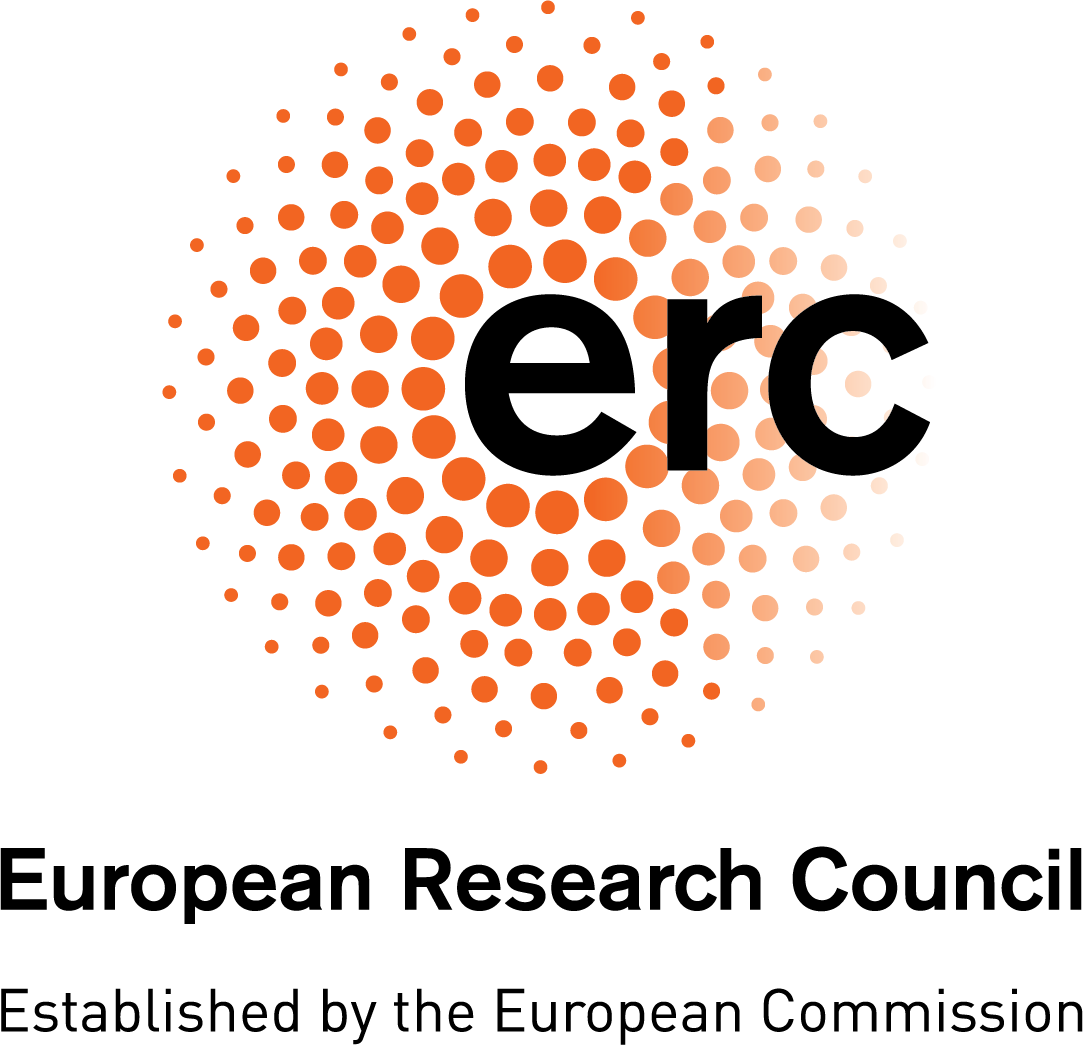}}%
}

\usepackage{tcolorbox}
\tcbset{colback=white,boxrule=0.5pt,arc=2pt,left=4pt,right=4pt,top=4pt,bottom=4pt}

\theoremstyle{plain}
\newtheorem{theorem}{Theorem}[section]
\newtheorem{lemma}[theorem]{Lemma}
\newtheorem{proposition}[theorem]{Proposition}
\newtheorem{corollary}[theorem]{Corollary}
\theoremstyle{definition} % upright font, not italic
\theoremstyle{definition}
\newtheorem{definition}{Definition}
\newtheorem{assumption}{Assumption}

\theoremstyle{remark}
\newtheorem{remark}{Remark}

\usepackage{authblk}

\let\oldproof\proof
\let\oldendproof\endproof

  \newlist{claimlist}{enumerate}{1}
  \setlist[claimlist]{label= (\arabic*)}
  \crefname{claimlisti}{Claim}{Claims}
  \Crefname{claimlisti}{Claim}{Claims}
\crefname{theorem}{Theorem}{Theorems}
\crefname{lemma}{Lemma}{Lemmas}
\crefname{corollary}{Corollary}{Corollaries}
\crefname{section}{Section}{Sections}
\crefname{proposition}{Proposition}{Propositions}
\crefalias{lemma}{lemma}
\crefalias{corollary}{corollary}
\crefalias{proposition}{proposition}
\crefname{definition}{Definition}{Definitions}
\crefname{remark}{Remark}{Remarks}
\crefname{conjecture}{Conjecture}{Conjectures}
\crefname{example}{Example}{Examples}
\crefname{annoyingproof}{Proof}{Proofs}
\crefname{assumption}{Assumption}{Assumptions}

\let\ORGhypersetup\hypersetup
\protected\def\hypersetup{\ORGhypersetup}
\title{Newtonian Dynamics in  Metric Spaces}

\author[1,2]{Massimo Fornasier\thanks{\texttt{massimo.fornasier@cit.tum.de}}}
\author[1,2]{Jona Klemenc\thanks{\texttt{jona.klemenc@tum.de}}}
\author[1,2]{Alessandro Scagliotti\thanks{\texttt{scag@ma.tum.de}}}

\affil[1]{Department of Mathematics, CIT School, Technical University of Munich, Munich, Germany}
\affil[2]{Munich Center for Machine Learning (MCML), Munich, Germany}

\date{\today}
\begin{document}

\maketitle
\begin{abstract}
 
We develop a variational framework for second-order, frictionless (Newtonian) dynamics in metric spaces, extending the metric theory of curves of maximal slope  from gradient flows to evolutions with inertia. The construction rests on three pillars, mirroring the first-order theory. First, we introduce a class of second-order absolutely continuous curves, characterized by a geometric--variational inequality on families of tangent geodesics, and define a metric notion of acceleration. 
Second, we define \emph{curves of Newtonian motion} as the curves achieving equality in an energy inequality involving an upper gradient--Hessian pair, which in Hilbert spaces is equivalent to Newton's equation $\ddot\gamma = -\nabla E(\gamma)$. Third, we propose a two-step predictor--corrector minimizing movement scheme and show that its discrete interpolants satisfy approximate versions of these inequalities. Under extendability and interceptor conditions on geodesics---which are satisfied, for instance, by geodesically uniformly complete, locally compact spaces with curvature bounded below---and under suitable regularity assumptions on the energy, we prove that the scheme converges, up to subsequences, to a curve of Newtonian motion. As an illustration, we recover the linear wave equation as a curve of Newtonian motion driven by the Dirichlet energy. This paper is complemented by the companion paper \cite{fornasierLionsLiftCompanion}, which specifies and adapts the results to the Wasserstein metric space $(\Pp_2,W_2)$ through Lions lift, which explains Vlasov equations as Newtonian flows on spaces of probability measures.
\end{abstract}



\noindent
\textit{AI use disclosure:} 
All concepts, research questions and proofs were developed by the authors. 
AI tools were employed for bug checking, copy-editing tasks and for the coding of the figures.
For more details, see \cref{tab:ai-disclosure}.

\tableofcontents

\section*{Introduction}
\renewcommand{\theequation}{I.\arabic{equation}}
In a Hilbert space $\mathcal{H}$, the gradient flow equation
\begin{equation} \label{eq:intro_gf}
	\dot\gamma_t = -\nabla E \big( \gamma_t \big)
\end{equation}
describes the trajectory $t\mapsto \gamma_t$ of an inertia-free particle $\gamma$ descending along the steepest slope of a smooth energy landscape $E\colon \mathcal{H} \to \R$. 
Its metric generalization in the sense of \emph{curves of maximal slope} (see the textbook \cite{ambrosioGradientFlows2008}) through a variational framework has unlocked a wide range of applications, ranging from optimization in metric spaces to the characterization of solutions to dissipative PDEs (\cite{jordan1998variational,OttoGeometry2001}).
Likewise, its higher-order sister equation, i.e., the Newtonian flow equation
\begin{equation}\label{eq:intro_newtonian_eq}
	\ddot\gamma = - \nabla E \big( \gamma_t \big),
\end{equation}
describes a \emph{massive} particle (i.e., with inertia) following Newton's second law 
%($F = ma$)
in a force field $F\colon \mathcal{H} \to \mathcal{H}$ induced by the potential energy $E$ (i.e., $F = - \nabla E$). 
\emph{Its} metric generalization through a variational framework promises progress on a further slate of applications, from inertia-driven optimization---such as heavy-ball or Nesterov-like methods---to momentum-driven PDEs---such as the wave equation, or Vlasov kinetic equations.

In this direction, there exist multiple strains of research.
An example which stays close in spirit to the program of curves of maximal slope is the construction of solutions to the anisotropic Euler equation proposed in \cite{westdickenbergProjectionsConeOptimal2010}, even though it focuses on a single application and does not come with a variational characterization of the resulting limit curves. 
A related, though structurally different, approach is the theory of Hamiltonian ODEs in Wasserstein spaces developed in \cite{ambrosio2008hamiltonian}, where conservative dynamics are formulated as first-order evolutions on $\mathcal{P}_2(\R^d\times \R^d)$, i.e., the space of probability measures on the phase space (see also the recent contribution \cite{brigatiKineticOptimalTransport2025}). 
In \cite{ambrosio2008hamiltonian} the authors relied on the position--velocity structure of the phase space $\R^d\times \R^d$. By contrast, our aim here is to develop an \emph{intrinsic second-order variational framework} directly on the metric space $(X,d_X)$.

In the present work, we attempt a more direct lift of the framework of \cite{ambrosioGradientFlows2008} to second-order evolutions. 
To this end, we first distill and recall three pillars on which the (so far) classical construction of curves of maximal slope rests, and which we aim to replicate for second-order evolutions.
\begin{enumerate}
	\item \textbf{A metric generalization of absolutely continuous curves.} A variational characterization yields the class of \emph{absolutely continuous} (or $\AC$-) curves with values in a metric space, and provides the metric definition of the scalar velocity $|\gamma'|$, for every $\AC$-curve $\gamma$. Moreover, this class includes curves of maximal slope. 
    More precisely, we recall that, given a metric space $(X, d_X)$, $\AC^2$-curves are defined to be those curves $\gamma\colon (a, b) \to X$ (for $a, b \in \R$ and $a < b$) for which there exists $m \in L^2((a, b), \R)$ such that
	\begin{equation}\label{eq:intro_ac_ineq}
		d_X(\gamma_s, \gamma_t) \leq \int_s^t m(r) \, \diff r, \qquad \forall s, t \in (a, b) \text{ s.t. } s \leq t.
	\end{equation}
	For every $\gamma \in \AC^2((a,b),X)$, the scalar velocity is defined as $|\gamma'|(t) = \lim_{s \to t}\frac{d_X(\gamma_s, \gamma_t)}{|s-t|}$, where the right-hand side is guaranteed to exist for almost every $t\in (a,b)$. It is then true that, for every $\gamma$ and $m$ as above, we have that $m(t) \geq |\gamma'|(t)$ for almost every $t$, and that $m \coloneqq |\gamma'|$ is admissible to make \cref{eq:intro_ac_ineq} true.
    
	\item \textbf{A metric generalization of gradient flow lines.} As before, this class of curves is characterized by a variational inequality. Namely, for a metric space $(X,d_X)$ and an energy $E\colon X\to \R\cup \{+\infty\}$, we say that $\gamma \in \AC^2((a,b),X)$ is a \emph{curve of maximal slope} if it fulfills the inequality
	\begin{equation} \label{eq:intro_coms_ineq}
		E(\gamma_s) - E(\gamma_t) \geq \int_s^t 
        \left(\frac{1}{2} |\gamma'|^2(r) + \frac{1}{2} g^2 \big(\gamma_r \big) \right) \, \diff r, \qquad \forall s, t \in (a, b) \text{ s.t. } s \leq t
	\end{equation}
	for an \emph{upper gradient} $g\colon X \to [0, +\infty]$, i.e., a proper and measurable function such that $E(\eta_s) - E(\eta_t) \leq \int_s^t |\eta'|(r)g(\eta_r) \, \diff r$ for every $\eta \in \AC^2((a,b),X)$. 
    We first notice that, whenever $g\colon X \to [0, +\infty]$ is an upper gradient, the inverse inequality of \eqref{eq:intro_coms_ineq} holds automatically for every $\eta \in \AC^2((a,b),X)$.
    Moreover, it is important to recall that, in the case of a smooth energy function $E$ defined on a Hilbert space $(\mathcal H, |\cdot |)$, the function $x\mapsto g(x) \coloneqq |\nabla E (x)|$ is an upper gradient. 
    In the general case, the \emph{metric slope} $|\partial E| \colon X \to [0, +\infty] $ is the natural generalization of $x\mapsto |\nabla E (x)|$, while the fact that (the lower semi-continuous envelope of) $|\partial E|$ is an upper gradient has to be proven on a case-by-case basis.

	\item \textbf{A discrete, variational scheme to generate and exhibit curves of maximal slope}. This scheme is called \emph{minimizing movement scheme} and is defined recursively by
	\begin{equation}
		U_{\tau, i+1} \in J_\tau[U_{\tau, i}]
	\end{equation}
	where $\tau > 0$ is the step size and 
	\begin{equation}
		J_\tau[x] \coloneqq \argmin_{y \in X} \Big\{E(y) + \frac{1}{2\tau}d_X(x, y)^2\Big\}.
	\end{equation}
	The piecewise constant interpolations $u_\tau\colon [0, \infty) \to X$ between the points $(U_{\tau, i})_{i \in \NN}$ fulfill approximate versions of the inequalities \eqref{eq:intro_ac_ineq} and \eqref{eq:intro_coms_ineq}.
\end{enumerate}
Finally, we mention that, under suitable assumptions, a compactness argument allows one to obtain a limit curve $\bar{u} \colon [0, \infty) \to X$ by selecting a sequence $(u_{\tau_n})_{n \in \NN} \subset (u_{\tau})_{\tau>0}$ with $\tau_n \to 0$ from the piecewise constant interpolations, and this binds together Pillar 1, Pillar 2 and Pillar 3.
%Binding together Pillar 1, Pillar 2 and Pillar 3 is the fact that under suitable assumptions and using a compactness argument, we obtain a limit curve $\bar{u}$ of the interpolating curves $u_\tau$ as $\tau \to 0$. 
Indeed, by taking the limit of the approximate inequalities fulfilled for every $\tau>0$ by $u_\tau$, it is possible to deduce that $\bar{u}$ is a curve of maximal slope.

In the present paper, we aim to build analogous pillars for the second-order conservative (i.e., frictionless) dynamics.
The main assumptions on the metric space $(X,d_X)$, besides completeness, local compactness and being geodesic, are an extendability and an interceptor condition on its geodesics (see \cref{ass:interceptor}). 
On the one hand, the interceptor condition implies a non-branching condition and is implied by a bound from below of the curvature (CBB). 
On the other hand, the extendability assumption excludes a direct application to the Wasserstein space $W_2$, where the framework of curves of maximal slope has seen tremendous success. 
We circumvent this issue by taking advantage of the Lagrangian route in the companion paper \cite{fornasierLionsLiftCompanion}. More tame spaces, such as compact manifolds, are covered by the assumptions of the present paper.
More generally, geodesically uniformly complete\footnote{By this, we mean that there exists a uniform horizon such that any unit-speed geodesic can be extended to this uniform horizon, see \cref{ass:interceptor}.}, locally compact CBB spaces are covered as well.

We are now in a position to illustrate the three pillars that support our construction.\\

\textbf{Pillar 1: A geometric-variational characterization of second-order absolutely continuous curves in metric spaces.} To this end, we define the space of geodesics\footnote{By \emph{geodesics in $X$}, we mean continuous curves $\xi\colon [a, b]\to X$ (with $a<b$) such that, for every $s, t \in [a, b]$, it holds that $	d_X(\xi_s, \xi_t) = |s - t| \frac{1}{b-a} d_X(\xi_a, \xi_b)$.}
\begin{equation} \label{eq:intro_geod_space}
	\Tang[\varepsilon] \coloneqq \{\xi \in AC([t, t + \varepsilon], X) \mid \xi \text{ is a constant-speed geodesic in } X, t \in \R\}, \quad \varepsilon > 0.
\end{equation}
Then, we say that a curve $\gamma\colon (a, b) \to X$ is \emph{\ACt} if there exists $\varepsilon_\gamma > 0$, $\mathfrak{E}_\gamma \geq 0$, $m \in L^2((a, b), \R)$, and a map $V^\gamma\colon (a, b)\to \Tang[\varepsilon_\gamma]$ with $V^\gamma_t  \in AC([t, t + \varepsilon], X)$ and $V^\gamma_t(t) = \gamma_t$ such that 
\begin{equation} \label{eq:intro_ac2}
			d_X \big( V_s^\gamma(T), V_t^\gamma(T) \big) \leq \int_s^t (T-r) \, m(r) \big(1 + \mathfrak{E}_\gamma  (T-r)^2\big) \, \diff r
\end{equation}
for every $s, t \in (a, b)$ and $T \in \R$ with $s \leq t \leq T \leq s + \varepsilon_\gamma$. 
If $\gamma \in \ACt((a,b),X)$ and $(X,d_X)$ satisfies the interceptor condition in \cref{ass:interceptor}, then the function $t\mapsto V^\gamma_t$ is uniquely determined up to a negligible set of instants (see \cref{rmk:tangents_recoverable}). For (almost every) instant $t$, $V^\gamma_t$ should be interpreted as the geodesic with scalar speed $|\gamma'_t|$ and \emph{tangent} to the curve $\gamma$ at $\gamma_t$.
To see this, we note that, in a Hilbert space $(\mathcal H, |\cdot |)$, if $\gamma \in W^{2,2}((a,b),\mathcal{H})$, setting $V^\gamma_t(T) \coloneqq \gamma_t + (T - t) \dot\gamma_t$, we have that $\big|V_s^\gamma(T) - V_t^\gamma(T) \big| \leq \int_s^t (T-r) \, m(r) \, \diff r$, and so \cref{eq:intro_ac2} is fulfilled with $m(t) = | \ddot \gamma_t |$, $\mathfrak{E}_\gamma = 0$ and for every $\varepsilon_\gamma > 0$.
We allow for $\mathfrak{E}_\gamma > 0$ to account for a possible curvature in the space $X$\footnote{When constructing curves of Newtonian motion, $\mathfrak{E}_\gamma$ can in practice be chosen to equal $C^2 \mathfrak{E}$, where $\mathfrak{E}$ comes from the \EpsE-interceptor bound (see \cref{def:interceptor_bound,ass:interceptor}), and $C$ is an upper bound on $|\gamma'|$. 
As shown in \cref{lem:cbb_interceptor_bound}, if $(X,d_X)$ is CBB (i.e., it has curvature bounded from below), then the \EpsE-interceptor bound holds.}.
Moreover, in a general metric space $(X,d_X)$, an upper bound on $\varepsilon_\gamma$ is needed as length-minimizing geodesics may not be arbitrarily extended.
Similarly as in the first-order case, the quantity
\begin{equation} \label{eq:intro_second_derivative_bound}
		|\gamma''|(s) \coloneqq \limsupTauTS\left(\frac{d_X \big(V^\gamma_s(T) V^\gamma_t(T)\big)}{(T-s)(t-s)}\right)
\end{equation}
plays the role of a \emph{metric acceleration}, and it turns out that $|\gamma''|(s) \leq m(s)$ for almost every $s \in (a, b)$. 
In addition, under a further curvature assumption on $X$ (a converse assumption to the interceptor bound, which is implied by a CAT\footnote{CAT stands for ``Élie \textbf{C}artan, Aleksandr Danilovich \textbf{A}leksandrov and Victor Andreevich \textbf{T}oponogov'', and denotes the property of having curvature bounded above---i.e., the inverse of CBB (\textbf{C}urvature \textbf{B}ounded from \textbf{B}elow).} property, together with a symmetry assumption, see \cref{prop:second_derivative_is_minimal}), the choice $m \coloneqq |\gamma''|$ for $\gamma \in \ACt((a,b),X)$ is admissible to make \cref{eq:intro_ac2} true.\\

\textbf{Pillar 2: A geometric, variational generalization of Newtonian flows.} To this end, we notice that in a Hilbert space $(\mathcal{H}, |\cdot|)$ and a smooth energy $E \colon \mathcal{H} \to \R$, the inequality
\begin{equation}\label{eq:intro_inequality}
	\langle \dot \gamma_s, \nabla E(\gamma_s) \rangle - \langle \dot \gamma_t, \nabla E(\gamma_t) \rangle \leq \int_s^t \frac{1}{2} |\ddot \gamma_r|^2  + \frac{1}{2} |\nabla E(\gamma_r)|^2 - \nabla^2 E(\gamma_r)(\dot \gamma_r, \dot \gamma_r) \, \diff r
\end{equation}
is fulfilled automatically for every $\gamma \in W^{2,2}((a,b),\mathcal{H})$ and every $s,t,\in(a,b)$ with $s\leq t$. Indeed, it is a direct consequence of the Young inequality yielding $-\langle \ddot \gamma_r, \nabla E(\gamma_r) \rangle \leq \frac{1}{2} |\ddot \gamma_r|^2 + \frac{1}{2} |\nabla E(\gamma_r)|^2$, and equality in \cref{eq:intro_inequality} implies that $\gamma$ fulfills \cref{eq:intro_newtonian_eq}. 
To generalize this variational characterization to metric spaces, we perform the following replacements:
\begin{equation*}
    \begin{split}
        & \langle \dot\gamma_s, \nabla E(\gamma_s)\rangle \longleftrightarrow \frac{\diff}{\diff \sigma} E \big(V_s^\gamma(\sigma) \big)\Big|_{\sigma = s}, \\
        & |\ddot \gamma_t| \longleftrightarrow |\gamma''|(t) \qquad \mbox{ (see \cref{eq:intro_second_derivative_bound})},
    \end{split}
\end{equation*}
and finally we replace the pair $\big(|\nabla E(\gamma_t)|, \nabla^2 E(\gamma_t)(\dot \gamma_t, \dot \gamma_t)\big)$ by a proper pair $(g, H)$, where $g\colon X \to [0, +\infty)$ and $H\colon 
\bigcup_{\varepsilon > 0}\Tang[\varepsilon] \to \R$ (see \cref{eq:intro_geod_space} for the definition of $\Tang[\varepsilon]$). 
We call such a replacement pair an \emph{upper gradient-Hessian pair} if the resulting inequality is fulfilled automatically for any $\eta \in \ACt((a,b),X)$, and we call such an \ACt-curve a \emph{curve of Newtonian motion} if equality is achieved. 
We thus obtain the characterization of curves of Newtonian motion as those curves $\gamma \in \ACt((a,b),X)$ that fulfill
\begin{equation} \label{eq:intro_curves_of_newtonian_motion}
		\frac{\diff}{\diff \sigma} E\big( V_s^\gamma(\sigma) \big) \Big|_{\sigma = s} - \frac{\diff}{\diff \sigma} E\big( V_t^\gamma(\sigma) \big) \Big|_{\sigma = t} 
        \geq \int_s^t \frac{1}{2}  |\gamma''|^2(r)\, \diff r + \int_s^t \frac{1}{2} g^2(\gamma_r)\, \diff r - \int_s^t H(V_r^\gamma)\, \diff r
\end{equation}
for every $s,t \in(a,b)$ with $s\leq t$.
When dealing with a Hilbert space $\mathcal{H}$ and a smooth energy $E\colon \mathcal{H} \to \R$, we can provide a more explicit description for these abstract objects.
Namely, recalling that any geodesic $\xi \in \bigcup_{\varepsilon > 0}\mathcal T_\varepsilon \mathcal{H}$ has the form $[t_\xi, t_\xi+\varepsilon] \ni \sigma \mapsto x_\xi + \sigma v_\xi$ for some $x_\xi , v_\xi \in \mathcal{H}$ and some $\varepsilon>0$, the pair
\begin{equation}
	\Big(g\colon x \mapsto |\nabla E(x)|, \quad H\colon \xi \mapsto \nabla^2 E(x_\xi) [v_\xi,v_\xi] \Big)
\end{equation}
is an upper gradient-Hessian pair, and a Sobolev curve $\gamma \in W^{2,2}((a,b), \mathcal{H})$ is a curve of Newtonian motion if and only if it fulfills \cref{eq:intro_newtonian_eq} for almost every $t$.
Finally, in \cref{rmk:friction_extension} we provide the variational characterization of solutions of the \emph{damped Newtonian motion} $\ddot \gamma_t = -\nabla E(\gamma_t)-\alpha\dot \gamma_t$, for $\alpha\in \R$. 
This fact paves the way to the study of Newtonian motion with viscous friction, even though this extension goes beyond the scope of the present paper and it is left for future work. \\

\textbf{Pillar 3: A discrete, variational scheme to construct and exhibit curves of Newtonian motion}. 
Here, we propose a two-step predictor-corrector scheme: For a fixed step size $\tau > 0$ and given the points $U_{\tau, i-1}$ and $U_{\tau, i}$ for $i \in \NN$, we recursively define
\begin{equation} \label{eq:definition_scheme_intro}
	\begin{split}
	Z_{\tau, i} &\coloneqq \argmin_{z \in X} \Big\{2d_X^2(U_{\tau, i}, z) - d_X^2(U_{\tau, i-1}, z) \Big\}, \\
	U_{\tau, i+1} &\in \argmin_{z \in X} \Big\{E(z) + \frac{1}{2\tau^2} d_X^2(Z_{\tau, i}, z)\Big\}.
	\end{split}
\end{equation}
We notice that the temporary point $Z_{\tau, i}$ is the metric analogue of $U_{\tau, i} + \big( U_{\tau, i} - U_{\tau, i-1} \big)$, and the extendability of geodesics is needed here to deduce that, if $U_{\tau, i}, U_{\tau, i-1}$ are sufficiently close, then $Z_{\tau, i}$ is such that
\begin{equation*}
    d_{X}(U_{\tau, i-1}, U_{\tau, i}) = d_{X}(U_{\tau, i}, Z_{\tau, i}) = \frac12 d_{X}(U_{\tau, i-1}, Z_{\tau, i}).
\end{equation*}
Moreover, the interceptor condition on geodesics gives uniqueness in the predictor step.
We show that the piecewise constant interpolations $u_\tau$ between the points $(U_{\tau, i})_{i \in \NN}$ together with ``approximate tangents'' $\mathcal{V}^{\tau}$ fulfill approximate versions of the inequalities \eqref{eq:intro_ac2} and \eqref{eq:intro_curves_of_newtonian_motion}. To initialize the scheme, we use a prescribed geodesic $\xi \colon [-\varepsilon, \varepsilon] \to X$ serving as an initial velocity vector, and set $U_{\tau, 0} = \xi_0$ and $U_{\tau, -1} = \xi_{-\tau}$ for $\tau \leq \varepsilon$.

Obtaining a limit curve requires assumptions on $E$: lower semicontinuity, coercivity, existence and regularity of the derivatives $\frac{\diff}{\diff s} E(\xi_s)$ for geodesics $\xi$\footnote{The main regularity assumption stems from the fact that the proxy $H$ for the Hessian in \cref{eq:intro_curves_of_newtonian_motion} has no sign, and thus getting meaningful a-priori estimates requires additional regularity. We expect that one can replace this assumption by ad-hoc arguments in applications of interest, as we did in the wave equation example.}; furthermore, we require a relaxation of\footnote{The precise pair we require to be an upper gradient-Hessian pair is $(|\partial^-E|, H_0^{E, +})$, defined in \cref{eq:def_lowerslope,eq:def_H0_plus}.} $(|\nabla E|, \xi \mapsto \frac{\diff^2}{\diff\sigma^2} E(\xi_\sigma)|_{\sigma=t_\xi})$ to be an upper gradient-Hessian pair.
Under these assumptions, we can obtain a limit \ACt-curve $\bar{u}\colon (a,b)\to X$ by extracting a sequence $(u_{\tau_n})_{n\in \NN} \subset (u_{\tau})_{\tau>0} $ with $\tau_n\to 0$ from the family of interpolating curves.
This construction consists of two steps: First, we use a compactness argument in the space of geodesics $\Tang[\varepsilon]$ to obtain a limit $V\colon (a,b) \to \Tang[\varepsilon]$ computed along a sequence $\{ \mathcal{V}^{\tau_n} \}_{n\in \NN} \subset \{ \mathcal{V}^{\tau} \}_{\tau>0}$ of the ``approximate tangents''. 
In a second step, recalling that for every $t\in (a,b)$ we have a geodesic $V_t\colon [t,t+\varepsilon]\to X$, we read off the limit curve $\bar{u}_t \coloneqq V_t(t)$.
Finally, we take the limit as $\tau_n\to 0$ in the approximate inequalities fulfilled for every $\tau>0$ by $u_{\tau}$ and $\mathcal{V}^\tau$, and we deduce that $\bar{u}$ is an \ACt-curve---with $V^{\bar{u}} = V$---as well as a curve of Newtonian motion.

We conclude the paper by providing an example usage of the second-order framework: The derivation of the linear wave equation as a curve of Newtonian motion along the driving Dirichlet energy.
This construction shows that the three pillars can be used as a flexible framework rather than a rigid set of theorems. Indeed, even though the Dirichlet energy lacks the regularity required for the direct application of the estimates derived in the abstract metric setting, those can readily be replaced by ad-hoc arguments.
Thus, while this final section does not yield a new existence theorem or weaken the classical regularity assumptions for the linear wave equation---in fact, our full existence statement requires $u_0 \in W^{3,2}(\R^d, \R)$ and $\dot u_0 \in W^{3,2}(\R^d, \R)$, both with compact support---it illustrates that the three pillars can be seen as guiding principles beyond the direct application of the abstract theory.\\

The example of the derivation of the wave equation as Newtonian dynamics is by no means the only application of the theory. 
Provided the well-known application of the metric gradient flow theory to Wasserstein spaces $(\Pp_p,W_p)$ and its implications for the descriptions of dissipative PDEs, a natural question is whether the abstract theory of the present paper can be equally applied in this context. In the companion paper \cite{fornasierLionsLiftCompanion} we 
explore this question, which is rendered particularly difficult by the fact that the geodesics in the Wasserstein space $(\Pp_2,W_2)$ are not extendable in general and  \cref{ass:interceptor} does not hold. The way to circumvent this technical difficulty is to identify $(\Pp_2,W_2)$ with the quotient space of $L^2(\Omega,\mathbb R^d)$ with respect to the law function $\mathcal L(X)= X_\#\mathbb P$ (we will use specifically $\Omega=[0,1]^{\mathbb R}$ endowed with the product Lebesgue measure). This identification, known as Lions lift, allows us to leverage the Hilbert space structure of $L^2(\Omega,\mathbb R^d)$ and to derive henceforth Newtonian curves on $(\Pp_2,W_2)$ by  projection. As a consequence, we are able to explain PDE of Vlasov type with suitable integral and interaction potentials as Newtonian curves on the Wasserstein space.

\subsection*{General notations}
In this section, we collect a few notations that we will use throughout the manuscript. 
We will denote by $(X, d_X)$ a metric space. By \emph{geodesics in $X$}, we mean continuous curves $\xi\colon [a, b]\to X$---for some $a,b\in \R$ with $a<b$---such that, for every $s, t \in [0, T]$, it holds that
\begin{equation*}
	d_X(\xi_s, \xi_t) = |s - t| \, \frac{d_X(\xi_a, \xi_b)}{b-a}.
\end{equation*}
It follows from the definition that geodesics are absolutely continuous in the sense of \cite{ambrosioGradientFlows2008}, and that for every such geodesic $\xi$, we have that $|\xi'|(t) = \frac{d_X(\xi_a, \xi_b)}{b-a}$ for every $t \in [a, b]$. 
If it is clear that some curve $\xi$ is a geodesic, we will sometimes simply write $|\xi'|$ instead of $\frac{d_X(\xi_a, \xi_b)}{b-a}$.

 % !TEX root = main_metric.tex
\begingroup
\setlength{\LTpre}{4pt}
\setlength{\LTpost}{4pt}
\renewcommand{\arraystretch}{1.1}

\begin{longtable}{@{} p{0.30\textwidth} p{0.65\textwidth} @{}}
\toprule
\textbf{Symbol} & \textbf{Description} \\
\midrule
\endfirsthead
\multicolumn{2}{l}{\emph{(continued)}} \\
\toprule
\textbf{Symbol} & \textbf{Description} \\
\midrule
\endhead
\midrule
\multicolumn{2}{r}{\emph{(continued on next page)}} \\
\endfoot
\bottomrule
\endlastfoot

%------------------------------------------------------------------
\multicolumn{2}{@{}l}{\textit{Spaces and sets}} \\[1pt]
$(X, d_X)$
  & Metric space with distance $d_X$. \\
$(\mathcal{H}, |\cdot|)$
  & Hilbert space with norm $|\cdot|$. \\
$\mathbb{R}^d$
  & $d$-dimensional Euclidean space. \\
$\mathbb{M}_2(\kappa)$
  & Two-dimensional model space of constant sectional curvature $\kappa$. \\
$AC([a,b], X)$
  & Absolutely continuous curves in $X$ on $[a,b]$. \\
\ACt{}
  & Second-order absolutely continuous curves; see \cref{def:AC2_curve}. \\
$\Tang[t,\varepsilon]$
  & Space of geodesics in $X$ starting at time $t$ with domain $[t,t+\varepsilon]$;
    $\Tang[\varepsilon] \coloneqq \bigcup_{t}\Tang[t,\varepsilon]$;
    $\mathcal{T}(X) \coloneqq \bigcup_{\varepsilon>0}\mathcal{T}_\varepsilon(X)$;
    see \cref{def:space_of_geodesics}. \\[4pt]

%------------------------------------------------------------------
\multicolumn{2}{@{}l}{\textit{Curves and their derivatives}} \\[1pt]
$\gamma,\, \eta$
  & Generic curves in $X$. \\
$\xi,\, \zeta$
  & Geodesics in $X$. \\
$|\gamma'|(t)$
  & Metric derivative: $\lim_{s\to t}d_X(\gamma_s,\gamma_t)/|s-t|$. \\
$|\gamma''|(s)$
  & Second-order metric derivative:
    $\limsup_{\substack{T\searrow s,\,t\searrow s\\s<t\leq T}}\frac{d_X(V^\gamma_s(T),V^\gamma_t(T))}{(T-s)(t-s)}$;
    see \cref{eq:metric_second_derivative_bound}. \\
$V^\gamma_t$
  & Tangent geodesic to an \ACt{}-curve $\gamma$ at $t$;
    $V^\gamma_t\in\Tang[t,\varepsilon_\gamma]$, $V^\gamma_t(t)=\gamma_t$;
    see \cref{def:AC2_curve}. \\
$\varepsilon_\gamma$, $\mathfrak{E}_\gamma$
  & Geodesic horizon and curvature bound of an \ACt{}-curve $\gamma$;
    see \cref{def:AC2_curve}. \\[4pt]

%------------------------------------------------------------------
\multicolumn{2}{@{}l}{\textit{Energies, slopes, and upper gradient-Hessian pairs}} \\[1pt]
$E$
  & Energy functional $E\colon X\to\mathbb{R}\cup\{+\infty\}$ on a metric space. \\
$|\partial E|(x)$
  & (Upper) slope: $\limsup_{y\to x}(E(x)-E(y))_+/d_X(x,y)$;
    see \cref{eq:def_slope}. \\
$|\partial^- E|(x)$
  & Lower slope (lower semicontinuous envelope of $|\partial E|$);
    see \cref{eq:def_lowerslope}. \\
$H^E_\tau(\xi_t,\xi_{t+\tau})$
  & Discrete Hessian:
    $\tau^{-2}\bigl(E(\xi_t)-2E(\xi_{t+\tau})+E(\xi_{t+2\tau})\bigr)$;
    see \cref{eq:def_H}. \\
$H^E_0(\xi)$, $H^{E,+}_0(\xi)$
  & Limit Hessian and its upper-semicontinuous regularisation;
    see \cref{eq:def_H0,eq:def_H0_plus}. \\
$(g,H)$
  & Upper gradient-Hessian pair for $E$; see \cref{def:upper_gradient}. \\
$J_\tau[x]$
  & Proximal map:
    $\argmin_{y\in X}\bigl(E(y)+\tfrac{1}{2\tau}d_X(x,y)^2\bigr)$;
    see \cref{eq:def_proximal}. \\[4pt]

%------------------------------------------------------------------
\multicolumn{2}{@{}l}{\textit{Discrete scheme}} \\[1pt]
$\tau$
  & Time-discretisation step size; see \cref{def:newtonian_sequence}. \\
$U_{\tau,k}$
  & $k$-th iterate; see \cref{def:newtonian_sequence}. \\
$U'_{\tau,k}$
  & Discrete velocity in the Hilbert-space setting (in particular, the wave application):
    $U'_{\tau,k}\coloneqq(U_{\tau,k}-U_{\tau,k-1})/\tau$. \\
$Z_{\tau,i}$
  & Predictor point, the metric analogue of $U_{\tau,i}+(U_{\tau,i}-U_{\tau,i-1})$:
    $Z_{\tau,i}\coloneqq\argmin_{z\in X}\bigl(2d_X^2(U_{\tau,i},z) - d_X^2(U_{\tau,i-1},z) \bigr)$;
    see \cref{def:newtonian_sequence}. \\[4pt]

%------------------------------------------------------------------
\multicolumn{2}{@{}l}{\textit{Miscellaneous}} \\[1pt]
$\Phi(\xi,\zeta)$
  & Angle between two geodesics $\xi,\zeta$ at their common starting point. \\
$\EpsE$
  & $\epsilon$-$\mathfrak{E}$ interceptor bound pair; see \cref{def:interceptor_bound}. \\
$\diff$
  & Differential operator ($\diff t$, $\diff x$, etc.). \\
$\mathbb{N},\,\mathbb{R}, \, \R_+$
  & Natural numbers $\mathbb{N}=\{0,1,2,\ldots\}$, real numbers $\R$, non-negative real numbers $\R_+ = [0,+\infty)$. \\
\end{longtable}
\endgroup

\renewcommand{\theequation}{\thesection.\arabic{equation}} \setcounter{equation}{0}

\section{First Pillar: Second-Order Absolute Continuity} \label{sec:defining_second-order_absolute_continuity}

In this section, we define the notion of \ACt-curves in metric spaces, which is the second-order analogue of the metric concept of \AC-curves \footnote{We stress the following distinction: An \ACt-curve is one which fulfills our \cref{def:AC2_curve}. A $\mathrm{AC}^p$-curve $\eta\colon (a, b) \to X$ is one for which there exists $m \in L^p((a, b))$ such that for all $s, t \in (a, b)$ with $s < t$, we have that $d_X(\eta_s, \eta_t) \leq \int_s^t m(r) \, \diff r$. In the latter definition (see, e.g., \cite[Def. 1.1.1]{ambrosioGradientFlows2008}), the superscript $p$ refers to the integrability class of $m$, while in the former, the index $\mathrm{II}$ refers to the second-order nature of the definition.}. 
Our definition has the form of a variational inequality, and it involves a family of geodesics defined on time-intervals whose time-lengths need to be uniformly positive.
This fact, in practice, restricts the choice of metric spaces where our definition can be used. For example, in Wasserstein spaces geodesics are in general not extendable, and this can create obstructions to the construction of such a family of geodesics. The case of Wasserstein spaces is directly addressed in the companion paper \cite{fornasierLionsLiftCompanion}.

The spaces described in the next definition play a pivotal role in our construction.

\begin{definition}[Space of geodesics] \label{def:space_of_geodesics}
	For every $\varepsilon > 0$ and $t \in \R$, we denote the \emph{space of geodesics with domain $[t, t+\varepsilon]$} with
	\begin{equation*}
		\Tang[t, \varepsilon] \coloneqq \big\{\xi \in \AC([t, t + \varepsilon], X) \mid \xi \text{ is a geodesic in } X \big\}.
	\end{equation*}
	Moreover, we define the \emph{space of geodesics whose domain has length $\varepsilon$} as
    $$\Tang[\varepsilon] \coloneqq \bigcup_{t \in \R} \Tang[t, \varepsilon],$$ 
    and finally the \emph{space of geodesics} as
    $$\mathcal{T} X = \bigcup_{\varepsilon > 0} \Tang[\varepsilon].$$
    Furthermore, we equip the space $\Tang[t, \varepsilon]$ with the supremum distance
	\begin{equation*}
		d_{\Tang[t, \varepsilon]}(\xi, \zeta) \coloneqq \sup_{s \in [t, t + \varepsilon]} d_X(\xi_s, \zeta_s).
	\end{equation*}
\end{definition}

\begin{remark}
    The following object may look closer to the usual notion of tangent space in differential geometry:
    \begin{equation*}
        T_{\varepsilon,x}X \coloneqq \big\{ \xi \in \AC([0, \varepsilon], X) \mid \varepsilon>0, \ \xi \text{ is a geodesic in } X, \ \xi(0)=x \big\}, 
    \end{equation*}
    and we notice that $\Tang[0, \varepsilon] = \bigcup_{x\in X} T_{\varepsilon,x}X$. However, in our construction, we are not actually making use of $ T_{\varepsilon,x}X$.
\end{remark}

We are now in a position to introduce the notion of \ACt-curve.

\begin{definition}[\ACt-curves] \label{def:AC2_curve}
	We say that a curve $\gamma\colon (a, b)\to X$, for $a, b \in \R$ and $a < b$, is \emph{second-order absolutely continuous} (or is an \emph{\ACt-curve}, or $\gamma \in \ACt((a,b),X)$), if there exist
	\begin{enumerate}
		\item A geodesic horizon $\varepsilon_\gamma > 0$;
		\item A function $V^\gamma\colon (a, b)\to \Tang[\varepsilon_\gamma]$ with $V^\gamma \colon t\mapsto V^\gamma_t$, such that for every $t \in (a, b)$, we have $V_t^\gamma \in \Tang[t, \varepsilon_\gamma]$ and it holds that $V_t^\gamma(t) = \gamma_t$;
		\item A function $m \in L^2\big((a, b), \R_+ \big)$ \footnote{The integrability regularity $L^2$ is required in the variational characterization of curves of Newtonian motion provided in \cref{sec:curves_of_newtonian_motion}. However, the results contained in the present section hold even if we require $m \in L^p\big((a, b), \R_+ \big)$ with $p\in (1,\infty)$. In the current version of the manuscript, we restrict our attention to the case  $m \in L^2\big((a, b), \R_+ \big)$.};
		\item A curvature bound $\mathfrak{E}_\gamma \geq 0$;
	\end{enumerate}
	such that the following holds:
	For all $s, t \in (a, b)$, $T \in \R$, with $s \leq t \leq T \leq s+
     \varepsilon_\gamma$, we have that
	\begin{equation} \label{eq:AC2_definition}
		d_X \big( V_s^\gamma(T), V_t^\gamma(T) \big) \leq \int_s^t (T-r) \ m(r) \big(1 + \mathfrak{E}_\gamma  (T-r)^2\big) \, \diff r.
	\end{equation}
    
\end{definition}

\begin{remark} \label{rmk:uniqueness_velocity_1}
    Without further assumptions on the metric space $(X,d_X)$, it is not possible to guarantee the uniqueness (up to a negligible set of instants) of the family of geodesics $V^\gamma\colon (a, b)\to \Tang[\varepsilon_\gamma]$ involved in \cref{def:AC2_curve}. Later in this section, we provide the \EpsE-interceptor bound condition (see \cref{def:interceptor_bound}), which holds for spaces with curvature bounded from below (see \cref{lem:cbb_interceptor_bound}) and which implies uniqueness of $V^\gamma$ up to a null-measure set of times (see \cref{rmk:tangents_recoverable}).

    We insist on the fact that \cref{lem:second_derivative_bound,lem:ac2_to_ac1} do not rely in any way on the uniqueness of the function $V^\gamma\colon (a, b)\to \Tang[\varepsilon_\gamma]$. 
\end{remark}

\begin{figure}[H]
	\centering
	\input{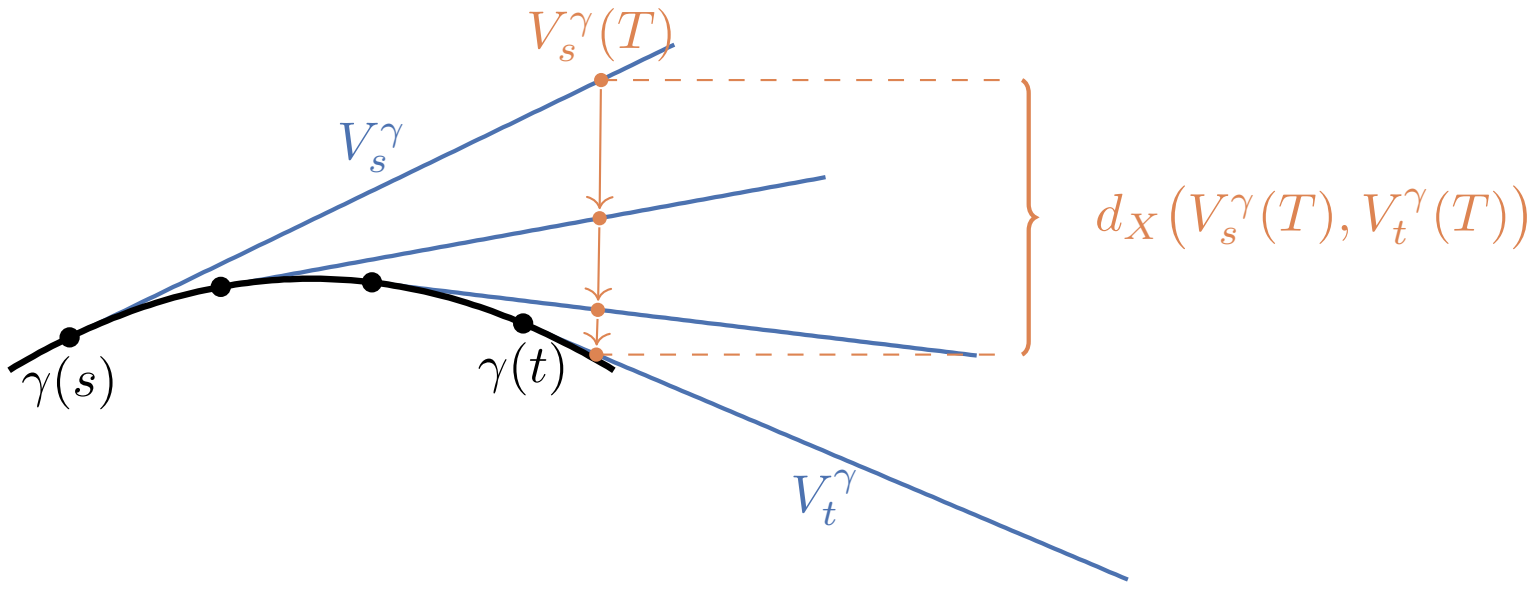}
	\caption{An example of an \ACt-curve $\gamma$ with its tangent geodesics $V_s^\gamma$ and $V_t^\gamma$. The distance $d_X\big(V_s^\gamma(T), V_t^\gamma(T) \big)$ is controlled by the integral on the right-hand side of \cref{eq:AC2_definition}.}
    \label{fig:AC2_interpr}
\end{figure}

\begin{remark} \label{rmk:First_Law}
    An interpretation of \cref{def:AC2_curve} is provided by Newton's first law of motion: \emph{``An object remains at rest or continues to move at a constant velocity (i.e., along a geodesic) unless acted upon by an external force''}.
    Hence, if we fix $s$ and $\Delta t > 0$ and set $t=T=s+\Delta t$ and if we understand $V_s^\gamma$ as the geodesic passing through $\gamma_s$ and `tangent' to $\gamma$, then the distance
    $$ d_X \big( V_s^\gamma(T), V_t^\gamma(T) \big) = d_X \big( V_s^\gamma(s+\Delta t), \gamma_{s + \Delta t} \big) \geq 0$$
    quantifies how the curve $\gamma$ \emph{accelerates}, i.e.,  how it deviates from being a geodesic in the interval $(s,t)$ (see \cref{fig:AC2_interpr}). 
    In this sense, the function $m \in L^2\big((a, b), \R_+\big)$ encodes an upper bound on the acceleration of the curve $\gamma$.
\end{remark}

\begin{remark} \label{rmk:long-short-realtions}
		For a complementary interpretation, we consider the limit $t \searrow s$. In this case, the bound in \cref{eq:AC2_definition} ensures both a \emph{short-term relation}---that $V_s^\gamma$ is in fact `tangent' to $\gamma$---and a \emph{long-term relation}---that the velocities of $V_s^\gamma$ are absolutely continuous.
        
        For the short-term relation, setting $T=t$, we see from \cref{eq:AC2_definition} that 
        $$ d_X \big( V_s^\gamma(t), V_t^\gamma(t) \big) = d_X \big( V_s^\gamma(t), \gamma_t \big) = O\big((t-s)^2\big) \qquad \mbox{as } t \searrow s$$
        for almost every $s \in (a, b)$.
        
        For the long-term relation, fixing $T > s$ while letting $t \searrow s$, for $T-s \gg t-s$, the evaluation of the geodesics at $T$ morally act as a magnifying glass on their directions. To better visualize this, let us argue in a Hilbert space, i.e., $(X,d_X)=(\mathcal{H}, |\cdot|)$, and let us assume that the curve $t\mapsto \gamma_t$ is smooth. 
        Owing to the short-term relation discussed above, we can write
        \begin{equation*}
            V_s^\gamma(T) = \gamma_s + \dot\gamma_s (T -s), \qquad V_t^\gamma(T) = \gamma_t + \dot\gamma_t (T -t),
        \end{equation*}
        so that the left-hand side of \cref{eq:AC2_definition} reads as
        \begin{equation} \label{eq:heuristic_long_term1}
        \begin{split}
            d_X \big( V_s^\gamma(T), V_t^\gamma(T) \big) &= \big| V_s^\gamma(T) -  V_t^\gamma(T)  \big| = \Big| \gamma_s -\gamma_t + (s-t)\dot\gamma_t + (T-s) (\dot\gamma_s - \dot\gamma_t) \Big| \\
            & = (T-s)\big| \dot\gamma_s - \dot\gamma_t \big| + O\big( (t-s)^2 \big)  \qquad \mbox{as } t \searrow s,
        \end{split}
        \end{equation}
        where we used the fact that $\gamma_s -\gamma_t + (s-t)\dot\gamma_t = O\big( (t-s)^2 \big)$ as $ t \searrow s$.
        Moreover, from the right-hand side of \cref{eq:AC2_definition}, we notice that
        \begin{equation} \label{eq:heuristic_long_term2}
        \begin{split}
            \int_s^t (T-r) \ m(r) \big(1 + &\mathfrak{E}_\gamma  (T-r)^2\big) \, \diff r\\
            &= 
            (T-s)\int_s^t m(r) \big(1 + \mathfrak{E}_\gamma  (T-r)^2\big) \, \diff r + o(t-s) \quad \mbox{as } t \searrow s.
        \end{split}
        \end{equation}
        By plugging \cref{eq:heuristic_long_term1,eq:heuristic_long_term2} into \cref{eq:AC2_definition} and neglecting the contributions of order $o(t-s)$ and $O\big( (t-s)^2 \big)$, we obtain that
        \begin{equation*}
            (T-s)\big| \dot\gamma_s - \dot\gamma_t \big| \leq (T-s)\int_s^t m(r) \big(1 + \mathfrak{E}_\gamma  (T-r)^2\big) \, \diff r,
        \end{equation*}
        which, up to the ``magnifying'' factor $T-s$ appearing on the two sides, yields a bound on the \emph{difference of the velocities} reminiscent of the variational definition of absolutely continuity.
        These informal observations are made perfectly rigorous in the Hilbertian setting in the proof of \cref{lem:euclidean_ac2_equivalence}.
\end{remark}

In the next result, we relate the function $m \in L^2 ((a, b), \R_+ )$ appearing at the right-hand side of \cref{eq:AC2_definition} with a point-wise notion of \emph{metric acceleration}.

\begin{lemma} \label{lem:second_derivative_bound}
	Let $\gamma\colon (a, b)\to X$ be an \ACt-curve. Then, we have that, for almost every $s \in (a, b)$,
	\begin{equation} \label{eq:metric_second_derivative_bound}
		|\gamma''|(s) \coloneqq \limsupTauTS\left(\frac{d_X \big( V^\gamma_s(T), V^\gamma_t(T) \big)}{(T-s)(t-s)}\right) \leq m(s).
	\end{equation}
	In particular, the function $|\gamma''|$---that we call \emph{metric acceleration}---is in $L^2((a, b), \R_+)$.
\end{lemma}

\begin{remark}
    Notice that the metric acceleration $|\gamma''|$ directly involves the function $V^\gamma\colon (a, b)\to \Tang[\varepsilon_\gamma]$.
    Therefore, if $V^\gamma$ is not uniquely determined (up to null-measure sets), we cannot expect a unique metric acceleration associated to the \ACt-curve $\gamma$.
\end{remark}

\begin{proof}[Proof of \cref{lem:second_derivative_bound}]
	To show \cref{eq:metric_second_derivative_bound}, by the Lebesgue theorem, it holds that $m(s) = \lim_{h \searrow 0} \frac{1}{h} \int_s^{s+h} m(u) \, \diff u$ for almost all $s \in (a, b)$.
    Then, for such $s$, the combination of \cref{eq:metric_second_derivative_bound,eq:AC2_definition} yields
    \begin{equation} \label{eq:Lebesgue_theorem_limsup}
        \begin{split}
            |\gamma''|(s) & \leq \limsupTauTS\frac{1}{(T-s) (t - s)}  \int_s^t (T-r) \ m(r)\big(1 + \mathfrak{E}_\gamma  (T-r)^2\big) \, \diff r\\
            &  \leq \limsupTauTS\frac{1}{(T-s) (t - s)}  \int_s^t (T-r) \ m(r) \, \diff r 
            \\ & \qquad \qquad  +  \limsupTauTS\frac{\mathfrak{E}_\gamma}{(T-s) (t - s)}  \int_s^t (T-r)^3 \ m(r) \, \diff r.
        \end{split}
    \end{equation}
    We first address the first term in the last inequality. As $r\geq s$, we have that
	\begin{equation} \label{eq:Lebesgue_theorem_2}
		\begin{split}
			\limsupTauTS\frac{1}{(T-s)(t - s)} & \int_s^t (T-r) \ m(r)\,\diff r  \leq  \limsupTauTS\frac{1}{t - s} \int_s^t  m(r)\,\diff r 
			 = m(s).
		\end{split}
	\end{equation}
	Finally, using \eqref{eq:Lebesgue_theorem_2}, we deduce that
	\begin{equation}\label{eq:Lebesgue_theorem_3}
		\limsupTauTS \frac{1}{(T-s) (t - s)} \int_s^t (T-r)^3 \ m(r) \, \diff r \leq \limsupTauTS \frac{T-s}{t - s} \int_s^t  m(r) \, \diff r = 0.
	\end{equation}
	Combining \cref{eq:Lebesgue_theorem_2,eq:Lebesgue_theorem_3} with \cref{eq:Lebesgue_theorem_limsup}, we obtain \cref{eq:metric_second_derivative_bound}. To see that $|\gamma''|$ is in $L^2((a, b))$, we are left to prove its measurability. 
    We preliminarily observe that, for any fixed $t$, the map $T\mapsto V^\gamma_t(T)$ is continuous. Furthermore, for any fixed $T$, the map $t\mapsto V^\gamma_t(T)$ is continuous as well by virtue of \cref{eq:AC2_definition}.
    Therefore, the $\limsup$ in \cref{eq:metric_second_derivative_bound} can be taken over rational $t$ and $T$. 
    Recalling the relations $s, t \in (a, b)$ with $s \leq t \leq T \leq s+
     \varepsilon_\gamma$, since for fixed $t$ and $T$ the argument inside the $\limsup$ is measurable as $s$ runs through its admissible range, we conclude---by setting the quotient to zero outside the admissible $s$-range---measurability in $s$ of the function $$s \mapsto \sup_{T \in \Q \cap (s, s +\varepsilon), t \in \Q \cap(s, T)}\left(\frac{d_X \big( V^\gamma_s(T), V^\gamma_t(T) \big)}{(T-s)(t-s)}\right)$$ for every $\varepsilon > 0$. 
     Taking the pointwise limit $\varepsilon \to 0$ shows measurability. 
\end{proof}

In the next result, we relate the notion of \ACt-curve to the classical one of \AC -curve taking values in a metric space.

\begin{lemma}\label{lem:ac2_to_ac1}
	Let $\gamma\colon (a, b)\to X$ be an \ACt-curve. Then, $\gamma$ is also an AC-curve. Furthermore, we have that $|\gamma'|(s) = |(V^\gamma_s)'|$ for almost every $s \in (a, b)$. Finally, we have that $s \mapsto |\gamma'|(s)$ has an absolutely continuous representative and that, for almost all $s \in (a, b)$,
	\begin{equation*}
		\left| \frac{\diff}{\diff s} |\gamma'|(s) \right| \leq |\gamma''|(s),
	\end{equation*}
    where the \emph{metric acceleration} $|\gamma''|$ is defined as in \cref{eq:metric_second_derivative_bound}.
\end{lemma}

\begin{remark} \label{rmk:acceleration_vs_der_speed}
    The inequality (in place of an equality) in the previous result is natural. Namely, it captures the fact that the acceleration, in general, is not merely due to the variation of the metric speed of a curve.
\end{remark}

Before we start with the proof, we introduce some notation that will be useful in the following: For an \ACt-curve $\gamma\colon (a, b)\to X$ with geodesic horizon $\varepsilon_\gamma > 0$, for every $s_1, s_2 \in (a, b)$ with $s_1 < s_2 <s_1+\frac{\varepsilon_\gamma}{2}$, we define the two functions
\begin{equation} \label{eq:def_G_cal_G}
	\begin{alignedat}{3}
		 & G^{s_1, s_2}\colon  [s_1, s_2]  \to \Tang[0, s_1+\varepsilon_\gamma - s_2], & \qquad \tau & \mapsto G^{s_1, s_2}_{\tau} \coloneqq \bigg[\sigma \mapsto V^\gamma_{\tau}(s_2+\sigma)\bigg], \\
		 & \mathcal{G}^{s_1, s_2}\colon  [0, s_1 + \varepsilon_\gamma - s_2]  \to X^{[s_1, s_2]}, & \qquad \sigma & \mapsto \mathcal{G}^{s_1, s_2}_{\sigma} \coloneqq \bigg[\tau \mapsto V^\gamma_{\tau}(s_2 + \sigma)\bigg],
	\end{alignedat}
\end{equation}
where $X^{[s_1, s_2]}$ denotes the set of all functions from $[s_1, s_2]$ to $X$.
\begin{figure}[H]
	\centering
	\input{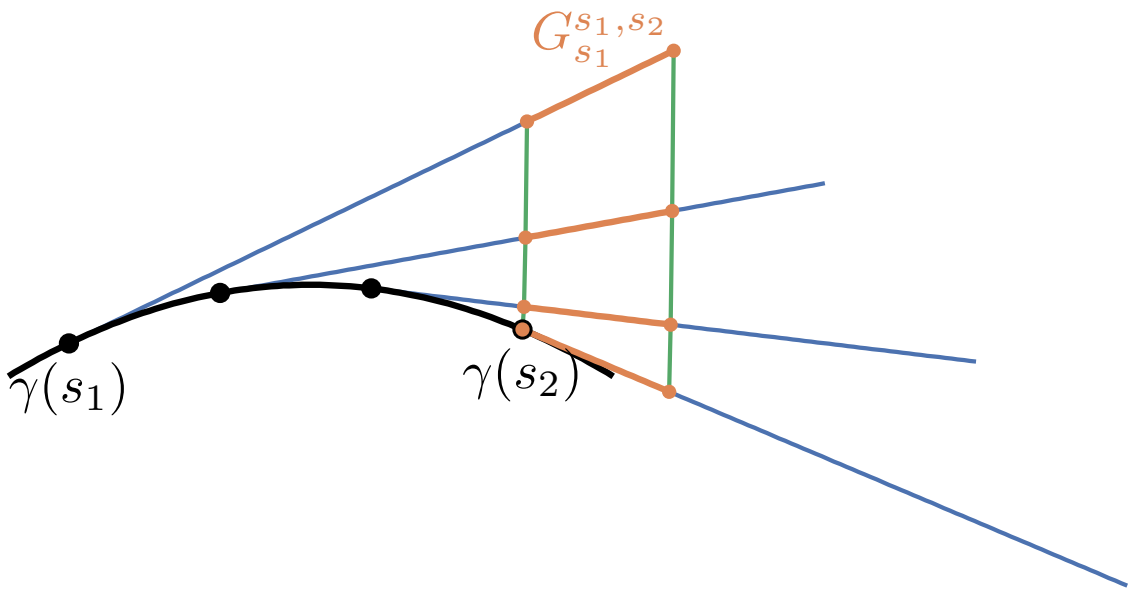}
	\caption{The construction of the functions $G^{s_1, s_2}_{\cdot}$. In black, we see the curve $\gamma$, with four points $\gamma_{s_1}, \gamma_{s'}, \gamma_{s''}$ and $\gamma_{s_2}$ (left to right) marked with a black dot. In blue, we see the tangent geodesics $V_{s_1}^\gamma$, $V_{s'}^\gamma$, $V_{s''}^\gamma$ and $V_{s_2}^\gamma$ (top to bottom). In orange, we see the geodesics $G^{s_1, s_2}_{s_1}, G^{s_1, s_2}_{s'}$, $G^{s_1, s_2}_{s''}$ and $G^{s_1, s_2}_{s_2}$ (top to bottom). Finally, in green, we see the paths $\mathcal{G}^{s_1, s_2}_0$ and $\mathcal{G}^{s_1, s_2}_{s_1 + \varepsilon_\gamma - s_2}$ (left to right).}
\end{figure}
Note that we have, for all allowed arguments $s_1 < s_2$, $\tau  \in [s_1, s_2]$ and $\sigma \in [0, s_1 + \varepsilon_\gamma - s_2]$, the two relations:
\begin{alignat}{2}
	G^{s_1, s_2}_{\tau}(\sigma) & = \mathcal{G}^{s_1, s_2}_{\sigma}(\tau); \label{eq:two_G_function_equivalence} \\
	|(G^{s_1, s_2}_{\tau})'| & = |(V^\gamma_{\tau})'|. \label{eq:velocity_g_equals_velocity_v}
\end{alignat}
Furthermore, using the triangle inequality (in fact, `quadrilateral inequality'), we get, for every $s_1< s_2$, every $\tau_1, \tau_2 \in [s_1,s_2]$ with $\tau_1<\tau_2$, and every  $\sigma_1, \sigma_2 \in [0, s_1 + \varepsilon_\gamma - s_2]$ with $\sigma_1 < \sigma_2$,
\begin{equation} \label{eq:estimate_change_of_pace_g}
	\begin{split}
		\Big| |(G^{s_1, s_2}_{\tau_1})'| - &|(G^{s_1, s_2}_{\tau_2})'| \Big|  = \frac{1}{\sigma_2 - \sigma_1}\Big| d_X\big(G^{s_1, s_2}_{\tau_1}(\sigma_1), G^{s_1, s_2}_{\tau_1}(\sigma_2)\big) - d_X\big(G^{s_1, s_2}_{\tau_2}(\sigma_1), G^{s_1, s_2}_{\tau_2}(\sigma_2)\big) \Big| \\
%		& \leq \frac{1}{\sigma_2 - \sigma_1} \Big| \Big(    d_X\big(G^{s_1, s_2}_{\tau_1}(\sigma_1), G^{s_1, s_2}_{\tau_2}(\sigma_1)\big) + d_X\big(G^{s_1, s_2}_{\tau_2}(\sigma_1), G^{s_1, s_2}_{\tau_1}(\sigma_2)\big) \Big) \\
%        & \qquad\qquad\qquad - \Big(  d_X\big(G^{s_1, s_2}_{\tau_2}(\sigma_1), G^{s_1, s_2}_{\tau_1}(\sigma_2)\big) -  d_X\big(G^{s_1, s_2}_{\tau_1}(\sigma_2), G^{s_1, s_2}_{\tau_2}(\sigma_2)\big) \Big) \Big|
        \\ & \leq \frac{1}{\sigma_2 - \sigma_1}  \Big(d_X\big(G^{s_1, s_2}_{\tau_1}(\sigma_1), G^{s_1, s_2}_{\tau_2}(\sigma_1)\big)+ d_X\big(G^{s_1, s_2}_{\tau_1}(\sigma_2), G^{s_1, s_2}_{\tau_2}(\sigma_2)\big)\Big)\\
        & \quad = \frac{1}{\sigma_2 - \sigma_1}  \Big(d_X\big(\mathcal G^{s_1, s_2}_{\sigma_1}(\tau_1), \mathcal G^{s_1, s_2}_{\sigma_1}(\tau_2)\big)+ d_X\big(\mathcal G^{s_1, s_2}_{\sigma_2}(\tau_1), \mathcal G^{s_1, s_2}_{\sigma_2}(\tau_2)\big)\Big).
	\end{split}
\end{equation}

\begin{proof}[Proof of \cref{lem:ac2_to_ac1}]
	We will prove the following four claims, which together imply the statement of the lemma:
	\begin{claimlist}
		\item\label{eq:ac2_to_ac1_1} The curve $(a, b) \ni s \mapsto |(V^\gamma_s)'|$ is an $\mathrm{AC}^2$-curve.
		\item\label{eq:ac2_to_ac1_2} $\Big|\frac{\diff}{\diff s} |(V^\gamma_s)'|\Big| \leq |\gamma''|(s)$ for almost every $s \in (a, b)$.
		\item\label{eq:ac2_to_ac1_3} The curve $\gamma$ is an $\mathrm{AC}^2$-curve.
		\item\label{eq:ac2_to_ac1_4} $|\gamma'|(s) = |(V^\gamma_s)'|$ for almost every $s \in (a, b)$.
	\end{claimlist}
	\textbf{Proof of \cref{eq:ac2_to_ac1_1}:}
	For any $s_1, s_2, \tau_1, \tau_2$ as in \cref{eq:estimate_change_of_pace_g} and $\sigma \in [0, s_1 + \varepsilon_\gamma - s_2]$, we can use \cref{eq:AC2_definition} with $s=\tau_1$, $t=\tau_2$ and $T = s_2 + \sigma$, to obtain
	\begin{equation*}
		\begin{split}
			d_X\big(G^{s_1, s_2}_{\tau_1}(\sigma), G^{s_1, s_2}_{\tau_2}(\sigma)\big) & \leq \int_{\tau_1}^{\tau_2} (s_2+\sigma -r) \ m(r) \big(1 + \mathfrak{E}_\gamma \ (s_2+\sigma - r)^2\big) \, \diff r \\
			& \leq \int_{\tau_1}^{\tau_2} \varepsilon_\gamma \ m(r) \ (1 + \mathfrak{E}_\gamma \ \varepsilon_\gamma^2) \, \diff r.
		\end{split}
	\end{equation*}
	By applying this estimate in \cref{eq:estimate_change_of_pace_g} with $\sigma_1 = 0$ and $\sigma_2 = s_1 + \varepsilon_\gamma - s_2$, we obtain that 
    \begin{equation*}
    \begin{split}
        \Big| |(G^{s_1, s_2}_{\tau_1})'| - |(G^{s_1, s_2}_{\tau_2})'| \Big|  
        & \leq \frac1{s_1 + \varepsilon_\gamma - s_2} \Big(d_X\big(G^{s_1, s_2}_{\tau_1}(0), G^{s_1, s_2}_{\tau_2}(0)\big)\\
        & \qquad \qquad\qquad\qquad + d_X\big(G^{s_1, s_2}_{\tau_1}(s_1 + \varepsilon_\gamma - s_2), G^{s_1, s_2}_{\tau_2}(s_1 + \varepsilon_\gamma - s_2)\big)\Big)\\
        & \leq \frac{2 \varepsilon_\gamma (1 + \mathfrak{E}_\gamma \ \varepsilon_\gamma^2) }{s_1 + \varepsilon_\gamma - s_2}  \int_{\tau_1}^{\tau_2}  \ m(r)  \, \diff r,
    \end{split}
    \end{equation*}
    which implies that
    $\tau \mapsto |(G^{s_1, s_2}_{\tau})'|$ is $\mathrm{AC}^2$, with
	\begin{equation*}
		\Big|\frac{\diff}{\diff \tau}|(G^{s_1, s_2}_{\tau})'|\Big| \leq \frac{2 \varepsilon_\gamma (1 + \mathfrak{E}_\gamma \ \varepsilon_\gamma^2) }{s_1 + \varepsilon_\gamma - s_2}  \ m(\tau) \leq 4  (1 + \mathfrak{E}_\gamma \ \varepsilon_\gamma^2) \ m(\tau) ,
	\end{equation*} 
	where we used that $s_2 - s_1 \leq \frac{\varepsilon_\gamma}{2}$.
	Thus, using \cref{eq:velocity_g_equals_velocity_v}, we see that also $s \mapsto |(V^\gamma_s)'|$ is $\mathrm{AC}^2$, where $\Big|\frac{\diff}{\diff s}|(V^\gamma_s)'|\Big|\leq 4  (1 + \mathfrak{E}_\gamma \ \varepsilon_\gamma^2) \ m(s)$. \\
	\textbf{Proof of \cref{eq:ac2_to_ac1_2}:}
	We will, again, rely on \cref{eq:estimate_change_of_pace_g}: For any curve $\eta$ taking values in $(X,d_X)$, we denote the right-sided metric derivative by
	\begin{equation*} %\label{eq:def_right_sided_metric_derivative}
		| \eta_+ |(t) \coloneqq \limsup_{h \searrow 0} \frac{d_X(\eta(t+h), \eta(t))}{h}.
	\end{equation*}
	 From \cref{eq:estimate_change_of_pace_g}, it follows that, for any $s_1,  s , s_2 \in (a, b)$ with $s_1 <  s < s_2$ and $\sigma_1 , \sigma_2 \in [0, s_1+ \varepsilon_\gamma - s_2)$ with $\sigma_1<\sigma_2$, we have that
	\begin{equation} \label{eq:estimate_change_of_pace_g_final}
					\Big| \big[\underline{\tau} \mapsto |(G^{s_1, s_2}_{\underline{\tau}})'|\big]_+ \Big|(s) \leq \frac{1}{\sigma_2 - \sigma_1} \Big(|(\mathcal{G}^{s_1, s_2}_{\sigma_1})_+|(s) + |(\mathcal{G}^{s_1, s_2}_{\sigma_2})_+|(s)\Big),
	\end{equation}
    where we used $\tau_1= s $ and we let $\tau_2 \searrow s$ in \cref{eq:estimate_change_of_pace_g}.
	We now estimate the right-hand side of the last inequality.
	By the definition of $|\gamma''|$ (see \cref{eq:metric_second_derivative_bound}), for any $s \in (s_1, s_2)$ and any $\alpha>0$, we may find some $\underline{T}^\star \in (s, s_2]$ such that for all $\underline{T} \in [s, \underline{T}^\star)$ and for all $\underline{t} \in (s, \underline{T})$, it holds that
	\begin{equation}\label{eq:use_def_second_metric_derivative}
		d_X\big(V^\gamma_s(\underline{T}), V^\gamma_{\underline{t}}(\underline{T})\big) \leq \big(|\gamma''|(s) + \alpha\big) \, (\underline{T}-s)\, (\underline{t}-s).
	\end{equation}
	Thus, recalling \cref{eq:def_G_cal_G}, we obtain the two estimates
	\begin{equation} \label{eq:estimate_metric_derivative_g}
		\begin{alignedat}{2}
			|(\mathcal{G}^{s_1, \underline{T}}_{0})_+|(s) & = \limsup_{h \searrow 0} \frac{d_X\big(V^\gamma_s(\underline{T}), V^\gamma_{s+h}(\underline{T})\big)}{h} & & \leq \big(|\gamma''|(s) + \alpha\big) \, (\underline{T} - s), \\
			|(\mathcal{G}^{s_1, \underline{T}}_{\underline{T}^\star - \underline{T}})_+|(s) & = \limsup_{h \searrow 0} \frac{d_X\big(V^\gamma_s(\underline{T}^\star), V^\gamma_{s+h}(\underline{T}^\star)\big)}{h} & & \leq \big(|\gamma''|(s) + \alpha\big) \, (\underline{T}^\star - s).
		\end{alignedat}
	\end{equation}
	Combining \cref{eq:velocity_g_equals_velocity_v,eq:estimate_metric_derivative_g,eq:estimate_change_of_pace_g_final} and setting $s_2=\underline{T}$, $\sigma_1=0$ and $\sigma_2 = \underline{T}^\star - \underline{T}$ (recall that $s<\underline{T}<\underline{T}^\star$), we obtain
	\begin{equation} \label{eq:right_sided_derivative_V_estimate}
    \begin{split}
		\big| |(V^\gamma_s)'|_+ \big| = \Big| \big[\tau \mapsto |(G^{s_1, \underline{T}}_{\tau})'|\big]_+ \big|(s) &
        \leq \frac{|(\mathcal{G}^{s_1, \underline{T}}_{0})_+|(s) + |(\mathcal{G}^{s_1, \underline{T}}_{\underline{T}^\star - \underline{T}})_+|(s)}{\underline{T}^\star - \underline{T}} \\
        &\leq \frac{(\underline{T}^\star - s) + (\underline{T} - s)}{\underline{T}^\star - \underline{T}} \, \big(|\gamma''|(s) + \alpha\big);
    \end{split}
	\end{equation}
	letting first $\underline{T} \searrow s$ and then $\alpha \searrow 0$, we obtain that $\big| |(V^\gamma_s)'|_+ \big| \leq |\gamma''|(s)$. By \cref{eq:ac2_to_ac1_1}, we know that the right-sided derivative coincides with the derivative almost everywhere, which concludes the proof of \cref{eq:ac2_to_ac1_2}.\\
	\textbf{Proof of \cref{eq:ac2_to_ac1_3}:}
	By \cref{eq:ac2_to_ac1_1} and boundedness of $(a, b)$, the curve $(a, b) \ni s \mapsto |(V^\gamma_s)'|$ is in particular an AC-curve, and so we can find a $V_\text{max}$ such that for all $s \in (a, b)$, $|(V^\gamma_s)'| \leq V_\text{max}$. By the definition of an \ACt-curve, we have, for all $s, t \in (a, b)$ such that $s < t \leq s + \varepsilon_\gamma$,
	\begin{equation}\label{eq:ac2_to_ac1_3_intermediate}
		\begin{split}
			d_X(\gamma_s, \gamma_t) & \leq d_X \big(\gamma_s, V_s^\gamma(t)\big) + d_X\big( V_s^\gamma(t), \gamma_t\big) \\
			& \leq V_\text{max} \, (t - s) + \int_s^t (t - r) \ m(r) \big(1 + \mathfrak{E}_\gamma \ (t - r)^2 \big) \, \diff r \\
			& \leq \int_s^t \Big(V_\text{max} + (b-a) \ m(r) \big(1 + \mathfrak{E}_\gamma \ (b-a)^2 \big)\Big) \, \diff r.
		\end{split}
	\end{equation}
	By partitioning the domain $(a, b)$ and using the triangle inequality, we see that \cref{eq:ac2_to_ac1_3_intermediate} extends for generic $s, t \in (a, b)$ such that $s < t$, which shows that $\gamma$ is an $\mathrm{AC}^2$-curve.\\
	\textbf{Proof of \cref{eq:ac2_to_ac1_4}:}
	Let us first show that we have that $|\gamma_+|(s) = |(V^\gamma_s)'|$ for every $s \in (a, b)$: To this end, by a simple application of the triangle inequality, we obtain
	\begin{equation}\label{eq:difference_one_sided_derivative_derivative_v}
		\big||\gamma_+|(s) - |(V^\gamma_s)'|\big| \leq \limsup_{h \searrow 0} \frac{d_X\big(\gamma_{s+h}, V_s^\gamma(s+h)\big)}{h} =0,
	\end{equation}
	where the $\limsup$ vanishes owing to the inequality \eqref{eq:AC2_definition} that characterizes \ACt-curves. To conclude the claim, note that by \cref{eq:ac2_to_ac1_1}, the function $s \mapsto |(V^\gamma_s)'|$ is continuous. 
    To this end, we use a standard argument\footnote{We provide a slight adaptation; the standard argument shows that a real-valued function with continuous right-derivative is differentiable everywhere (see~\cite{continuous_right_derivative_implies_differentiability}).} to show that, for any AC curve $\eta\colon [a, b]\to X$ such that the function
	\begin{equation*}
		s \mapsto |\eta_+|(s) \coloneqq \limsup_{h \searrow 0} \frac{d_X(\eta_{s+h}, \eta_s)}{h}
	\end{equation*}
	is continuous on $[a,b)$, it holds that $|\eta'|(s) = |\eta_+|(s)$ for almost all $s \in (a, b)$. To prove this, we first show that if $|\eta_+|(s) \leq M$ for all $s \in [a, b)$, then $\frac1{s-a} d_X(\eta_a, \eta_s) \leq M$ for all $s \in (a, b]$. 
    To see this, we define the function
	\begin{equation*}
		L(a) \coloneqq  |\eta_+|(a), \qquad \text{and} \qquad L(s) \coloneqq  \frac{d_X(\eta_a, \eta_s)}{s - a}, \quad \forall s \in (a, b],
	\end{equation*}
	and the set
	\begin{equation*}
		E_\delta \coloneqq \{s \in (a, b] \mid \forall y \in (a, s]: L(y) \leq M + \delta\},
	\end{equation*}
	for some $\delta > 0$. Note that $E_\delta$ is non-empty and, by continuity of $L$ on $(a, b]$, it has the form $(a, c_\delta]$ for some $c_\delta \in (a, b]$. However, if $c_\delta < b$, then by the definition of $|\eta_+|$, there exists some $\varepsilon > 0$ such that for all $h \in (0, \varepsilon)$, it holds that $\frac1h d_X(\eta_{c_\delta}, \eta_{c_\delta+h}) \leq M + \delta$, which, with an application of the triangle inequality to the three points $\eta_a$, $\eta_{c_\delta}$ and $\eta_{c_\delta + h}$, contradicts the maximality of $c_\delta$.
	Thus, we have that $c_\delta = b$ for all $\delta > 0$, which shows that $\frac1{r-a} d_X(\eta_a, \eta_r) \leq M$ for all $r \in (a, b]$. \\
	Now let us show that $|\eta'|(s) = |\eta_+|(s)$ for almost all $s \in (a, b)$. To this end, we fix some $s \in (a, b)$ and we set $M_h \coloneqq \max_{r \in [s-h, s]} |\eta_+|(r)$ for $h > 0$ sufficiently small. 
    By the previous argument, we have that for all such $h > 0$ and all $r \in (s-h, s)$, it holds that $\frac1{s - r} {d_X(\eta_r, \eta_s)} \leq M_h$. We thus have that
	\begin{equation*}
		\limsup_{r \nearrow s} \frac{d_X(\eta_r, \eta_s)}{s - r} \leq \lim_{h \searrow 0} M_h = |\eta_+|(s),
	\end{equation*}
	which proves that $|\eta_+|(s) = \limsup_{r \to s} \frac{d_X(\eta_r, \eta_s)}{|r - s|} = |\eta'|(s)$ for almost all $s \in (a, b)$.
\end{proof}
Let us show that our definition of \ACt-curves is consistent with the usual definition of twice absolutely continuous curves in Euclidean and Hilbert spaces.

\begin{lemma} \label{lem:euclidean_ac2_equivalence}
	Let $(X, d_X) = (\mathcal{H}, |\cdot|)$ be a Hilbert space and let $\gamma\colon (a, b)\to \mathcal{H}$, for $a, b \in \R$ and $a < b$, be a curve. Then, the following two statements are equivalent:
	\begin{enumerate}
		\item\label{item:our_def_twice_ac} $\gamma$ is an \ACt-curve as per \cref{def:AC2_curve}.
		\item\label{item:usual_def_twice_ac} $\gamma \in C^1((a, b), \mathcal{H})$, and the velocity $s \mapsto \dot{\gamma}(s)$ is an \AC-curve with metric speed in $L^2$.
	\end{enumerate}
	Furthermore, in this case, we have that $|\gamma ''|(r) = |\ddot \gamma_r|$ for almost every $r \in (a, b)$.
\end{lemma}
\begin{proof}
	\textbf{\ref{item:our_def_twice_ac} $\implies$~\ref{item:usual_def_twice_ac}:}
	The proof works analogously to the proof of \cref{lem:ac2_to_ac1}: In $\mathcal{H}$, we may replace \cref{eq:estimate_change_of_pace_g} with the following stronger estimate, for all allowed arguments $s_1, s_2$, $\tau_1, \tau_2$ and any $\sigma_1, \sigma_2 \in [0, s_1 + \varepsilon_\gamma - s_2]$ with $\sigma_1 < \sigma_2$:
	\begin{equation*}
		\Big| \dot{G}^{s_1, s_2}_{\tau_1}(0) - \dot{G}^{s_1, s_2}_{\tau_2}(0) \Big| \leq \frac{1}{\sigma_2 - \sigma_1}  \Big(d_X\big(G^{s_1, s_2}_{\tau_1}(\sigma_1), G^{s_1, s_2}_{\tau_2}(\sigma_1)\big) + d_X\big(G^{s_1, s_2}_{\tau_1}(\sigma_2), G^{s_1, s_2}_{\tau_2}(\sigma_2)\big)\Big).
	\end{equation*}
	We note that in $\mathcal{H}$, $G^{s_1, s_2}_{\tau_1}$ and $V_s^\gamma$ are affine functions whose derivatives we denote by $\dot{G}^{s_1, s_2}_{\tau_1}$ and $\dot{V}^\gamma_s$, respectively.
	In place of \cref{eq:ac2_to_ac1_1}, this allows us to conclude that
	$\tau \mapsto \dot{G}^{s_1, s_2}_{\tau}$ is absolutely continuous and thus so is $s \mapsto \dot{V}^\gamma_s$. In place of \cref{eq:ac2_to_ac1_2}, we can apply the same modification to \cref{eq:right_sided_derivative_V_estimate} to conclude that $|\frac{\diff}{\diff s} \dot{V}^\gamma_s| \leq |\gamma''|(s)$ for almost every $s \in (a, b)$. 
    By \cref{lem:ac2_to_ac1}, $\gamma$ is absolutely continuous.
    Moreover, the fact that $\dot{\gamma}_s = \dot{V}^\gamma_s$ for almost every $s$ follows from the definition of an \ACt-curve, together with an estimate similar to \cref{eq:difference_one_sided_derivative_derivative_v}, which yields $|\dot\gamma_s - \dot{V}^\gamma_s| \leq \limsup_{h \searrow 0} \frac{|\gamma_{s+h} - V_s^\gamma(s+h)|}{h} = 0$.
    Finally, by continuity of $s \mapsto \dot{V}^\gamma_s$, we have $\gamma \in C^1$. 
    Thus, we see \ref{item:our_def_twice_ac} implies \ref{item:usual_def_twice_ac} with $| \gamma''|(r) \geq |\ddot \gamma_r|$ for almost every $r \in (a, b)$.\\
	\textbf{\ref{item:usual_def_twice_ac} $\implies$~\ref{item:our_def_twice_ac}:}
	Let $\gamma$ be as in~\ref{item:usual_def_twice_ac}. 
    Then, setting $V_t^\gamma(T) \coloneqq \gamma_t + (T-t) \ \dot{\gamma}_t$, we have, for any $s, t \in (a, b)$, $T \in \R$, with $s < t < T$,
	that the function $r \mapsto \gamma_r + (T - r) \ \dot{\gamma}_r$ is absolutely continuous on $[s, t]$ with derivative $r \mapsto (T-r) \ \ddot{\gamma}_r$ almost everywhere. Thus, we have that
	\begin{equation*}
		\begin{split}
			d_X \big(V_s^\gamma(T), V_t^\gamma(T) \big) & =|\gamma_s + (T-s) \ \dot{\gamma}_s - \gamma_t - (T-t) \ \dot{\gamma}_t| \\
			& =\left| \int_s^t (T-r)  \ \ddot{\gamma}_r \, \diff r \right| \leq \int_s^t (T-r) \ |\ddot{\gamma}_r| \, \diff r,
		\end{split}
	\end{equation*}
	from which, by virtue of \cref{lem:second_derivative_bound}, we conclude that \ref{item:usual_def_twice_ac} implies \ref{item:our_def_twice_ac} with $| \gamma''|(r) \leq |\ddot \gamma_r|$ for almost every $r \in (a, b)$.
\end{proof}

%Furthermore, to be able to prove that certain curves fulfill the definition, it is useful to assume that these geodesics are well-behaved; this last requirement is captured in the following definition.
In the proof of \cref{lem:euclidean_ac2_equivalence} (more precisely, of the implication \textbf{\ref{item:our_def_twice_ac} $\implies$~\ref{item:usual_def_twice_ac}}) we show that, for an \ACt-curve $\gamma$ in a Hilbert space, $\dot \gamma_s = \dot V^\gamma_s$ for almost every $s\in (a,b)$. 
We provide below a geometric condition that guarantees that the function $s\mapsto V^\gamma_s$ is (up to a negligible set of instants) uniquely determined by the \ACt-curve $\gamma$.
Moreover, this condition plays a pivotal role in \cref{sec:curves_of_newtonian_motion,sec:discrete_scheme}.

\begin{definition}[\EpsE-Interceptor Bound] \label{def:interceptor_bound}
	We say that a metric space $X$ fulfills the \emph{\EpsE-interceptor bound} if, for every two geodesics $\xi, \zeta\colon [0, T]\to X$ such that $\xi(0) = \zeta(0)$ and $\max\{|\xi'|, |\zeta'|\} \leq 1$, the following holds: For every $s, t \in (0, T]$ such that $0 < s \leq t \leq\epsilon$, we have that
	\begin{equation*}
		d_X(\xi_t, \zeta_t) \leq \frac{d_X(\xi_s, \zeta_s)}{s} \, ( t + t^3 \,  \mathfrak{E}).
	\end{equation*}
\end{definition}
Informally, the \EpsE-interceptor bound ensures that geodesics do not diverge too quickly from each other.

\begin{remark} \label{rmk:interceptor_Hilbert}
    We notice that a Hilbert space $(\mathcal{H}, |\cdot|)$ satisfies \cref{def:interceptor_bound} with $\mathfrak{E}=0$ and $\epsilon=T$.
\end{remark}

\begin{remark} \label{rem:interceptor_bound_scaling}
	The upper bound $\max\{|\xi'|, |\zeta'|\} \leq 1$ on the velocities of the geodesics in \cref{def:interceptor_bound} is arbitrarily chosen, in the following sense: For any $V > 0$, the space $X$ fulfills the \EpsE-interceptor bound for $\epsilon, \mathfrak{E} > 0$ if and only if, for every two geodesics $\xi, \zeta\colon [0, T]\to X$ such that $\xi(0) = \zeta(0)$ and $\max\{|\xi'|, |\zeta'|\} \leq V$, the following holds: For every $s, t \in (0, T]$ such that $0 < s \leq t \leq\frac{\epsilon}{V}$, we have that
	\begin{equation*}
		d_X(\xi_t, \zeta_t) \leq \frac{d_X(\xi_s, \zeta_s)}{s} \, ( t + t^3 \, V^2 \, \mathfrak{E}).
	\end{equation*}
\end{remark}
Perhaps unsurprisingly, the \EpsE-interceptor bound is fulfilled by spaces with curvature bounded from below (CBB).
We postpone the proof of this result to \cref{sec:proof-of-CBB}.
\begin{lemma} \label{lem:cbb_interceptor_bound}
	Let $X$ be a complete length $CBB(\kappa)$-space for some $\kappa \in \mathbb{R}$. Then $X$ fulfills the \emph{\EpsE-interceptor bound} for some $\epsilon, \mathfrak{E} > 0$ depending only on $\kappa$.
\end{lemma}

We are now in a position to establish the (essential) uniqueness of the function $V^\gamma\colon (a, b)\to \Tang[\varepsilon_\gamma]$ associated to an \ACt-curve $\gamma$.

\begin{proposition} \label{rmk:tangents_recoverable}
	Let us assume that $(X,d_X)$ fulfills the \EpsE-interceptor bound (see \cref{def:interceptor_bound}) and that $\gamma\colon (a, b) \to X$ is an \ACt-curve as per \cref{def:AC2_curve}.
    Then, the function $V^\gamma\colon (a, b)\to \Tang[\varepsilon_\gamma]$ is uniquely determined up to a negligible set of instants.
\end{proposition}
    \begin{proof}
    To see this, let us assume for a moment that \cref{def:AC2_curve} is fulfilled for the same curve $\gamma$ with two different functions $V^\gamma$ and $W^\gamma$, with their respective acceleration bounds $m_1, m_2 \in L^2((a, b), \R_+)$ and geodesic horizons $\varepsilon_\gamma^1,\varepsilon_\gamma^2>0$.
    Then, after setting $\varepsilon_\gamma\coloneqq \min \{ \varepsilon_\gamma^1,\varepsilon_\gamma^2 \} >0$, for every $s,t\in (a,b)$ with $s\leq t\leq s+\varepsilon_\gamma$, recalling that by definition $ W^\gamma_t(t)  = V^\gamma_t(t) = \gamma_t$, by the triangle inequality we compute
    \begin{equation*}
        \begin{split}
            \frac{d_X \big( V^\gamma_s(t), W^\gamma_s(t) \big)}{(t-s)^2} 
            \leq \frac{d_X \big( V^\gamma_s(t), \gamma_t \big)}{(t-s)^2} + \frac{d_X \big( W^\gamma_s(t), \gamma_t \big)}{(t-s)^2}.
        \end{split}
    \end{equation*}
    Then, it follows from \cref{lem:second_derivative_bound} that 
    \begin{equation} \label{eq:uniqueness_tangents_1}
        \limsup_{t \searrow s} \frac{d_X \big( V^\gamma_s(t), W^\gamma_s(t) \big)}{(t-s)^2} \leq m_1(s) + m_2(s) < \infty
    \end{equation}
    for almost every $s \in (a, b)$. 
    Let $\epsilon, \mathfrak{E}$ be the constants provided by the \EpsE-interceptor bound, and let us set $t\mapsto \xi(t)\coloneqq W^\gamma_s(t)$ and $t\mapsto \zeta(t)\coloneqq V^\gamma_s(t)$.
    If we define $\epsilon' \coloneqq \min \left\{  \varepsilon_\gamma, \epsilon/| \xi'|, \epsilon/| \zeta'|  \right\}$ (see also \cref{rem:interceptor_bound_scaling}) and we fix $T\in (s, s+\epsilon')$, we can apply the \EpsE-interceptor bound to $\xi$ and $\zeta$ (recall once again that $ \xi_s = W^\gamma_s(s)  = \gamma_s = V^\gamma_s(s) = \zeta_s$) to get
    \begin{equation} \label{eq:uniqueness_tangents_2}
    d_X ( \xi_T, \zeta_T ) =
        d_X \big( V^\gamma_s(T), W^\gamma_s(T) \big) \leq \frac{d_X \big( V^\gamma_s(t), W^\gamma_s(t) \big)}{t-s} \Big( (T-s) + (T-s)^3 \, \max\{ |\xi'|^2, |\zeta'|^2 \} \, \mathfrak{E} \Big)
    \end{equation}
    for almost every $s\in (a,b)$ and for every $t\in (s,T]$.
    Since in \cref{eq:uniqueness_tangents_2} we can choose $t$ arbitrarily close to $s$, by virtue of \cref{eq:uniqueness_tangents_1} we deduce that $d_X \big( V^\gamma_s(T), W^\gamma_s(T) \big) = 0$ for every $T\in (s,s+\epsilon')$, i.e., that $V^\gamma_s,W^\gamma_s \in \Tang[s, \varepsilon_\gamma]$ coincide for almost every $s\in (a,b)$.
    %For this reason, as our main interest lies in spaces which fulfill the \EpsE-interceptor bound, we will treat $V^\gamma$ as recoverable from $\gamma$.
\end{proof}

To work with \ACt-curves, we will often require the space $X$ to fulfill the following assumption:
\begin{assumption} \label{ass:interceptor}
	The metric space $X$ fulfills the $\epsilon$-$\mathfrak{E}$-interceptor bound for some $\epsilon, \mathfrak{E} > 0$ (see \cref{def:interceptor_bound}). Furthermore, any geodesic $\xi\colon [0, T] \to X$ with $|\xi'| \leq 1$ and $T < \epsilon$ can be extended to a geodesic $\zeta\colon [0, \epsilon] \to X$ such that $\zeta|_{[0, T]} = \xi$. 
    Lastly, $X$ is a geodesic space, i.e., every two points can be connected by a geodesic.
\end{assumption}

The next proposition allows us to better specify the constants $\mathfrak{E}_\gamma, \varepsilon_\gamma$ appearing in the definition of \ACt-curves under further geometrical assumptions on $(X,d_X)$.
Moreover, under the proposition hypotheses, we can set $m=|\gamma''|$ (see  \cref{eq:metric_second_derivative_bound} for the definition of the metric acceleration $|\gamma''|$) in the integral at the right-hand side of the inequality \eqref{eq:AC2_definition}.

\begin{proposition} \label{prop:second_derivative_is_minimal}
	Let us assume that $X$ fulfills \cref{ass:interceptor}, that $X$ is $\text{CAT}(\kappa)$ for some $\kappa \in \R$, and furthermore that it fulfills the following condition: There exists some $C_X > 0, l_X > 0$ such that for every $l \in (0, l_X]$ and for every two geodesics $\xi, \zeta\colon [-l, l]\to X$ such that $\xi_0 = \zeta_0$ and such that $|\xi'| \leq 1$ and $|\zeta'| \leq 1$, it holds that $d_X(\xi_{-l}, \zeta_{-l}) \leq C_X d_X (\xi_l, \zeta_l)$. \\
    Then, for any \ACt-curve $\gamma$, the function $|\gamma''|$
	is a minimal choice of the function $m$ in the definition of \ACt-curves (see \cref{eq:AC2_definition}). 
    More precisely, setting $V_{\max} \coloneqq \sup_{t \in (a, b)} |(V_t^\gamma)'|$, we have that $\gamma$ fulfills the \cref{def:AC2_curve} of an \ACt-curve with the choice $m = |\gamma''|$, $\mathfrak{E}_\gamma = V_{\max}^2 \, \mathfrak{E}$ and some geodesic horizon $\varepsilon_\gamma \in (0, \frac{\epsilon}{2 V_{\max} + 1}]$.
\end{proposition}
\begin{proof}
	Since $\gamma$ is an \ACt-curve, there exist $\mathfrak{E}^{\mathrm{old}}_\gamma, \varepsilon_\gamma^\mathrm{old} > 0$ such that 
    $$d_X \big(V_s^\gamma(T), V_t^\gamma(T)\big) \leq \int_s^t (T-r)\ m(r)\big(1+\mathfrak{E}^\mathrm{old}_\gamma(T-r)^2\big)\,\diff r$$
    for all $s, t \in (a, b)$, $T \in \R$, such that $s \leq t \leq T \leq s + \varepsilon_\gamma^\mathrm{old}$.
	We are going to show that, under the assumptions of the proposition, for all $s, t \in (a, b)$, $T \in \R$, with $s < t < T$ and such that $T - s \leq \min\{\frac{\epsilon}{V_{\max} + 1}, \varepsilon_\gamma^\mathrm{old}\}$, it holds that
	\begin{equation}\label{eq:goal_AC2_proof}
		d_X \big( V_s^\gamma(T), V_t^\gamma(T) \big) \leq \int_s^t (T - r) \ |\gamma''|(r) \big(1 + V_{\max}^2 \ \mathfrak{E} \ (T - r)^2 \big) \,\diff r.
	\end{equation}
	To this end, let us first show that for any such $s, T$, we have that
	\begin{equation} \label{eq:intermediate_AC2_proof}
		\limsup_{h \searrow 0} \frac{d_X \big( V_s^\gamma(T), V_{s+h}^\gamma(T) \big)}{h} \leq (T - s) \ |\gamma''|(s) \ \big( 1 + V_{\max}^2 \ \mathfrak{E} \ (T - s)^2 \big).
	\end{equation}
	We fix such $s, T$ and we pick an arbitrary $\alpha > 0$. 
    By the definition of $|\gamma''|(s)$ in \cref{eq:metric_second_derivative_bound}, there exists some $\underline{T}^\star \in (s, s + \varepsilon_\gamma^{\mathrm{old}}]$ such that for all $\underline{T} \in (s, \underline{T}^\star]$, $\underline{t} \in (s, \underline{T}] \cap (a, b)$, we have that
	\begin{equation}\label{eq:use_def_second_metric_derivative_AC2_proof}
		\frac{d_X\big(V^\gamma_s(\underline{T}), V^\gamma_{\underline{t}}(\underline{T})\big)}{(\underline{T} - s)(\underline{t}-s)} \leq |\gamma''|(s) + \alpha;
	\end{equation}
	we can further assume that $\frac{l_X}{\underline{T}^\star - s} \geq V_{\max} \coloneqq \sup_{t \in (a, b)} |(V_t^\gamma)'|$.
	From \cref{eq:use_def_second_metric_derivative_AC2_proof}, we derive for such $\underline{T}$ and $\underline{t}$, the two inequalities
	\begin{align}
		d_X\big(V^\gamma_s(\underline{t}), \gamma_{\underline{t}}\big) & \leq (|\gamma''|(s) + \alpha) \ (\underline{t} - s)^2 = O((\underline{t} - s)^2) \label{eq:intermediate_AC_proof_11} \\
		d_X\big(V^\gamma_s(\underline{T}), V^\gamma_{\underline{t}}(\underline{T})\big) & \leq (|\gamma''|(s) + \alpha) \ (\underline{T} - s) \ (\underline{t} - s) = O((\underline{t} - s)). \label{eq:intermediate_AC_proof_12}
	\end{align}
	We fix such a $\underline{T}$ such that $\underline{T} \leq \frac{T + s}{2}$.
    Then, we have, for every $\underline{t} \in (s, \underline{T}] \cap (a, b)$, that 
    \begin{equation*}
        d_X(\gamma_s, V_{\underline{t}}^\gamma(\underline{T})) \leq d_X(\gamma_s, \gamma_{\underline{t}}) + d_X(\gamma_{\underline{t}}, V_{\underline{t}}^\gamma(\underline{T})) \leq (\underline{T} - s) \ V_{\max}.    
    \end{equation*}
    Since $(\underline{T} - s) V_{\max} \leq (T - s) V_{\max} \leq \epsilon$, we can apply \cref{ass:interceptor} to choose, for every such $\underline{t}$, a geodesic $\xi^{\underline{t}} \colon [s, \underline{T}] \to X$ such that $\xi^{\underline{t}}_s = \gamma_s$ and $\xi^{\underline{t}}_{\underline{T}} = V_{\underline{t}}^\gamma(\underline{T})$. Furthermore, we can deduce that $|(\xi^{\underline{t}})'| \leq V_{\max}$. 
    Thus, we can use the extensibility assumption\footnote{From \cref{ass:interceptor}, one obtains immediately that geodesics with any velocity can be extended up to \emph{length} $\epsilon$.} in \cref{ass:interceptor} to extend $\xi$ (without renaming) to the domain $[s, T]$ by using that $2 (\underline{T} - s) V_{\max} \leq (T-s) V_{\max} \leq \epsilon$. By the inverse \EpsE-interceptor bound stemming from the $\text{CAT}(\kappa)$ assumption\footnote{By ``inverse \EpsE-interceptor bound'', we mean the following: \emph{There exist $\epsilon, \mathfrak{E} > 0$ depending only on $\kappa$ such that, for every two geodesics $\xi, \zeta\colon [0, T]\to X$ with $\xi(0) = \zeta(0)$ and $\max\{|\xi'|, |\zeta'|\} \leq 1$, the following holds: For every $s, t \in (0, T]$ such that $0 < s \leq t \leq\epsilon$, we have that
	$d_X(\xi_s, \zeta_s) \leq \frac{d_X(\xi_t, \zeta_t)}{t} \cdot ( s + s \cdot t^2 \cdot \mathfrak{E})$.} The proof that the inverse \EpsE-interceptor bound is implied by the $\text{CAT}(\kappa)$ assumption follows the same structure as the proof of \cref{lem:cbb_interceptor_bound}, and we omit it in the present version of the manuscript.}, we get from \cref{eq:intermediate_AC_proof_12} that, for fixed $\underline{T}$ small enough\footnote{I.e., such that $(\underline{T} - s) V_{\max} \leq \epsilon_{\text{inv}}$, where $\epsilon_{\text{inv}}$ is the length bound of the inverse \EpsE-interceptor bound.}
	\begin{equation*}
		d_X \big(\xi^{\underline{t}}_{\underline{t}}, V_{s}^\gamma(\underline{t})\big) \leq
        d_X\big(\xi^{\underline{t}}_{\underline{T}}, V_{s}^\gamma(\underline{T})\big) \ O((\underline{t} - s)) = d_X \big(V^\gamma_{\underline{t}}(\underline{T}), V^\gamma_s(\underline{T})\big) \ O((\underline{t} - s)) \leq O\big((\underline{t} - s)^2\big).
	\end{equation*}
	Since $(\underline{T} - \underline{t}) \ V_{\max} \leq l_X$, we can use the symmetry assumption\footnote{By a simple rescaling step, the inequality from the symmetry assumption extends to those geodesics whose \emph{length} in space is smaller than $l_X$.} on the space $X$ together with \cref{eq:intermediate_AC_proof_11} to further obtain (see \cref{fig:cross_figure})
	\begin{equation}\label{eq:xi_cross_estimate}
		\begin{split}
			d_X\Big(\xi^{\underline{t}}_{2 \underline{T} - \underline{t}}, V_{\underline{t}}^\gamma\big(2 \underline{T} - \underline{t}\big)\Big) & \leq C_X \ d_X \big(\xi^{\underline{t}}_{\underline{t}}, \gamma_{\underline{t}}\big) = O\Big( d_X \big(\xi^{\underline{t}}_{\underline{t}}, \gamma_{\underline{t}} \big) \Big) \\
			& \leq O \Big( d_X\big( \xi^{\underline{t}}_{\underline{t}}, V^\gamma_s(\underline{t}) \big) +d_X \big( V^\gamma_s(\underline{t}), \gamma_{\underline{t}} \big)\Big)= O((\underline{t} - s)^2).
		\end{split}
	\end{equation}
	\begin{figure}[h]
		\centering
		\definecolor{fig_color_1}{HTML}{000000}
\definecolor{fig_color_2}{HTML}{4C72B0}
\definecolor{fig_color_3}{HTML}{DD8452}
\definecolor{fig_color_4}{HTML}{55A868}
\definecolor{fig_color_5}{HTML}{C44E52}
\definecolor{fig_color_6}{HTML}{8172B2}
\begin{tikzpicture}[scale=2]
  \draw [fig_color_2, thick] (0.300,-0.148) -- (0.320,-0.140) -- (0.340,-0.131) -- (0.361,-0.123) -- (0.381,-0.114) -- (0.401,-0.105) -- (0.421,-0.097) -- (0.441,-0.088) -- (0.462,-0.080) -- (0.482,-0.071) -- (0.502,-0.063) -- (0.522,-0.054) -- (0.542,-0.046) -- (0.563,-0.037) -- (0.583,-0.028) -- (0.603,-0.020) -- (0.623,-0.011) -- (0.643,-0.003) -- (0.664,0.006) -- (0.684,0.014) -- (0.704,0.023) -- (0.724,0.031) -- (0.744,0.040) -- (0.765,0.049) -- (0.785,0.057) -- (0.805,0.066) -- (0.825,0.074) -- (0.845,0.083) -- (0.866,0.091) -- (0.886,0.100) -- (0.906,0.108) -- (0.926,0.117) -- (0.946,0.126) -- (0.967,0.134) -- (0.987,0.143) -- (1.007,0.151) -- (1.027,0.160) -- (1.047,0.168) -- (1.068,0.177) -- (1.088,0.185) -- (1.108,0.194) -- (1.128,0.203) -- (1.148,0.211) -- (1.169,0.220) -- (1.189,0.228) -- (1.209,0.237) -- (1.229,0.245) -- (1.249,0.254) -- (1.270,0.262) -- (1.290,0.271) -- (1.310,0.280) -- (1.330,0.288) -- (1.351,0.297) -- (1.371,0.305) -- (1.391,0.314) -- (1.411,0.322) -- (1.431,0.331) -- (1.452,0.339) -- (1.472,0.348) -- (1.492,0.357) -- (1.512,0.365) -- (1.532,0.374) -- (1.553,0.382) -- (1.573,0.391) -- (1.593,0.399) -- (1.613,0.408) -- (1.633,0.416) -- (1.654,0.425) -- (1.674,0.434) -- (1.694,0.442) -- (1.714,0.451) -- (1.734,0.459) -- (1.755,0.468) -- (1.775,0.476) -- (1.795,0.485) -- (1.815,0.493) -- (1.835,0.502) -- (1.856,0.511) -- (1.876,0.519) -- (1.896,0.528) -- (1.916,0.536) -- (1.936,0.545) -- (1.957,0.553) -- (1.977,0.562) -- (1.997,0.570) -- (2.017,0.579) -- (2.037,0.588) -- (2.058,0.596) -- (2.078,0.605) -- (2.098,0.613) -- (2.118,0.622) -- (2.138,0.630) -- (2.159,0.639) -- (2.179,0.647) -- (2.199,0.656) -- (2.219,0.665) -- (2.239,0.673) -- (2.260,0.682) -- (2.280,0.690) -- (2.300,0.699);
  \draw [fig_color_2, thick] (1.200,-0.012) -- (1.220,-0.015) -- (1.240,-0.017) -- (1.261,-0.019) -- (1.281,-0.022) -- (1.301,-0.024) -- (1.321,-0.027) -- (1.341,-0.029) -- (1.362,-0.032) -- (1.382,-0.034) -- (1.402,-0.037) -- (1.422,-0.039) -- (1.442,-0.041) -- (1.463,-0.044) -- (1.483,-0.046) -- (1.503,-0.049) -- (1.523,-0.051) -- (1.543,-0.054) -- (1.564,-0.056) -- (1.584,-0.059) -- (1.604,-0.061) -- (1.624,-0.063) -- (1.644,-0.066) -- (1.665,-0.068) -- (1.685,-0.071) -- (1.705,-0.073) -- (1.725,-0.076) -- (1.745,-0.078) -- (1.766,-0.081) -- (1.786,-0.083) -- (1.806,-0.085) -- (1.826,-0.088) -- (1.846,-0.090) -- (1.867,-0.093) -- (1.887,-0.095) -- (1.907,-0.098) -- (1.927,-0.100) -- (1.947,-0.103) -- (1.968,-0.105) -- (1.988,-0.107) -- (2.008,-0.110) -- (2.028,-0.112) -- (2.048,-0.115) -- (2.069,-0.117) -- (2.089,-0.120) -- (2.109,-0.122) -- (2.129,-0.125) -- (2.149,-0.127) -- (2.170,-0.129) -- (2.190,-0.132) -- (2.210,-0.134) -- (2.230,-0.137) -- (2.251,-0.139) -- (2.271,-0.142) -- (2.291,-0.144) -- (2.311,-0.147) -- (2.331,-0.149) -- (2.352,-0.151) -- (2.372,-0.154) -- (2.392,-0.156) -- (2.412,-0.159) -- (2.432,-0.161) -- (2.453,-0.164) -- (2.473,-0.166) -- (2.493,-0.169) -- (2.513,-0.171) -- (2.533,-0.173) -- (2.554,-0.176) -- (2.574,-0.178) -- (2.594,-0.181) -- (2.614,-0.183) -- (2.634,-0.186) -- (2.655,-0.188) -- (2.675,-0.191) -- (2.695,-0.193) -- (2.715,-0.195) -- (2.735,-0.198) -- (2.756,-0.200) -- (2.776,-0.203) -- (2.796,-0.205) -- (2.816,-0.208) -- (2.836,-0.210) -- (2.857,-0.213) -- (2.877,-0.215) -- (2.897,-0.217) -- (2.917,-0.220) -- (2.937,-0.222) -- (2.958,-0.225) -- (2.978,-0.227) -- (2.998,-0.230) -- (3.018,-0.232) -- (3.038,-0.235) -- (3.059,-0.237) -- (3.079,-0.239) -- (3.099,-0.242) -- (3.119,-0.244) -- (3.139,-0.247) -- (3.160,-0.249) -- (3.180,-0.252) -- (3.200,-0.254);
  \draw [fig_color_1, very thick] (0.000,-0.303) -- (0.020,-0.290) -- (0.040,-0.279) -- (0.061,-0.267) -- (0.081,-0.256) -- (0.101,-0.244) -- (0.121,-0.234) -- (0.141,-0.223) -- (0.162,-0.213) -- (0.182,-0.203) -- (0.202,-0.193) -- (0.222,-0.183) -- (0.242,-0.174) -- (0.263,-0.164) -- (0.283,-0.156) -- (0.303,-0.147) -- (0.323,-0.139) -- (0.343,-0.130) -- (0.364,-0.123) -- (0.384,-0.115) -- (0.404,-0.107) -- (0.424,-0.100) -- (0.444,-0.093) -- (0.465,-0.087) -- (0.485,-0.080) -- (0.505,-0.074) -- (0.525,-0.068) -- (0.545,-0.062) -- (0.566,-0.057) -- (0.586,-0.052) -- (0.606,-0.047) -- (0.626,-0.042) -- (0.646,-0.038) -- (0.667,-0.034) -- (0.687,-0.030) -- (0.707,-0.026) -- (0.727,-0.023) -- (0.747,-0.019) -- (0.768,-0.016) -- (0.788,-0.014) -- (0.808,-0.011) -- (0.828,-0.009) -- (0.848,-0.007) -- (0.869,-0.005) -- (0.889,-0.004) -- (0.909,-0.002) -- (0.929,-0.002) -- (0.949,-0.001) -- (0.970,-0.000) -- (0.990,-0.000) -- (1.010,-0.000) -- (1.030,-0.000) -- (1.051,-0.001) -- (1.071,-0.002) -- (1.091,-0.003) -- (1.111,-0.004) -- (1.131,-0.005) -- (1.152,-0.007) -- (1.172,-0.009) -- (1.192,-0.011) -- (1.212,-0.014) -- (1.232,-0.016) -- (1.253,-0.019) -- (1.273,-0.023) -- (1.293,-0.026) -- (1.313,-0.030) -- (1.333,-0.034) -- (1.354,-0.038) -- (1.374,-0.042) -- (1.394,-0.047) -- (1.414,-0.052) -- (1.434,-0.057) -- (1.455,-0.062) -- (1.475,-0.068) -- (1.495,-0.074) -- (1.515,-0.080) -- (1.535,-0.087) -- (1.556,-0.093) -- (1.576,-0.100) -- (1.596,-0.107) -- (1.616,-0.115) -- (1.636,-0.123) -- (1.657,-0.130) -- (1.677,-0.139) -- (1.697,-0.147) -- (1.717,-0.156) -- (1.737,-0.164) -- (1.758,-0.174) -- (1.778,-0.183) -- (1.798,-0.193) -- (1.818,-0.203) -- (1.838,-0.213) -- (1.859,-0.223) -- (1.879,-0.234) -- (1.899,-0.244) -- (1.919,-0.256) -- (1.939,-0.267) -- (1.960,-0.279) -- (1.980,-0.290) -- (2.000,-0.303);
  \draw [fig_color_4, thick] (0.300,-0.148) -- (0.320,-0.147) -- (0.341,-0.146) -- (0.361,-0.144) -- (0.381,-0.143) -- (0.401,-0.141) -- (0.422,-0.140) -- (0.442,-0.139) -- (0.462,-0.137) -- (0.482,-0.136) -- (0.503,-0.135) -- (0.523,-0.133) -- (0.543,-0.132) -- (0.563,-0.131) -- (0.584,-0.129) -- (0.604,-0.128) -- (0.624,-0.127) -- (0.645,-0.125) -- (0.665,-0.124) -- (0.685,-0.122) -- (0.705,-0.121) -- (0.726,-0.120) -- (0.746,-0.118) -- (0.766,-0.117) -- (0.786,-0.116) -- (0.807,-0.114) -- (0.827,-0.113) -- (0.847,-0.112) -- (0.867,-0.110) -- (0.888,-0.109) -- (0.908,-0.108) -- (0.928,-0.106) -- (0.948,-0.105) -- (0.969,-0.103) -- (0.989,-0.102) -- (1.009,-0.101) -- (1.030,-0.099) -- (1.050,-0.098) -- (1.070,-0.097) -- (1.090,-0.095) -- (1.111,-0.094) -- (1.131,-0.093) -- (1.151,-0.091) -- (1.171,-0.090) -- (1.192,-0.089) -- (1.212,-0.087) -- (1.232,-0.086) -- (1.252,-0.085) -- (1.273,-0.083) -- (1.293,-0.082) -- (1.313,-0.080) -- (1.334,-0.079) -- (1.354,-0.078) -- (1.374,-0.076) -- (1.394,-0.075) -- (1.415,-0.074) -- (1.435,-0.072) -- (1.455,-0.071) -- (1.475,-0.070) -- (1.496,-0.068) -- (1.516,-0.067) -- (1.536,-0.066) -- (1.556,-0.064) -- (1.577,-0.063) -- (1.597,-0.061) -- (1.617,-0.060) -- (1.637,-0.059) -- (1.658,-0.057) -- (1.678,-0.056) -- (1.698,-0.055) -- (1.719,-0.053) -- (1.739,-0.052) -- (1.759,-0.051) -- (1.779,-0.049) -- (1.800,-0.048) -- (1.820,-0.047) -- (1.840,-0.045) -- (1.860,-0.044) -- (1.881,-0.042) -- (1.901,-0.041) -- (1.921,-0.040) -- (1.941,-0.038) -- (1.962,-0.037) -- (1.982,-0.036) -- (2.002,-0.034) -- (2.023,-0.033) -- (2.043,-0.032) -- (2.063,-0.030) -- (2.083,-0.029) -- (2.104,-0.028) -- (2.124,-0.026) -- (2.144,-0.025) -- (2.164,-0.024) -- (2.185,-0.022) -- (2.205,-0.021) -- (2.225,-0.019) -- (2.245,-0.018) -- (2.266,-0.017) -- (2.286,-0.015) -- (2.306,-0.014);
  \filldraw [fig_color_1] (0.300,-0.148) circle (0.02);
  \filldraw [fig_color_3] (2.300,0.699) circle (0.02);
  \filldraw [fig_color_1] (1.200,-0.012) circle (0.02);
  \filldraw [fig_color_4] (1.604,-0.061) circle (0.02);
  \filldraw [fig_color_5] (1.988,-0.107) circle (0.015);
  \filldraw [fig_color_5] (1.200,-0.012) circle (0.015);
  \filldraw [fig_color_3] (2.271,-0.142) circle (0.02);
  \filldraw [fig_color_5] (2.002,-0.034) circle (0.015);
  \filldraw [fig_color_5] (1.171,-0.090) circle (0.015);
  \filldraw [fig_color_3] (2.286,-0.015) circle (0.02);
  \draw [fig_color_5] (1.200,-0.012) -- (1.171,-0.090);
  \draw [fig_color_5] (1.988,-0.107) -- (2.002,-0.034);
  \draw [fig_color_3] (2.286,-0.015) -- (2.300,0.699);
  \draw [fig_color_3] (2.271,-0.142) -- (2.286,-0.015);
  \draw [fig_color_1] (0.300, -0.148)  node[ below] {$\gamma(s)$};
  \draw [fig_color_2] (1.310, 0.280)  node[ above] {$V_s^\gamma$};
  \draw [fig_color_2] (2.533, -0.173)  node[ below] {$V_{\underline{t}}^\gamma$};
  \draw [fig_color_4] (0.969, -0.103)  node[ below] {$\xi^{\underline{t}}$};
  \draw [fig_color_4] (1.604, -0.061)  node[ above] {$\xi^{\underline{t}}(\underline{T})$};
  \draw [fig_color_3, dashed] (2.300, 0.699) -- (3.271, 0.699);
  \draw [fig_color_3, dashed] (2.271, -0.142) -- (3.271, -0.142);
  \draw [fig_color_3, decorate, decoration={brace}, thick] (3.321, 0.699) -- (3.321, -0.142)  node[midway, right] {$\quad d_X\big(V^\gamma_{s}(T), V^\gamma_{\underline{t}}(T)\big)$};
\end{tikzpicture}
		\caption{The estimates of $d(V_s^\gamma(T), V_{\underline{t}}^\gamma(T))$ using $\xi^{\underline{t}}$. In black, we see the curve $\gamma$, in blue, the geodesics $V_s^\gamma$ and $V_{\underline{t}}^\gamma$ and in green the geodesic $\xi^{\underline{t}}$, which crosses $V_{\underline{t}}$ in $\xi^{\underline{t}}(\underline{T}) = V_{\underline{t}}^\gamma(\underline{T})$. In red, we see the four points (clockwise, starting top right) $\xi^{\underline{t}}_{2\underline{T}-\underline{t}}$, $V_{\underline{t}}^\gamma(2 \underline{T} - \underline{t})$, $\xi^{\underline{t}}_{\underline{t}}$, and $V_{\underline{t}}^\gamma(\underline{t}) = \gamma(\underline{t})$, used in \cref{eq:xi_cross_estimate}. In orange, we see the three points (top to bottom) $V_s^\gamma(T)$, $\xi^{\underline{t}}(T)$ and $V_{\underline{t}}^\gamma(T)$, used in \cref{eq:intermediate_AC2_proof_2}.}\label{fig:cross_figure}
	\end{figure}
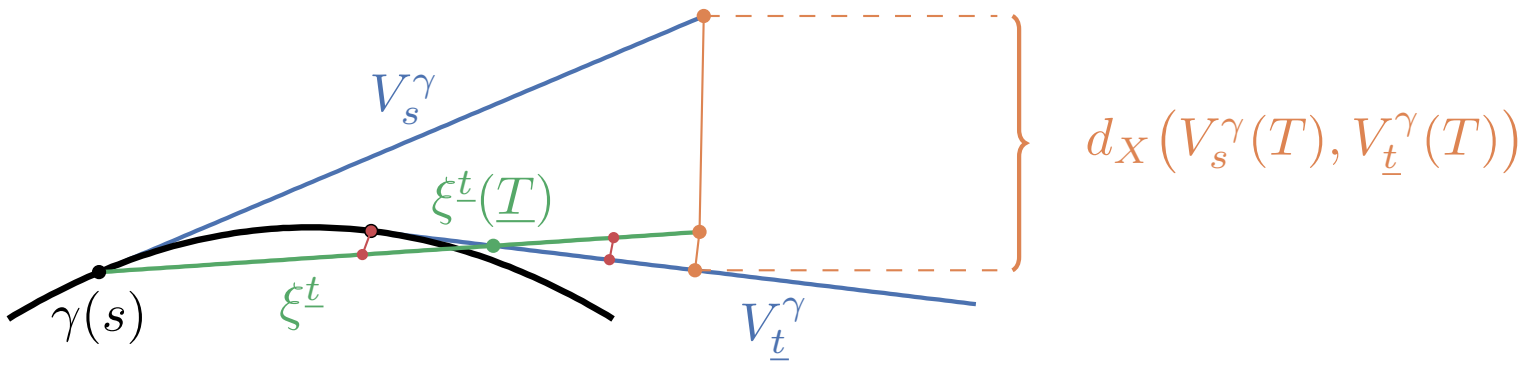
	Using the \EpsE-interceptor bound, we get $d_X \big( \xi^{\underline{t}}_T, V_{\underline{t}}^\gamma(T) \big) \leq O((\underline{t} - s)^2)$.
	By leveraging the triangle inequality once again, we thus obtain
	\begin{equation}\label{eq:intermediate_AC2_proof_2}
		\begin{split}
			\frac{d_X \big(V_s^\gamma(T), V_{\underline{t}}^\gamma(T) \big)}{\underline{t} - s} & \leq \frac{d_X \big(V_s^\gamma(T), \xi^{\underline{t}}_T \big)}{ \underline{t} - s} + \frac{d_X \big( \xi^{\underline{t}}_T, V_{\underline{t}}^\gamma(T) \big)}{\underline{t} - s} = \frac{d_X \big( V_s^\gamma(T), \xi^{\underline{t}}_T \big)}{ \underline{t} - s} + O(\underline{t} - s).
		\end{split}
	\end{equation}
	Motivated by the above calculation, we analyze $d_X \big(V_s^\gamma(T), \xi^{\underline{t}}_T\big)$: By the \EpsE-interceptor bound and \cref{eq:intermediate_AC_proof_12}, we have that
	\begin{equation}\label{eq:intermediate_AC2_proof_3}
		\begin{split}
			d_X \big(V_s^\gamma(T), \xi^{\underline{t}}_T \big) & \leq d_X(V_s^\gamma(\underline{T}), \xi^{\underline{t}}_{\underline{T}}) \ \frac{T-s}{\underline{T}-s} \ \big(1 + V_{\max}^2 \ \mathfrak{E} \ (T-s)^2 \big) \\
			& \leq (|\gamma''|(s) + \alpha) \ (\underline{t} - s)\ (T-s) \ \big(1 + V_{\max}^2 \ \mathfrak{E} \ (T-s)^2 \big),
		\end{split}
	\end{equation}
	where we used \cref{eq:intermediate_AC_proof_12} and the fact that $\xi^{\underline{t}}_{\underline{T}} = V_{\underline{t}}^\gamma(\underline{T})$ in the last inequality.
	Combining \cref{eq:intermediate_AC2_proof_2,eq:intermediate_AC2_proof_3}, we obtain
	\begin{equation*}
		\limsup_{t \searrow s} \frac{d_X \big( V_s^\gamma(T), V_t^\gamma(T) \big)}{ (t - s)} \leq (|\gamma''|(s) + \alpha) \ (T-s) \ (1 + V_{\max}^2 \ \mathfrak{E} \ (T-s)^2);
	\end{equation*}
	since $\alpha > 0$ was arbitrary, we obtain \cref{eq:intermediate_AC2_proof}.\\
	Finally, let us show how \cref{eq:intermediate_AC2_proof} implies \cref{eq:goal_AC2_proof}. 
    For a fixed $s < T$, we define the function $G_{s, T}(t) \coloneqq d_X\big(V^\gamma_s(T), V^\gamma_t(T)\big)$ for $t \in [s, T] \cap (a, b)$. 
    From $\gamma$ being an \ACt-curve, using \cref{eq:AC2_definition} and the reverse triangle inequality, we get
    \begin{equation*}
        |G_{s,T}(t_1)-G_{s,T}(t_2)| \leq d_X\big(V_{t_1}^\gamma(T), V_{t_2}^\gamma(T)\big) \leq \int_{t_1}^{t_2} (T-r)\ m(r)\big(1+\mathfrak{E}^{\mathrm{old}}_\gamma(T-r)^2\big)\,\diff r
    \end{equation*}
    for $t_1 \leq t_2$, from which we deduce that $G_{s, T}(\cdot)$ is locally absolutely continuous. 
    Furthermore, using \cref{eq:intermediate_AC2_proof} together with the triangle inequality, we get, for almost every $r \in (s, T) \cap (a, b)$,
	\begin{equation*}
		\begin{split}
			G_{s, T}'(r) = \limsup_{h \searrow 0} \frac{G_{s, T}(r + h) - G_{s, T}(r)}{h} & \leq \limsup_{h \searrow 0} \frac{d_X\big(V_r^\gamma(T), V_{r+h}^\gamma(T)\big)}{h} \\
			& \leq (T-r) \ |\gamma''|(r) \ (1 + V_{\max}^2 \ \mathfrak{E} \ (T-r)^2),
		\end{split}
	\end{equation*}
	which proves \cref{eq:goal_AC2_proof} by integrating from $s$ to $t$ and thus finishes the proof of the proposition.
	% To this end, we fix $s$ and $\tau$ as above and define the function $G_{s, \tau}(t) \coloneqq d_X\big(V_s^\gamma(\tau), V_t^\gamma(\tau - t + s)\big)$ for $t \in [s, s + \tau]$. 
\end{proof}

\section{Second Pillar: Variational Characterization of Curves of Newtonian Motion} \label{sec:curves_of_newtonian_motion}

In this section, we work with a driving energy $E\colon X \to (-\infty,+\infty]$, and we construct a metric variational characterization of curves of Newtonian motion for that driving energy. 
To this end, we will often require some assumptions on $E$.

\begin{assumption}[Lower Semicontinuity, cf.\ 2.1a,~\cite{ambrosioGradientFlows2008}]\label{ass:lsc}
	$E$ is lower semicontinuous.  %on $d_X$-bounded sets.
\end{assumption}
The most basic remaining assumption ensures that the energy does not tend to minus infinity faster than `quadratically'.

\begin{assumption}[Coercivity; cf.\ 2.1b,~\cite{ambrosioGradientFlows2008}]\label{ass:coercivity}
	There exists some $\tau_\star > 0$ and some $x_\star \in X$ such that
	\begin{equation*}
		M^\star \coloneqq \inf_{x \in X} \left\{ E(x) + \frac{1}{2\tau_\star} d_X(x, x_\star)^2 \right\} > -\infty.
	\end{equation*}
\end{assumption}

To be able to construct curves of Newtonian motion, we require a regularity assumption on $E$, for which we differentiate between a weak and a strong one:
{\renewcommand{\theassumption}{\arabic{assumption}.w}
\begin{assumption}[Existence of the derivative]\label{ass:existence_derivative}
	The functional $E$ fulfills the following regularity condition: For every $\varepsilon > 0$ and every geodesic $\xi \in \Tang[0, \varepsilon]$, $\frac{\diff}{\diff s} E(\xi_s)|_{s=t}$ exists for all $t \in [0, \varepsilon]$ and $t\mapsto \frac{\diff}{\diff s} E(\xi_s)|_{s=t}$ is absolutely continuous.
\end{assumption}}
\addtocounter{assumption}{-1}
{\renewcommand{\theassumption}{\arabic{assumption}.s}
\begin{assumption}[Lipschitz continuity of the derivative]\label{ass:lipschitz}
	The functional $E$ fulfills the following regularity condition: there is a global constant $L$ such that, for every $\varepsilon > 0$ and for every geodesic $\xi \in \Tang[0, \varepsilon]$ with $|\xi'| \leq 1$, $\frac{\diff}{\diff s} E(\xi_s)|_{s=t}$ exists for all $t \in [0, \varepsilon]$ and $t\mapsto \frac{\diff}{\diff s} E(\xi_s)|_{s=t}$ is $L$-Lipschitz continuous.
\end{assumption}
}
\begin{remark}\label{rmk:rescaling_lipschitz}
	By rescaling any geodesic $\xi \in \Tang[0, \varepsilon]$, we immediatley see if $E$ fulfills \cref{ass:lipschitz}, then $\frac{\diff}{\diff s} E(\xi_s)|_{s=t}$ exists for all $t \in [0, \varepsilon]$ and $t\mapsto \frac{\diff}{\diff s} E(\xi_s)|_{s=t}$ is Lipschitz with constant $|\xi'|^2 L$.
\end{remark}
\begin{assumption}\label{ass:functional}
	Let $\{\xi^{(i)}\}_{i\in \NN} \subset \Tang[0, \varepsilon]$ be a sequence of geodesics such that $\xi^{(i)} \to \xi$ uniformly as $i \to \infty$. Then, we have that
	\begin{equation*}
		\lim_{i \to \infty} \frac{\diff}{\diff t} E(\xi^{(i)}_t)|_{t=0} = \frac{\diff}{\diff t} E(\xi_t)|_{t=0}.
	\end{equation*}
\end{assumption}
\begin{assumption}\label{ass:local_compactness}
	The space $X$ is complete and locally compact.
\end{assumption}

\begin{remark}
    \cref{ass:lipschitz} and \cref{ass:functional} are both relatively strong ones, requiring in practice a certain smoothness of the energy $E$. Yet, in a second stage analysis and depending on the specific application context, one can extend results of existence of Newtonian curves also for nonsmooth $E$, by means of {\it ad hoc}  arguments, see, for example, Section \ref{sec:linear_wave_equation} for the case of the Dirichlet energy and the wave equation. Also in the companion paper \cite{fornasierLionsLiftCompanion} we can obtain Newtonian curves by {\it ad hoc} surgery also for certain nonsmooth interaction energies on the Wasserstein space.
\end{remark}

%\subsection{Defining Curves of Newtonian Motion} %\label{sec:defining_curves_of_newtonian_motion}

The challenge that we tackle in this section is the variational characterization of the \emph{curves of Newtonian motion} in metric spaces, where the classical tangent-vector calculus is not available. 
Similarly as done for the notion of \emph{curves of maximal slope} in metric spaces in \cite{ambrosioGradientFlows2008}, we will define them through an integral inequality involving exclusively metric objects such as the metric acceleration and an \emph{upper gradient-Hessian pair} (see the definition below). The inequality in the definition of \emph{curves of Newtonian motion} is meant to enforce the identity, as the opposite one is always provided for every \ACt-curve.

\begin{definition}[Upper gradient-Hessian pair] \label{def:upper_gradient}
    Let us assume that the metric space $(X,d_X)$ complies with \cref{ass:interceptor}.
	A pair of functions $g\colon X \to [0, +\infty]$, $H\colon \mathcal{T}(X) \to [-\infty,+\infty]$, is called an \emph{upper gradient-Hessian pair} for the functional $E\colon X \to (-\infty,+\infty]$ if, for every \ACt-curve $\gamma\colon (a, b)\to X$ such that
    \begin{itemize}
        \item $t \mapsto \frac{\diff}{\diff \sigma} E(V_t^\gamma(\sigma))\big|_{\sigma=t}$ exists as a finite real number for all $t \in (a, b)$ and has finite total variation on every compact subinterval of $(a, b)$,
        \item $t \mapsto H(V_t^\gamma)$ is in $L^1_{\mathrm{loc}}((a, b), \R)$,
    \end{itemize}
     we have that $g \circ \gamma \colon (a,b) \to [0,+\infty]$ is Borel-measurable, and that, for almost every $s , t \in (a, b)$  with $s<t$, the following inequality holds
	\begin{equation} \label{eq:def_upper_gradient}
		\frac{\diff}{\diff \sigma} E(V_s^\gamma(\sigma)) \Big|_{\sigma = s} 
        - \frac{\diff}{\diff \sigma}  E(V_t^\gamma(\sigma)) \Big|_{\sigma = t} 
        \leq \int_s^t \frac{1}{2}  |\gamma''|^2(r)\, \diff r + \int_s^t \frac{1}{2} g^2(\gamma_r)\, \diff r - \int_s^t H(V_r^\gamma)\, \diff r,
	\end{equation}
	where $|\gamma''|$ is defined in \cref{eq:metric_second_derivative_bound}, and $V_s^\gamma$, $V_t^\gamma$ are defined according to \cref{def:AC2_curve}.
\end{definition}

\begin{remark}
    We reiterate the fact that $(X,d_X)$ fulfilling \cref{ass:interceptor} guarantees that the geodesics $V_t^\gamma \in \mathcal{T}(X)$ involved in the definition of \ACt-curves are uniquely determined for almost every instant $t$ (see \cref{rmk:tangents_recoverable}). 
\end{remark}

We are now in a position to define curves of Newtonian motion in metric spaces.
\begin{definition}[Curve of Newtonian motion] \label{def:curve_newtonian_motion}
    Let us assume that the metric space $(X,d_X)$ complies with \cref{ass:interceptor}.
	Let $E\colon X \to (-\infty,+\infty]$ be a functional and let $(g, H)$ be an upper gradient-Hessian pair for $E$. \\
    An \ACt-curve $\gamma\colon (a, b)\to X$ is called a \emph{curve of Newtonian motion} for the functional $E$ with respect to the upper gradient-Hessian pair $(g, H)$ if
    \begin{itemize}
        \item $t \mapsto \frac{\diff}{\diff \sigma} E(V_t^\gamma(\sigma))\big|_{\sigma=t}$ exists as a finite real number for every $t\in(a,b)$ and has finite total variation on every compact subinterval of $(a, b)$,
        \item $t \mapsto H(V_t^\gamma)$ is in $L^1_{\mathrm{loc}}((a, b))$,
        \item for almost every $s , t \in (a, b)$  with $s<t$, we have 
	\begin{equation} \label{eq:def_curve_newtonian_motion}
    \begin{split}
		\frac{\diff}{\diff \sigma} E(V_s^\gamma(\sigma))\Big|_{\sigma = s} - &\frac{\diff}{\diff \sigma} E(V_t^\gamma(\sigma))\Big|_{\sigma = t} \\
        & \geq \int_s^t \frac{1}{2} |\gamma''|^2(r)\, \diff r + \int_s^t \frac{1}{2}  g^2(\gamma_r)\, \diff r - \int_s^t H(V_r^\gamma)\, \diff r.
    \end{split}
	\end{equation}
    \end{itemize}
\end{definition}

Below we show that our definition extends the existing notions of curve of Newtonian motion in Hilbert spaces (see \cref{ex:hilbert_space_newtonian_motion}) and in Riemannian manifolds (see \cref{ex:Riemannian_newtonian_motion}). 
Before proceeding, we introduce some notation that we require in the rest of the paper.

\begin{definition} \label{def:slopes_and_H}
	Given a functional $E\colon X \to (-\infty,+\infty]$, write $D(E)=\{x\in X:E(x)<+\infty\}$. We define the slopes $|\partial E|,|\partial^-E|:X\to[0,+\infty]$ and the following objects:
	\begin{align}
		|\partial E|(x) & \coloneqq \begin{cases}\limsup_{y\to x,\,y\ne x}\frac{(E(x)-E(y))^+}{d_X(x,y)},&x\in D(E),\\+\infty,&x\notin D(E),\end{cases} \label{eq:def_slope} \\
		|\partial^- E|(x) & \coloneqq \liminf_{y \to x} |\partial E|(y), \label{eq:def_lowerslope} \\
		E_\tau(x) & \coloneqq \inf_{y \in X} \left\{ E(y) + \frac{1}{2\tau} d_X(x, y)^2 \right\},\qquad \tau > 0, \label{eq:def_Yosida} \\
		J_\tau[x] & \coloneqq \argmin_{y \in X} \left\{ E(y) + \frac{1}{2\tau} d_X(x, y)^2 \right\},\qquad \tau > 0, \qquad J_0[x] = \{x\} \label{eq:def_proximal} \\
		K(x, y, z) & \coloneqq 2 d_X(x, y)^2 + 2 d_X(y, z)^2 - d_X(x, z)^2,\qquad x, y, z \in X \label{eq:def_K} \\
		K^\star(x, y) & \coloneqq \argmin_{z \in X} K(x, y, z). \label{eq:def_Kstar}
	\end{align}
	Note that under \cref{ass:interceptor} and for $d_X(x, y) \leq \frac{\epsilon}{2}$, we can extend the geodesic $\xi \in \mathcal{T}_{\frac{\epsilon}{2}}X, \xi_0=x, \xi_{\frac{\epsilon}{2}} = y$ to a geodesic $\overline{\xi} \in \Tang$. Then, setting $z = \overline{\xi}_{\epsilon}$, we have that $K(x, y, z) = 0$. Furthermore, the \EpsE-interceptor bound ensures that no other point $z' \neq z$ achieves $K(x, y, z') = 0$, as any such point $z'$ gives rise to a geodesic $\zeta\colon [0, \epsilon] \to X$ such that $\zeta_0 = x = \overline{\xi}_0$, $\zeta(\frac{\epsilon}{2}) = y = \overline{\xi}_{\frac{\epsilon}{2}}$ and $\zeta_\epsilon = z'$. Justified by this observation and under a slight abuse of notation, we set $K^\star(x, y)$ to be the unique element $z$ for which $K(x, y, z)=0$ whenever $d_X(x, y) \leq \frac{\epsilon}{2}$. We further define:
	\begin{align}
		H^E_\tau(x, y) & \coloneqq \frac{1}{\tau^2}\Big(E(x) - 2 E(y) + E(K^\star(x, y))\Big), \qquad & x,y\in X,\ \tau>0,\ d_X(x,y)\leq\frac\epsilon2, \label{eq:def_H} \\
		H^E_0(\xi) & \coloneqq \limsup_{\tau \searrow 0} H^E_\tau(\xi_t, \xi_{t+\tau}) \nonumber\\ &= \limsup_{\tau \searrow 0} \frac{E(\xi_t) - 2 E(\xi_{t + \tau}) + E(\xi_{t + 2 \tau})}{\tau^2}, \qquad & \xi \in \Tang[t, \varepsilon], \varepsilon > 0, t \in \R, \label{eq:def_H0} \\
		H^{E, +}_0(\xi) & \coloneqq \limsup_{\Tang[t, \varepsilon] \ni \zeta \to \xi} H^E_0(\zeta), \qquad & \xi \in \Tang[t, \varepsilon], \varepsilon > 0, t \in \R, \label{eq:def_H0_plus} \\
		N_\tau(x, y, z) & \coloneqq \frac{d_X\big(K^\star(x, y), z\big)}{\tau^2}, \qquad & x, y, z \in X, \tau > 0, d_X(x, y) \leq \frac{\epsilon}{2}, \label{eq:def_N}
	\end{align}
	where the topology on $\Tang[t, \varepsilon]$ in the definition of $H^{E, +}_0$ is given by uniform convergence on $[t, t+\varepsilon]$. In the definition of $H^E_\tau$, we follow the convention that $+\infty - \infty = \infty$.
\end{definition}

We first prove that, in a Hilbert space, the notion of curve of Newtonian motion proposed in \cref{def:curve_newtonian_motion} is equivalent to be a solution of $\ddot\gamma = - \nabla E(\gamma)$.

\begin{proposition} \label{ex:hilbert_space_newtonian_motion}
	Let $(X, d_X) = (\mathcal{H}, |\cdot|)$ be a Hilbert space and let $E \in C^2(\mathcal{H})$. Then, we have that
	\begin{equation} \label{eq:example_Euclidean_H_and_slope}
		\begin{alignedat}{2}
			H_0^{E, +}(\xi) & = H_0^E(\xi) = \nabla^2E(\xi_t)[\dot\xi_t, \dot\xi_t], \qquad & & \forall \xi \in \Tang[t, \varepsilon] \\
			|\partial^- E|(x) & = |\partial E|(x) = |\nabla E|(x), \qquad & & \forall x \in X.
		\end{alignedat}
	\end{equation}
	Furthermore, the pair $(|\partial^- E|, H_0^{E, +})$ is an upper gradient-Hessian pair for $E$, and $E$ fulfills \cref{ass:functional}. Finally, an \ACt-curve $\gamma\colon (a, b) \to X$ is a curve of Newtonian motion for $E$ with respect to $(|\partial^- E|, H_0^{E, +})$ (as per \cref{def:curve_newtonian_motion}) if and only if the identity
	\begin{equation} \label{eq:newtonian_flow_ex_1}
		\ddot\gamma_t = - \nabla E(\gamma_t).
	\end{equation}
    holds for almost all $t \in (a, b)$.
\end{proposition}
\begin{proof}
	We start by showing \cref{eq:example_Euclidean_H_and_slope}. For every geodesic $\xi\colon [0, \varepsilon] \to X$, we have
	\begin{equation*}
		\begin{split}
			H^E_0(\xi) & = \limsup_{\tau \searrow 0} H^E_\tau(\xi_t, \xi_{t+\tau}) = \limsup_{\tau \searrow 0}\frac{1}{\tau^2}\Big(E(\xi_t) - 2 E(\xi_{t+\tau}) + E(\xi_{t+2\tau})\Big) = \frac{\diff^2}{\diff s^2} E(\xi_s)_{s=t} \\
			&=\nabla^2E(\xi_t)[\dot\xi_t, \dot\xi_t].
		\end{split}
	\end{equation*}
	Furthermore, if $\{\zeta_i\}_{i=1}^\infty$ is a sequence in $\Tang[0, \varepsilon]$ converging uniformly to $\xi$, then we can write $\zeta_i \colon t \mapsto x_i + t v_i$ and $\xi \colon t \mapsto x + t v$, and the uniform convergence implies that $x_i \to x$ and $v_i \to v$ as $i \to \infty$. Thus, by continuity of $\nabla^2 E$, we have that $H_0^E(\zeta_i) \to H_0^E(\xi)$ as $i \to \infty$, which shows that $H_0^{E, +}(\xi) = H_0^E(\xi)$ for all $\xi \in \Tang[0, \varepsilon]$. The second equality in \cref{eq:example_Euclidean_H_and_slope} is a similar continuity argument. \\
	Next, we show that the pair $(|\partial^- E|, H_0^{E, +})$ is an upper gradient-Hessian pair for $E$, and that an \ACt-curve $\gamma$ is a curve of Newtonian motion if and only if it fulfills \cref{eq:newtonian_flow_ex_1}. To this end, let $\gamma\colon (a, b) \to \mathcal{H}$ be an \ACt-curve. Then, by \cref{lem:euclidean_ac2_equivalence}, $\gamma$ is a usual twice absolutely continuous curve in $\mathcal{H}$, and $|\gamma''|(t) = |\ddot \gamma_t|$; since $E \in C^2$, we also have that $t \mapsto \langle \nabla E(\gamma_t), \dot \gamma_t\rangle$ is locally absolutely continuous (and thus has locally finite total variation) and $t \mapsto \nabla^2 E(\gamma_t)[\dot \gamma_t, \dot \gamma_t]$ is continuous (and thus is in $L^1_{\mathrm{loc}}$),  as demanded by \cref{def:curve_newtonian_motion}. We have, for all $s, t \in (a, b)$ with $s<t$,
	\begin{equation} \label{eq:comput_hilbert_conservative}
		\begin{split}
			\frac{\diff}{\diff \sigma} E(V_s^\gamma(\sigma))\Big|_{\sigma = s} & - \frac{\diff}{\diff \sigma} E(V_t^\gamma(\sigma))\Big|_{\sigma = t}
			 = \langle \nabla E(\gamma_s), \dot{\gamma}_s \rangle - \langle \nabla E(\gamma_t), \dot{\gamma}_t \rangle \\
             %= - \int_s^t \frac{\diff}{\diff r} \langle \nabla E(\gamma_r), \dot{\gamma}_r \rangle \, \diff r
			& = - \int_s^t \langle D^2 E(\gamma_r) \dot{\gamma}_r, \dot{\gamma}_r \rangle + \langle \nabla E(\gamma_r), \ddot{\gamma}_r \rangle \, \diff r \\
			& \leq \int_s^t \frac{1}{2} |\ddot{\gamma}_r|^2 + \frac{1}{2} |\nabla E(\gamma_r)|^2 - \nabla^2E(\gamma_r)[\dot\gamma_r, \dot\gamma_r] \, \diff r \\
			& = \int_s^t \frac{1}{2} |\gamma''|^2(r)\, \diff r + \int_s^t \frac{1}{2} |\partial^- E|^2(\gamma_r)\, \diff r - \int_s^t H_0^{E, +}(V_r^\gamma)\, \diff r,
		\end{split}
	\end{equation}
	which shows \cref{eq:def_upper_gradient}, where we have equality in the third line if and only if $\ddot{\gamma}_r = - \nabla E(\gamma_r)$ for almost all $r \in (a, b)$, which also shows the last statement of the proposition. \\
	Finally, the fact that $E$ fulfills \cref{ass:functional} follows from the continuity of $\nabla E$ combined with an argument similar to the one used to show $H_0^{E, +} = H_0^E$.
\end{proof}

Leveraging on the differentiable structure of the space, in Riemannian manifolds curves of Newtonian motion are classically the solutions of the second-order equation $\nabla_{\dot \gamma} \dot\gamma = - \nabla E(\gamma)$.
In the next result, we prove that---when $(X,d_X)$ is a Riemannian manifold---the latter identity and \cref{def:curve_newtonian_motion} are actually equivalent.

\begin{proposition} \label{ex:Riemannian_newtonian_motion}
	Let $(X,d_X)$ be a Riemannian manifold with bounded sectional curvature, and let $E\colon X \to \R$ be a $C^2$-functional. Then, we have that
	\begin{equation}\label{eq:example_Riemannian_H_and_slope}
		\begin{alignedat}{2}
			H_0^{E, +}(\xi) & = H_0^E(\xi) = \mathrm{Hess}\, E(\dot \xi_t, \dot \xi_t) \coloneqq \langle\nabla_{\dot \xi_t} \nabla E, \dot \xi_t \rangle_g \qquad & & \forall \xi \in \Tang[t, \varepsilon], \\
			|\partial^- E|(x) & = |\partial E|(x) = |\nabla E|(x) \qquad & & \forall x \in X.
		\end{alignedat}
	\end{equation}
	Furthermore, the pair $(|\partial^- E|, H_0^{E, +})$ is an upper gradient-Hessian pair for $E$, and $E$ fulfills \cref{ass:functional}. Finally, an \ACt-curve $\gamma\colon (a, b) \to X$ is a curve of Newtonian motion for $E$ with respect to $(|\partial^- E|, H_0^{E, +})$ if and only if the identity
	\begin{equation*}
		\nabla_{\dot \gamma_t} \dot\gamma_t = - \nabla E(\gamma_t)
	\end{equation*}
    holds for almost all $t \in (a, b)$.
\end{proposition}
\begin{proof}
	To show \cref{eq:example_Riemannian_H_and_slope}, we fix a geodesic $\xi\colon [t, t+\varepsilon] \to X$ and choose a normal chart $\varphi\colon U \subseteq X \to \R^n$ centered at $\xi_t$. 
    We may assume that $\varepsilon$ is small enough such that $\mathrm{Im}(\xi) \subseteq U$. 
    If we set $\bar{\xi} \coloneqq \varphi \circ \xi \colon [t,t+\varepsilon] \to \R^n$ and $\bar{E} \coloneqq E \circ \varphi^{-1} \colon \varphi(U)\to \R$, then, owing to the fact that we work with a normal chart, we have that
    \begin{equation*}
     s \mapsto \bar{\xi}(s) =  \varphi\Big(\exp_{\xi_t} \big( (s-t) \, \dot{\xi}_t \big)\Big)     
    \end{equation*}
    is a geodesic in $\R^n$, and $\langle\nabla_{\dot \xi_t} \nabla E, \dot \xi_t \rangle_g = \nabla^2 \bar{E}(\bar{\xi_t}) [\dot{\bar{\xi}}_t, \dot{\bar{\xi}}_t]$. Thus, $H_0^E(\xi) = \langle\nabla_{\dot \xi_t} \nabla E, \dot \xi_t \rangle_g$. Furthermore, for any sequence of geodesics $\zeta_i: [t, t+\varepsilon] \to X$ converging uniformly to $\xi$, we may again assume that $\mathrm{Im}(\zeta_i) \subseteq U$ for all $i$ large enough. If we set again $\bar{\zeta}_i \coloneqq \varphi \circ \zeta_i$, we have that $\bar{\zeta}_i(t) \to \bar{\xi}(t)$ and $\dot{\bar{\zeta}}_i(t) \to \dot{\bar{\xi}}(t)$ as $i \to \infty$\footnote{For fixed $0<h\leq\varepsilon$, the uniform convergence provides a bound for $|\zeta_i'|=d_X(\zeta_i(t),\zeta_i(t+h))/h$; local equivalence of metric and coordinate norms therefore bounds $\dot{\bar\zeta}_i(t)$. It then suffices to argue by contradiction: assume that, after picking a subsequence, $\dot{\bar{\zeta}}_{i, t} \to v \neq \dot{\bar{\xi}}_t$. Then, using continuity of the geodesic ODE in the initial data, we have that $\bar{\zeta}_i \to \big[s \mapsto (s-t) \, v]$ uniformly.}. Furthermore, the difference between $\langle\nabla_{\dot \zeta_i(t)} \nabla E, \dot \zeta_i(t) \rangle_g$ and $\nabla^2 \bar{E}(\bar{\zeta}_i(t))\big[\dot{\bar{\zeta}}_i(t), \dot{\bar{\zeta}}_i(t)\big]$ goes to zero as $\zeta_i(t) \to \xi(t)$\footnote{This follows from the fact that we can express $\langle\nabla_{\dot \zeta_i(t)} \nabla E, \dot \zeta_i(t) \rangle_g$ in coordinates using the Christoffel symbols, that this expression reduces to $\nabla^2 \bar{E}(\bar{\zeta}_i(t))\big[\dot{\bar{\zeta}}_i(t), \dot{\bar{\zeta}}_i(t)\big]$ whenever the Christoffel symbols vanish, that this expression is continuous in the Christoffel symbols and from the fact that the Christoffel symbols vanish at $\xi_t$.}, which shows the first line in \cref{eq:example_Riemannian_H_and_slope}. The second line follows from a similar argument using normal charts. \\
	Next, note that for any \ACt-curve $\gamma\colon (a, b) \to X$ locally on a time interval whose image and short tangent segments lie in a relatively compact subset of a chart $(U, \varphi)$, also $\varphi \circ \gamma$ is an \ACt-curve in $\R^d$\footnote{Here, it is necessary to see that we can find a bound on $(\varphi \circ \xi)''$, uniformly over all geodesics $\xi$ in $U$ with $|\xi'| \leq \sup_{s \in (a, b)} |(V^\gamma_s)'|$, and that, adding this bound to a constant multiple of $m$ in \cref{def:AC2_curve}, the integral requirement is fulfilled for the modified geodesics $\bar{V}_s^\gamma\colon s + \sigma \mapsto \varphi(\gamma_s) + (\varphi \circ V_s^\gamma)'_s \, \sigma $.}; using \cref{lem:euclidean_ac2_equivalence}, $\varphi \circ \gamma$ is thus a usual twice absolutely continuous curve in $\R^d$. As a next step, let us show that $|\nabla_{\dot{\gamma}_r} \dot{\gamma}_r| = |\gamma''|(r)$ for almost every $r$, where the latter term is defined in \cref{eq:metric_second_derivative_bound}. To this end, we fix a point $r$ where the coordinate second derivative and the Euclidean acceleration identity hold, and pick a normal chart $(U, \phi)$ centered at $\gamma_r$. Then $|\nabla_{\dot\gamma_r} \dot\gamma_r| = |(\phi \circ \gamma)_r''|$; and the defining limit in \cref{eq:metric_second_derivative_bound} can be calculated in either $M$ or $U$: For this, we calculate for $T \in \R$ such that $0 < T - r$ small and $h \in (0, T-r]$:
	\begin{align*}
		\frac{d_M\big(V^\gamma_r(T), V^\gamma_{r+h}(T)\big)}{(T-r)h} &= \frac{|\phi(V^\gamma_r(T)) - \phi(V^\gamma_{r+h}(T))|}{(T-r)h} \cdot (1 + o(1)) \\
		&= \frac{|V^{\phi \circ \gamma}_r(T) - V^{\phi \circ \gamma}_{r+h}(T)|}{(T-r)h} \cdot (1 + o(1)) + O\left(\frac{|\phi(V^\gamma_{r+h}(T)) - V^{\phi \circ \gamma}_{r+h}(T)|}{(T-r)h}\right),
	\end{align*}
	where we take the $o(1)$ term to mean that this term goes to zero as $T \searrow r$, and where we used that $V^{\phi \circ \gamma}_r(T) = \phi(V^\gamma_r(T)) = (T - r) (\phi \circ \gamma)'_r$ because $(U, \phi)$ is centered at $\gamma_r$. The $o(1)$ term above vanishes in the limit; what remains to show is that the last error term vanishes as well. To this end, it suffices to see that the function $G(x, v, t) \coloneqq \phi\big(\exp(t (\phi^{-1})'_x(v))\big) - x - tv$ is smooth, vanishes at $x=0$ and vanishes quadratically at $t=0$. Thus, $G(x, v, t) = O(|x| t^2)$ locally uniformly over bounded $v$. Applying this estimate with $x = \phi(\gamma_{r+h})$, $v = (\phi \circ \gamma)'_{r+h}$, and $t = T - (r+h)$, we see that the last error term vanishes, since $O(|\phi(\gamma_{r+h})|) = O(h)$.

	Since $E \in C^2$, using the local coordinate regularity above, we also have that $t \mapsto \langle \nabla E(\gamma_t), \dot \gamma_t\rangle_g$ is locally absolutely continuous (and thus has locally finite total variation) and $t \mapsto \mathrm{Hess}\, E(\dot \gamma_t, \dot \gamma_t)$ is continuous (and thus is in $L^1_{\mathrm{loc}}$), as demanded by \cref{def:curve_newtonian_motion}.
	Thus, $|\nabla_{\dot{\gamma}_r} \dot{\gamma}_r| = |\gamma''|(r)$ for almost every $r$, and showing \cref{eq:def_upper_gradient} is equivalent to showing that, for all $s, t \in (a, b)$ such that $s < t$ and all such $\gamma$,
	\begin{equation}\label{eq:riemannian_manifold_proof_0}
		\frac{\diff}{\diff \sigma} E(\gamma_\sigma)\Big|_{\sigma = s} - \frac{\diff}{\diff \sigma} E(\gamma_\sigma)\Big|_{\sigma = t} \leq \int_s^t \frac{1}{2}  |\nabla_{\dot \gamma_r} \dot\gamma_r|^2\, \diff r + \int_s^t \frac{1}{2}  |\nabla E|^2(\gamma_r)\, \diff r - \int_s^t \text{Hess } E(\dot \gamma_r, \dot \gamma_r)\, \diff r.
	\end{equation}
	By a telescoping sum argument, it suffices to show the \cref{eq:riemannian_manifold_proof_0} for small $t - s$. Because $\dot \gamma$ is bounded,
	we thus may assume that we can work on a single chart $(U, \varphi)$ containing $\gamma_{[s, t]}$. Note that \cref{eq:riemannian_manifold_proof_0} would follow automatically if $\gamma$ and $E$ were smooth: In this case, we would have, as in \cref{ex:hilbert_space_newtonian_motion},
	\begin{equation}
    \begin{split} \label{eq:riemannian_manifold_proof_1}
		\frac{\diff}{\diff \sigma} E(\gamma_\sigma)\Big|_{\sigma = s} - \frac{\diff}{\diff \sigma} E(\gamma_\sigma)\Big|_{\sigma = t} &= - \int_s^t \frac{\diff^2}{\diff \sigma^2} E(\gamma_\sigma) \Big|_{\sigma=r} \diff r = - \int_s^t \frac{\diff}{\diff \sigma} \langle \nabla E\big|_{\gamma_\sigma}, \dot{\gamma}_\sigma\rangle_g \Big|_{\sigma=r} \diff r \\ 
		&= - \int_s^t  \langle \nabla E, \nabla_{\dot{\gamma}} \dot{\gamma} \rangle_g\Big|_r + \langle\nabla_{\dot \gamma_r} \nabla E, \dot \gamma_r \rangle_g \, \diff r \\
		&\leq \int_s^t \frac{1}{2}\Big(|\nabla E|_g^2 + |\nabla_{\dot{\gamma}_r} \dot{\gamma}_r|_g^2\Big) - \mathrm{Hess}\, E(\dot \gamma_r, \dot \gamma_r) \, \diff r,
	\end{split}
    \end{equation}
	with equality if $\nabla_{\dot{\gamma}_r} \dot{\gamma}_r = -\nabla E(\gamma_r)$ for almost every $r$, where we used the metric compatibility of the Levi-Civita connection in \cref{eq:riemannian_manifold_proof_1}. \\
	To show \cref{eq:riemannian_manifold_proof_0} in the general case, we use the fact that $\bar{\gamma} \coloneqq \varphi \circ \gamma$ can be strongly approximated in $H^2$ by a smooth curve $\bar{\eta}$, and, by the Sobolev embedding, $\dot{\bar{\eta}}$ and $\bar{\eta}$ approximate $\dot{\bar{\gamma}}$ and $\bar{\gamma}$ even uniformly. Likewise, $\bar{E} \coloneqq E \circ \varphi^{-1}$ can be approximated by a smooth energy $\bar{F}$ in $C^2$ locally. By writing $\nabla_{\dot\gamma_r} \dot\gamma_r$ and $\mathrm{Hess}\, E(\dot{\gamma}_r, \dot{\gamma}_r)$ in coordinates and setting $\eta \coloneqq \varphi^{-1} \circ \bar{\eta}$, $F \coloneqq \bar{F} \circ \varphi$, we thus have that the right hand side of \cref{eq:riemannian_manifold_proof_0} is approximated by
	\begin{equation*}
		\int_s^t \frac{1}{2}  |\nabla_{\dot \eta_r} \dot\eta_r|^2\, \diff r + \int_s^t \frac{1}{2}  |\nabla F|^2(\eta_r)\, \diff r - \int_s^t \text{Hess } F(\dot \eta_r, \dot \eta_r)\, \diff r;
	\end{equation*}
	likewise, the left hand side is approximated by
	\begin{equation*}
		\frac{\diff}{\diff \sigma} F(\eta_\sigma)\Big|_{\sigma = s} - \frac{\diff}{\diff \sigma} F(\eta_\sigma)\Big|_{\sigma = t}.
	\end{equation*}
	Thus, the general case follows from the smooth case. \\
	Finally, the fact that $E$ fulfills \cref{ass:functional} follows from the continuity of $\nabla E$ and from another argument using local charts as for proving $H_0^{E, +} = H_0^E$.
\end{proof}

\begin{remark} \label{rmk:friction_extension}
    The variational definition of curve of Newtonian motion provided in \cref{def:curve_newtonian_motion} can be promptly extended to the solutions of \emph{damped Newton equation} with linear and isotropic viscous friction.
    To see this, let us assume that $(X, d_X) = (\mathcal{H}, |\cdot|)$ is a Hilbert space and that $E \in C^2(\mathcal{H})$.
    %and that the curve $\gamma \colon (a,b) \to \mathcal{H}$ is an \ACt-curve that solves
%    \begin{equation} \label{eq:damped_Newtonian}
%        \ddot \gamma_t = - \nabla E(\gamma_t) - \alpha \dot \gamma_t
%    \end{equation}
%    for some $\alpha \in \R$ (notice that, strictly speaking, the viscous friction requires $\alpha>0$). 
     Keeping in mind the computations done in \cref{eq:comput_hilbert_conservative}, we notice that, \emph{for every} \ACt-curve  $\gamma \colon (a,b) \to \mathcal{H}$ and for every $s,t\in (a,b)$ with $s<t$ and every $\alpha \in \R$, we have
    \begin{equation} \label{eq:comput_hilbert_damped_1}
		\begin{split}
			\langle \nabla E(\gamma_s), \dot{\gamma}_s \rangle - & \langle \nabla E(\gamma_t), \dot{\gamma}_t \rangle 
            - \alpha \big( E(\gamma_t) -  E(\gamma_s) \big) \\
             %= - \int_s^t \frac{\diff}{\diff r} \langle \nabla E(\gamma_r), \dot{\gamma}_r \rangle \, \diff r
			 &= - \int_s^t \langle D^2 E(\gamma_r) \dot{\gamma}_r, \dot{\gamma}_r \rangle + \langle \nabla E(\gamma_r), \ddot{\gamma}_r + \alpha \dot \gamma_r \rangle \, \diff r \\
			& \leq \int_s^t \frac{1}{2} |\ddot{\gamma}_r + \alpha \dot \gamma_r|^2 + \frac{1}{2} |\nabla E(\gamma_r)|^2 - \nabla^2E(\gamma_r)[\dot\gamma_r, \dot\gamma_r] \, \diff r\\
            & =  \int_s^t \frac{1}{2} |\ddot{\gamma}_r|^2 + \frac12 \alpha^2 |\dot \gamma_r|^2 + \alpha \langle \ddot \gamma_r, \dot \gamma_r  \rangle+ \frac{1}{2} |\nabla E(\gamma_r)|^2 - \nabla^2E(\gamma_r)[\dot\gamma_r, \dot\gamma_r] \, \diff r ,
		\end{split}
	\end{equation}
    and observing that  $\int _s^t \langle \ddot \gamma_r, \dot \gamma_r \rangle \, \diff r  = \frac12 |\dot \gamma_t|^2 - \frac12 |\dot \gamma_s|^2$, we can rearrange \cref{eq:comput_hilbert_damped_1} as follows:
    \begin{equation} \label{eq:comput_hilbert_damped_2}
		\begin{split}
			\langle \nabla E(\gamma_s), \dot{\gamma}_s \rangle - & \langle \nabla E(\gamma_t), \dot{\gamma}_t \rangle 
            - \alpha \left( E(\gamma_t) + \frac12|\dot \gamma_t|^2 -  E(\gamma_s)  - \frac12 |\dot \gamma_s|^2 \right) \\
            & \leq  \int_s^t \frac{1}{2} |\ddot{\gamma}_r|^2 + \frac12 \alpha^2 |\dot \gamma_r|^2 + \frac{1}{2} |\nabla E(\gamma_r)|^2 - \nabla^2E(\gamma_r)[\dot\gamma_r, \dot\gamma_r] \, \diff r.
		\end{split}
	\end{equation}
    We observe that the equality in \cref{eq:comput_hilbert_damped_2} holds if and only if $t\mapsto \gamma_t$ solves 
    \begin{equation*} %\label{eq:damped_Newtonian}
        \ddot \gamma_t = - \nabla E(\gamma_t) - \alpha \dot \gamma_t,
    \end{equation*}
    where $\alpha$ is the viscous friction parameter (notice that, strictly speaking, the viscous friction requires $\alpha>0$). Moreover, \cref{eq:comput_hilbert_damped_2} can be expressed in terms of purely metric objects, as its conservative counterpart \cref{eq:def_curve_newtonian_motion}.
    The study of damped Newtonian motion goes beyond the scope of the present paper, and it is left for future work.
\end{remark}

In the present section, we have characterized curves of Newtonian motion in metric spaces, and we have shown that our notion extends the one already available in the smooth setting. 
The final goal pursued is to construct them using the \emph{purely metric tools} that we develop in \cref{sec:discrete_scheme}.

\section{Third Pillar: Discrete Variational Scheme} \label{sec:discrete_scheme}

After having established a notion of curves of Newtonian motion, in the present section we define a discrete scheme yielding approximations of such curves. 
In \cref{sec:compactness_limiting} we will use a compactness argument to show that we can take the limit of the step size going to zero to obtain continuous curves of Newtonian motion.

\begin{definition}[Predictor-corrector scheme] \label{def:newtonian_sequence}
	Let $(X,d_X)$ fulfill \cref{ass:interceptor} and let $\epsilon>0$ be the corresponding geodesic horizon.
	Let us consider $\tau > 0$ and a functional $E\colon X\to \R\cup \{+\infty\}$. 
    We say a sequence $U_{\tau, -1}, \dots, U_{\tau, N}$ is a \emph{Newtonian sequence of stepsize $\tau$ for the energy $E$} whenever we have, for all $k =1,\ldots, N$, that
	\begin{equation}\label{eq:def_newtonian_sequence}
		d_X(U_{\tau, k-2}, U_{\tau, k - 1}) \leq \frac{\epsilon}{3}, \qquad U_{\tau, k} \in J_{\tau^2}[K^\star(U_{\tau, k-2}, U_{\tau, k-1})],
	\end{equation}
	where $K^\star$ is defined in \cref{eq:def_Kstar}, $J_{\tau^2}$ is defined in \cref{eq:def_proximal}.
\end{definition}

\begin{remark}
    \cref{def:newtonian_sequence} introduces a two-step \emph{predictor-corrector scheme}: For a fixed stepsize $\tau > 0$ and given the points $U_{\tau, k-2}$ and $U_{\tau, k-1}$ for $k =1,\ldots,N$, we recursively define
\begin{equation} \label{eq:remark_scheme}
	\begin{split}
	Z_{\tau, k-1} \coloneqq \argmin_{z \in X} \Big\{2\, d_X^2(U_{\tau, k-1}, z) - d_X^2(U_{\tau, k-2}, z) \Big\} \qquad &\mbox{(predictor)}, \\
	U_{\tau, k} \in \argmin_{z \in X} \Big\{E(z) + \frac{1}{2\tau^2} d_X^2(Z_{\tau, k-1}, z) \Big\} \qquad &\mbox{(corrector)}.
	\end{split}
\end{equation}
We notice that the temporary point $Z_{\tau, k-1}$ is the metric analogue of $U_{\tau, k-1} + \big( U_{\tau, k-1} - U_{\tau, k-2} \big)$, and that the geodesic extendability and the \EpsE-interceptor bound give, respectively, existence---whenever $d_X(U_{\tau, k-1},U_{\tau, k-2})\leq \epsilon/2$---and uniqueness in the predictor step.
\end{remark}

\begin{remark} \label{rmk:scheme_initialization}
    As we may expect due to the second-order nature of Newtonian motion, the construction of $U_{\tau,k}$ involves the previous two points $U_{\tau,k-1}$ and $U_{\tau,k-2}$. 
    This consideration immediately yields the question about the initialization of the scheme, as we need to properly set $U_{\tau,0}$ and $U_{\tau,-1}$.
    In the case $(X, d_X)= (\mathcal{H}, |\cdot|)$ is a Hilbert space and $E$ is a smooth energy, the equation $\ddot \gamma_t = -\nabla E(\gamma_t)$ requires a Cauchy datum in the form of an initial point $x_0$ and an initial velocity $v_0$.
    In our setting of a metric space $(X, d_X)$ complying with \cref{ass:interceptor}, the Cauchy datum is given in the form of an \emph{initial tangent geodesic} $\xi\colon [-\epsilon', \epsilon'] \to X$ such that $\xi(0)=x_0$, where $\epsilon'$ is such that $2 |\xi'|\epsilon' \leq \epsilon$, with $\epsilon$ provided by \cref{ass:interceptor}. 
    Therefore, for every $\tau \in (0, \epsilon')$ such that $\tau \, |\xi'|\leq \epsilon/3$, we set
    \begin{equation} \label{eq:scheme_initalization}
        U_{\tau,0} \coloneqq\xi_0=x_0, \qquad U_{\tau,-1}\coloneqq \xi_{-\tau}.
    \end{equation}
    We emphasize that, in the predictor phase yielding the construction of $Z_{\tau,k-1}$ (see \cref{eq:remark_scheme}), the stepsize $\tau$ does not appear \emph{explicitly}. This is due to the fact that the information on the stepsize is \emph{implicitly} contained in the distance $d_X(U_{\tau, k-1}, U_{\tau, k-2})$, as is apparent from \cref{eq:scheme_initalization} for $k=1$.
\end{remark}

\begin{remark} \label{rmk:corrector_stepsize}
    From \cref{eq:remark_scheme}, we observe that the corrector phase consists of a step of the minimizing movement scheme \emph{with stepsize $\tau^2$}.
    The latter choice is motivated by the following fact. Let us consider $(X,d_X)=(\R, |\cdot|)$, and the energies $E^1,E^2\colon \R\to \R$ defined respectively as $E^1(x) \coloneqq a_1 x $ and $E^2(x) \coloneqq \frac12 a_2 x^2$. We fix the starting point $x_0\in \R$ and the initial velocity $v_0\in \R$. 
    Then, we take $\tau>0$, we set
    \begin{equation*}
        U^\ell_{\tau, 0} \coloneqq  x_0, \quad U_{\tau,-1}^\ell \coloneqq x_0-\tau v_0 \qquad \mbox{for } \ \ell=1,2,
    \end{equation*}
    and we denote with $ U^\ell_{\tau, 1}$ the point obtained by applying the predictor-corrector scheme as per \cref{def:newtonian_sequence} (see also \cref{eq:remark_scheme}) to the energy $E^\ell$, for $\ell=1,2$.
    Finally, if we introduce $x^\ell\colon \R \to \R$ as exact solutions of
    \begin{equation*}
        \ddot x^\ell(t) = -\nabla E^\ell\big( x^\ell(t) \big), \quad x^\ell(0) = x_0, \quad \dot x^\ell(0)=v_0, \qquad \mbox{for } \ \ell=1,2,
    \end{equation*}
    a direct computation shows that
    \begin{equation*}
        \big| U^\ell_{\tau,1} - x^\ell(\tau) \big| = O\big(\tau^2 \big) \qquad \mbox{as $\tau \to 0$}, \qquad \mbox{for $\ell=1,2$},
    \end{equation*}
    as well as
    \begin{equation*}
        \left| \frac{U^\ell_{\tau,1} - U^\ell_{\tau,0}}{\tau} - \dot x^\ell(\tau) \right| = O\big(\tau^2 \big) \qquad \mbox{as $\tau \to 0$}, \qquad \mbox{for $\ell=1,2$}.
    \end{equation*}
%    are, respectively,
%    \begin{equation*}
%        x^1(t) = x_0 +v_0\, t + \frac{a}2t^2, \qquad  x^2(t) = x_0\cos\big( \sqrt{a} t \big) + \frac{v_0}{\sqrt{a}}\sin\big( \sqrt{a} t \big)
%    \end{equation*}
\end{remark}

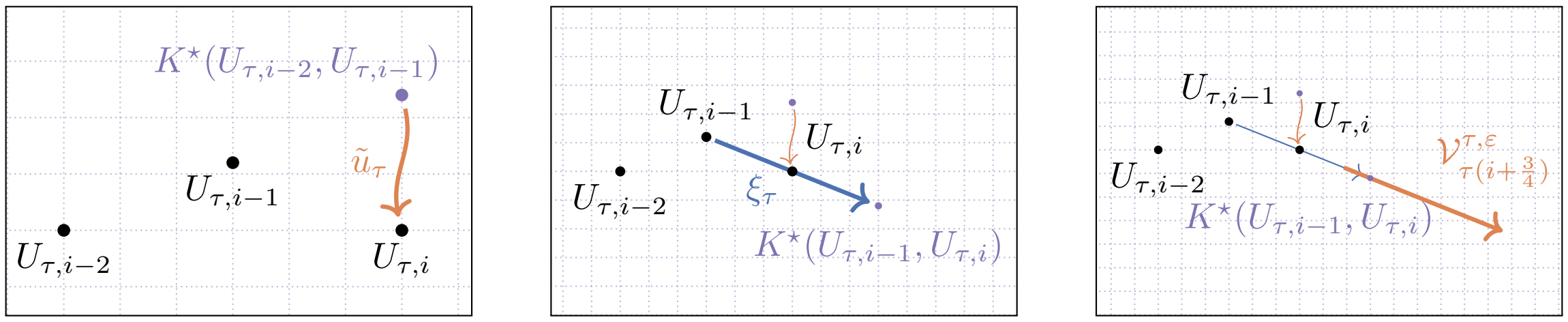
\begin{figure}[ht]
	\centering
	\begin{subfigure}[t]{0.3\textwidth}
		\centering
		\begin{adjustbox}{clip,margin=0pt,trim=0 0 0 0,fbox,width=\linewidth} % trim behaves like in includegraphics
			\begin{minipage}[t]{\linewidth}
				\definecolor{fig_color_1}{HTML}{000000}
\definecolor{fig_color_2}{HTML}{4C72B0}
\definecolor{fig_color_3}{HTML}{DD8452}
\definecolor{fig_color_4}{HTML}{55A868}
\definecolor{fig_color_5}{HTML}{C44E52}
\definecolor{fig_color_6}{HTML}{8172B2}
\begin{tikzpicture}[scale=0.550]
                  \path[use as bounding box] (-1,-1.500) rectangle (7,3.955);
  \draw[step=1,fig_color_6!50,densely dotted] (-10,-10) grid (10,10);
  \draw [->, fig_color_3, very thick] (6.059,2.158) -- (6.064,2.133) -- (6.069,2.109) -- (6.073,2.085) -- (6.078,2.061) -- (6.081,2.036) -- (6.085,2.012) -- (6.088,1.988) -- (6.091,1.964) -- (6.093,1.939) -- (6.095,1.915) -- (6.097,1.891) -- (6.098,1.867) -- (6.099,1.842) -- (6.100,1.818) -- (6.100,1.794) -- (6.100,1.770) -- (6.099,1.745) -- (6.098,1.721) -- (6.096,1.697) -- (6.095,1.673) -- (6.092,1.648) -- (6.090,1.624) -- (6.087,1.600) -- (6.083,1.576) -- (6.080,1.552) -- (6.076,1.527) -- (6.071,1.503) -- (6.067,1.479) -- (6.062,1.455) -- (6.057,1.430) -- (6.051,1.406) -- (6.046,1.382) -- (6.040,1.358) -- (6.034,1.333) -- (6.028,1.309) -- (6.022,1.285) -- (6.016,1.261) -- (6.010,1.236) -- (6.003,1.212) -- (5.997,1.188) -- (5.990,1.164) -- (5.984,1.139) -- (5.978,1.115) -- (5.972,1.091) -- (5.966,1.067) -- (5.960,1.042) -- (5.954,1.018) -- (5.949,0.994) -- (5.943,0.970) -- (5.938,0.945) -- (5.933,0.921) -- (5.929,0.897) -- (5.924,0.873) -- (5.920,0.848) -- (5.917,0.824) -- (5.913,0.800) -- (5.910,0.776) -- (5.908,0.752) -- (5.905,0.727) -- (5.904,0.703) -- (5.902,0.679) -- (5.901,0.655) -- (5.900,0.630) -- (5.900,0.606) -- (5.900,0.582) -- (5.901,0.558) -- (5.902,0.533) -- (5.903,0.509) -- (5.905,0.485) -- (5.907,0.461) -- (5.909,0.436) -- (5.912,0.412) -- (5.915,0.388) -- (5.919,0.364) -- (5.922,0.339) -- (5.927,0.315) -- (5.931,0.291) -- (5.936,0.267) -- (5.941,0.242);
  \filldraw [fig_color_1] (0.000,0.000) circle (0.1) node[below]  {$U_{\tau, i-2}$};
  \filldraw [fig_color_1] (3.000,1.200) circle (0.1) node[below]  {$U_{\tau, i-1}$};
  \filldraw [fig_color_6] (6.000,2.400) circle (0.1) node[anchor=south east, xshift=5mm]  {$K^\star(U_{\tau, i-2}, U_{\tau, i-1})$};
  \filldraw [fig_color_1] (6.000,0.000) circle (0.1) node[below]  {$U_{\tau, i}$};
  \draw [fig_color_3] (6.000, 1.200)  node[anchor=east] {$\tilde{u}_\tau$};
\end{tikzpicture}
				\centering
			\end{minipage}
		\end{adjustbox}
		\caption{While it is in general \emph{not} necessarily piecewise continuous, we can still imagine $\tilde{u}_\tau(t)$ to \emph{interpolate} between $K^\star(U_{\tau, i-1}, U_{\tau, i})$ and $U_{\tau, i+1}$ for $t \in [\tau i, \tau (i+1))$.}\label{fig:interp:1}
	\end{subfigure}
	\hfill
	\begin{subfigure}[t]{0.3\textwidth}
		\centering
		\begin{adjustbox}{clip,margin=0pt,trim=0 0 0 0,fbox,width=\linewidth} % trim behaves like in includegraphics
			\begin{minipage}[t]{\linewidth}
				\centering
				\definecolor{fig_color_1}{HTML}{000000}
\definecolor{fig_color_2}{HTML}{4C72B0}
\definecolor{fig_color_3}{HTML}{DD8452}
\definecolor{fig_color_4}{HTML}{55A868}
\definecolor{fig_color_5}{HTML}{C44E52}
\definecolor{fig_color_6}{HTML}{8172B2}
\begin{tikzpicture}[scale=0.280]
                  \path[use as bounding box] (1,-5.000) rectangle (10,5.714);
  \draw[step=1,fig_color_6!50,densely dotted] (-10,-10) grid (20,10);
  \draw [->, fig_color_2, very thick] (3.303,1.079) -- (3.364,1.055) -- (3.424,1.030) -- (3.485,1.006) -- (3.545,0.982) -- (3.606,0.958) -- (3.667,0.933) -- (3.727,0.909) -- (3.788,0.885) -- (3.848,0.861) -- (3.909,0.836) -- (3.970,0.812) -- (4.030,0.788) -- (4.091,0.764) -- (4.152,0.739) -- (4.212,0.715) -- (4.273,0.691) -- (4.333,0.667) -- (4.394,0.642) -- (4.455,0.618) -- (4.515,0.594) -- (4.576,0.570) -- (4.636,0.545) -- (4.697,0.521) -- (4.758,0.497) -- (4.818,0.473) -- (4.879,0.448) -- (4.939,0.424) -- (5.000,0.400) -- (5.061,0.376) -- (5.121,0.352) -- (5.182,0.327) -- (5.242,0.303) -- (5.303,0.279) -- (5.364,0.255) -- (5.424,0.230) -- (5.485,0.206) -- (5.545,0.182) -- (5.606,0.158) -- (5.667,0.133) -- (5.727,0.109) -- (5.788,0.085) -- (5.848,0.061) -- (5.909,0.036) -- (5.970,0.012) -- (6.030,-0.012) -- (6.091,-0.036) -- (6.152,-0.061) -- (6.212,-0.085) -- (6.273,-0.109) -- (6.333,-0.133) -- (6.394,-0.158) -- (6.455,-0.182) -- (6.515,-0.206) -- (6.576,-0.230) -- (6.636,-0.255) -- (6.697,-0.279) -- (6.758,-0.303) -- (6.818,-0.327) -- (6.879,-0.352) -- (6.939,-0.376) -- (7.000,-0.400) -- (7.061,-0.424) -- (7.121,-0.448) -- (7.182,-0.473) -- (7.242,-0.497) -- (7.303,-0.521) -- (7.364,-0.545) -- (7.424,-0.570) -- (7.485,-0.594) -- (7.545,-0.618) -- (7.606,-0.642) -- (7.667,-0.667) -- (7.727,-0.691) -- (7.788,-0.715) -- (7.848,-0.739) -- (7.909,-0.764) -- (7.970,-0.788) -- (8.030,-0.812) -- (8.091,-0.836) -- (8.152,-0.861) -- (8.212,-0.885) -- (8.273,-0.909) -- (8.333,-0.933) -- (8.394,-0.958) -- (8.455,-0.982) -- (8.515,-1.006) -- (8.576,-1.030) -- (8.636,-1.055) -- (8.697,-1.079);
  \draw [->, fig_color_3] (6.059,2.158) -- (6.064,2.133) -- (6.069,2.109) -- (6.073,2.085) -- (6.078,2.061) -- (6.081,2.036) -- (6.085,2.012) -- (6.088,1.988) -- (6.091,1.964) -- (6.093,1.939) -- (6.095,1.915) -- (6.097,1.891) -- (6.098,1.867) -- (6.099,1.842) -- (6.100,1.818) -- (6.100,1.794) -- (6.100,1.770) -- (6.099,1.745) -- (6.098,1.721) -- (6.096,1.697) -- (6.095,1.673) -- (6.092,1.648) -- (6.090,1.624) -- (6.087,1.600) -- (6.083,1.576) -- (6.080,1.552) -- (6.076,1.527) -- (6.071,1.503) -- (6.067,1.479) -- (6.062,1.455) -- (6.057,1.430) -- (6.051,1.406) -- (6.046,1.382) -- (6.040,1.358) -- (6.034,1.333) -- (6.028,1.309) -- (6.022,1.285) -- (6.016,1.261) -- (6.010,1.236) -- (6.003,1.212) -- (5.997,1.188) -- (5.990,1.164) -- (5.984,1.139) -- (5.978,1.115) -- (5.972,1.091) -- (5.966,1.067) -- (5.960,1.042) -- (5.954,1.018) -- (5.949,0.994) -- (5.943,0.970) -- (5.938,0.945) -- (5.933,0.921) -- (5.929,0.897) -- (5.924,0.873) -- (5.920,0.848) -- (5.917,0.824) -- (5.913,0.800) -- (5.910,0.776) -- (5.908,0.752) -- (5.905,0.727) -- (5.904,0.703) -- (5.902,0.679) -- (5.901,0.655) -- (5.900,0.630) -- (5.900,0.606) -- (5.900,0.582) -- (5.901,0.558) -- (5.902,0.533) -- (5.903,0.509) -- (5.905,0.485) -- (5.907,0.461) -- (5.909,0.436) -- (5.912,0.412) -- (5.915,0.388) -- (5.919,0.364) -- (5.922,0.339) -- (5.927,0.315) -- (5.931,0.291) -- (5.936,0.267) -- (5.941,0.242);
  \filldraw [fig_color_1] (0.000,0.000) circle (0.15) node[below]  {$U_{\tau, i-2}$};
  \filldraw [fig_color_1] (3.000,1.200) circle (0.15) node[above]  {$U_{\tau, i-1}$};
  \filldraw [fig_color_1] (6.000,0.000) circle (0.15) node[anchor=south west]  {$U_{\tau, i}$};
  \filldraw [fig_color_6] (9.000,-1.200) circle (0.1) node[anchor=north,yshift=-1mm]  {$K^\star(U_{\tau, i-1}, U_{\tau, i})$};
  \filldraw [fig_color_6] (6.000,2.400) circle (0.1);
  \draw [fig_color_2] (6.000, -0.000)  node[anchor=north east,yshift=1mm] {$\xi_\tau$};
\end{tikzpicture}
			\end{minipage}
		\end{adjustbox}
		\caption{The geodesic $\xi^{i,\tau}$ starts in $U_{\tau, i-1}$ and ends in $K^\star(U_{\tau, i-1}, U_{\tau, i})$. Note that $|(\xi^{i,\tau})'| = \frac{d(U_{\tau, i-1}, U_{\tau, i})}{\tau}$.}\label{fig:interp:2}
	\end{subfigure}
	\hfill
	\begin{subfigure}[t]{0.3\textwidth}
		\centering
		\begin{adjustbox}{clip,margin=0pt,trim=0 0 0 0,fbox,width=\linewidth} % trim behaves like in includegraphics
			\begin{minipage}[t]{\linewidth}
				\centering
				\definecolor{fig_color_1}{HTML}{000000}
\definecolor{fig_color_2}{HTML}{4C72B0}
\definecolor{fig_color_3}{HTML}{DD8452}
\definecolor{fig_color_4}{HTML}{55A868}
\definecolor{fig_color_5}{HTML}{C44E52}
\definecolor{fig_color_6}{HTML}{8172B2}
\begin{tikzpicture}[scale=0.230]
                  \path[use as bounding box] (4,-7.000) rectangle (10,6.043);
  \draw[step=1,fig_color_6!50,densely dotted] (-10,-10) grid (20,10);
  \draw [->, fig_color_2] (3.303,1.079) -- (3.364,1.055) -- (3.424,1.030) -- (3.485,1.006) -- (3.545,0.982) -- (3.606,0.958) -- (3.667,0.933) -- (3.727,0.909) -- (3.788,0.885) -- (3.848,0.861) -- (3.909,0.836) -- (3.970,0.812) -- (4.030,0.788) -- (4.091,0.764) -- (4.152,0.739) -- (4.212,0.715) -- (4.273,0.691) -- (4.333,0.667) -- (4.394,0.642) -- (4.455,0.618) -- (4.515,0.594) -- (4.576,0.570) -- (4.636,0.545) -- (4.697,0.521) -- (4.758,0.497) -- (4.818,0.473) -- (4.879,0.448) -- (4.939,0.424) -- (5.000,0.400) -- (5.061,0.376) -- (5.121,0.352) -- (5.182,0.327) -- (5.242,0.303) -- (5.303,0.279) -- (5.364,0.255) -- (5.424,0.230) -- (5.485,0.206) -- (5.545,0.182) -- (5.606,0.158) -- (5.667,0.133) -- (5.727,0.109) -- (5.788,0.085) -- (5.848,0.061) -- (5.909,0.036) -- (5.970,0.012) -- (6.030,-0.012) -- (6.091,-0.036) -- (6.152,-0.061) -- (6.212,-0.085) -- (6.273,-0.109) -- (6.333,-0.133) -- (6.394,-0.158) -- (6.455,-0.182) -- (6.515,-0.206) -- (6.576,-0.230) -- (6.636,-0.255) -- (6.697,-0.279) -- (6.758,-0.303) -- (6.818,-0.327) -- (6.879,-0.352) -- (6.939,-0.376) -- (7.000,-0.400) -- (7.061,-0.424) -- (7.121,-0.448) -- (7.182,-0.473) -- (7.242,-0.497) -- (7.303,-0.521) -- (7.364,-0.545) -- (7.424,-0.570) -- (7.485,-0.594) -- (7.545,-0.618) -- (7.606,-0.642) -- (7.667,-0.667) -- (7.727,-0.691) -- (7.788,-0.715) -- (7.848,-0.739) -- (7.909,-0.764) -- (7.970,-0.788) -- (8.030,-0.812) -- (8.091,-0.836) -- (8.152,-0.861) -- (8.212,-0.885) -- (8.273,-0.909) -- (8.333,-0.933) -- (8.394,-0.958) -- (8.455,-0.982) -- (8.515,-1.006) -- (8.576,-1.030) -- (8.636,-1.055) -- (8.697,-1.079);
  \draw [->, fig_color_3] (6.059,2.158) -- (6.064,2.133) -- (6.069,2.109) -- (6.073,2.085) -- (6.078,2.061) -- (6.081,2.036) -- (6.085,2.012) -- (6.088,1.988) -- (6.091,1.964) -- (6.093,1.939) -- (6.095,1.915) -- (6.097,1.891) -- (6.098,1.867) -- (6.099,1.842) -- (6.100,1.818) -- (6.100,1.794) -- (6.100,1.770) -- (6.099,1.745) -- (6.098,1.721) -- (6.096,1.697) -- (6.095,1.673) -- (6.092,1.648) -- (6.090,1.624) -- (6.087,1.600) -- (6.083,1.576) -- (6.080,1.552) -- (6.076,1.527) -- (6.071,1.503) -- (6.067,1.479) -- (6.062,1.455) -- (6.057,1.430) -- (6.051,1.406) -- (6.046,1.382) -- (6.040,1.358) -- (6.034,1.333) -- (6.028,1.309) -- (6.022,1.285) -- (6.016,1.261) -- (6.010,1.236) -- (6.003,1.212) -- (5.997,1.188) -- (5.990,1.164) -- (5.984,1.139) -- (5.978,1.115) -- (5.972,1.091) -- (5.966,1.067) -- (5.960,1.042) -- (5.954,1.018) -- (5.949,0.994) -- (5.943,0.970) -- (5.938,0.945) -- (5.933,0.921) -- (5.929,0.897) -- (5.924,0.873) -- (5.920,0.848) -- (5.917,0.824) -- (5.913,0.800) -- (5.910,0.776) -- (5.908,0.752) -- (5.905,0.727) -- (5.904,0.703) -- (5.902,0.679) -- (5.901,0.655) -- (5.900,0.630) -- (5.900,0.606) -- (5.900,0.582) -- (5.901,0.558) -- (5.902,0.533) -- (5.903,0.509) -- (5.905,0.485) -- (5.907,0.461) -- (5.909,0.436) -- (5.912,0.412) -- (5.915,0.388) -- (5.919,0.364) -- (5.922,0.339) -- (5.927,0.315) -- (5.931,0.291) -- (5.936,0.267) -- (5.941,0.242);
  \draw [->, fig_color_3, very thick] (7.879,-0.752) -- (7.955,-0.782) -- (8.030,-0.812) -- (8.106,-0.842) -- (8.182,-0.873) -- (8.258,-0.903) -- (8.333,-0.933) -- (8.409,-0.964) -- (8.485,-0.994) -- (8.561,-1.024) -- (8.636,-1.055) -- (8.712,-1.085) -- (8.788,-1.115) -- (8.864,-1.145) -- (8.939,-1.176) -- (9.015,-1.206) -- (9.091,-1.236) -- (9.167,-1.267) -- (9.242,-1.297) -- (9.318,-1.327) -- (9.394,-1.358) -- (9.470,-1.388) -- (9.545,-1.418) -- (9.621,-1.448) -- (9.697,-1.479) -- (9.773,-1.509) -- (9.848,-1.539) -- (9.924,-1.570) -- (10.000,-1.600) -- (10.076,-1.630) -- (10.152,-1.661) -- (10.227,-1.691) -- (10.303,-1.721) -- (10.379,-1.752) -- (10.455,-1.782) -- (10.530,-1.812) -- (10.606,-1.842) -- (10.682,-1.873) -- (10.758,-1.903) -- (10.833,-1.933) -- (10.909,-1.964) -- (10.985,-1.994) -- (11.061,-2.024) -- (11.136,-2.055) -- (11.212,-2.085) -- (11.288,-2.115) -- (11.364,-2.145) -- (11.439,-2.176) -- (11.515,-2.206) -- (11.591,-2.236) -- (11.667,-2.267) -- (11.742,-2.297) -- (11.818,-2.327) -- (11.894,-2.358) -- (11.970,-2.388) -- (12.045,-2.418) -- (12.121,-2.448) -- (12.197,-2.479) -- (12.273,-2.509) -- (12.348,-2.539) -- (12.424,-2.570) -- (12.500,-2.600) -- (12.576,-2.630) -- (12.652,-2.661) -- (12.727,-2.691) -- (12.803,-2.721) -- (12.879,-2.752) -- (12.955,-2.782) -- (13.030,-2.812) -- (13.106,-2.842) -- (13.182,-2.873) -- (13.258,-2.903) -- (13.333,-2.933) -- (13.409,-2.964) -- (13.485,-2.994) -- (13.561,-3.024) -- (13.636,-3.055) -- (13.712,-3.085) -- (13.788,-3.115) -- (13.864,-3.145) -- (13.939,-3.176) -- (14.015,-3.206) -- (14.091,-3.236) -- (14.167,-3.267) -- (14.242,-3.297) -- (14.318,-3.327) -- (14.394,-3.358) -- (14.470,-3.388) -- (14.545,-3.418) -- (14.621,-3.448);
  \filldraw [fig_color_1] (0.000,0.000) circle (0.15) node[below]  {$U_{\tau, i-2}$};
  \filldraw [fig_color_1] (3.000,1.200) circle (0.15) node[above]  {$U_{\tau, i-1}$};
  \filldraw [fig_color_1] (6.000,0.000) circle (0.15) node[anchor=south west]  {$U_{\tau, i}$};
  \filldraw [fig_color_6] (9.000,-1.200) circle (0.1) node[anchor=north,xshift=-6mm,yshift=-1mm]  {$K^\star(U_{\tau, i-1}, U_{\tau, i})$};
  \filldraw [fig_color_6] (6.000,2.400) circle (0.1);
  \draw [fig_color_3] (11.250, -2.100)  node[anchor=south west] {$\mathcal{V}^{\tau, \varepsilon}_{\tau (i + \frac{3}{4})}$};
\end{tikzpicture}
			\end{minipage}
		\end{adjustbox}
		\caption{The geodesic $\mathcal{V}_{\tau(i+\frac{3}{4})}^{\tau, \varepsilon}$ starts between the points $U_{\tau, i}$ and $K^\star(U_{\tau, i-1}, U_{\tau, i})$ and extends---for large enough $\varepsilon$---beyond $K^\star(U_{\tau, i-1}, U_{\tau, i})$. Note that $|(\mathcal{V}_{\tau(i+\frac{3}{4})}^{\tau,\varepsilon})'| = |(\xi^{i,\tau})'|$.}\label{fig:interp:3}
	\end{subfigure}

	\caption{The objects introduced in \cref{def:interpolations} in a prototypical example.} \label{fig:interpolation_multifig}
\end{figure}

\begin{figure}[ht]
	\centering
	\begin{subfigure}[t]{0.45\textwidth}
		\centering
		\begin{adjustbox}{clip,margin=0pt,trim=0 0 0 0,width=\linewidth} % trim behaves like in includegraphics
			\begin{minipage}[t]{\linewidth}
				\definecolor{fig_color_1}{HTML}{000000}
\definecolor{fig_color_2}{HTML}{4C72B0}
\definecolor{fig_color_3}{HTML}{DD8452}
\definecolor{fig_color_4}{HTML}{55A868}
\definecolor{fig_color_5}{HTML}{C44E52}
\definecolor{fig_color_6}{HTML}{8172B2}
\begin{tikzpicture}[xscale=2,yscale=1.4]
        \path[use as bounding box] (-1,-0.5) rectangle (1.3,1.75);
        \draw [fig_color_4, very thick] (0,0) -- (0.010,0.071) -- (0.020,0.100) -- (0.030,0.122) -- (0.040,0.141) -- (0.050,0.158) -- (0.060,0.173) -- (0.070,0.187) -- (0.080,0.200) -- (0.090,0.212) -- (0.100,0.224) -- (0.110,0.235) -- (0.120,0.245) -- (0.130,0.255) -- (0.140,0.265) -- (0.150,0.274) -- (0.160,0.283) -- (0.170,0.292) -- (0.180,0.300) -- (0.190,0.308) -- (0.200,0.316) -- (0.210,0.324) -- (0.220,0.332) -- (0.230,0.339) -- (0.240,0.346) -- (0.250,0.354) -- (0.260,0.361) -- (0.270,0.367) -- (0.280,0.374) -- (0.290,0.381) -- (0.300,0.387) -- (0.310,0.394) -- (0.320,0.400) -- (0.330,0.406) -- (0.340,0.412) -- (0.350,0.418) -- (0.360,0.424) -- (0.370,0.430) -- (0.380,0.436) -- (0.390,0.442) -- (0.400,0.447) -- (0.410,0.453) -- (0.420,0.458) -- (0.430,0.464) -- (0.440,0.469) -- (0.450,0.474) -- (0.460,0.480) -- (0.470,0.485) -- (0.480,0.490) -- (0.490,0.495) -- (0.500,0.500) -- (0.500,0.500) -- (0.510,0.505) -- (0.520,0.510) -- (0.530,0.515) -- (0.540,0.520) -- (0.550,0.526) -- (0.560,0.531) -- (0.570,0.536) -- (0.580,0.542) -- (0.590,0.547) -- (0.600,0.553) -- (0.610,0.558) -- (0.620,0.564) -- (0.630,0.570) -- (0.640,0.576) -- (0.650,0.582) -- (0.660,0.588) -- (0.670,0.594) -- (0.680,0.600) -- (0.690,0.606) -- (0.700,0.613) -- (0.710,0.619) -- (0.720,0.626) -- (0.730,0.633) -- (0.740,0.639) -- (0.750,0.646) -- (0.760,0.654) -- (0.770,0.661) -- (0.780,0.668) -- (0.790,0.676) -- (0.800,0.684) -- (0.810,0.692) -- (0.820,0.700) -- (0.830,0.708) -- (0.840,0.717) -- (0.850,0.726) -- (0.860,0.735) -- (0.870,0.745) -- (0.880,0.755) -- (0.890,0.765) -- (0.900,0.776) -- (0.910,0.788) -- (0.920,0.800) -- (0.930,0.813) -- (0.940,0.827) -- (0.950,0.842) -- (0.960,0.859) -- (0.970,0.878) -- (0.980,0.900) -- (0.990,0.929)-- (1,1);
        \draw[->] (-0.2, 0) -- (1.2,0) node[right] {$\delta$};
        \draw[->] (0, -0.2) -- (0,1.2) node[above,yshift=1mm] {$\psi(\delta)$};
        \draw [fig_color_6] (0.5,0.04) -- (0.5,-0.04) node[below,font=\small] {$\frac{\tau}{2}$};
        \draw [fig_color_6] (1,0.04) -- (1,-0.04) node[below,font=\small] {$\tau$};
        \draw [fig_color_6] (0.04,0.5) -- (-0.04,0.5) node[left,font=\small] {$\frac{\tau}{2}$};
        \draw [fig_color_6] (0.04,1) -- (-0.04,1) node[left,font=\small] {$\tau$};
        \draw[dotted] (0,1) -- (1, 1);
        \draw[dotted] (1,0) -- (1, 1);
        \draw[dotted] (0,0.5) -- (0.5, 0.5);
        \draw[dotted] (0.5,0) -- (0.5, 0.5);
        \end{tikzpicture}
        
				\centering
			\end{minipage}
		\end{adjustbox}
		\caption{$\psi$ is a strictly monotone function which bijectively maps $[0, \tau]$ to $[0, \tau]$.}
	\end{subfigure}
	\hfill
	\begin{subfigure}[t]{0.45\textwidth}
		\centering
		\begin{adjustbox}{clip,margin=0pt,trim=0 0 0 0,width=\linewidth} % trim behaves like in includegraphics
			\begin{minipage}[t]{\linewidth}
				\centering
				\definecolor{fig_color_1}{HTML}{000000}
\definecolor{fig_color_2}{HTML}{4C72B0}
\definecolor{fig_color_3}{HTML}{DD8452}
\definecolor{fig_color_4}{HTML}{55A868}
\definecolor{fig_color_5}{HTML}{C44E52}
\definecolor{fig_color_6}{HTML}{8172B2}
\begin{tikzpicture}[xscale=2,yscale=0.7]
        \path[use as bounding box] (0,-1) rectangle (1.3,3.5);
        \draw [fig_color_3, very thick] (0,0) -- (0.071,0.283) -- (0.100,0.400) -- (0.122,0.490) -- (0.141,0.566) -- (0.158,0.632) -- (0.173,0.693) -- (0.187,0.748) -- (0.200,0.800) -- (0.212,0.849) -- (0.224,0.894) -- (0.235,0.938) -- (0.245,0.980) -- (0.255,1.020) -- (0.265,1.058) -- (0.274,1.095) -- (0.283,1.131) -- (0.292,1.166) -- (0.300,1.200) -- (0.308,1.233) -- (0.316,1.265) -- (0.324,1.296) -- (0.332,1.327) -- (0.339,1.356) -- (0.346,1.386) -- (0.354,1.414) -- (0.361,1.442) -- (0.367,1.470) -- (0.374,1.497) -- (0.381,1.523) -- (0.387,1.549) -- (0.394,1.575) -- (0.400,1.600) -- (0.406,1.625) -- (0.412,1.649) -- (0.418,1.673) -- (0.424,1.697) -- (0.430,1.720) -- (0.436,1.744) -- (0.442,1.766) -- (0.447,1.789) -- (0.453,1.811) -- (0.458,1.833) -- (0.464,1.855) -- (0.469,1.876) -- (0.474,1.897) -- (0.480,1.918) -- (0.485,1.939) -- (0.490,1.960) -- (0.495,1.980) -- (0.500,2.000) -- (0.500,2.000) -- (0.505,1.980) -- (0.510,1.960) -- (0.515,1.939) -- (0.520,1.918) -- (0.526,1.897) -- (0.531,1.876) -- (0.536,1.855) -- (0.542,1.833) -- (0.547,1.811) -- (0.553,1.789) -- (0.558,1.766) -- (0.564,1.744) -- (0.570,1.720) -- (0.576,1.697) -- (0.582,1.673) -- (0.588,1.649) -- (0.594,1.625) -- (0.600,1.600) -- (0.606,1.575) -- (0.613,1.549) -- (0.619,1.523) -- (0.626,1.497) -- (0.633,1.470) -- (0.639,1.442) -- (0.646,1.414) -- (0.654,1.386) -- (0.661,1.356) -- (0.668,1.327) -- (0.676,1.296) -- (0.684,1.265) -- (0.692,1.233) -- (0.700,1.200) -- (0.708,1.166) -- (0.717,1.131) -- (0.726,1.095) -- (0.735,1.058) -- (0.745,1.020) -- (0.755,0.980) -- (0.765,0.938) -- (0.776,0.894) -- (0.788,0.849) -- (0.800,0.800) -- (0.813,0.748) -- (0.827,0.693) -- (0.842,0.632) -- (0.859,0.566) -- (0.878,0.490) -- (0.900,0.400) -- (0.929,0.283)-- (1,0);
        \draw[->] (-0.2, 0) -- (1.2,0) node[right] {$\psi(\delta)$};
        \draw[->] (0, -0.2) -- (0,2.3) node[above,yshift=1mm] {$\frac{1}{\psi'(\delta)}$};
        \draw [fig_color_6] (0.5,0.08) -- (0.5,-0.08) node[below,font=\small] {$\frac{\tau}{2}$};
        \draw [fig_color_6] (1,0.08) -- (1,-0.08) node[below,font=\small] {$\tau$};
        %\draw [fig_color_6] (0.04,1) -- (-0.04,1) node[left,font=\small] {$1$};
        \draw [fig_color_6] (0.04,2) -- (-0.04,2) node[left,font=\small] {$2$};
        \draw[dotted] (0,2) -- (1, 2);
        \draw[dotted] (1,0) -- (1, 2);
        \draw[dotted] (0.5,0) -- (0.5, 2);
        \end{tikzpicture}
        
			\end{minipage}
		\end{adjustbox}
		\caption{$\frac{1}{\psi'(\psi^{-1}(\delta))}$ is piecewise linear.}
	\end{subfigure}
	\hfill
	\caption{The function $\psi$ is designed to solve the ODE $\dot \psi = \frac{D}{\frac{\tau}{2} - |\frac{\tau}{2} - \psi|}$ for some constant $D$. This property will be useful in the proof of \cref{prop:discrete_newtonian_motion}.}\label{fig:phi_function}
\end{figure}
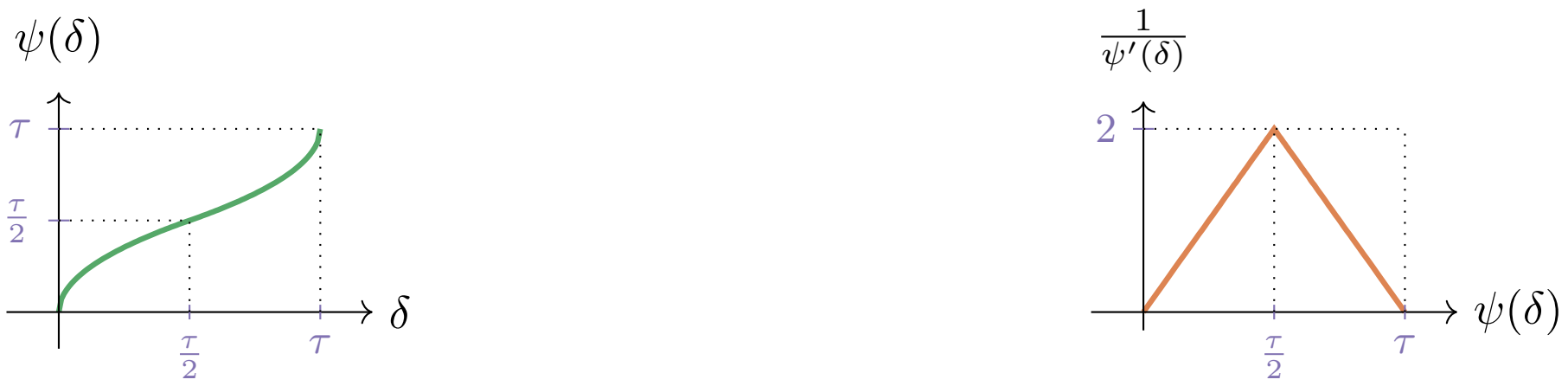

We start by defining two curves arising from a Newtonian sequence:
\begin{definition}[Interpolations]\label{def:interpolations}
	Let $(U_{\tau, k})_{k=-1}^N$ be a Newtonian sequence.
	Recalling the definition of $J_{\cdot}$ from \cref{eq:def_proximal}, we denote by $u$, $\tilde{u}$---see \cref{fig:interp:1}---any choice of curves satisfying
	\begin{align}
		u_\tau & \colon [0, \tau N] \to X,\qquad u_\tau(\tau i + \delta) = U_{\tau, i} \quad & i \in \{0,\ldots,N-1\}, \delta \in [0, \tau), \label{eq:piecewise_constant_interpolation} \\
		\tilde{u}_\tau & \colon [0, \tau N) \to X,\qquad \tilde{u}_\tau(\tau i + \delta) \in J_{\delta \cdot \tau}[K^\star(U_{\tau, i-1}, U_{\tau, i})], \quad & i \in \{0,\ldots,N-1\}, \delta \in [0, \tau). \label{eq:de_giorgi_interpolation}
	\end{align}
	We set $u_\tau(\tau N)=U_{\tau,N}$. Note that \cref{eq:piecewise_constant_interpolation} defines $u_\tau$ uniquely, while \cref{eq:de_giorgi_interpolation} may define $\tilde{u}_\tau$ non-uniquely.
	Furthermore, we temporarily set $\xi^{i, \tau}\colon [0, 2 \tau] \to X$ to be the geodesic with $\xi^{i, \tau}_0 = U_{\tau, i-1}$, $\xi^{i, \tau}_\tau = U_{\tau, i}$---see \cref{fig:interp:2}---and define, for $\varepsilon >0$, the functions $\mathcal{V}^{\tau, \varepsilon}, \tilde{\mathcal{V}}^{\tau, \varepsilon}, \bar{\mathcal{V}}^{\tau, \varepsilon} \colon [0, \tau N) \to \Tang[\varepsilon]$ by
	\begin{alignat}{3}
		\mathcal{V}^{\tau, \varepsilon}_{\tau i + \delta} & \coloneqq \bigg[\tau i + \delta + \sigma \mapsto \xi^{i, \tau}_{2\delta + \sigma}\bigg]\quad & i \in \{0,\ldots,N-1\}, \delta \in [0, \tau), \sigma \in [0, 2 (\tau - \delta)], \label{eq:velocity_interpolation} \\
		\mathcal{V}^{\tau, \varepsilon}_{t} & \in \Tang[t, \varepsilon], & t \in [0, \tau N),\nonumber \\
		\tilde{\mathcal{V}}^{\tau, \varepsilon}_{\tau i + \delta} &\coloneqq \bigg[\tau i + \delta + \sigma \mapsto \xi^{i, \tau}_{2\psi(\delta) + \sigma}\bigg] \quad & i \in \{0,\ldots,N-1\}, \delta \in [0, \tau), \sigma \in [0, 2 (\tau - \psi(\delta))], \label{eq:fancy_velocity_interpolation} \\
		\tilde{\mathcal{V}}^{\tau, \varepsilon}_{t} & \in \Tang[t, \varepsilon], & t \in [0, \tau N),\nonumber \\
		\bar{\mathcal{V}}^{\tau, \varepsilon}_{\tau i + \delta} & \coloneqq \bigg[\tau i + \delta + \sigma \mapsto \xi^{i, \tau}_\sigma\bigg] \quad & i \in \{0,\ldots,N-1\}, \delta \in [0, \tau), \sigma \in [0, 2 \tau], \label{eq:bar_velocity_interpolation} \\
		\bar{\mathcal{V}}^{\tau, \varepsilon}_{t} & \in \Tang[t, \varepsilon], & t \in [0, \tau N),\nonumber
	\end{alignat}
	where we implicitly restrict and extend the geodesics $\mathcal{V}^{\tau, \varepsilon}_t$, $\tilde{\mathcal{V}}^{\tau, \varepsilon}_t$ and $\bar{\mathcal{V}}^{\tau, \varepsilon}_t$ such that their domain is $[t, t + \varepsilon]$---see \cref{fig:interp:3}---, and where $\mathcal{V}^{\tau, \varepsilon}$, $\tilde{\mathcal{V}}^{\tau, \varepsilon}$ and $\bar{\mathcal{V}}^{\tau, \varepsilon}$ are only defined whenever such an extension is uniquely possible for all $t \in [0, \tau N)$; moreover, $\psi \in C^0([0, \tau)) \cap C^1((0, \tau))$---see \cref{fig:phi_function}---is defined as
	\begin{equation}\label{eq:def_phi}
		\psi(\delta) \coloneqq \begin{cases}
			\phantom{\tau - {}}\sqrt{\frac{\delta \cdot \tau}{2}}, \qquad & \delta \leq \frac{\tau}{2}, \\
			\tau - \sqrt{\frac{\tau^2 - \delta \cdot \tau}{2}}, \qquad & \delta > \frac{\tau}{2},
		\end{cases}
		\quad \text{thus} \quad \psi'(\delta) \coloneqq \begin{cases}
			\frac{\tau}{4} \frac{1}{\psi(\delta)}, \qquad & \delta \leq \frac{\tau}{2}, \\
			\frac{\tau}{4}\frac{1}{\tau - \psi(\delta)}, \qquad & \delta > \frac{\tau}{2}.
		\end{cases}
	\end{equation}
    Below, $N_\tau \colon X\times X\times X \to \R $ is defined as in \cref{eq:def_N}, and $d^+_{\cdot}, d^-_{\cdot}$ are defined as
	\begin{equation} \label{eq:def_d_pm}
		d^+_{r}(x) \coloneqq \sup_{y \in J_{r}[x]} d_X(x, y), \qquad d^-_{r}(x) \coloneqq \inf_{y \in J_{r}[x]} d_X(x, y),
	\end{equation}
	with $r>0$ and $J_{r}[x]$ introduced in \cref{eq:def_proximal}.
	Finally, we define the four functions \footnote{We let $\delta$ run through the left-open interval $(0, \tau]$ in the definition of $|G_\tau|$ to avoid division by zero; in the rest of this paper, only its equivalence class w.r.t.~equality almost everywhere matters  (and not its pointwise values).} for every $i=1,\ldots,N-1$ and every $\delta \in [0, \tau)$:
	\begin{align}
		|u_\tau^{(1)}| & \colon [0, \tau N) \to \R_+, \qquad |u_\tau^{(1)}|(\tau i + \delta) \coloneqq \frac{1}{\tau}d_X(U_{\tau, i}, U_{\tau, i+1}), \label{eq:velocity_norm_interpolation} \\
		|u_\tau^{(2)}| & \colon [0, \tau N) \to \R_+, \qquad |u_\tau^{(2)}|(\tau i + \delta) \coloneqq N_\tau(U_{\tau, i-1}, U_{\tau, i}, U_{\tau, i+1}), \label{eq:acceleration_norm_interpolation} \\
		d(Eu_\tau) & \colon [0, \tau N] \to \R, \qquad d(Eu_\tau)(\tau i + \delta) \coloneqq \frac{1}{\tau} (E(U_{\tau, i}) - E(U_{\tau, i-1})), \label{eq:energy_derivative_interpolation} \\
		|G_\tau| & \colon [0, \tau N] \to \R_+, \qquad |G_\tau|(\tau i + \delta) \coloneqq \frac{d^{+}_{\delta \ \tau}(K^\star(U_{\tau, i-1}, U_{\tau, i}))}{\delta \ \tau}, \label{eq:def_G_tau}
	\end{align}
	and we set 
    \begin{equation*}
        d(Eu_\tau)(\tau N) \coloneqq \frac{E(U_{\tau,N})-E(U_{\tau,N-1})}\tau, \qquad  |G_\tau|(0)=0.
    \end{equation*}
\end{definition}
\begin{remark}\label{rmk:varphi_tau}
We note that the function $\psi$ introduced in \cref{eq:def_phi} depends on $\tau$. Strictly speaking, it should therefore be denoted by $\psi^\tau$. However, to avoid overly cumbersome notation, we simply write $\psi$, unless it is necessary to emphasize the dependence on $\tau$.
\end{remark}

\newlength{\interpheight}
\setlength{\interpheight}{3cm}
\setlength{\fboxsep}{0pt}

We state in \cref{lem:helper_lemma_ambrosio} two auxiliary results from~\cite{ambrosioGradientFlows2008} that will be required in the following.

\begin{lemma}[\cite{ambrosioGradientFlows2008} Lemma 3.1.3, Lemma 3.1.4] \label{lem:helper_lemma_ambrosio}
	Let $E$ be proper and let $\varepsilon > 0$, $x \in X$. If $y \in J_\varepsilon[x]$, then $|\partial E|(y) \leq \frac{1}{\varepsilon} d_X(x, y)$.
	Furthermore, whenever $E$ fulfills \cref{ass:lsc}, $J_\tau[x] \neq \emptyset$ for all $0 \leq \tau \leq \varepsilon$, and $E(x) < +\infty$, we have
	\begin{equation}
		\frac{d_X^2(x, y)}{2 \varepsilon} + \int_0^\varepsilon \frac{(d^+_s(x))^2}{2 s^2} \, \diff s = E(x) - E(y).
	\end{equation}
\end{lemma}

We derive below a discrete version of \cref{eq:def_curve_newtonian_motion} for Newtonian sequences.

\begin{proposition}[Discrete identity, cf.\ \cref{eq:def_curve_newtonian_motion}]\label{prop:discrete_newtonian_motion}
	We assume that $(X, d_X)$ and $E \colon X \to \R \cup \{ +\infty\}$ fulfill \cref{ass:interceptor,ass:existence_derivative,ass:lsc}. Let $(U_{\tau, k})_{k=-1}^N$ be a Newtonian sequence, and let $|u_\tau^{(2)}|$, $d(Eu_\tau)$, and $\tilde{\mathcal{V}}^{\tau, \varepsilon}$ be as in \cref{eq:acceleration_norm_interpolation,eq:energy_derivative_interpolation,eq:fancy_velocity_interpolation}, respectively, where $\varepsilon > 0$ is assumed to be small enough such that $\tilde{\mathcal{V}}^{\tau, \varepsilon}$ exists. Setting $t_\tau^j \coloneqq \tau \cdot j$ for all $j \in \NN$, we have, for all $i, j \in \NN, 0 \leq i < j \leq N$,
	\begin{equation}\label{eq:discrete_newtonian_motion_inequality}
		d(Eu_\tau)(t_\tau^i) - d(Eu_\tau)(t_\tau^j) \geq \frac{1}{2}\int_{t_\tau^i}^{t_\tau^j}  |u_\tau^{(2)}|^2(r)\, \diff r - \int_{t_\tau^i}^{t_\tau^j} H^E_{0}(\tilde{\mathcal{V}}^{\tau, \varepsilon}_{r}) \, \diff r,
	\end{equation}
	If furthermore $J_s[K^\star(U_{\tau,k-1},U_{\tau,k})]\neq\emptyset$ for every $i\leq k<j$ and $0<s\leq\tau^2$ \footnote{E.g., for sufficiently small $\tau$, if $X$ and $E$ fulfill \cref{ass:coercivity,ass:lsc,ass:local_compactness,ass:interceptor}.}, let $\tilde u_\tau$ and $|G_\tau|$ be as in \cref{eq:de_giorgi_interpolation,eq:def_G_tau}. 
    Then, we have the identity
	\begin{equation} \label{eq:discrete_newtonian_motion}
		d(Eu_\tau)(t_\tau^i) - d(Eu_\tau)(t_\tau^j) = \frac{1}{2}\int_{t_\tau^i}^{t_\tau^j}  |u_\tau^{(2)}|^2(r)\, \diff r + \frac{1}{2} \int_{t_\tau^i}^{t_\tau^j} |G_\tau|^2(r)\, \diff r - \int_{t_\tau^i}^{t_\tau^j} H^E_{0}(\tilde{\mathcal{V}}^{\tau, \varepsilon}_{r}) \, \diff r,
	\end{equation}
	where $H_0^E$ has been introduced in \cref{eq:def_H0}.
	Under the same additional assumption, for almost all $t \in (t_\tau^i,t_\tau^j)$, we have
	\begin{equation}\label{eq:discrete_newtonian_motion_estimate_G}
		|\partial E|(\tilde{u}_\tau(t)) \leq |G_\tau|(t).
	\end{equation}
\end{proposition}
\begin{proof}
	We prove \cref{eq:discrete_newtonian_motion} and comment on how to obtain \cref{eq:discrete_newtonian_motion_inequality} at the end of the proof. By a telescoping sum argument, it suffices to show \cref{eq:discrete_newtonian_motion} for $j = i + 1$. To this end, let $i = 0, \ldots,N-1$. In the rest of this proof, we set $Z \coloneqq K^{\star}(U_{\tau, i-1}, U_{\tau, i})$. First, by applying \cref{lem:helper_lemma_ambrosio} with $x=Z$ and $y = U_{\tau, i+1}$---note that by \cref{eq:def_newtonian_sequence}, $U_{\tau, i+1} \in J_{\tau^2}(Z)$---, we have
	\begin{equation} \label{eq:discrete_newtonian_motion_proof_1}
		E(Z) - E(U_{\tau, i+1}) = \frac{d^2_X(Z, U_{\tau, i+1})}{2 \tau^2} + \int_0^{\tau^2} \frac{(d^+_s(Z))^2}{2 s^2} \, \diff s 
        = \frac{d^2_X(Z, U_{\tau, i+1})}{2 \tau^2} + \frac{\tau}{2} \int_0^{\tau} \frac{(d^+_{\delta \, \tau}(Z))^2}{(\delta \ \tau)^2} \, \diff \delta.
	\end{equation}
	Recalling \cref{eq:acceleration_norm_interpolation,eq:def_G_tau}, since for every $s \in [t_\tau^i, t_\tau^{i+1})$ we have that
    \begin{equation*}
        |u^{(2)}_\tau|(s) = N_\tau(U_{\tau, i-1}, U_{\tau, i}, U_{\tau, i+1}) = \frac{d_X(Z, U_{\tau, i+1})}{\tau^2}, \qquad  |G_\tau|(s) = \frac{d^+_{(s - t_\tau^i)\, \tau}(Z)}{(s - t_\tau^i) \ \tau},
    \end{equation*}
    we can divide \cref{eq:discrete_newtonian_motion_proof_1} by $\tau$ to obtain
	\begin{equation} \label{eq:proof_discrete_newtonian_motion_1}
		\frac{1}{\tau} (E(Z) - E(U_{\tau, i+1})) = \frac{1}{2}\int_{t_\tau^i}^{t_\tau^{i+1}}  |u_\tau^{(2)}|^2(r)\, \diff r + \frac{1}{2}\int_{t_\tau^i}^{t_\tau^{i+1}}  |G_\tau|^2(r)\, \diff r.
	\end{equation}
	Next, we analyze the last integral on the right hand side in \cref{eq:discrete_newtonian_motion}.
	We define the geodesic $\zeta\colon [t_\tau^i, t_\tau^{i} + \tau] \to X$ by $\zeta_{t_\tau^i} = U_{\tau, i-1}$, $\zeta_{t_\tau^{i+1}} = Z$. % and note by \cref{ass:lipschitz} that $\frac{\diff}{\diff s} E(\zeta_s)|_{s=t}$ exists for all $t \in [t_\tau^i, t_\tau^{i+1})$, is Lipschitz continuous---and thus, in particular, absolutely continuous.
	Furthermore, for $t \in (t_\tau^i,t_\tau^{i+1})$, setting $\xi = \mathcal{V}^{\tau, \varepsilon}_{t}$, by the uniqueness of geodesics following from the \EpsE-interceptor bound, we have that $\zeta_{t + \sigma} = \xi_{t + 2 \sigma}$ for sufficiently small $\sigma \geq 0$.
    Thus, for almost every $t \in (t_\tau^i, t_\tau^{i+1})$, it follows that
	\begin{equation*}
		\begin{split}
			\frac{\diff^2}{\diff s^2} E(\zeta_s)|_{s=t} & = \lim_{\sigma \searrow 0} \frac{E(\zeta_t) - 2 E(\zeta_{t+\sigma}) + E(\zeta_{t+2\sigma})}{\sigma^2} = \lim_{\sigma \searrow 0} \frac{E(\xi_{t}) - 2 E(\xi_{t + 2\sigma}) + E(\xi_{t + 4 \sigma})}{\sigma^2} \\
			& = 4 \lim_{\sigma \searrow 0}\frac{E(\xi_{t}) - 2 E(\xi_{t + \sigma}) + E(\xi_{t + 2\sigma})}{\sigma^2} = 4 H_0^E(\mathcal{V}^{\tau, \varepsilon}_{t}),
		\end{split}
	\end{equation*}
	where $\mathcal{V}^{\tau, \varepsilon}_{t}$ is defined as in \cref{def:interpolations} and \cref{ass:existence_derivative} guarantees existence of these quantities almost everywhere. We further see from \cref{eq:def_H0,eq:velocity_interpolation,eq:fancy_velocity_interpolation} that $H_0^E(\tilde{\mathcal{V}}^{\tau, \varepsilon}_{t_i + \delta}) = H_0^E(\mathcal{V}^{\tau, \varepsilon}_{t_i + \psi(\delta)})$ for $\delta \in [0, \tau)$.
    Hence, we calculate
	\begin{equation} \label{eq:proof_discrete_newtonian_motion_2}
		\begin{split}
			\int_{t_\tau^i}^{t_\tau^{i+1}} H^E_0(\tilde{\mathcal{V}}^{\tau, \varepsilon}_{r}) \, \diff r & = \int_0^\tau H_0^E(\mathcal{V}^{\tau, \varepsilon}_{t_\tau^i + \psi(\delta)}) \, \diff \delta = \int_0^\tau \frac{1}{\psi'(\psi^{-1}(\delta))} H_0^E(\mathcal{V}^{\tau, \varepsilon}_{t_{\tau}^i + \delta}) \, \diff \delta \\
			&= \int_0^\tau \frac{1}{4}\frac{1}{\psi'(\psi^{-1}(\delta))} \frac{\diff^2}{\diff s^2} E(\zeta_s)\Big|_{s=t_{\tau}^i + \delta} \, \diff \delta = \frac{1}{\tau}\Big(E(U_{\tau, i-1}) - 2 E(U_{\tau, i}) + E(Z)\Big),
		\end{split}
	\end{equation}
	where we used \cref{eq:def_phi} to see that the second distributional derivative of the function
	\begin{equation*}
		\delta \mapsto \left\{
	\begin{array}{ll}
		\frac{1}{\psi'(\psi^{-1}(\delta))}, \quad & \text{if } \delta \in [0, \tau] \\
		0 \quad & \text{otherwise}
	\end{array}
	\right\} = \begin{cases}
		\frac{4}{\tau} \delta, \quad &\text{if } \delta \in [0, \frac{\tau}{2}], \\
		\frac{4}{\tau} (\tau - \delta) \quad &\text{if } \delta \in (\frac{\tau}{2}, \tau], \\
		0 \quad &\text{otherwise},
		\end{cases}
	\end{equation*}
	is given by $\frac{4}{\tau} (\boldsymbol{\delta}_{0} - 2 \boldsymbol{\delta}_{\frac{\tau}{2}} + \boldsymbol{\delta}_{\tau})$ (see also \cref{fig:phi_function}), where $\boldsymbol{\delta}$ denotes the Dirac delta distribution. Here, \cref{ass:existence_derivative} justifies the integration by parts.
	Subtracting \cref{eq:proof_discrete_newtonian_motion_2} from \cref{eq:proof_discrete_newtonian_motion_1} yields
	\begin{equation*}
		\begin{split}
			d(Eu_\tau)(t_\tau^i) - d(Eu_\tau)(t_\tau^{i+1}) & = \frac{1}{\tau} \big(E(U_{\tau, i}) - E(U_{\tau, i-1})\big) - \frac{1}{\tau}\big(E(U_{\tau, i+1}) - E(U_{\tau, i})\big) \\
			& = \int_{t_\tau^i}^{t_\tau^{i+1}} \frac{1}{2} |u_\tau^{(2)}|^2(r)\, \diff r + \int_{t_\tau^i}^{t_\tau^{i+1}} \frac{1}{2} |G_\tau|^2(r)\, \diff r - \int_{t_\tau^i}^{t_\tau^{i+1}} H^E_{0}(\tilde{\mathcal{V}}^{\tau, \varepsilon}_{r}) \, \diff r,
		\end{split}
	\end{equation*}
	which concludes the proof of \cref{eq:discrete_newtonian_motion}. \cref{eq:discrete_newtonian_motion_estimate_G} is a direct application of \cref{lem:helper_lemma_ambrosio}. \\
	To obtain \cref{eq:discrete_newtonian_motion_inequality}, one simply replaces \cref{eq:discrete_newtonian_motion_proof_1} by the inequality
	\begin{equation}
		E(Z) - E(U_{\tau, i+1}) \geq \frac{d^2_X(Z, U_{\tau, i+1})}{2 \tau^2},
	\end{equation}
	and repeats the above calculations mutatis mutandis.
\end{proof}

In view of the limiting construction detailed in \cref{sec:compactness_limiting}, we need to establish some \emph{a priori} bounds on discrete Newtonian sequences.

\begin{lemma}[A priori bounds I]\label{lem:apriori_bounds}
	We assume that $(X, d_X)$ and $E \colon X \to \R \cup \{ +\infty\}$ fulfill \cref{ass:interceptor,ass:coercivity,ass:lipschitz,ass:lsc}. 
    Then, for every $R, T > 0$, there exist a constant\footnote{Recall that $\tau_\star,x_\star$ are introduced in \cref{ass:coercivity}. Moreover, for brevity, we choose not to make the dependence on the fixed constants $L$ (coming from \cref{ass:lipschitz}) and $M^\star$ (coming from \cref{ass:coercivity}) explicit.} $C = C(x_\star, \tau_\star, R, T) \geq R$  and $\underline{\tau} > 0$ such that for every Newtonian sequence $(U_{\tau, k})_{k=-1}^N$ satisfying
	\begin{equation*}
		E(U_{\tau, 0}) \leq R, \quad d^2_X (U_{\tau, 0}, x_\star) \leq R,
        \quad N \, \tau \leq T,
        \quad d(Eu_\tau)(0) \leq R,
        \quad\frac{d_X\big(U_{\tau, -1}, U_{\tau, 0}\big)}{\tau} \leq R, \quad\tau \leq \underline{\tau},
	\end{equation*}
	where $d(Eu_{\tau})$ is defined as in \cref{eq:energy_derivative_interpolation}, we have
	\begin{gather}
		d^2_X(U_{\tau, n}, x_\star) \leq C,  \qquad\qquad \sum_{k=1}^n \frac{d^2_X(U_{\tau, k}, U_{\tau, k-1})}{\tau} \leq C, \label{eq:apriori_prop_1}
	\end{gather}
    for every $n =1, \ldots , N$, 
	where $t_\tau^n \coloneqq n \, \tau$.\\ 
    Furthermore, if $(X,d_X)$ and $E$ fulfill  \cref{ass:interceptor,ass:coercivity,ass:lipschitz,ass:lsc,ass:local_compactness}, we have that any such Newtonian sequence can be extended to a Newtonian sequence $(U_{\tau, k})_{k=-1}^{N'}$ with $\tau N' \geq T$. 
    In this case, we also have 
	\begin{gather}
		\int_0^{t^n_\tau} |u_\tau^{(2)}|^2(t) \diff t \leq C, \qquad \qquad 
        \frac{d_X(U_{\tau, n}, U_{\tau, n-1})}{\tau} \leq C, \label{eq:apriori_prop_1_5}\\
		 d(Eu_\tau)(t^n_\tau) \leq C, \qquad  \qquad d(Eu_\tau)(t^n_\tau) - d(Eu_\tau)(t^n_\tau-\tau) \leq C \, \tau. \label{eq:apriori_prop_2}
	\end{gather}
    for every $n=1, \ldots, \lfloor N'/2\rfloor$.
\end{lemma}
\begin{proof}
	We fix such a Newtonian sequence $(U_{\tau, k})_{k=-1}^N$ and set, for all $k = 0, \ldots , N$,
	\begin{equation*}
		\begin{split}
			\FDiff_{k} & \coloneqq d_X(U_{\tau, k}, U_{\tau, k-1}), \\
			\FDiffDiff_k & \coloneqq d_X(K^\star(U_{\tau, k-2}, U_{\tau, k-1}), U_{\tau, k}), \qquad 1\leq k\leq N, \\
			E_k & \coloneqq E(U_{\tau, k}).
		\end{split}
	\end{equation*}
	For $1\leq k\leq N$, we have $\FDiff_k=\tau|u_\tau^{(1)}|((k-1)\tau)$ and $\FDiffDiff_k=\tau^2|u_\tau^{(2)}|((k-1)\tau)$ (see \cref{eq:velocity_norm_interpolation,eq:acceleration_norm_interpolation}, respectively). 
    We set $\widehat{R} \coloneqq R+LR^2$ (with $L>0$ given by \cref{ass:lipschitz}), and we restrict our focus to $\tau \leq 1$.
    Let us start by collecting some inequalities. 
    First, by the triangle inequality we have that $\FDiff_k - \FDiff_{k-1} \leq \FDiffDiff_k$ for all $k = 1, \ldots , N$, and hence also 
    \begin{equation*}
        \FDiff_k^2 - \FDiff_{k-1}^2 \leq  2 \FDiff_k^2 - 2\FDiff_k \FDiff_{k-1} \leq 2 \FDiff_k \FDiffDiff_{k}.
    \end{equation*}
    Thus, we get
	\begin{equation}\label{eq:apriori_inequality_1}
		\frac{1}{\tau}\big(\FDiff^2_k - \FDiff^2_0\big) = \frac{1}{\tau}\Big(\sum_{j=1}^k \big(\FDiff_j^2 - \FDiff_{j-1}^2\big) \Big)\leq \frac{2}{\tau}\sum_{j=1}^k \FDiff_j \FDiffDiff_j \leq \frac{2}{\tau_\star} \sum_{j=1}^k \FDiff_j^2 + \frac{\tau_\star}{2 \tau^2} \sum_{j=1}^k \FDiffDiff_j^2,
	\end{equation}
	where $\tau_\star$ is as in \cref{ass:coercivity}.
	Second, we deduce from \cref{eq:discrete_newtonian_motion_inequality}---using the Lipschitz condition \cref{ass:lipschitz} as in \cref{rmk:rescaling_lipschitz}---that, for 
    $$\varepsilon \coloneqq \frac{\epsilon \tau}{1 + \max_{0\leq k\leq N} d_X(U_{\tau, k}, U_{\tau, k-1})},$$
    we have
	\begin{equation}
		\frac{1}{2} \tau \sum_{j=1}^k \frac{\FDiffDiff_j^2}{\tau^4} = \frac{1}{2} \int_0^{t_\tau^k} |u_\tau^{(2)}|^2(t) \, \diff t \leq \underbrace{d(Eu_\tau)(0)}_{\leq R} - \underbrace{d(Eu_\tau)(t_\tau^k)}_{= \frac{1}{\tau} \big(E_k - E_{k-1}\big)} + \int_0^{t_\tau^k} \underbrace{H^E_0(\tilde{\mathcal{V}}^{\tau, \varepsilon}_{t})}_{\leq L\FDiff_{\lfloor t/\tau\rfloor}^2/\tau^2\text{ a.e.}} \, \diff t,
	\end{equation}
	and hence that
	\begin{equation}\label{eq:apriori_inequality_2}
		\begin{split}
			\frac{1}{2 \tau^2} \sum_{j=1}^k \FDiffDiff_j^2 & \leq \tau R + E_{k-1} - E_k + L\sum_{j=0}^{k-1}\FDiff_j^2 \leq \tau\widehat{R} + E_{k-1} - E_k + L \sum_{j=1}^{k} \FDiff_j^2.
		\end{split}
	\end{equation}
	Third, from \cref{ass:coercivity}, we get that, for some $x_\star \in X$,
	\begin{equation}\label{eq:apriori_inequality_3}
		E_k \geq - \frac{1}{2\tau_\star} d^2_X(U_{\tau, k}, x_\star) + M^\star.
	\end{equation}
	Using \cref{eq:apriori_inequality_1,eq:apriori_inequality_2,eq:apriori_inequality_3}, we calculate
	\begin{equation}
		\begin{split}
			G(n) & \coloneqq \frac{1}{\tau} \sum_{k=1}^{n} \FDiff_k^2 \leq \frac{n}{\tau} \FDiff_0^2 + \mathlarger{\sum_{k=1}^n} \bigg(\frac{2}{\tau_\star} \sum_{j=1}^k \FDiff_j^2 + \frac{\tau_\star}{2 \tau^2} \sum_{j=1}^k \FDiffDiff_j^2\bigg) \\
			& \leq \frac{n}{\tau} \FDiff_0^2 + \mathlarger{\sum_{k=1}^n} \bigg(\frac{2}{\tau_\star} \sum_{j=1}^k \FDiff_j^2 + \tau_\star \tau\widehat{R} + \tau_\star E_{k-1} - \tau_\star E_k + L \tau_\star \sum_{j=1}^{k} \FDiff_j^2\bigg) \\
			& = \frac{n}{\tau} \FDiff_0^2 + n\tau \tau_\star\widehat{R} + \tau_\star E_0 - \tau_\star E_n + \tau \Big(\frac{2}{\tau_\star} + L \tau_\star\Big) \cdot \mathlarger{\sum_{k=1}^n} \bigg(\frac{1}{\tau}\sum_{j=1}^{k} \FDiff_j^2\bigg) \\
			& \leq \underbrace{T R^2 + \tau_\star(T \widehat{R}+R) - \tau_\star M^\star}_{\eqqcolon M'} + \frac{1}{2} d^2_X(U_{\tau, n}, x_\star) + \tau \underbrace{\Big(\frac{2}{\tau_\star} + L \tau_\star\Big)}_{\eqqcolon A'} \cdot \sum_{k=1}^n G(k).
		\end{split}
	\end{equation}
	On the other hand, arguing similarly as in \cref{eq:apriori_inequality_1}, we have
	\begin{equation}
		\begin{split}
			H(n) & \coloneqq d^2_X(U_{\tau, n}, x_\star) = d^2_X(U_{\tau, 0}, x_\star) + \sum_{k=1}^n \big(d^2_X(U_{\tau, k}, x_\star) - d^2_X(U_{\tau, k-1}, x_\star)\big) \\
			& \leq R + 2 \sum_{k=1}^n d_X(U_{\tau, k}, x_\star) \cdot \FDiff_k \leq R + 2\tau \sum_{k=1}^n d^2_X(U_{\tau, k}, x_\star) + \frac{1}{2\tau} \sum_{k=1}^n \FDiff_k^2 \\
			& = R + 2\tau \sum_{k=1}^n H(k) + \frac{1}{2}G(n).
		\end{split}
	\end{equation}
	We will finish by using a Gronwall-like argument of $G$ and $H$. Setting $M \coloneqq \max\{M', R\}$, $A \coloneqq \max\{A', 2\}$, and $F\colon [N] \mapsto \R^2$ to $F(n) \coloneqq \big(G(n), H(n) \big)$, we obtain the vector inequality
	\begin{equation*}
		\Bigg[F(n)\Bigg] \leq \Bigg[\begin{pmatrix}
			M \\
			M
		\end{pmatrix} + \begin{pmatrix}
			0 & \frac{1}{2} \\
			\frac{1}{2} & 0
		\end{pmatrix} F(n) + \tau A \sum_{k=1}^n F(k)\Bigg],  %\quad q \in \{1, 2\},
	\end{equation*}
	which, after taking the supremum norm and rearranging, yields 
    $$\|F(n)\|_\infty \leq 2 M + 2 \tau A \sum_{k=1}^n \|F(k)\|_\infty.$$
	Finally, applying a discrete version of a Gronwall Lemma (see, e.g., \cite[Lemma 3.2.4]{ambrosioGradientFlows2008}), we obtain, for $\tau \leq \underline{\tau} \coloneqq \min\{1,\frac{1}{4 A}\}$,
	\begin{equation} \label{eq:constant_1_priori_1}
		\|F(n)\|_\infty \leq \frac{2M}{1 - 2 \tau A} \cdot e^{\frac{2\tau A}{1- 2 \tau A} n} \leq 4 M e^{4 A T} \coloneqq C_1,
	\end{equation}
	where $C_1$ depends only on $x_\star, \tau_\star, R, T, L, M^\star$. This gives us uniform bounds on $G(n)$ and $H(n)$, which immediately yield the two inequalities in \cref{eq:apriori_prop_1} with the constant $C_1$. \\
	Before we show \cref{eq:apriori_prop_1_5,eq:apriori_prop_2}, we establish the extensibility claim at the end of \cref{lem:apriori_bounds} for
    \begin{equation*}
        \underline{\tau} \coloneqq \min\left\{ \frac{1}{4A}, \frac{\epsilon^2}{9 C_1}, \frac{\epsilon}{3R}, \frac{\sqrt{\tau_\star}}{2}, 1 \right\}
    \end{equation*} 
    To this end, for $\tau \leq \underline{\tau}$, we can always extend $(U_{\tau, k})_{k=-1}^{N}$ to $(U_{\tau, k})_{k=-1}^{N+1}$ by choosing 
    \begin{equation*}
        U_{\tau, N+1} \in J_{\tau^2}[K^\star(U_{\tau, N-1}, U_{\tau, N})];
    \end{equation*}
    note that the latter set is nonempty and well-defined by \cref{ass:interceptor,ass:lsc,ass:coercivity,ass:local_compactness} whenever $d_X(U_{\tau, N-1}, U_{\tau, N}) \leq \frac{\epsilon}{3}$. 
    For $N=0$, the initial bound and the choice of $\underline{\tau}$ give
    \begin{equation*}
        d_X(U_{\tau,-1},U_{\tau,0})\leq \tau R\leq \underline{\tau} R  \leq \frac\epsilon3 .
    \end{equation*}
    Thus, for $N\geq1$, the only way for the extensibility claim to fail is if, for some Newtonian sequence $(U_{\tau, k})_{k=-1}^{N}$ with $\tau \leq \underline{\tau}$, we have that $d_X\big(U_{\tau, N-1}, U_{\tau, N} \big) > {\epsilon}/{3}$ while $\tau N \leq T$. However in this case, we have, according to \cref{eq:apriori_prop_1},
		\begin{equation*}
			\sum_{k=1}^{N} \frac{d^2_X(U_{\tau, k}, U_{\tau, k-1})}{\tau} \leq C_1,
		\end{equation*}
		and in particular $\frac{d^2_X(U_{\tau, N-1}, U_{\tau, N})}{\tau} \leq C_1$ (recall that the constant $C_1$ was obtained in \cref{eq:constant_1_priori_1}).
		Rewriting, we get
		\begin{equation*}
			d_X(U_{\tau, N-1}, U_{\tau, N}) \leq \sqrt{C_1 \tau} \leq \sqrt{C_1 \underline{\tau}} \leq \frac{\epsilon}{3},
		\end{equation*}
		so this failure case cannot occur.
	
	To show \cref{eq:apriori_prop_1_5,eq:apriori_prop_2} we will, by virtue of this extensibility result\footnote{We can arrive at such a sequence by invoking the extensibility result with $2T$ in place of $T$.}, extend the Newtonian sequence to have length $N' \geq N$, such that $n \leq \frac{N'}{2}$ and $N' \geq \frac{T}{\tau}$.
	Let us show the bound on the integral of $|u_\tau^{(2)}|^2$; note that we can safely assume that $t^n_\tau \geq \frac{T}{4}$:
	\begin{equation*}
		\begin{split}
			\int_0^{t^n_\tau} |u_\tau^{(2)}|^2(t) \diff t & = \tau \sum_{k=1}^n \frac{\FDiffDiff_k^2}{\tau^4} \leq \frac{1}{n} \tau \mathlarger{\sum_{i=n}^{2n}}\bigg(\sum_{k=1}^i \frac{\FDiffDiff_k^2}{\tau^4} \bigg) \leq \frac{\tau^2}{t^n_\tau} \mathlarger{\sum_{i=1}^{2n}} \bigg(\sum_{k=1}^i \frac{\FDiffDiff_k^2}{\tau^4} \bigg) \\
			& \overset{(*)}{\leq} \frac{1}{t^n_\tau} \bigg(4 \tau\widehat{R} n + 2 \cdot (E_0 - E_{2n}) + 2 L \tau \sum_{i=1}^{2n} G(i)\bigg) \\
			& \leq \frac{1}{t^n_\tau} \bigg(4\widehat{R} t^n_\tau + 2 R + \frac{1}{\tau_\star} C_1 - 2 M^\star + 4 L t^n_\tau C_1\bigg) \\
			& \leq \frac{4}{T} \bigg(4\widehat{R} T + 2 R + \frac{1}{\tau_\star} C_1 + 2|M^\star| + 4 L T C_1\bigg) \eqqcolon C_2,
		\end{split}
	\end{equation*}
	where we repeatedly used \cref{eq:apriori_inequality_2} in $(*)$.
	Next, we show the uniform bound on $|u^{(1)}_\tau|$, by using \cref{eq:apriori_inequality_1}. Here, we use again that $n \leq \lfloor \frac{N'}{2} \rfloor$:
	\begin{align*}
		\frac{\FDiff^2_n}{\tau^2} &\leq R^2 + \frac{2}{\tau \tau_\star} \sum_{j=1}^n \FDiff_j^2 + \frac{\tau_\star}{2 \tau^3} \sum_{j=1}^n \FDiffDiff_j^2 = R^2 + \frac{2}{\tau_\star} G(n) + \frac{\tau_\star}{2} \int_0^{t_\tau^n} |u_\tau^{(2)}|^2(t) \diff t  \\
		&\leq R^2 + \frac{2}{\tau_\star} C_1 + \frac{\tau_\star}{2} C_2 \eqqcolon C_3.
	\end{align*}
	Finally, we show the bounds on $d(Eu_\tau)(s)$ and its increments for $1\leq n\leq\lfloor N'/2\rfloor$; the initial incoming speed is bounded by $R\leq\sqrt{C_3}$. First, using \cref{eq:discrete_newtonian_motion} and ignoring some non-positive terms, for 
    $$\varepsilon \coloneqq \frac{\epsilon \tau}{1 + \max_{0\leq k\leq N} d(U_{\tau, k}, U_{\tau, k-1})},$$
    we get the estimate
	\begin{equation*}
		d(Eu_\tau)(t^n_\tau) \leq d(Eu_\tau)(0) + \int_0^{t_\tau^n} H_0^E(\tilde{\mathcal{V}}^{\tau, \varepsilon}_{t}) \, \diff t \leq R + \frac{L}{\tau}\sum_{j=0}^{n-1}\FDiff_j^2 \leq R + L T C_3 \eqqcolon C_4.
	\end{equation*}
	Second, again by using \cref{eq:discrete_newtonian_motion}, we have
	\begin{equation*}
		\begin{split}
			d(Eu_\tau)(t^n_\tau) - d(Eu_\tau)(t^n_\tau - \tau) & \leq \int_{t_\tau^{n-1}}^{t_\tau^n} H_0^E(\tilde{\mathcal{V}}^{\tau, \varepsilon}_{t}) \diff t \leq \frac{L}{\tau}\FDiff_{n-1}^2 \leq \tau\underbrace{LC_3}_{\eqqcolon C_5}.
		\end{split}
	\end{equation*}
	Thus, setting $C \coloneqq \max\{C_1, C_2, \sqrt{C_3}, C_4, C_5, R\}$ concludes the proof.
\end{proof}

The next result provides other useful estimates. 

\begin{lemma}[A priori bounds II] \label{lem:apriori_bounds_2}
	We assume that $(X, d_X)$ and $E \colon X \to \R \cup \{ +\infty\}$ fulfill \cref{ass:interceptor,ass:coercivity,ass:lipschitz,ass:lsc}.
    Then, for any $C > 0$, $N \in \NN$ and any Newtonian sequence $(U_{\tau, k})_{k=-1}^N$ such that
	\begin{equation*}
		\frac{d_X(U_{\tau, i}, U_{\tau, i-1})}{\tau} \leq C  \qquad \forall \  i \in \{0, \dots, N\},
	\end{equation*}
	we have that, setting $\varepsilon \coloneqq \frac{1}{C}\, \epsilon$ and assuming $2\tau\leq\varepsilon$,
	\begin{gather}
		d_X\big(u_\tau(t), \tilde{u}_\tau(t)\big) \leq 3 \tau C, \qquad
        d_{\Tang[t, \varepsilon]}(\mathcal{V}^{\tau, \varepsilon}_t, \tilde{\mathcal{V}}^{\tau, \varepsilon}_t) \leq 2 \tau C, 
        \qquad d_{\Tang[t, \varepsilon]}(\mathcal{V}^{\tau, \varepsilon}_t, \bar{\mathcal{V}}^{\tau, \varepsilon}_t) \leq 2 \tau C, \label{eq:apriori_bounds_2}\\
		\left|d(E u_\tau)(t) - \frac{\diff}{\diff \sigma} E(\bar{\mathcal{V}}_t^{\tau, \varepsilon}(\sigma))\Big|_{\sigma = t}\right| \leq \tau C^2 L, \qquad \left|d(E u_\tau)(t) - \frac{\diff}{\diff \sigma} E(\tilde{\mathcal{V}}_t^{\tau, \varepsilon}(\sigma))\Big|_{\sigma = t}\right| \leq 2\tau C^2 L, \label{eq:apriori_bounds_2_2}\\
		\left|d(E u_\tau)(t) - \frac{\diff}{\diff \sigma} E(\mathcal{V}_t^{\tau, \varepsilon}(\sigma))\Big|_{\sigma = t}\right| \leq 2\tau C^2 L, \label{eq:apriori_bounds_2_3}
	\end{gather}
    for every $t \in [0, \tau N)$, where $u_\tau$, $\tilde{u}_\tau$, $\mathcal{V}^{\tau, \varepsilon}$, $\tilde{\mathcal{V}}^{\tau, \varepsilon}$, $\bar{\mathcal{V}}^{\tau, \varepsilon}$ and $d(E u_\tau)$ are defined as in \cref{eq:piecewise_constant_interpolation,eq:de_giorgi_interpolation,eq:velocity_interpolation,eq:fancy_velocity_interpolation,eq:bar_velocity_interpolation,eq:energy_derivative_interpolation}, and $d_{\Tang[t, \varepsilon]}$ is defined as in \cref{def:space_of_geodesics}.
\end{lemma}
\begin{proof}
	For $i \in \{0, \dots, N-1\}$, $t^i_\tau \coloneqq i \tau$ and $\delta \in [0, \tau)$, and defining $Z \coloneqq K^\star(U_{\tau, i-1}, U_{\tau, i})$, we have
	\begin{equation*}
		\begin{split}
		d_X\big(u_\tau(t^i_\tau + \delta), \tilde{u}_\tau(t^i_\tau + \delta)\big) &= d_X(U_{\tau, i}, J_{\delta \cdot \tau}[Z]) \leq d_X(U_{\tau, i}, Z) + d_X(Z, J_{\delta \cdot \tau}[Z]) \\
		&\leq C \tau + d^{-}_{\tau^2}(Z) \leq C \tau + d_X(Z, U_{\tau, i+1}) \\
		&\leq C \tau + d_X(Z, U_{\tau, i}) + d_X(U_{\tau, i}, U_{\tau, i+1}) \\
		&\leq 3 \tau C,
		\end{split}
	\end{equation*}
	where $d^{\pm}_{\cdot}$ is defined in~\cref{eq:def_d_pm}, and where we used that $d_X(Z, J_{\delta \cdot \tau}[Z]) \leq d^{+}_{\delta\tau}(Z) \leq d^{-}_{\tau^2}(Z)$ for every $\delta < \tau$ (see~\cite[Lemma~3.1.2]{ambrosioGradientFlows2008}, and~\cite{Tradeoff2026} for extensions). This proves the first bound in~\cref{eq:apriori_bounds_2}. 
    
    For the other two bounds, let us first show the following claim: Given two geodesics $\xi\colon I \to X$ and $\zeta\colon I \to X$, where $I$ is a closed interval in $\R$ of length $\varepsilon$, if there exists an interval $J \subseteq I$ of positive length and some $r \in \R$ such that $J + r \subseteq I$ and such that for all $\sigma \in J$, $\xi_\sigma = \zeta_{\sigma+r}$, we have that $|\xi'| = |\zeta'|$ and
	\begin{equation} \label{eq:geodesic_interpolation_bound}
		d_{\Tang[\varepsilon]}(\xi, \zeta) = |r| \ |\xi'|.
	\end{equation}
	First, note that for $\sigma_1, \sigma_2 \in J$ such that $\sigma_1 \neq \sigma_2$, we have that 
    \begin{equation*}
        |\xi'| = \frac{d_X(\xi_{\sigma_1}, \xi_{\sigma_2})}{|\sigma_1 - \sigma_2|} = \frac{d_X(\zeta_{\sigma_1 + r}, \zeta_{\sigma_2 + r})}{|(\sigma_1 + r) - (\sigma_2 + r)|} = |\zeta'|.
    \end{equation*}
	Next, let us note that \EpsE-interceptor bound implies that the set $J' \coloneqq \{\sigma \in I \cap I - r \mid \xi_\sigma = \zeta_{\sigma + r}\}$ is open in $I \cap I - r$; since it is also closed (by continuity), we have $J' = I \cap (I - r)$.
	We choose some arbitrary $t \in I$ and show that $d_X(\xi_t, \zeta_t) \leq |r| \ |\xi'|$. We can---by possibly switching the roles of $\xi$ and $\zeta$---assume that $r \geq 0$. 
    Then, since $J' = I \cap (I - r)$ is nonempty, there exists $r' \in [0, r]$ such that $t - r' \in J'$. Thus, we have 
    \begin{equation*}
        d_X(\xi_t, \zeta_t) \leq d_X(\xi_t, \xi_{t - r'}) + d_X(\zeta_{t - r' + r}, \zeta_{t}) = |r| \ |\xi'|.
    \end{equation*} 
    Running through all $t \in I$, this gives one inequality in \cref{eq:geodesic_interpolation_bound}; for the other, choose $t \in J'$ to see that $d_X(\xi_t, \zeta_t) = d_X(\zeta_{t+r},  \zeta_{t}) = |r| \ |\xi'|$.

	To use the claim related to \cref{eq:geodesic_interpolation_bound}, we recall that $\varepsilon \geq 2 \tau$ by the hypothesis and that $0 \leq \delta < \tau$. 
    Then, setting $I = [t^i_\tau + \delta, t^i_\tau + \delta + \varepsilon]$ and setting $r \coloneqq 2 (\delta - \psi(\delta))$ ($\psi$ being defined in \cref{eq:def_phi}) and $r' \coloneqq 2 \delta$, there exist $J \subseteq I$, $J' \subseteq I$ such that for all $\sigma \in J, \sigma' \in J'$, we have that $\sigma + r \in I$ and $\sigma' + r' \in I$. Furthermore, recalling the definition of $\xi^{i, \tau}$ from \cref{eq:velocity_interpolation,eq:fancy_velocity_interpolation,eq:bar_velocity_interpolation}, we have, for such $\sigma, \sigma'$, that
	\begin{align*}
		\tilde{\mathcal{V}}_{t^i_\tau + \delta}^{\tau, \varepsilon}(\sigma + r) &= \xi^{i, \tau}_{2 \psi(\delta) + \sigma + r - (t^i_\tau + \delta)} = \xi^{i, \tau}_{2 \delta + \sigma - (t^i_\tau + \delta)} = \mathcal{V}_{t^i_\tau + \delta}^{\tau, \varepsilon}(\sigma),\\
		 \bar{\mathcal{V}}_{t^i_\tau + \delta}^{\tau, \varepsilon}(\sigma' + r') &= \xi^{i, \tau}_{\sigma' + r'- (t^i_\tau + \delta)} = \xi^{i, \tau}_{\sigma' + 2 \delta - (t^i_\tau + \delta)} = \mathcal{V}_{t^i_\tau + \delta}^{\tau, \varepsilon}(\sigma').
	\end{align*}
	We can thus apply the previous claim to get
	\begin{equation*}
		d_{\Tang[\varepsilon]}(\mathcal{V}^{\tau, \varepsilon}_{t^i_\tau + \delta}, \tilde{\mathcal{V}}^{\tau, \varepsilon}_{t^i_\tau + \delta}) = |r| \ |(\mathcal{V}^{\tau, \varepsilon}_{t^i_\tau + \delta})'| \leq 2 \tau C, 
        \qquad d_{\Tang[\varepsilon]}(\mathcal{V}^{\tau, \varepsilon}_{t^i_\tau + \delta}, \bar{\mathcal{V}}^{\tau, \varepsilon}_{t^i_\tau + \delta}) = |r'| \ |(\mathcal{V}^{\tau, \varepsilon}_{t^i_\tau + \delta})'| \leq 2 \tau C,
	\end{equation*}
	finishing the proof of \cref{eq:apriori_bounds_2}. 
    For \cref{eq:apriori_bounds_2_2,eq:apriori_bounds_2_3}, we first consider the case of $\bar{\mathcal{V}}$. 
    In this case, we have, for $t \in [t^\tau_i, t^\tau_{i+1})$,
	\begin{align*}
		\left| d(E u_\tau)(t) - \frac{\diff}{\diff \sigma} E(\bar{\mathcal{V}}_t^{\tau, \varepsilon}(\sigma))\Big|_{\sigma = t}\right| &= \left| \frac{1}{\tau} \int_{t^\tau_i}^{t^\tau_{i+1}} \frac{\diff}{\diff \sigma} E(\bar{\mathcal{V}}_{t^\tau_i}^{\tau, \varepsilon}(\sigma))\Big|_{\sigma=s} \, \diff s - \frac{1}{\tau} \int_{t^\tau_i}^{t^\tau_{i+1}} \frac{\diff}{\diff \sigma} E(\bar{\mathcal{V}}_t^{\tau, \varepsilon}(\sigma))\Big|_{\sigma = t} \, \diff s \right| \\
		&\leq \frac{1}{\tau} \int_{t^\tau_i}^{t^\tau_{i+1}} \underbrace{\left| \frac{\diff}{\diff \sigma} E(\bar{\mathcal{V}}_{t^\tau_i}^{\tau, \varepsilon}(\sigma))\Big|_{\sigma=s} - \frac{\diff}{\diff \sigma} E(\bar{\mathcal{V}}_{t^\tau_i}^{\tau, \varepsilon}(\sigma))\Big|_{\sigma=t^\tau_i} \right|}_{\leq |s-t^\tau_i| C^2 L \text{ by \cref{ass:lipschitz} as in \cref{rmk:rescaling_lipschitz}.}} \, \diff s \\
		&\leq \tau \ C^2  \ L,
	\end{align*}
	where we used the definition of $\bar{\mathcal{V}}^{\tau, \varepsilon}$ to identify its derivative at $t$ with that at $t^\tau_i$; the other two inequlaties follow similarly.
\end{proof}

In the next proposition we show that the velocity interpolations $\bar{\mathcal{V}}^{\tau, \varepsilon}$ (see \cref{eq:bar_velocity_interpolation}) related to a Newtonian sequence $(U_{\tau, n})_{n=-1}^N$ satisfy an inequality reminiscent of the variational characterization of \ACt-curves provided in \cref{eq:AC2_definition}. 
For this reason, we say that \emph{Newtonian sequences are proto-\ACt}.

\begin{proposition}[Newtonian sequences are proto-\ACt]\label{prop:newtonian_sequences_are_proto_ac}
	We assume that $(X,d_X)$ fulfills \cref{ass:interceptor}.
    Let $(U_{\tau, n})_{n=-1}^N$ be a Newtonian sequence, and let $C>0$ be a constant such that the bounds in \cref{eq:apriori_prop_1,eq:apriori_prop_1_5,eq:apriori_prop_2} are fulfilled and, in addition, 
    $$\frac{d_X(U_{\tau,0},U_{\tau,-1})}{\tau}\leq C.$$ 
    We set $\varepsilon \coloneqq \frac{1}{C}\, \epsilon$ and let $\bar{\mathcal{V}}^{\tau, \varepsilon}$, $|u_\tau^{(2)}|$ be defined as in \cref{eq:bar_velocity_interpolation,eq:acceleration_norm_interpolation}. 
    Then, for every $i, j =1, \ldots,N-1$ with $ i \leq j$, and for any $T >0$ such that $\tau \leq T - t^j_\tau \leq T - t^i_\tau \leq \varepsilon$---where $t^k_\tau \coloneqq \tau \cdot k$---we have
	\begin{equation}\label{eq:newtonian_sequences_are_proto_ac}
		d_X \big( \bar{\mathcal{V}}^{\tau, \varepsilon}_{t^i_\tau}(T), \bar{\mathcal{V}}^{\tau, \varepsilon}_{t^j_\tau}(T) \big) \leq \int_{t^i_\tau}^{t^j_\tau} (T - r) \ |u_\tau^{(2)}|(r) \ (1 + C^2 \mathfrak{E}\ (T - r)^2) \, \diff r.
	\end{equation}
\end{proposition}
\begin{proof}
	Since $d_X \big( \bar{\mathcal{V}}^{\tau, \varepsilon}_{t^i_\tau}(T), \bar{\mathcal{V}}^{\tau, \varepsilon}_{t^j_\tau}(T) \big) \leq d_X \big( \bar{\mathcal{V}}^{\tau, \varepsilon}_{t^i_\tau}(T), \bar{\mathcal{V}}^{\tau, \varepsilon}_{t^{i+1}_\tau}(T) \big) + \cdots + d_X \big(\bar{\mathcal{V}}^{\tau, \varepsilon}_{t^{j-1}_\tau}(T), \bar{\mathcal{V}}^{\tau, \varepsilon}_{t^j_\tau}(T) \big)$, it suffices to prove the claim for the special case $j = i+1$. 
    We set $\underline{T} \coloneqq t^i_\tau + 2 \tau$ and note
	\begin{equation*}
		\begin{split}
			d_X \big(\bar{\mathcal{V}}^{\tau, \varepsilon}_{t^i_\tau}(\underline{T}), \bar{\mathcal{V}}^{\tau, \varepsilon}_{t^{i+1}_\tau}(\underline{T}) \big) & = d_X \big(K^\star(U_{\tau, i-1}, U_{\tau, i}), U_{\tau, i+1} \big) = \tau^2 N_\tau(U_{\tau, i-1},U_{\tau, i},U_{\tau, i+1}) \\
			& = \tau \int_{t^i_\tau}^{t^{i+1}_\tau} |u_\tau^{(2)}|(r) \diff r.
		\end{split}
	\end{equation*}
	Since $\bar{\mathcal{V}}^{\tau, \varepsilon}_{t^i_\tau}(t^{i+1}_\tau) = \bar{\mathcal{V}}^{\tau, \varepsilon}_{t^{i+1}_\tau}(t^{i+1}_\tau)$, we may use the \EpsE-interceptor bound as explained in \cref{rem:interceptor_bound_scaling} together with the fact that, for all $0 \leq j < N$, $|(\mathcal{V}_{t^j_\tau}^{\tau, \varepsilon})'| \leq C$,  to get
	\begin{equation*}
		\begin{split}
			d_X \big( \bar{\mathcal{V}}^{\tau, \varepsilon}_{t^i_\tau}(T), \bar{\mathcal{V}}^{\tau, \varepsilon}_{t^{i+1}_\tau}(T) \big) &
            \leq \frac{d_X\left(\bar{\mathcal{V}}^{\tau, \varepsilon}_{t^i_\tau}(\underline{T}), \bar{\mathcal{V}}^{\tau, \varepsilon}_{t^{i+1}_\tau}(\underline{T})\right)}{\tau}\big((T - t^{i+1}_\tau) + (T - t^{i+1}_\tau)^3 \ C^2 \ \mathfrak{E}\big) \\
			&= \int_{t^i_\tau}^{t^{i+1}_\tau} |u_\tau^{(2)}|(r) \, \diff r \ \big((T - t^{i+1}_\tau) + (T - t^{i+1}_\tau)^3 \ C^2 \ \mathfrak{E}\big) \\
			& \leq \int_{t^i_\tau}^{t^{i+1}_\tau} (T - r) |u_\tau^{(2)}|(r) \ \big(1 + (T - r)^2 \ C^2 \ \mathfrak{E}\big) \, \diff r,
		\end{split}
	\end{equation*}
    and this concludes the proof.
\end{proof}

\section{Compactness Argument and Limiting Procedure} \label{sec:compactness_limiting}

To prove convergence of the Newtonian sequences to an AC$^2$-curve, we will proceed in two steps:
\begin{enumerate}
	\item In \cref{prop:discrete_solutions_proto_ac2_geodesic_space}, we show that the Newtonian sequences, realized as curves into the space of geodesics, are proto-AC. This will allow us to use a compactness argument to obtain a limit curve of geodesics in~\cref{prop:compactness_discrete_solutions_geodesic_space}.
	\item In \cref{prop:compactness_discrete_solutions_ac2}, we show that this limit curve of geodesics gives rise to an \ACt-curve into $(X, d_X)$.
\end{enumerate}
To this end, we first need to define the shift operator on the space of geodesics.
\begin{definition} \label{def:Shift}
	For any $\varepsilon > 0$, recalling that $\Tang = \bigcup_{s\in \R} \Tang[s,\varepsilon]$, we define the shift operator $\mathrm{Shift}_t \colon \Tang \to \Tang[t, \varepsilon]$
	to be the operator which maps a geodesic $\xi \in \Tang[s,\varepsilon]$ to the geodesic $\mathrm{Shift}_t \xi \in \Tang[t, \varepsilon]$ defined by $\mathrm{Shift}_t \xi (r) = \xi(r - t + s)$ for all $r \in [t, t+\varepsilon]$.
\end{definition}

\begin{proposition}[Newtonian sequences are proto-AC into the space of geodesics]\label{prop:discrete_solutions_proto_ac2_geodesic_space}
	We assume that $(X,d_X)$ fulfills \cref{ass:interceptor}.
    Let $(U_{\tau, n})_{n=-1}^N$ be a Newtonian sequence, and let $C>0$ be a constant such that the bounds in \cref{eq:apriori_prop_1,eq:apriori_prop_1_5,eq:apriori_prop_2} are fulfilled and, in addition, 
    $$\frac{d_X(U_{\tau,0},U_{\tau,-1})}{\tau}\leq C.$$ 
    We set $\varepsilon \coloneqq \frac{1}{C}\, \epsilon$ and let $\bar{\mathcal{V}}^{\tau, \varepsilon}$, $|u_\tau^{(2)}|$ be defined as in \cref{eq:bar_velocity_interpolation,eq:acceleration_norm_interpolation}.
    Then, for every $i, j =1, \ldots,N-1$ with $ i \leq j$, setting $t^i_\tau \coloneqq \tau \, i$, we have that
	\begin{equation*}
		d_{\Tang[0, \varepsilon]}\Big(\mathrm{Shift}_0 \bar{\mathcal{V}}^{\tau, \varepsilon}_{t^i_\tau}, \mathrm{Shift}_0 \bar{\mathcal{V}}^{\tau, \varepsilon}_{t^j_\tau}\Big) \leq \int_{t^i_\tau}^{t^j_\tau} 3C \, \diff t + \varepsilon \int_{t^i_\tau}^{t^j_\tau} (1 + C^2 \mathfrak{E} \varepsilon^2)\ |u_\tau^{(2)}|(t) \, \diff t.
	\end{equation*}
\end{proposition}
\begin{proof}
	Again, by the triangle inequality, it suffices to consider the case $j = i+1$.
	We set 
    $$\xi \coloneqq \mathrm{Shift}_0 \bar{\mathcal{V}}^{\tau, \varepsilon}_{t^i_\tau}, 
    \qquad \zeta \coloneqq \mathrm{Shift}_0 \bar{\mathcal{V}}^{\tau, \varepsilon}_{t^{i+1}_\tau}.$$ 
	We will show that for any $s \in [0, \varepsilon]$, we have
	\begin{equation*}
		d_X\big(\xi(s), \zeta(s)\big) \leq \tau \ \Big(3C + \varepsilon \ (1 + C^2 \mathfrak{E} \varepsilon^2)\ |u_\tau^{(2)}|(t^i_\tau) \Big).
	\end{equation*}
	First, let us consider the case $s \leq 2\tau$. In this case, we have
	\begin{equation*}
		\begin{split}
			d_X\big(\xi(s), \zeta(s)\big) & \leq d_X(\xi(s), U_{\tau, i}) + d_X(U_{\tau, i}, \zeta(s)) \\
			& \leq \frac{d_X(U_{\tau,i}, U_{\tau,i-1})}{\tau} \ |(s - \tau)| + |u^{(1)}_\tau|(t^{i}_\tau) \ |s| \\
			& \leq \tau C + 2\tau C = 3 \tau C.
		\end{split}
	\end{equation*}
	Second, let us consider the case $s > 2\tau$, so that $s \in (2\tau, \varepsilon]$. Using \cref{prop:newtonian_sequences_are_proto_ac} with $T = t^i_\tau + s$---note in particular that in this case, $T \geq t^{i+1}_\tau + \tau$---, we have
	\begin{equation*}
		d_X \big(\xi(s), \zeta(s - \tau) \big) \leq \tau \Big(\varepsilon |u_\tau^{(2)}|(t^i_\tau)\ (1 + C^2 \mathfrak{E} \varepsilon^2) \Big).
	\end{equation*}
	Combining this with $d_X \big(\zeta(s - \tau), \zeta(s) \big) \leq |u^{(1)}_\tau|(t^i_\tau) \, \tau \leq C \, \tau$, we obtain the claim.
\end{proof}

In the proof of the next result, we use the fact that, under \cref{ass:interceptor,ass:local_compactness}, the space $(X, d_X)$ is \emph{proper}, i.e., every closed and bounded subset $K \subset X$ is compact.
More precisely, since $(X, d_X)$ being geodesic implies that $(X, d_X)$ is a length space and since we assume  that $(X, d_X)$ is locally compact, by virtue of the metric Hopf-Rinow Theorem (or  Hopf-Rinow-Cohn-Vossen Theorem, see \cite[Theorem~2.5.28]{BuragoMetric}) the completeness is equivalent to the properness.

\begin{proposition}[Compactness of discrete solutions into the space of geodesics]\label{prop:compactness_discrete_solutions_geodesic_space}
	We assume that $(X, d_X)$ and $E \colon X \to \R \cup \{ +\infty\}$ fulfill \cref{ass:interceptor,ass:local_compactness,ass:lsc,ass:coercivity,ass:lipschitz}. 
    Let $T > 0$ and $R > 0$, and let $(\tau_i)_{i=1}^\infty$ be a sequence of step sizes such that $\tau_i \to 0$ as $i \to \infty$. Furthermore, we assume that for each $i$, we are given a Newtonian sequence $(U_{\tau_i, n})_{n=-1}^{N_{\tau_i}}$ fulfilling
	\begin{equation*}
		N_{\tau_i} \tau_i > T, \quad E(U_{\tau_i, 0}) \leq R, \quad d^2_X(U_{\tau_i, 0}, x_\star) \leq R,\quad d(Eu_{\tau_i})(0) \leq R, \quad \frac{d_X(U_{\tau_i,-1},U_{\tau_i,0})}{\tau_i} \leq R.
	\end{equation*}
	Then, defining $C\coloneqq C(x_\star, \tau_\star, R, T)$ as in \cref{lem:apriori_bounds} and setting $\varepsilon = \frac{1}{C}\, \epsilon$, there exists---without relabeling---a subsequence $(\tau_i)_{i=1}^\infty$ and a curve $V\colon [0, T] \to \Tang$ such that for all $t$, $V_t \in \Tang[t, \varepsilon]$ and such that, for all $t \in [0, T]$,
	\begin{equation}\label{eq:compactness_in_geodesic_space}
		V_t = \lim_{i \to \infty} \bar{\mathcal{V}}^{\tau_i, \varepsilon}_{t} = \lim_{i \to \infty} \mathcal{V}^{\tau_i, \varepsilon}_{t} = \lim_{i \to \infty} \tilde{\mathcal{V}}^{\tau_i, \varepsilon}_{t} \qquad \text{ in } \Tang[t, \varepsilon].
	\end{equation}
\end{proposition}
\begin{proof}
	We will only prove the first convergence in \cref{eq:compactness_in_geodesic_space}; the others follow then from \cref{lem:apriori_bounds_2}. 
    In fact, we can instead prove existence of some $W\colon [0, T] \to \Tang[0, \varepsilon]$ such that, for all $t \in [0, T]$,
	\begin{equation*}
		W_t = \lim_{i \to \infty} \mathrm{Shift}_0 \bar{\mathcal{V}}^{\tau_i, \varepsilon}_{t} \qquad \text{ in } \Tang[0, \varepsilon],
	\end{equation*}
	since we can set $V_t \coloneqq \mathrm{Shift}_t W_t$.
	By \cref{lem:apriori_bounds}, there exists some $C = C(x_\star, \tau_\star, R, T)$ such that $|u^{(1)}_{\tau_i}|(t) \leq C$ for all $t \in [0, T]$ and all $i$, and that $\big\| |u_{\tau_i}^{(2)}| \big\|_{L^2([0, T])} \leq C$ for all $i$. Thus, we first can---without relabeling---extract a weakly convergent subsequence such that $|u_{\tau_i}^{(2)}| \rightharpoonup_{L^2} A$ for some $A \in L^2([0, T], \R_+)$. \\
	As in~\cite[Proof of Corollary~3.3.4]{ambrosioGradientFlows2008}, we fix $0 \leq s \leq t \leq T$ and define
	\begin{equation} \label{eq:compactness_discrete_solutions_geodesic_space_s_i_t_i}
    \begin{split}
		&s(i) \coloneqq \max\{r \mid r = \tau_i \ n \text{ for some } n \in \NN, \ r \leq s\}, \\
        &t(i) \coloneqq \max\{r \mid r = \tau_i \ n \text{ for some } n \in \NN, \ r \leq t\}.
    \end{split}
	\end{equation}
	By \cref{prop:discrete_solutions_proto_ac2_geodesic_space} and using that $\mathrm{Shift}_0 \bar{\mathcal{V}}^{\tau_i, \varepsilon}_{s} = \mathrm{Shift}_0 \bar{\mathcal{V}}^{\tau_i, \varepsilon}_{s(i)}$ and similarly for $t$ and $t(i)$, we have
	\begin{equation*}
		d_{\Tang[0, \varepsilon]}\Big(\mathrm{Shift}_0 \bar{\mathcal{V}}^{\tau_i, \varepsilon}_{s}, \mathrm{Shift}_0 \bar{\mathcal{V}}^{\tau_i, \varepsilon}_{t}\Big) 
        \leq \int_{s(i)}^{t(i)} \Big(3C + \varepsilon \ (1 + C^2 \mathfrak{E} \varepsilon^2)\ |u_{\tau_i}^{(2)}| (r) \Big) \, \diff r,
	\end{equation*}
	and hence
	\begin{equation*}
		\limsup_{i \to \infty} 
        d_{\Tang[0, \varepsilon]} \Big(\mathrm{Shift}_0 \bar{\mathcal{V}}^{\tau_i, \varepsilon}_{s}, \mathrm{Shift}_0 \bar{\mathcal{V}}^{\tau_i, \varepsilon}_{t}\Big) 
        \leq \int_{s}^{t} \Big(3C + \varepsilon \ (1 + C^2 \mathfrak{E} \varepsilon^2)\ A(r)  \Big)\,\diff r.
	\end{equation*}
	We would like to apply an Arzelà-Ascoli-type argument to obtain the existence of the limit curve $W\colon [0, T] \to \Tang[0, \varepsilon]$. To this end, we first need to make sure that there exists some compact $K_{\mathcal{T}} \subset \Tang[0, \varepsilon]$ such that for all $i$ and all $t \in [0, T]$, $\mathrm{Shift}_0 \bar{\mathcal{V}}^{\tau_i, \varepsilon}_{t} \in K_{\mathcal{T}}$. 
    To see this, we first note that there exists some compact $K_X \subseteq X$ such that for all $i$ and all $t \in [0, T]$ and for all $r \in [t, t + \varepsilon]$, $\bar{\mathcal{V}}^{\tau_i, \varepsilon}_{t}(r) \in K_X$; this follows from the uniform bound on $|u^{(1)}_{\tau_i}|$ and on $d_X(U_{\tau_i, k}, x_\star)$ provided by \cref{lem:apriori_bounds} combined with the fact that $(X, d_X)$ is proper. 
    Using this and again the bound on $|u^{(1)}_{\tau_i}|$, we find we can use Arzelà-Ascoli to see that any sequence $\{\mathrm{Shift}_0 \bar{\mathcal{V}}^{\tau_{i_j}, \varepsilon}_{t^j_\tau}\}_{j=1}^\infty$ has a convergent subsequence.
    Thus, the set $\{\mathrm{Shift}_0 \bar{\mathcal{V}}^{\tau_i, \varepsilon}_{t} \mid i \in \NN, t \in [0, T]\}$ is precompact, and we can define $K_{\mathcal{T}}$ as its closure.
    Moreover, $K_{\mathcal{T}}$ consists of geodesics, as the limit of a sequence of geodesics uniformly convergent is in turn a geodesic.\\
	We are now in a position to apply the version of the Arzelà-Ascoli Theorem~\cite[3.3.1]{ambrosioGradientFlows2008}, which yields the existence of a curve $W \colon [0, T] \to \Tang[0, \varepsilon]$ such that, without relabeling, $\mathrm{Shift}_0 \bar{\mathcal{V}}^{\tau_i, \varepsilon}_{t} \to W_t$ in $\Tang[0, \varepsilon]$ for all $t \in [0, T]$.
\end{proof}

Finally, we are in a position to obtain an \ACt-trajectory as the limit of Newtonian sequences. 

\begin{proposition}[Compactness in the space of \ACt-curves]\label{prop:compactness_discrete_solutions_ac2}
	We assume that $(X, d_X)$ and $E \colon X \to \R \cup \{ +\infty\}$ fulfill \cref{ass:interceptor,ass:functional,ass:local_compactness,ass:lsc,ass:lipschitz,ass:coercivity}. Let $R > 0$ and $T > 0$, and let $(\tau_i)_{i=1}^\infty$ be a sequence of step sizes such that $\tau_i \to 0$ as $i \to \infty$. Furthermore, we assume that for each $i$, we are given a Newtonian sequence $(U_{\tau_i, n})_{n=-1}^{N_{\tau_i}}$ fulfilling
	\begin{equation*}
		N_{\tau_i} \tau_i > T, \quad E(U_{\tau_i, 0}) \leq R, \quad d^2_X(U_{\tau_i, 0}, x_\star) \leq R,\quad d(Eu_{\tau_i})(0) \leq R, \quad \frac{d_X(U_{\tau_i,-1},U_{\tau_i,0})}{\tau_i} \leq R.
	\end{equation*}
	Then, defining $C\coloneqq C(x_\star, \tau_\star, R, T)$ as in \cref{lem:apriori_bounds} and setting $\varepsilon = \frac{1}{C}\, \epsilon$, there exists---without relabeling---a subsequence $(\tau_i)_{i=1}^\infty$, an \ACt-curve $\bar{u}\colon [0, T] \to X$, a function $\varphi\colon [0, T] \to \R$ of finite total variation and a function $A \in L^2_{loc}([0, T], \R)$ such that
	\begin{gather}
		V^{\bar{u}}_t = \lim_{i \to \infty} \bar{\mathcal{V}}^{\tau_i, \varepsilon}_{t} = \lim_{i \to \infty} \mathcal{V}^{\tau_i, \varepsilon}_{t} = \lim_{i \to \infty} \tilde{\mathcal{V}}^{\tau_i, \varepsilon}_{t} \quad \text{ in } \Tang[t, \varepsilon], \text{ for all } t \in [0, T], \label{eq:compactness_1}\\
		\bar{u}(t) = \lim_{i \to \infty} u_{\tau_i}(t) = \lim_{i \to \infty} \tilde{u}_{\tau_i}(t), \quad \text{ in } X, \text{ for all } t \in [0, T], \label{eq:compactness_2}\\
		\varphi(t) \coloneqq \lim_{i \to \infty} d(Eu_{\tau_i})(t) = \frac{\diff}{\diff \sigma} E(V^{\bar{u}}_t(\sigma)) \Big|_{\sigma=t} \quad \text{ for all } t \in [0, T], \label{eq:compactness_3}\\
		|u^{(2)}_{\tau_i}| \rightharpoonup A, \quad A(t) \geq |\bar{u}''|(t) \quad \text{ for a.e. } t \in [0, T], \label{eq:compactness_4}\\
		\liminf_{i \to \infty} |G_{\tau_i}|(t) \geq |\partial^- E|(\bar{u}(t)) \quad \text{ for almost all } t \in [0, T], \label{eq:compactness_5}
	\end{gather}
	where $|G_{\tau_i}|$ is defined as in \cref{eq:def_G_tau}\footnote{In particular, since $X$ and $E$ fulfill \cref{ass:coercivity,ass:lsc,ass:local_compactness,ass:interceptor}, $|G_{\tau_i}|$ exists for $\tau_i$ sufficiently small.}.
\end{proposition}
\begin{proof}
	First, from \cref{eq:apriori_prop_2}, it follows that we can write the functions $d(Eu_{\tau_i})(t)$ as the sum of (1.) a bounded monotone non-decreasing function and (2.) a monotone non-increasing function. Thus, we may apply Helly~\cite[Theorem 3.3.3]{ambrosioGradientFlows2008} to obtain the existence of a subsequence---without relabeling---and a function $\varphi\colon [0, T] \to \R \cup \{-\infty\}$ such that for all $t$, $d(Eu_{\tau_i})(t) \to \varphi(t)$ as $i \to \infty$. 
  Next, we can apply \cref{prop:compactness_discrete_solutions_geodesic_space} to obtain the existence of a further subsequence---without relabeling---and a curve $V\colon [0, T] \to \Tang$ such that for all $t$, $V_t \in \Tang[t, \varepsilon]$ and such that \cref{eq:compactness_1} holds.
	As in the previous proof, we further extract a subsequence such that $|u^{(2)}_{\tau_i}| \rightharpoonup_{L^2} A$  for some $A \in L^2_{loc}([0, T], \R_+)$. Defining $\bar{u}(t) \coloneqq V_t(t)$, \cref{eq:compactness_2} follows immediately from \cref{eq:compactness_1,lem:apriori_bounds_2}. \cref{eq:compactness_3} follows from \cref{ass:functional} together with \cref{eq:apriori_bounds_2_2}; it further implies that $\varphi$ takes values in $\R$. That $\varphi$ has bounded total variation comes from the fact that the increase of the functions $d(Eu_{\tau_i})(t)$ is uniformly bounded, a property that transfers to the limit $\varphi$ and implies finite total variation. \\
	Next, we want to show that $\bar u$ is an \ACt-curve and that $A \geq |\bar u''|$. By the definition of an \ACt-curve (\cref{def:AC2_curve}) and by \cref{lem:second_derivative_bound}, it suffices to show that for any $0 \leq s \leq t \leq T$ and $S \in \R$ such that $t < S < s+\frac{1}{C}\, \epsilon$ (the boundary cases follow by continuity), we have
    \begin{equation} \label{eq:proof_compactness_as_ac2}
		d_X \big( V_s(S), V_t(S) \big) \leq \int_s^t (S - r) \ A(r) \ \big(1 + C^2 \mathfrak{E} (S -r)^2 \big) \, \diff r.
	\end{equation}
	We fix such an $s, t, S$ and define $s(n)$ and $t(n)$ as in \cref{eq:compactness_discrete_solutions_geodesic_space_s_i_t_i}. For $n$ large enough, we have $S - s(n) \leq \frac{1}{C}\, \epsilon$ and $S - t(n) \geq \tau_n$. Thus, for $n$ large enough, we may apply \cref{prop:newtonian_sequences_are_proto_ac} to obtain
	\begin{equation*}
		d_X \big(\bar{\mathcal{V}}^{\tau_n, \varepsilon}_{s(n)}(S), \bar{\mathcal{V}}^{\tau_n, \varepsilon}_{t(n)}(S) \big) 
        \leq \int_{s(n)}^{t(n)} (S - r) |u_{\tau_n}^{(2)}|(r) \ \big(1 + C^2 \mathfrak{E} \ (S - r)^2\big) \, \diff r,
	\end{equation*}
	and therefore
	\begin{equation*}
		\begin{split}
			d_X \big( V_s(S), V_t(S) \big) & 
            = \lim_{n \to \infty} d_X \big( \bar{\mathcal{V}}^{\tau_n, \varepsilon}_{s(n)}(S), \bar{\mathcal{V}}^{\tau_n, \varepsilon}_{t(n)}(S) \big) \\
			& \leq \liminf_{n \to \infty} \int_{s(n)}^{t(n)} (S - r) |u_{\tau_n}^{(2)}|(r) \ \big(1 + C^2 \mathfrak{E} \  (S - r)^2\big) \, \diff r \\
			& \leq \int_{s}^{t} (S - r) A(r) \  \big(1 + C^2 \mathfrak{E}\ (S - r)^2 \big) \, \diff r,
		\end{split}
	\end{equation*}
	which shows \cref{eq:proof_compactness_as_ac2}, identifies $V$ as a tangent family $V^{\bar u}$, and yields \cref{eq:compactness_4}; note that we used $d(\bar{\mathcal{V}}^{\tau_n, \varepsilon}_{s(n)}(S), \bar{\mathcal{V}}^{\tau_n, \varepsilon}_{s}(S)) \leq C \tau_n$, and likewise for $t$ and $t(n)$. Finally, \cref{eq:compactness_5} follows, as in~\cite[Proof of Corollary 3.3.4]{ambrosioGradientFlows2008}, from \cref{eq:discrete_newtonian_motion_estimate_G}, which gives, for almost every $t \in [0, T]$,
	\begin{equation*}
		|\partial^- E|(\bar{u}(t)) \leq \liminf_{i \to \infty} |\partial E|(\tilde{u}_{\tau_i}(t)) \leq \liminf_{i \to \infty} |G_{\tau_i}|(t),
	\end{equation*}
    and concludes the proof.
\end{proof}

\begin{proposition}[Limit curves are curves of Newtonian motion] \label{prop:limit_curves_are_newtonian_curves}
	Under the same assumptions and notations as in \cref{prop:compactness_discrete_solutions_ac2},
	the limit \ACt-curve $\bar u\colon [0,T]\to X$ fulfills, for almost every $s, t \in [0, T]$ such that $s < t$, the following inequality:
	\begin{equation*}
		\frac{\diff}{\diff \sigma} E(V^{\bar{u}}_s(\sigma)) \Big|_{\sigma = s} - \frac{\diff}{\diff \sigma} E(V^{\bar{u}}_t(\sigma)) \Big|_{\sigma = t}
        \geq
        \frac{1}{2}\int_s^t  |\bar{u}''|^2(r) \diff r + \frac{1}{2}\int_s^t  |\partial^- E|^2(\bar{u}(r)) \diff r - \int_s^t H_0^{E, +}(V^{\bar{u}}_{r}) \diff r .
	\end{equation*}
	Furthermore, the map $t \mapsto \frac{\diff}{\diff \sigma} E(V^{\bar{u}}_t(\sigma))\big|_{\sigma=t}$ has finite total variation and $t \mapsto H_0^{E, +}(V^{\bar{u}}_{t}) \in L^1([0, T], \R)$.
	In particular, if $(|\partial^- E|, H_0^{E, +})$ is an upper gradient-Hessian pair (see \cref{def:upper_gradient}) for $E$, then $\bar{u}$ is a curve of Newtonian motion (see \cref{def:curve_newtonian_motion}).
\end{proposition}
\begin{proof}
	First, by \cref{ass:lipschitz}, for every geodesic $\xi$, we have that $|H_0^E(\xi)| \leq L \, |\xi'|^2$, and thus also $|H_0^{E, +}(\xi)| \leq L \, |\xi'|^2$.
    This fact, owing to \cref{eq:compactness_1} and to the uniform bound $C>0$ from \cref{lem:apriori_bounds}, ensures that $\big|H_0^E(\tilde{\mathcal{V}}_t^{\tau_i, \varepsilon}) \big| \leq C^2 L$, which in turn gives the boundedness of $H_0^{E, +}(V_s^{\bar{u}})$ and guarantees that $t \mapsto H_0^{E, +}(V^{\bar{u}}_{t}) $ belongs to $ L^1([0, T], \R)$. \\
	Moreover, we have, for any $s', s, t', t \in [0, T]$ such that $s' < s < t < t'$:
	\begin{equation}\label{eq:limit_curves_are_newtonian}
		\begin{split}
			& \int_s^t \frac{1}{2} |\bar{u}''|^2(r) \diff r + \int_s^t \frac{1}{2} |\partial^- E|^2(\bar{u}(r)) \diff r - \int_s^t H_0^{E, +}(V^{\bar{u}}_{r}) \, \diff r \\
			& \leq \liminf_{i \to \infty} \bigg(\int_s^t \frac{1}{2} |u_{\tau_i}^{(2)}|^2(r) \, \diff r\bigg) + \int_s^t \liminf_{i \to \infty} \frac{1}{2} |G_{\tau_i}|^2(r) \, \diff r - \limsup_{i \to \infty} \int_s^t H_0^{E}(\tilde{\mathcal{V}}^{\tau_i, \varepsilon}_{r}) \, \diff r \\
			& \leq \liminf_{i \to \infty} \bigg(\int_s^t \frac{1}{2} |u_{\tau_i}^{(2)}|^2(r) \, \diff r + \int_s^t \frac{1}{2} |G_{\tau_i}|^2(r) \, \diff r - \int_s^t H_0^{E}(\tilde{\mathcal{V}}^{\tau_i, \varepsilon}_{r})\, \diff r\bigg) \\
			& \leq \liminf_{i \to \infty} \big(d(Eu_{\tau_i})(s') - d(Eu_{\tau_i})(t')\big) + C (t' - t) + C (s - s') \\
			 &=  \frac{\diff}{\diff \sigma} E(V^{\bar{u}}_{s'}(\sigma)) \Big|_{\sigma = s'} - \frac{\diff}{\diff \sigma} E(V^{\bar{u}}_{t'}(\sigma)) \Big|_{\sigma = t'} + C (t' - t) + C (s - s'),
		\end{split}
	\end{equation}
	where we used Fatou's lemma in the second step---note that $|G_{\tau_i}|^2(s) \geq 0$ for all $s$ and $i$, and that $\big| H_0^{E}(\tilde{\mathcal{V}}^{\tau_i, \varepsilon}_{s}) \big| \leq C^2 \, L$  as discussed above---, and we used the boundedness of $H_0^E$ in the third passage, and \cref{eq:compactness_3} in the last step.
	Note that $t \mapsto \frac{\diff}{\diff \sigma} E(V^{\bar{u}}_t(\sigma))\big|_{\sigma=t} = \varphi(t)$ has finite variation by \cref{prop:compactness_discrete_solutions_ac2}. Thus, for almost every $s, t \in [0, T]$, $\varphi$ is continuous in $s$ and $t$, which allows us to take the limit in \cref{eq:limit_curves_are_newtonian} to finish the proof.
\end{proof}

	Finally, let us finish this section by assembling the hitherto proven results into an existence statement for curves of Newtonian motion.

	\begin{theorem} \label{thm:existence_newtonian_motion}
		We assume that $(X, d_X)$ and $E \colon X \to \R \cup \{ +\infty\}$ fulfill \cref{ass:interceptor,ass:functional,ass:local_compactness,ass:lipschitz,ass:coercivity,ass:lsc}, and that the pair $(|\partial^- E|, H_0^{E,+})$ is an upper gradient-Hessian pair for $E$ (\cref{def:upper_gradient}). Then, for each $T > 0$ and $\varepsilon' > 0$ and for each $\xi \in \Tang[-\varepsilon'/2,\varepsilon']$---we recall \cref{def:space_of_geodesics}---, there exists $\varepsilon > 0$ and a curve of Newtonian motion $\bar{u}\colon [0, T] \to X$ such that $V_0^{\bar{u}} = \xi|_{[0, \varepsilon]}$.
	\end{theorem}
	\begin{proof}
		By \cref{lem:apriori_bounds} applied with horizon $2T$ and for $\tau \leq \varepsilon'/2$ small enough, there exists $N_\tau \in \NN$ such that $N_\tau \tau > T$ and such that we can construct a Newtonian sequence $(U_{\tau, i})_{i=-1}^{N_\tau}$ such that $U_{\tau, 0} = \xi_0$ and $U_{\tau, -1} = \xi_{-\tau}$.
        One can immediately verify the uniform bounds with an appropriate $R$ (in particular, 
        $$\frac{d_X(U_{\tau, -1}, U_{\tau, 0})}{\tau} = |\xi'|,$$ 
        and 
        $$d(E u_\tau)(0) = \frac{E(\xi_0) - E(\xi_{-\tau})}{\tau} \leq \frac{\diff}{\diff \sigma} E(\xi_\sigma)|_{\sigma=0} + \frac{1}{2} \tau L |\xi'|^2,$$
        where the derivative exists according to \cref{ass:existence_derivative} and with $L>0$ chosen according to \cref{ass:lipschitz}).
		Let $C>0$ be the corresponding uniform a priori constant. Setting $\varepsilon = \min\{\epsilon/C,\varepsilon'/2\}$ and restricting the tangent curves to this duration, we have by \cref{def:interpolations} that $\mathcal{V}^{\tau, \varepsilon}_0(\sigma) = \xi_{\sigma - \tau}$ for $\sigma \in [0, \varepsilon]$, and thus $\mathcal{V}^{\tau, \varepsilon}_0(\sigma) \to \xi(\sigma)$ uniformly over $\sigma \in [0, \varepsilon]$ as $\tau \to 0$. Thus, \cref{prop:compactness_discrete_solutions_ac2}---in particular, \cref{eq:compactness_1}---together with \cref{prop:limit_curves_are_newtonian_curves} yields the thesis, where we use the fact that $(|\partial^- E|, H_0^{E,+})$ is an upper gradient-Hessian pair for $E$.
	\end{proof}

\section{Linear wave equation as a curve of Newtonian motion} \label{sec:linear_wave_equation}

In this section, we study the solutions of the wave equation in $\R^d$ as curves of Newtonian motion in $L^2 (\R^d,\R)$ driven by the Dirichlet energy.
We recall that, in this setting, the Dirichlet energy is the functional $E\colon L^2 (\R^d,\R) \to \R \cup \{+\infty\} $ defined as 
\begin{equation*}
    E(u) \coloneqq \begin{cases}
        \frac12 \int_{\R^d} |\nabla u(x)|^2 \, \diff x & u\in H^1(\R^d,\R),\\
        +\infty & \mbox{otherwise}.
    \end{cases}
\end{equation*}
The aim of this section is twofold. First, we plan to prove the existence of solutions of the linear wave equation with regular initial data. Second, as the space $L^2(\R^d,\R)$ is not locally compact---we recall that this hypothesis is required in \cref{ass:local_compactness}---we also want to show that the metric-Newtonian machinery developed before is flexible enough to be adapted to milder settings with \textit{ad hoc} arguments.
We shall consider the linear wave equation 
\begin{equation} \label{eq:linear_wave}
 \begin{cases}
     \partial_{tt}^2 u_t - \Delta u_t =0,\\
     u_{t=0}= u_0 \in H^3_c(\R^d,\R),\\
     \partial_t u_{t=0}= \dot{u}_0 \in H^2_c(\R^d,\R),
 \end{cases} 
\end{equation}
where $H^k_c(\R^d,\R)$ denotes the subset of $H^k(\R^d,\R)$ consisting of functions with compact support.
The novel contribution of this section is the \emph{integral variational characterization} of the solution $u_\cdot\colon [0,T]\to L^2(\R^d,\R)$ of \cref{eq:linear_wave}, which reads as
	\begin{equation*} 
		DE(u_s)[\partial_t u_s] - DE(u_t)[\partial_t u_t] = \int_{s}^{t} \frac{1}{2} \|\partial^2_{tt} u_r\|^2_{L^2} + \frac{1}{2} \|\Delta u_r\|^2_{L^2} - \|\nabla \partial_t u_r\|^2_{L^2} \, \diff r
	\end{equation*}
for every $s, t \in [0, T)$ with $s < t$. 
We emphasize that the last identity is not implied by the well-known total energy conservation
\begin{equation*}
   \frac12 \| \partial_t u_t \|_{L^2}^2 + E(u_t) =
   \frac12 \| \dot u_0 \|_{L^2}^2 + E(u_0) \qquad \forall t \in (0,T),
\end{equation*}
which holds for regular solutions of \cref{eq:linear_wave} (see, e.g., \cite[Section~2.4.3]{evans2022partial}).

\begin{remark}\label{rmk:gradient_sol_wave}
    Let us assume that $u_\cdot \in C^1([0,T],H^2(\R^d,\R))$ is a weak solution of \cref{eq:linear_wave}. Then, setting $ \nabla u_t \eqqcolon v_t =(v_t^1,\ldots,v_t^d) \in H^1(\R^d,\R^d)$ for every $t\in [0,T]$, we observe that the trajectory $v_\cdot^j \in C^1([0,T],H^1(\R^d,\R))$ is in turn a weak solution of the linear wave equation $\partial^2_{tt}v_t^j - \Delta v_t^j=0$ for $j=1,\ldots,d$.
\end{remark}

\begin{remark}\label{rmk:geometry_L2}
    Since $L^2$ is a Hilbert space, given two geodesics $\xi, \zeta\colon [0, T]\to X \coloneqq L^2(\R^d,\R)$ such that $\xi(0) = \zeta(0)$ and $\max\{|\xi'|, |\zeta'|\} \leq 1$, we have that 
    \begin{equation*}
        \|\xi_t - \zeta_t\|_{L^2} =\frac{\|\xi_s - \zeta_s\|_{L^2}}{s}t
    \end{equation*}
    for every $t,s\in (0,T]$. In particular, this means that the \EpsE-interceptor bound holds with $\mathfrak{E}=0$ and for every $\epsilon>0$.
\end{remark}
\begin{remark}\label{rmk:objects_here_1}
    Leveraging the vector space structure of $X\coloneqq L^2(\R^d,\R)$, we have that
    \begin{equation} \label{eq:K_star_Dir}
        K^\star(x,y) = y+(y-x)=2y -x
    \end{equation}
    for every $x,y\in X$, where we recall that the operator $K^\star$ has been defined in \cref{eq:def_Kstar}.
\end{remark}
\begin{remark}\label{rmk:gradient_syst}
    Before proceeding, we compute the expression of the slope $|\partial E|$. Namely, for every $u\in H^1(\R^d,\R)$ we have
    \begin{equation*}
        \begin{split}
            |\partial E|(u)&= \limsup_{v\to 0} \frac{\big( E(u) - E(u+v) \big)^+}{\| v \|_{L^2}} = \limsup_{v\to 0} \frac{\big( -\langle \nabla u , \nabla v \rangle - \frac12 \|\nabla v \|_{L^2}^2  \big)^+}{\| v \|_{L^2}} \leq \limsup_{v\to 0} \frac{|\langle \nabla u , \nabla v \rangle|}{\| v \|_{L^2}},
        \end{split}
    \end{equation*}
    where we used that $E(u+v) = +\infty$ if $v\not \in H^1$.
    Thus, on the one hand, if $u\in H^2$, we have $\langle \nabla u , \nabla v \rangle = - \langle \Delta u ,  v \rangle$, and by choosing $(w_n)_n \subset H^1$ and $(r_n)_n \subset \R$ such that $w_n\to_{L^2} \Delta u$ and $r_n\searrow0$ such that $r_n \|\nabla w_n\|_{L^2}^2 \to 0$ as $n\to\infty$---and using $v_n \coloneqq r_n w_n$ as a test vector---, it follows that $|\partial E|(u) = \| \Delta u \|_{L^2}$. On the other hand, if $u \not \in H^2$, we can similarly deduce that $|\partial E|(u)=+\infty$.
    Moreover, for every $t\in [0,T]$ and $\varepsilon>0$, we have that, for every geodesic $\xi\in \Tang[t, \varepsilon]$, there exist $u,v\in X$ such that $\xi\colon [t,t+\varepsilon]\ni \sigma \mapsto u + \sigma v$.
    Therefore, a direct computation shows that
    \begin{equation} \label{eq:def_H0_Dir}
        H_0^{E}(\xi) = \int_{\R^d} |\nabla v|^2 \diff x = 2 E(v)
    \end{equation}
    whenever $u, v \in H^1(\R^d, \R)$, where $H_0^E$ (for a generic driving energy $E$) has been originally introduced in \cref{eq:def_H0}. To use a replacement of the Hessian of $E$ on the whole of $L^2$, we extend this formula:
    
	\begin{definition} \label{def:hessian_wave}
		We define, for $t \in \R, \varepsilon \in (0, +\infty)$ and for any geodesic $\xi\colon [t, t + \varepsilon] \ni \sigma \mapsto u + \sigma v$ for $u, v \in L^2$, the quantity
			\begin{equation}
				\bar{H}_0^E(\xi) \coloneqq \begin{cases}
					2 E(v) \qquad & \text{if } u, v \in H^1(\R^d, \R), \\
					+\infty \qquad &\text{otherwise.}
				\end{cases}
			\end{equation}
			We note that for $u, v \in H^1$, we have that $\bar{H}_0^E(\xi) = H_0^E(\xi)$.
	\end{definition}
        
    Finally, we also recall that, for every $u, v \in H^1$, the directional derivative of $E$ at $u$ in the direction of $v$ is $DE(u)[v] = \langle \nabla u, \nabla v \rangle_{L^2}$.
\end{remark}
\subsection{Upper gradient-Hessian pair for the Dirichlet energy}
Before we proceed to construct a curve of Newtonian motion for the Dirichlet energy, we motivate this construction by (1.)\ providing an upper gradient-Hessian pair for the Dirichlet energy, and (2.)\ showing that any curve of Newtonian motion fulfills the wave equation~\eqref{eq:linear_wave}.
\begin{proposition}\label{prop:ugs_wave}
	The pair $(|\partial E|, \bar{H}^E_0)$, with $\bar{H}^E_0$ defined as in \cref{def:hessian_wave}, is an upper gradient-Hessian pair for $E$, according to \cref{def:upper_gradient}. Furthermore, any curve of Newtonian motion $u$ according to \cref{def:curve_newtonian_motion} fulfills the wave equation
	\begin{align*}
		u_t \in H^2(\R^d, \R), \,\partial_t u_t \in H^1(\R^d, \R), \qquad &\text{for a.e. } t, \\
		\partial_{tt}^2 u_t - \Delta u_t = 0	, \qquad &\text{for a.e. } t.
	\end{align*}
\end{proposition}
\begin{proof}
	By \cref{def:upper_gradient}, and using \cref{lem:euclidean_ac2_equivalence}, we have to show that for any $u_{\cdot}\colon (a, b) \to X = L^2(\R^d, \R)$, for $a, b \in \R$ and $a < b$, such that $u \in C^1((a, b), X)$, $\partial_t u$ is an AC-curve, and such that $s \mapsto DE(u_s)[\partial_t u_s]$ has finite total variation on every compact subinterval of $(a,b)$, $|\partial E| \circ u$ is Borel, and the following inequality holds for almost all $s, t \in (a, b)$ such that $s < t$:
	\begin{equation}\label{eq:proof_of_ugs_wave_1}
		DE(u_s)[\partial_t u_s] - DE(u_t)[\partial_t u_t] \leq \int_{s}^{t} \frac{1}{2} \|\partial^2_{tt} u_r\|^2_{L^2} + \frac{1}{2} |\partial E|^2(u_r) - \bar{H}_0^E(V_r^u) \, \diff r
	\end{equation}
	whenever $\bar{H}_0^E(V_r^u) \in L^1_{\mathrm{loc}}$, where $DE(u_r)$ denotes the directional derivative of $E$ at $u_r$; we also recall the identities
	\begin{align}
		|\partial E|(u_t) &= \begin{cases}
			\| \Delta u_t \|_{L^2}, \qquad &\text{ if } u_t \in H^2, \\
			+\infty \qquad &\text{ otherwise,}
		\end{cases} \label{eq:proof_of_ugs_wave_2}\\
		\bar{H}_0^E(V_t^u) &= \begin{cases}
			\|\nabla \partial_t u_t\|^2_{L^2}\qquad &\text{ if } u_t \in H^1, \partial_t u_t \in H^1, \\
			+\infty \qquad &\text{ otherwise.}
		\end{cases}
	\end{align}
	from \cref{rmk:gradient_syst,def:hessian_wave}. The Borel claim follows immediately, since any such curve is in particular continuous with regards to the strong topology on $L^2$ (and thus Borel) and $|\partial E|$ is seen from \cref{eq:proof_of_ugs_wave_2} to be Borel as well. First, we would like to assume integrability of the integrand in \cref{eq:proof_of_ugs_wave_1}. To this end, we set $g(r) \coloneqq \frac{1}{2} \|\partial^2_{tt} u_r\|^2_{L^2} + \frac{1}{2} |\partial E|^2(u_r) - \bar{H}_0^E(V_r^u)$ and set
	\begin{equation*}
		S \coloneqq \{r \in (a, b) \mid g \text{ is integrable in some neighbourhood around } r\},
	\end{equation*}
	noting that $S$ is open and thus consists of countably many disjoint open intervals. For almost all $s, t \in (a, b)$\footnote{The nullset of pairs $(s, t)$ which do not fall in either of those cases is contained in those pairs $(s, t)$ where either $s$ is the left endpoint of a constituting interval of $S$, or $t$ is the right endpoint of a constituting interval of $S$.}, there are two possibilities: (1.) $(s, t) \setminus S \neq \emptyset$, in which case the right hand side of \cref{eq:proof_of_ugs_wave_1} is $+\infty$---$g \not\in L^1((s, t))$, but $\int_s^t - \bar{H}_0^E(V_r^u) \, \diff r > -\infty$ because the integrand is in $L^1_{\mathrm{loc}}((a, b))$---and \cref{eq:proof_of_ugs_wave_1} is trivial; or (2.) $[s, t]$ is inside of one of the disjoint open intervals constituting $S$. In case (2.), there exists $a', b' \in (a, b) \cap \Q$ such that $[a', b'] \subseteq S$ and $a' < s < t < b'$, and hence $g \in L^1((a', b'))$ by compactness of $[a', b']$. It therefore suffices to treat case (2.), so we will proceed to fix $a', b' \in (a, b) \cap \Q$ such that $a' < b'$ and such that $g \in L^1((a', b'))$, and show that \cref{eq:proof_of_ugs_wave_1} is true for almost every $s, t \in (a', b')$. This suffices, since then the pairs $s, t$ for which \cref{eq:proof_of_ugs_wave_1} fails are either among the nullset excluded above, or are part of the exceptional set for some $a', b' \in (a, b) \cap \Q$; a countable union of null sets. Thus, we have that $u \in L^2((a', b'); H^2(\R^d))$ (using $\|\nabla^2 f\|_{L^2} = \|\Delta f\|_{L^2}$ on $\R^d$ and boundedness of $\|u_r\|_{L^2}$ on $[a', b'] \subseteq (a, b)$), that $\partial_t u \in L^2((a', b'); H^1(\R^d))$ and that $\partial^2_{tt} u \in L^2((a', b'); L^2(\R^d))$; we can also rewrite \cref{eq:proof_of_ugs_wave_1} (for almost every $s, t \in (a', b')$) as
	\begin{equation} \label{eq:proof_of_ugs_wave_1b}
		\langle \nabla u_s, \nabla \partial_t u_s \rangle_{L^2} - \langle \nabla u_t, \nabla \partial_t u_t \rangle_{L^2} \leq \int_{s}^{t} \frac{1}{2} \|\partial^2_{tt} u_r\|^2_{L^2} + \frac{1}{2} \| \Delta u_r \|^2_{L^2} - \|\nabla \partial_t u_r\|^2_{L^2} \, \diff r.
	\end{equation}
	We assume the following claim, which we will show at the end of the proof:
	\begin{equation}\label{eq:claim_in_proof_ugs_wave}
		\nabla u_t = \nabla u_s + \int_s^t \nabla \partial_t u_r \, \diff r, \qquad \text{for almost every } s, t \in (a', b') \text{ s.t. } s < t.
	\end{equation}
	Next, we set $v \coloneqq \nabla u$ and $w \coloneqq \partial_t u$---so that $\partial_t v = \nabla w$ by \cref{eq:claim_in_proof_ugs_wave}---and we aim to show that for almost every $s, t \in (a', b')$ such that $s < t$, we have that
	\begin{equation}\label{eq:product_rule_ugs_wave_proof}
		\langle v_s, \nabla w_s \rangle - \langle v_t, \nabla w_t \rangle = \int_s^t  \langle \nabla \cdot v_r, \partial_t w_r \rangle - \langle \partial_t v_r, \nabla w_r \rangle \, \diff r.
	\end{equation}
	Note that \cref{eq:product_rule_ugs_wave_proof} would be automatically true if $v$ and $\nabla w$ were absolutely continuous (as functions into $L^2(\R^d)$) using the product rule and integrating by parts. However, we cannot guarantee this regularity, so we argue by convolution instead: We set $u^\varepsilon \coloneqq u \ast_t \rho_\varepsilon$, $v^\varepsilon \coloneqq \nabla u^\varepsilon = v \ast_t \rho_\varepsilon$ and $w^\varepsilon \coloneqq w \ast_t \rho_\varepsilon$ for a smooth convolution kernel $\rho_\varepsilon$ with $\operatorname{supp} \rho_\varepsilon \subseteq (-\varepsilon, \varepsilon)$ and $\int \rho_\varepsilon = 1$. Then $v^\varepsilon$ and $\nabla w^\varepsilon$ are smooth functions into $L^2(\R^d)$, and so we have, for $\varepsilon$ small enough that $[s, t] \subseteq (a' + \varepsilon, b' - \varepsilon)$ and integrating by parts in $x$ (which is admissible since $u^\varepsilon_r \in H^2$ and $\partial_t w^\varepsilon_r = (\partial^2_{tt} u) \ast_t \rho_\varepsilon|_r \in H^1$),
	\begin{equation}\label{eq:product_rule_convoluted_ugs_wave_proof}
		\langle v^\varepsilon_s, \nabla w^\varepsilon_s \rangle - \langle v^\varepsilon_t, \nabla w^\varepsilon_t \rangle = \int_s^t  \langle \nabla \cdot v^\varepsilon_r, \partial_t w^\varepsilon_r \rangle - \langle \partial_t v^\varepsilon_r, \nabla w^\varepsilon_r \rangle \, \diff r.
	\end{equation}
	Furthermore, for $z \in \{v, \nabla \cdot v, \nabla w, \partial_t v, \partial_t w \}$, we have that $z^\varepsilon \to z$ strongly in $L^2_{\mathrm{loc}}((a', b'); L^2(\R^d))$ as $\varepsilon \to 0$. Thus, the right hand side of \cref{eq:product_rule_convoluted_ugs_wave_proof} converges to the right hand side of \cref{eq:product_rule_ugs_wave_proof} for \emph{all} $s, t \in (a', b')$ and any nullsequence $\{\varepsilon_n\}_{n \in \NN}$, while the left hand side of \cref{eq:product_rule_convoluted_ugs_wave_proof} converges to the left hand side of \cref{eq:product_rule_ugs_wave_proof} for \emph{almost all} $s, t \in (a', b')$ along some common subsequence $\{\varepsilon_{n_k}\}_{k \in \NN}$. Thus, \cref{eq:product_rule_ugs_wave_proof} holds for almost all $s, t \in (a', b')$ such that $s < t$. \\
	To prove the proposition from \cref{eq:product_rule_ugs_wave_proof}, note that we can now expand
	\begin{align}\label{eq:ugs_wave_end}
		\langle \nabla u_s, \nabla \partial_t u_s \rangle_{L^2} - \langle \nabla u_t, \nabla \partial_t u_t \rangle_{L^2} &= \int_s^t \langle \Delta u_r, \partial_{tt}^2 u_r \rangle_{L^2} - \langle \nabla \partial_t u_r, \nabla \partial_t u_r \rangle_{L^2} \, \diff r \\
		&\overset{(*)}{\leq} \int_s^t \frac{1}{2} \| \partial_{tt}^2 u_r \|_{L^2}^2 + \frac{1}{2} \| \Delta u_r \|_{L^2}^2 - \| \nabla \partial_t u_r \|^2_{L^2} \, \diff r,
	\end{align}
	with equality in $(*)$ if and only if $\partial^2_{tt} u_r = \Delta u_r$ for a.e. $r \in (s, t)$. If $u$ is a curve of Newtonian motion, the inverse inequality is true for almost all $s, t \in (a, b)$; combining this with \cref{eq:ugs_wave_end}, we have that for such a curve, for almost all $s, t \in (a, b)$, and for almost all $r \in (s, t)$, we have that $\partial^2_{tt} u_r = \Delta u_r$. Saturating $(a, b)$ with countably many such intervals $(s, t)$, we see that $\partial^2_{tt} u_r = \Delta u_r$ holds true for almost every $r \in (a, b)$. Lastly, for such curves of Newtonian motion the right hand side in \cref{eq:proof_of_ugs_wave_1} is finite for almost every $s, t \in (a, b)$, and thus $u_r \in H^2(\R^d, \R)$, $\partial_t u_r \in H^1(\R^d, \R)$ for a.e. $r \in (a, b)$. \\
	What remains to be proven is the claim \eqref{eq:claim_in_proof_ugs_wave}. To this end, note that the Bochner integral on the right hand side of \eqref{eq:claim_in_proof_ugs_wave} exists because $\|\nabla \partial_t u\|$ is in $L^1((a', b'))$ and $\nabla \partial_t u$ is strongly measurable by Pettis (weak measurability following from $u \in C^1((a, b), X)$; see, e.g.,~\cite[Theorem 2.104, Theorem 2.108]{fonsecaModernMethodsCalculus2007}). Furthermore, \cref{eq:claim_in_proof_ugs_wave} is true if tested against $C^\infty_c(\R^d; \R^d)$-functions---this follows from $u_t, \partial_t u_t \in H^1$, and an integration by parts argument on the equality $u_t = u_s + \int_s^t \partial_t u_r \, \diff r$ following from $u \in \mathrm{AC}$---yielding, by density of $C^\infty_c(\R^d; \R^d)$ in $L^2(\R^d; \R^d)$, \eqref{eq:claim_in_proof_ugs_wave}.
\end{proof}
\subsection{Discrete construction}
As we did for the general case, for every evolution horizon $T>0$ and every step-size $\tau>0$, we introduce the discrete trajectory $(U_{\tau,i})_{i=-1}^{N_\tau}$, with $N_\tau\coloneqq \big \lceil \frac{T}{\tau} \big\rceil$; as before, we write $N_\tau$ in place of the general $N$ to make the step-size dependence explicit.
Namely, we set
\begin{equation} \label{eq:wave_initial_states}
    U_{\tau,0} \coloneqq u_0, \qquad U_{\tau,-1} \coloneqq u_0 -\tau \dot{u}_0,
\end{equation}
and, for every $i\geq 0$, recalling \cref{rmk:objects_here_1}, we assign
\begin{equation} \label{eq:update_wave}
    Z_{\tau,i} = K^\star(U_{\tau, i-1}, U_{\tau, i})= 2U_{\tau,i} - U_{\tau,i-1} ,
    \quad U_{\tau, i+1} \in \argmin_{w \in L^2(\R^d,\R)} \left( \frac{1}{2\tau^2}\|w -  Z_{\tau,i} \|_{L^2}^2 + E(w)  \right).
\end{equation}
Using classical arguments of calculus of variations (i.e., the direct method and the Euler-Lagrange equations), it turns out that the minimization problem in \cref{eq:update_wave} has a unique minimizer $u^\star \in H^1(\R^d,\R)$, which is a weak solution of 
\begin{equation} \label{eq:Euler_Lagrange_update}
        -\tau^2 \Delta u^\star + u^\star = Z_{\tau,i},
\end{equation}
and which can be proved to belong to $H^2(\R^d,\R)$. We recall that, in accordance with \cref{eq:def_proximal}, for every $Y\in X$ and every $\delta>0$, $J_{\delta\cdot \tau}[Y]$ denotes the (unique) solution of the equation $-\delta \tau \Delta u^\star + u^\star = Y$.

\begin{remark}\label{rmk:gradient_sol_MMS}
    In analogy with what was observed in \cref{rmk:gradient_sol_wave}, if $Z_{\tau, i}\in H^1(\R^d, \R)$, then, for every $j=1,\ldots,d$, we have that
    \begin{equation}
        -\tau^2 \Delta (\partial_{x_j} u^\star) + \partial_{x_j} u^\star = \partial_{x_j} Z_{\tau, i},
    \end{equation}
    where we adopted the same notations as in \cref{eq:update_wave,eq:Euler_Lagrange_update}.
    In particular, this implies that, if $U_{\tau,0}=u_0$ and $U_{\tau,-1} = u_0 - \tau \dot{u}_0$ are in $H^1$, then, for every $j=1,\ldots,d$, $(\partial_{x_j}U_{\tau, i})_{i=-1}^{N_\tau}$ is also a Newtonian sequence for the Dirichlet energy $E$.
\end{remark}

\begin{remark}\label{rmk:objects_here_2}
As previously done in \cref{rmk:objects_here_1}, we provide the explicit expression of some objects of interest in the present framework.
We recall that the curves $u_\tau$ and $\tilde u_\tau$ are given by
\begin{align}
		u_\tau & \colon [0, \tau N_\tau] \to X,\qquad u_\tau(\tau i + \delta) = U_{\tau, i} \qquad & i \in \NN, \delta \in [0, \tau), \label{eq:piecewise_constant_interpolation_wave} \\
		\tilde{u}_\tau & \colon [0, \tau N_\tau] \to X,\qquad \tilde{u}_\tau(\tau i + \delta) = J_{\delta \cdot \tau}[K^\star(U_{\tau, i-1}, U_{\tau, i})] = J_{\delta \cdot \tau}[Z_{\tau,i}], \qquad & i \in \NN, \delta \in [0, \tau). \label{eq:de_giorgi_interpolation_wave}
\end{align}
Moreover, as regards the family of curves $\mathcal{V}^{\tau, \varepsilon}, \tilde{\mathcal{V}}^{\tau, \varepsilon}, \bar{\mathcal{V}}^{\tau, \varepsilon}$ originally defined in \cref{eq:velocity_interpolation,eq:fancy_velocity_interpolation,eq:bar_velocity_interpolation}, we have
	\begin{alignat}{3}
		&\mathcal{V}^{\tau, \varepsilon} \colon [0, \tau N_\tau)  \to \Tang[\varepsilon], \quad \mathcal{V}^{\tau, \varepsilon}_{t}  \in \Tang[t, \varepsilon], \quad  t \in [0, \tau N_\tau),\nonumber \\
        &\mathcal{V}^{\tau, \varepsilon}_{\tau i + \delta}  \coloneqq \bigg[\tau i + \delta + \sigma \mapsto U_{\tau,i-1} + \frac{2\delta + \sigma}{\tau} \big( U_{\tau,i}- U_{\tau,i-1} \big) \bigg]\quad & i \in \NN, \, \delta \in [0, \tau), \, \sigma \in [0, \varepsilon], \label{eq:velocity_interpolation_Dir} \\
        &\tilde{\mathcal{V}}^{\tau, \varepsilon} \colon [0, \tau N_\tau)  \to \Tang[\varepsilon], \quad \tilde{\mathcal{V}}^{\tau, \varepsilon}_t  \in \Tang[t, \varepsilon], \quad  t \in [0, \tau N_\tau),\nonumber \\
        &\tilde{\mathcal{V}}^{\tau, \varepsilon}_{\tau i + \delta}  \coloneqq \bigg[\tau i + \delta + \sigma \mapsto U_{\tau,i-1} + \frac{2\psi(\delta) + \sigma}{\tau} \big( U_{\tau,i}- U_{\tau,i-1} \big) \bigg]\quad & i \in \NN,\, \delta \in [0, \tau),\, \sigma \in [0, \varepsilon], \label{eq:fancy_velocity_interpolation_Dir} \\
        &\bar{\mathcal{V}}^{\tau, \varepsilon} \colon [0, \tau N_\tau)  \to \Tang[\varepsilon], \quad \bar{\mathcal{V}}^{\tau, \varepsilon}_t  \in \Tang[t, \varepsilon], \quad  t \in [0, \tau N_\tau),\nonumber \\
        &\bar{\mathcal{V}}^{\tau, \varepsilon}_{\tau i + \delta}  \coloneqq \bigg[\tau i + \delta + \sigma \mapsto U_{\tau,i-1} + \frac{\sigma}{\tau} \big( U_{\tau,i}- U_{\tau,i-1} \big) \bigg]\quad & i \in \NN, \, \delta \in [0, \tau), \, \sigma \in [0, \varepsilon]. \label{eq:bar_velocity_interpolation_Dir}
	\end{alignat}
    We observe that, in our case, geodesics can be arbitrarily extended. Then, we recall that
    \begin{equation} \label{eq:def_N_wave}
        N_\tau(x, y, z)  \coloneqq \frac{\|z - 2y +x\|_{L^2}}{\tau^2}, \qquad  x, y, z \in X, \  \tau > 0, 
    \end{equation}
    and that, for every $i \in \NN$ and $ \delta \in [0, \tau)$,
    	\begin{align}
		|u_\tau^{(1)}| & \colon [0, \tau N_\tau) \to \R_+, \qquad |u_\tau^{(1)}|(\tau i + \delta) = \frac{1}{\tau} \|U_{\tau, i} - U_{\tau, i+1}\|_{L^2}, \label{eq:velocity_norm_interpolation_wave} \\
		|u_\tau^{(2)}| & \colon [0, \tau N_\tau) \to \R_+, \qquad |u_\tau^{(2)}|(\tau i + \delta) = N_\tau(U_{\tau, i-1}, U_{\tau, i}, U_{\tau, i+1})= \| \Delta U_{\tau,i+1} \|_{L^2},  \label{eq:acceleration_norm_interpolation_wave} \\
		d(Eu_\tau) & \colon [0, \tau N_\tau] \to \R, \qquad d(Eu_\tau)(\tau i + \delta) = \frac{1}{\tau} (E(U_{\tau, i}) - E(U_{\tau, i-1})), \label{eq:energy_derivative_interpolation_wave}
	\end{align}
	and for every $i \in \NN$ and $ \delta \in [0, \tau)$\footnote{As in \cref{eq:def_G_tau}, we let $\delta$ run through the left-open interval $(0, \tau]$ in the definition of $|G_\tau|$ to avoid division by zero; in the rest of this paper, only its equivalent class in $L^1$ (and not its pointwise values) matter.},
	\begin{align}
		|G_\tau| & \colon [0, \tau N_\tau] \to \R_+, \qquad |G_\tau|(\tau i + \delta) = \frac{  \| J_{\delta \cdot \tau}[Z_{\tau,i}] - Z_{\tau, i} \|_{L^2} }{\delta \tau} = \|\Delta J_{\delta \cdot \tau}[Z_{\tau,i}] \|_{L^2}, \label{eq:def_G_tau_wave}
	\end{align}
    where $N_\tau$ is defined as in \cref{eq:def_N,eq:def_N_wave}, and $Z_{\tau,i}=K^\star(U_{\tau, i-1}, U_{\tau, i}) = 2 U_{\tau, i} - U_{\tau, i-1}$.
\end{remark}

We are in a position to establish a variational identity for the discrete trajectory.

\begin{proposition}[Discrete identity for the Dirichlet energy, cf.~\cref{prop:discrete_newtonian_motion}]\label{prop:discrete_newtonian_motion_Dir}
	Let $(U_{\tau, k})_{k=-1}^{N_\tau}$ be a Newtonian sequence in $X=L^2(\R^d, \R)$ constructed using \cref{eq:update_wave}, with $U_{\tau,-1},U_{\tau,0}\in H^1$, and let $\tilde{u}_\tau, \tilde{\mathcal{V}}^{\tau, \varepsilon}, |u_\tau^{(2)}|$, $d(Eu_\tau)$, $|G_\tau|$ be as in \eqref{eq:de_giorgi_interpolation_wave}, 
	\eqref{eq:fancy_velocity_interpolation_Dir}, \eqref{eq:acceleration_norm_interpolation_wave}, \eqref{eq:energy_derivative_interpolation_wave} and \eqref{eq:def_G_tau_wave}, respectively.
    Setting $t_\tau^j \coloneqq \tau \cdot j$ for all $j \in \NN$, we have, for all $i, j \in \NN$ such that $0 \leq i < j \leq N_\tau$,
	\begin{equation} \label{eq:discrete_newtonian_motion_Dir}
		d(Eu_\tau)(t_\tau^i) - d(Eu_\tau)(t_\tau^j) =  \frac{1}{2} \int_{t_\tau^i}^{t_\tau^j}   |u_\tau^{(2)}|^2(r)\, \diff r + \frac{1}{2} \int_{t_\tau^i}^{t_\tau^j} |G_\tau|^2(r)\, \diff r - \int_{t_\tau^i}^{t_\tau^j} \bar{H}^E_{0}(\tilde{\mathcal{V}}^{\tau, \varepsilon}_{r}) \, \diff r,
	\end{equation}
	where $\bar{H}_0^E$ is as in \cref{def:hessian_wave}.
	Furthermore, for almost all $t \in [0, T]$, we have
	\begin{equation} \label{eq:discrete_newtonian_motion_estimate_G_dir}
		|\partial E|(\tilde{u}_\tau(t)) = |G_\tau|(t) = \| \Delta \tilde u_\tau(t) \|_{L^2}.
	\end{equation}
\end{proposition}
\begin{proof}
	First, we note that for any such Newtonian sequence, we have that $E(U_{\tau, i}) < +\infty$ for all $-1 \leq i \leq N_\tau$, thus $U_{\tau, i} \in H^1$ and hence $\bar{H}_0^E$ as in \cref{def:hessian_wave} coincides with $H_0^E$ as defined in \cref{eq:def_H0} for all relevant geodesics.
    As observed for \cref{prop:discrete_newtonian_motion}, it is not restrictive to show the result for $j=i+1$.
    On the one hand, since the arguments used in the proof of \cref{prop:discrete_newtonian_motion} to deduce \cref{eq:proof_discrete_newtonian_motion_1} are completely general, we have that
	\begin{equation} \label{eq:der_balance_wave_1}
		\frac{1}{\tau} (E(Z) - E(U_{\tau, i+1})) = \frac{1}{2}\int_{t_\tau^i}^{t_\tau^{i+1}}  |u_\tau^{(2)}|^2(r)\, \diff r + \frac{1}{2}\int_{t_\tau^i}^{t_\tau^{i+1}}  |G_\tau|^2(r)\, \diff r,
	\end{equation}
    where $Z=K^\star(U_{\tau, i-1}, U_{\tau, i}) = 2 U_{\tau, i} - U_{\tau, i-1}$.
    On the other hand, since for every $r\in (t^i_\tau , t^{i+1}_\tau)$ the geodesic $\sigma \mapsto \tilde{\mathcal{V}}^{\tau, \varepsilon}_{r}(\sigma)$ is a straight line whose velocity is $\frac{\diff}{\diff \sigma} \tilde{\mathcal{V}}^{\tau, \varepsilon}_{r}(\sigma) = \frac1\tau \big( U_{\tau, i} - U_{\tau, i-1} \big)$ (see \cref{eq:fancy_velocity_interpolation_Dir}), from \cref{eq:def_H0_Dir} it follows that
    \begin{equation*}
        H_0^E \big( \tilde{\mathcal{V}}^{\tau, \varepsilon}_{r} \big) = \frac1{\tau^2} \| \nabla (U_{\tau, i} - U_{\tau, i-1}) \|_{L^2}^2
    \end{equation*}
    for every $r\in (t^i_\tau , t^{i+1}_\tau)$, yielding
    \begin{equation} \label{eq:der_balance_wave_2}
        \int_{t_\tau^i}^{t_\tau^{i+1}} H_0^E \big( \tilde{\mathcal{V}}^{\tau, \varepsilon}_{r} \big) \, \diff r = \frac1\tau \| \nabla (U_{\tau, i} - U_{\tau, i-1}) \|_{L^2}^2 = \frac1\tau
        \big( E(U_{\tau, i-1}) - 2E(U_{\tau,i}) + E(Z) \big).
    \end{equation}
    By combining \cref{eq:der_balance_wave_1,eq:der_balance_wave_2}, we obtain \cref{eq:discrete_newtonian_motion_Dir}. Finally, \cref{eq:discrete_newtonian_motion_estimate_G_dir} is a direct combination of \cref{rmk:gradient_syst} with \cref{eq:de_giorgi_interpolation_wave,eq:def_G_tau_wave}.
\end{proof}

In the next result, we show that the sum of the driving energy and the discrete kinetic energy is non-increasing along the discrete Newtonian motion. The proof does not actually rely on the specific form of the energy $E$, but rather on the fact that it is convex.

\begin{lemma}\label{lem:nonincreas_total_energy}
    Let $X$ be a Hilbert space, and let $E\colon X \to \R \cup \{+\infty\}$ be a proper and lower semi-continuous convex functional. 
    For every $\tau>0$ and $T>0$, given $U_{\tau,-1}, U_{\tau,0} \in X$, let $(U_{\tau,i})_{i=-1}^{N_\tau}$ (with $N_\tau\coloneqq \big \lceil {T}/{\tau} \big\rceil$) be the Newtonian sequence defined according to \cref{eq:update_wave}, and let us set $t_\tau^i\coloneqq \tau i$.
    Then, for every $i\geq 0$ we have that
    \begin{equation} \label{eq:nonincreas_total_energy}
        E\big(u_\tau(t_\tau^i)\big) + \frac12 \big| u^{(1)}_\tau \big|^2(t_\tau^{i-1}) \geq E\big(u_\tau(t_\tau^{i+1})\big) + \frac12 \big| u^{(1)}_\tau \big|^2(t_\tau^{i}),
    \end{equation}
    where $u_\tau$ and $|u^{(1)}_\tau|$ are defined in \cref{eq:piecewise_constant_interpolation_wave,eq:velocity_norm_interpolation_wave}, respectively.
    In particular, if $E$ is bounded from below and if $E\big(u_\tau(t_\tau^0)\big) + \frac12 \big| u^{(1)}_\tau \big|^2(t_\tau^{-1})<\infty$, then $t\mapsto E\big(u_\tau(t)\big)$ and $t\mapsto \frac12 \big| u^{(1)}_\tau \big|^2(t)$ are uniformly bounded on $[0,T]$.
\end{lemma}
\begin{proof}
    We recall that, if $x^\star \in X$ is such that $x^\star \in \arg \min_{x\in X} \left\{ \frac1{2\tau^2} \|x-y\|^2_X + E(x) \right\}$, then, by virtue of the \emph{sum rule for subdifferentials} (see, e.g.,~\cite[Theorem~23.8]{Rockafellar1970}) we have $0\in (x^\star - y) + \tau^2 \partial E(x^\star)$, which implies that $\frac1{\tau^2}(y - x^\star)\in \partial E(x^\star)$. Therefore, recalling \cref{eq:update_wave}, we deduce that 
    \begin{equation*}
    \begin{split}
        E(U_{\tau, i}) & \geq  E(U_{\tau, i+1}) + \frac1{\tau^2} \langle Z_{\tau, i} - U_{\tau,i+1}, U_{\tau,i} - U_{\tau,i+1} \rangle \\
         & \quad = E(U_{\tau, i+1}) + \frac1{2\tau^2} \left( \big\| U_{\tau, i+1} - U_{\tau, i} \big\|_X^2 + \big\| Z_{\tau, i} - U_{\tau, i+1} \big\|_X^2 - \big\| U_{\tau, i} - Z_{\tau, i} \big\|_X^2 \right),
    \end{split}
    \end{equation*}
    where we used the parallelogram identity for the second passage. Recalling that $U_{\tau, i} - Z_{\tau, i} = U_{\tau, i-1} - U_{\tau, i}$ (see again \cref{eq:update_wave}) and neglecting the contribution of $\big\| Z_{\tau, i} - U_{\tau, i+1} \big\|_X^2$ in the last inequality, we obtain
    \begin{equation*}
        E(U_{\tau, i}) + \frac1{2\tau^2} \big\| U_{\tau, i} - U_{\tau, i-1} \big\|_X^2 \geq
        E(U_{\tau, i+1}) + \frac1{2 \tau^2}\big\| U_{\tau, i+1} - U_{\tau, i} \big\|_X^2,
    \end{equation*}
    which, combined with the definitions in \cref{eq:piecewise_constant_interpolation_wave,eq:velocity_norm_interpolation_wave}, yields \cref{eq:nonincreas_total_energy}. The second part of the statement follows directly.
\end{proof}

In the case where $E$ is the Dirichlet energy on $X=L^2 (\R^d, \R)$, the previous result can be used to deduce further \textit{a priori} estimates.

\begin{corollary}[Energy dissipation inequality]\label{cor:nonincreas_total_energy}
    Let $E\colon X \to \R \cup \{+\infty\}$ be the Dirichlet energy, with $X=L^2 (\R^d, \R)$.
    For every $\tau>0$ and $T>0$, given $U_{\tau,-1}, U_{\tau,0} \in H^1(\R^d,\R)$, let $(U_{\tau,i})_{i=-1}^{N_\tau}$ (with $N_\tau\coloneqq \big \lceil {T}/{\tau} \big\rceil$) be the Newtonian sequence defined according to \cref{eq:update_wave}.
    Then, for every $i\geq 0$ and for every $j=1,\ldots,d$ we have that
    \begin{equation} \label{eq:nonincreas_total_energy_deriv}
        E\big( \partial_{x_j} U_{\tau,i} \big) + \frac1{2\tau^2} \big\| \partial_{x_j} U_{\tau,i} - \partial_{x_j} U_{\tau,i-1} \big\|^2_{L^2} \geq E\big( \partial_{x_j} U_{\tau,i+1} \big) + \frac1{2\tau^2} \big\| \partial_{x_j} U_{\tau,i+1} - \partial_{x_j} U_{\tau,i} \big\|^2_{L^2}.
    \end{equation}
    Moreover, if $U_{\tau,0} \in H^2(\R^d,\R)$, we have that
    \begin{equation} \label{eq:estimate_grad_laplacian}
        \max_{i=1,\ldots,N} \left\{ | U_{\tau, i} |_{H^2}^2 \, , \frac1{\tau^2} \| \nabla U_{\tau, i} - \nabla U_{\tau, i-1} \|_{L^2}^2 \right\} \leq
        | U_{\tau, 0} |_{H^2}^2 + \frac1{\tau^2} \| \nabla U_{\tau, 0} - \nabla U_{\tau, -1} \|_{L^2}^2,
    \end{equation}
    where $| u |_{H^2}^2 \coloneqq \sum_{1\leq j_1,j_2\leq d} \|\partial_{x_{j_1},x_{j_2}}^2 u \|_{L^2}^2 = \sum_{j=1}^d \| \nabla (\partial_{x_j} u)  \|_{L^2}^2$.
\end{corollary}

\begin{proof}
    As observed in \cref{rmk:gradient_sol_MMS}, as $U_{\tau,-1}, U_{\tau,0} \in H^1(\R^d,\R)$, it turns out that $(\partial_{x_j} U_{\tau,i})_{i=-1}^N$ is a Newtonian sequence in $L^2$ with respect to the Dirichlet energy, for every $j=1,\ldots, d$. Hence, by virtue of \cref{lem:nonincreas_total_energy}, we deduce the relation in \cref{eq:nonincreas_total_energy_deriv}. Moreover, for every $j=1,\ldots, d$ and for every $i\geq 1$, we further have that
    \begin{equation*}
    \begin{split}
        \frac{\big\|\nabla( \partial_{x_j} U_{\tau,0} ) \big\|_{L^2}^2}{2} + \frac {\big\| \partial_{x_j} U_{\tau,0} - \partial_{x_j} U_{\tau,-1} \big\|^2_{L^2}}{2\tau^2} & \geq 
         \frac{\big\|\nabla( \partial_{x_j} U_{\tau,i} ) \big\|_{L^2}^2}2  + \frac{\big\| \partial_{x_j} U_{\tau,i} - \partial_{x_j} U_{\tau,i-1} \big\|^2_{L^2}}{2\tau^2} .
    \end{split}
    \end{equation*}
    Then, summing over the index $j$ in the previous formula, we deduce \cref{eq:estimate_grad_laplacian}.
\end{proof}

We are now in a position to establish the first block of the \textit{a priori} bounds (see \cref{lem:apriori_bounds}).

\begin{lemma}[A priori bounds I for the Dirichlet energy, cf.~\cref{lem:apriori_bounds}]\label{lem:apriori_bounds_wave}
    Let $u_0 \in H^2(\R^d,\R)$ and $\dot{u}_0\in H^1(\R^d,\R)$, both with compact support, and let $T > 0$. For every $\tau >0$, we set $U_{\tau, 0} \coloneqq u_0$ and $U_{\tau,-1} \coloneqq u_0 -\tau \dot{u}_0$, and let $(U_{\tau,i})_{i=-1}^{N_\tau}$ (with $N_\tau\coloneqq \big \lceil {T}/{\tau} \big\rceil$) be the Newtonian sequence constructed according to \cref{eq:update_wave}. Finally, let $u_\tau, |u_\tau^{(1)}|, |u_\tau^{(2)}|, d(Eu_\tau)$ be the functions defined as in \cref{eq:piecewise_constant_interpolation_wave,eq:velocity_norm_interpolation_wave,eq:acceleration_norm_interpolation_wave,eq:energy_derivative_interpolation_wave}.
    Then, there exists a constant $C>0$ such that, for every $\tau>0$ small enough, for every $0\leq i < {N_\tau}$ and for every $t \in [0,N_\tau\tau)$, we have
    \begin{equation} \label{eq:apriori_1_wave_a}
        \frac1\tau \| U_{\tau,i} - U_{\tau,i+1}\|_{L^2} \leq C, \quad 
        |u^{(1)}_\tau|(t) \leq C, \qquad  E( U_{\tau, i} ) \leq \frac{C^2}{2}, \quad E\big(u_\tau(t)\big) \leq \frac{C^2}{2}, 
    \end{equation}
    \begin{equation} \label{eq:apriori_1_wave_b}
        \| \Delta U_{\tau,i} \|_{L^2} \leq C, \quad |u^{(2)}_\tau|(t) \leq C, \qquad
       \left| \frac{E( U_{\tau, i} )-E( U_{\tau, i-1} )}{\tau} \right| \leq C^2,  \quad  |d(Eu_\tau)(t) | \leq C^2.
    \end{equation}
    For $1 \leq i \leq N_\tau$ and $t \in [\tau, \tau N_\tau]$, we additionally have
    \begin{equation} \label{eq:apriori_1_wave_c}
        \frac{E( U_{\tau, i} )-2E( U_{\tau, i-1}) +E( U_{\tau, i-2})}{\tau}\leq C^2 \tau, \quad  d(Eu_\tau)(t) - d(Eu_\tau)(t -\tau) \leq C^2 \tau.   
    \end{equation}
    Moreover, for every $\tau>0$ small enough, for every $0\leq i < {N_\tau}$ and for every $t \in [0,T]$, we have
    \begin{equation} \label{eq:apriori_1_wave_d}
        \| U_{\tau,i} - u_0 \|_{L^2} \leq C \tau i, \quad \| u_{\tau}(t) - u_0 \|_{L^2} \leq C t.
    \end{equation}
\end{lemma}
\begin{proof}
    Owing to \cref{lem:nonincreas_total_energy}---recall also the definitions in \cref{eq:piecewise_constant_interpolation_wave,eq:velocity_norm_interpolation_wave}---, we deduce that \cref{eq:apriori_1_wave_a} holds with the constant $C_1 = \sqrt{2} \sqrt{E(u_0) + \frac{1}{2} \|\dot{u}_0\|_{L^2}^2}$.
    Then, we observe that $\| U_{\tau,i} - u_0 \|_{L^2} \leq \sum_{k=1}^i\| U_{\tau,k} - U_{\tau,k-1} \|_{L^2} \leq C_1\tau i$, and we establish \cref{eq:apriori_1_wave_d}.
    Moreover, by using that $u_0 \in H^2, \dot{u}_0\in H^1$ and by applying \cref{cor:nonincreas_total_energy}, we deduce the first two inequalities in \cref{eq:apriori_1_wave_b} with the constant $C_2 = \sqrt{|u_0|_{H^2}^2+\|\nabla \dot{u}_0\|_{L^2}^2 }$ (see also \cref{eq:acceleration_norm_interpolation_wave}). Then, for $1\leq i\leq N_\tau$, we observe that
    \begin{equation*}
        \begin{split}
            \left| \frac{E( U_{\tau, i} )-E( U_{\tau, i-1} )}{\tau} \right| &
            %= \frac{1}{2\tau} \big| \langle \nabla U_{\tau, i} - \nabla U_{\tau, i-1}, \nabla U_{\tau, i} + \nabla U_{\tau, i-1} \rangle \big|
            \leq \frac{1}{2\tau} \big\| \nabla U_{\tau, i} - \nabla U_{\tau, i-1} \big\|_{L^2} \big\| \nabla U_{\tau, i} + \nabla U_{\tau, i-1} \big\|_{L^2} \\
            & \leq \frac{C_2}{\sqrt 2}\left( \sqrt{E(U_{\tau, i})} + \sqrt{E(U_{\tau,i-1})}  \right) \leq C_1C_2,
        \end{split}
    \end{equation*}
    For $i=0$ and $\tau\leq1$, initialization gives $|d(Eu_\tau)(0)|\leq\|\nabla u_0\|_{L^2}\|\nabla\dot u_0\|_{L^2}+\frac12\|\nabla\dot u_0\|_{L^2}^2\leq C_1C_2+\frac12 C_2^2$. Together, these establish the last two inequalities in \cref{eq:apriori_1_wave_b} (see also \cref{eq:energy_derivative_interpolation_wave}).
    Then, by virtue of \cref{eq:discrete_newtonian_motion_Dir}---where we neglect the negative terms---and \cref{eq:der_balance_wave_2}, we compute, for $1\leq i\leq N_\tau$,
    \begin{equation*}
        \begin{split}
            \frac1\tau \big( E( U_{\tau, i} )-2E( U_{\tau, i-1}) + E( U_{\tau, i-2}) \big) &= 
            \big( d(E u_\tau)(i\tau) - d(E u_\tau)((i-1)\tau) \big)\\
            &  \leq \frac1\tau \| \nabla (U_{\tau, i-1} - U_{\tau, i-2}) \|_{L^2}^2 \leq C_2^2 \tau,
        \end{split}
    \end{equation*}
    where we used \cref{cor:nonincreas_total_energy} for $i\geq2$ and $U_{\tau,0}-U_{\tau,-1}=\tau\dot u_0$ for $i=1$ in the last passage. 
    The second bound in~\eqref{eq:apriori_1_wave_c} follows from the last inequality and from the definition in~\cref{eq:energy_derivative_interpolation_wave}.
    Finally, we set $C\coloneqq C_1+C_2$.
\end{proof}

In the next result, we provide the second block of auxiliary estimates. Before proceeding, we recall that, in the present framework, given $\xi^1, \xi^2 \in \Tang[t, \varepsilon]$, we have
\begin{equation*}
    d_{\Tang[t, \varepsilon]} (\xi^1, \xi^2 ) = \sup_{\sigma\in [0,\varepsilon]} \| \xi^1_{t+\sigma} - \xi^2_{t+\sigma} \|_{L^2}.
\end{equation*}

\begin{lemma}[A priori bounds II for the Dirichlet energy, cf.~\cref{lem:apriori_bounds_2}]\label{lem:apriori_bounds_2_wave}
    For every $\tau >0$, let us consider $u_0 \in H^2(\R^d,\R)$ and $\dot{u}_0\in H^1(\R^d,\R)$, let us set $U_{\tau, 0} \coloneqq u_0$ and $U_{\tau,-1} \coloneqq u_0 -\tau \dot{u}_0$, and let $(U_{\tau,i})_{i=-1}^{N_\tau}$ be the Newtonian sequence constructed according to \cref{eq:update_wave}.
    Finally, let $u_\tau, \tilde u_\tau, \mathcal{V}^{\tau, \varepsilon}, \tilde{\mathcal{V}}^{\tau, \varepsilon}, \bar{\mathcal{V}}^{\tau, \varepsilon}$ be the functions defined, respectively, in \cref{eq:piecewise_constant_interpolation_wave,eq:de_giorgi_interpolation_wave,eq:velocity_interpolation_Dir,eq:fancy_velocity_interpolation_Dir,eq:bar_velocity_interpolation_Dir}.
    Then, there exists a constant $C>0$ such that, for every $\tau>0$, $\varepsilon>0$ and $t\geq 0$ we have
    \begin{equation*}
        \| u_\tau(t) - \tilde u_\tau(t) \|_{L^2} \leq C\tau (1 + \tau), \qquad 
        d_{\Tang[t, \varepsilon]}(\mathcal{V}^{\tau, \varepsilon}_t, \tilde{\mathcal{V}}^{\tau, \varepsilon}_t) \leq \frac{C}{4}\tau,
        \qquad d_{\Tang[t, \varepsilon]}(\mathcal{V}^{\tau, \varepsilon}_t, \bar{\mathcal{V}}^{\tau, \varepsilon}_t) \leq 2C\tau.
    \end{equation*}
\end{lemma}
\begin{proof}
    Given $\tau>0$ and $t\geq 0$, let us consider $i\geq 0$ and $\delta\in [0,\tau)$ such that $t=i\tau +\delta$. 
    Then, recalling \cref{eq:update_wave,eq:de_giorgi_interpolation_wave}, we have
    \begin{equation*}
        \begin{split}
            \| u_\tau(t) - \tilde u_\tau(t) \|_{L^2} &= \big\| U_{\tau, i} -  J_{\delta\tau}[Z_{\tau,i}] \big\|_{L^2} 
            \leq \big\| U_{\tau, i} -  Z_{\tau,i} \big\|_{L^2} + \big\| Z_{\tau,i} -  J_{\delta\tau}[Z_{\tau,i}] \big\|_{L^2} \label{eq:apriori_1_1_wave}\\
            & \leq \big\| U_{\tau, i} -  U_{\tau, i-1} \big\|_{L^2} + \big\| Z_{\tau,i} -  J_{\tau^2}[Z_{\tau,i}] \big\|_{L^2} \leq C\tau (1 + \tau),
        \end{split}
    \end{equation*}
    where we used that $J_{\tau^2}[Z_{\tau,i}] - Z_{\tau,i} = \tau^2 \Delta U_{\tau,i+1}$ and the estimate established in \cref{lem:apriori_bounds_wave}. \\
    Then, recalling \cref{eq:velocity_interpolation_Dir,eq:fancy_velocity_interpolation_Dir} and writing $t=i\tau +\delta$ as above, we observe that
    \begin{equation*}
        \begin{split}
            d_{\Tang[t, \varepsilon]}(\mathcal{V}^{\tau, \varepsilon}_t, \tilde{\mathcal{V}}^{\tau, \varepsilon}_t) = |\psi^\tau(\delta)-\delta|\frac{2}{\tau} \big\| U_{\tau, i} -  U_{\tau, i-1} \big\|_{L^2} \leq \frac{C}{4}\tau
        \end{split}
    \end{equation*}
    for every $\varepsilon > 0$, where we used \cref{lem:apriori_bounds_wave} and the fact that $\max_{\delta\in [0,\tau]}|\psi^\tau(\delta)-\delta| = \tau/8$ (see \cref{eq:def_phi} and \cref{rmk:varphi_tau}).\\
    Finally, from \cref{eq:bar_velocity_interpolation_Dir}, as $\delta\in [0,\tau]$, we deduce that
    \begin{equation*}
        \begin{split}
            d_{\Tang[t, \varepsilon]}(\mathcal{V}^{\tau, \varepsilon}_t, \bar{\mathcal{V}}^{\tau, \varepsilon}_t) = \frac{2 \delta}{\tau} \big\| U_{\tau, i} -  U_{\tau, i-1} \big\|_{L^2} \leq 2 C \tau
        \end{split}
    \end{equation*}
    for every $\varepsilon > 0$.
\end{proof}

\begin{lemma}[A priori bounds III for the Dirichlet energy]\label{lem:apriori_3_wave}
	For every $\tau \in (0,1)$, let us consider $u_0 \in H^2(\R^d,\R)$ and $\dot{u}_0\in H^1(\R^d,\R)$, let us set $U_{\tau, 0} \coloneqq u_0$ and $U_{\tau,-1} \coloneqq u_0 -\tau \dot{u}_0$, and let $(U_{\tau,i})_{i=-1}^{N_\tau}$ be the Newtonian sequence constructed according to \cref{eq:update_wave}.
	Let $A \subset \R^d$ be such that $\supp(u_0) \cup \supp(\dot{u}_0) \subseteq A$.
	Then, for each $T > 0$, there exist $L^0$ and $v > 0$, independent of $\tau$, such that, denoting for a set $A \subseteq \R^d$ its enlargement by $r$ as $A_r \coloneqq \{x \in \R^d \mid \dist(x, A) \leq r \}$, we have that
	\begin{equation} \label{eq:sol_wave_Equation_decay}
		\int_{A_{v (t + \tau k)}^c} | U_{\tau, \lfloor t / \tau \rfloor }|^2(x) \, \diff x  \leq L^0 2^{-\lfloor t / \tau \rfloor  - k}, \qquad \forall t \in [0, T],\ k \in \NN.
	\end{equation}
	In particular, if $A$ is compact, then, for each $T > 0$, there exists $M_T > 0$, independent of $\tau$, such that
	\begin{equation} \label{eq:sol_wave_equation_in_weighted_l2}
		\int_{\R^d} |U_{\tau, i}|^2(x) (1 + |x|^2) \, \diff x \leq M_T, \qquad \int_{\R^d} |U'_{\tau, i}|^2(x) (1 + |x|^2) \, \diff x \leq M_T
	\end{equation}
	for $\tau$ small enough and $ i \leq N_\tau$, where we set $U'_{\tau, i} \coloneqq \frac{U_{\tau, i} - U_{\tau, i-1}}{\tau}$.
\end{lemma}
The proof is deferred to \cref{sec:proof-of-apriori-3-wave}.

The next result concerns a bound on the velocity interpolation of the Newtonian sequences.

\begin{proposition}[Newtonian sequences are proto-\ACt, cf.~\cref{prop:newtonian_sequences_are_proto_ac,prop:compactness_discrete_solutions_geodesic_space}]\label{prop:newtonian_sequences_are_proto_ac_wave}
	    For every $\tau >0$, let us consider $u_0 \in H^2(\R^d,\R)$ and $\dot{u}_0\in H^1(\R^d,\R)$, let us set $U_{\tau, 0} \coloneqq u_0$ and $U_{\tau,-1} \coloneqq u_0 -\tau \dot{u}_0$, and let $(U_{\tau,i})_{i=-1}^{N_\tau}$ be the Newtonian sequence constructed according to \cref{eq:update_wave}.
        Finally, let $\bar{\mathcal{V}}^{\tau, \varepsilon}, |u_\tau^{(2)}|$ be defined as in \cref{eq:bar_velocity_interpolation_Dir,eq:acceleration_norm_interpolation_wave}. \\
        Then, for every $i, j \in \NN$ such that $1 \leq i \leq j < N_\tau$, for every $\varepsilon>0$, and for any $T >0$ such that $\tau \leq T - t^j_\tau \leq T - t^i_\tau$, we have
	\begin{equation*}
		\big\| \bar{\mathcal{V}}^{\tau, \varepsilon}_{t^i_\tau}(T) - \bar{\mathcal{V}}^{\tau, \varepsilon}_{t^j_\tau}(T) \big\|_{L^2} \leq \int_{t^{i}_\tau}^{t^{j}_\tau} |u_\tau^{(2)}|(r)(T-r) \, \diff r ,
	\end{equation*}
    and
    \begin{equation*}
		d_{\Tang[0, \varepsilon]}\Big(\mathrm{Shift}_0 \bar{\mathcal{V}}^{\tau, \varepsilon}_{t^i_\tau}, \mathrm{Shift}_0 \bar{\mathcal{V}}^{\tau, \varepsilon}_{t^j_\tau}\Big) \leq \int_{t_\tau^{i}}^{t_\tau^{j}} C + \varepsilon |u_\tau^{(2)}|(r)  \,\diff r,
	\end{equation*}
    where we set $t^k_\tau \coloneqq  k \tau$, the operator $\mathrm{Shift}_0$ was introduced  in \cref{def:Shift}, and $C>0$ is provided by \cref{lem:apriori_bounds_wave}.
\end{proposition}
\begin{proof}
    As usual, we prove the statement for $j=i+1$, and the general case follows by applying the triangle inequality.
    We observe that
    \begin{equation*}
        \begin{split}
            \bar{\mathcal{V}}^{\tau, \varepsilon}_{t^{i+1}_\tau}(T) - \bar{\mathcal{V}}^{\tau, \varepsilon}_{t^i_\tau}(T) &= U_{\tau, i} + \frac{T-t_\tau^{i+1}}{\tau} \big( U_{\tau, i+1} -U_{\tau, i} \big) - U_{\tau, i-1}- \frac{T-t_\tau^{i}}{\tau} \big( U_{\tau, i} -U_{\tau, i-1} \big) \\
            &= \frac{T-t_\tau^{i+1}}{\tau} \big( U_{\tau, i+1} - 2 U_{\tau, i}  +U_{\tau, i-1} \big),
        \end{split}
    \end{equation*}
    so that
    \begin{equation*}
        \begin{split}
            \big\| \bar{\mathcal{V}}^{\tau, \varepsilon}_{t^{i+1}_\tau}(T) - \bar{\mathcal{V}}^{\tau, \varepsilon}_{t^{i}_\tau}(T) \big\|_{L^2}  & = |u_\tau^{(2)}|(t^{i}_\tau)(T-t_\tau^{i+1})\tau = \int_{t^{i}_\tau}^{t^{i+1}_\tau} |u_\tau^{(2)}|(r)(T-t_\tau^{i+1}) \, \diff r  \\
            &\leq \int_{t^{i}_\tau}^{t^{i+1}_\tau} |u_\tau^{(2)}|(r)(T-r) \, \diff r .
        \end{split}
    \end{equation*}
    Concerning the second part of the statement, we compute
    \begin{equation*}
            \mathrm{Shift}_0 \bar{\mathcal{V}}^{\tau, \varepsilon}_{t^{i+1}_\tau}(\sigma) - \mathrm{Shift}_0 \bar{\mathcal{V}}^{\tau, \varepsilon}_{t^{i}_\tau}(\sigma) =
            U_{\tau, i} - U_{\tau, i-1} + \frac{\sigma}{\tau} \big( U_{\tau, i+1} -2 U_{\tau, i} + U_{\tau, i-1}\big),
    \end{equation*}
    and we get, for every $\sigma\in [0,\varepsilon]$,
    \begin{equation*}
            \big\| \mathrm{Shift}_0 \bar{\mathcal{V}}^{\tau, \varepsilon}_{t^{i+1}_\tau}(\sigma) - \mathrm{Shift}_0 \bar{\mathcal{V}}^{\tau, \varepsilon}_{t^{i}_\tau}(\sigma) \big\|_{L^2} \leq  |u_\tau^{(1)}|(t_\tau^{i-1}) \tau
            + \sigma \tau |u_\tau^{(2)}|(t_\tau^{i})
            \leq \int_{t_\tau^{i}}^{t_\tau^{i+1}} C + \varepsilon |u_\tau^{(2)}|(r)  \,\diff r,
    \end{equation*}
    where we used \cref{lem:apriori_bounds_wave} to bound $ |u_\tau^{(1)}|(t_\tau^{i-1}) \leq C$.
\end{proof}
\begin{proposition}[Compactness of discrete solutions into space of geodesics, cf.~\cref{prop:compactness_discrete_solutions_geodesic_space}]\label{prop:compactness_discrete_solutions_geodesic_space_wave}
	Let us consider $k \in \NN$, and let $(\tau_i)_{i=1}^\infty$ be a sequence of step sizes such that $\tau_i \to 0$ as $i \to \infty$. Let furthermore $u_0 \in H^{k+2}(\R^d,\R)$ and $\dot{u}_0\in H^{k+1}(\R^d,\R)$, both with compact support. For every $\tau_i >0$, we set $U_{\tau, 0} \coloneqq u_0$ and $U_{\tau,-1} \coloneqq u_0 -\tau \dot{u}_0$, and let $(U_{\tau,i})_{i=-1}^{\infty}$ be the Newtonian sequence constructed according to \cref{eq:update_wave}. 
    Finally, let $\bar{\mathcal{V}}^{\tau_i, \varepsilon}, {\mathcal{V}}^{\tau_i, \varepsilon}, \tilde{\mathcal{V}}^{\tau_i, \varepsilon}$ be the functions defined as in \cref{eq:bar_velocity_interpolation_Dir,eq:fancy_velocity_interpolation_Dir,eq:velocity_interpolation_Dir}.
	Then, for any $\varepsilon > 0$, there exists---without relabeling---a subsequence $(\tau_i)_{i=1}^\infty$ and a curve $V\colon [0, \infty) \to \Tang$ such that for all $t$, $V_t \in \Tang[t, \varepsilon]$ and such that, for all $t \in [0, \infty)$,
	\begin{equation} \label{eq:compactness_in_geodesic_space_wave}
		V_t = \lim_{i \to \infty} \bar{\mathcal{V}}^{\tau_i, \varepsilon}_{t} = \lim_{i \to \infty} \mathcal{V}^{\tau_i, \varepsilon}_{t} = \lim_{i \to \infty} \tilde{\mathcal{V}}^{\tau_i, \varepsilon}_{t} \qquad \text{ in } \mathcal{T}_{t, \varepsilon} H^k.
	\end{equation}
\end{proposition}
\begin{proof}
	To show the convergences in \cref{eq:compactness_in_geodesic_space_wave}, it suffices to show that for each multiindex $\alpha$ such that $|\alpha| \leq k$, we have that
	\begin{equation} \label{eq:compactness_in_geodesic_space_wave_transformed}
		D^\alpha V_t = \lim_{i \to \infty}  D^\alpha\bar{\mathcal{V}}^{\tau_i, \varepsilon}_{t} = \lim_{i \to \infty}  D^\alpha \mathcal{V}^{\tau_i, \varepsilon}_{t} = \lim_{i \to \infty} D^\alpha \tilde{\mathcal{V}}^{\tau_i, \varepsilon}_{t} \qquad \text{ in } \mathcal{T}_{t, \varepsilon} L^2,
	\end{equation}
	where, for a geodesic $\eta\colon \sigma \mapsto u + \sigma v$ in $L^2$ such that $u, v \in H^k$, we set $D^\alpha \eta\colon \sigma \mapsto D^\alpha u + \sigma D^\alpha v$.
	By repeatedly applying \cref{rmk:gradient_sol_MMS}, we can reduce \cref{eq:compactness_in_geodesic_space_wave_transformed} to the case where $|\alpha| = 0$. In other words, in the rest of the proof, we may assume without loss of generality that $k=0$, and we set $X=H^0=L^2$. \\
	From here on, the proof proceeds similarly as in \cref{prop:compactness_discrete_solutions_geodesic_space}; we will only sketch it. We only prove the first convergence in \cref{eq:compactness_in_geodesic_space_wave}; the others follow then from \cref{lem:apriori_bounds_2_wave}. By a diagonal argument, it suffices to show such a convergence for each time interval $[0, T]$; we let such a $T$ be fixed and set $N_\tau\coloneqq \big \lceil \frac{T}{\tau} \big\rceil + 1$. Furthermore, we may instead prove existence of some $W\colon [0, T] \to \Tang[0, \varepsilon]$ such that, for all $t \in [0, T]$,
	\begin{equation*}
		W_t = \lim_{i \to \infty} \mathrm{Shift}_0 \bar{\mathcal{V}}^{\tau_i, \varepsilon}_{t} \qquad \text{ in } \Tang[0, \varepsilon],
	\end{equation*}
	since then we may set $V_t = \mathrm{Shift}_t W_t$. As in the proof of \cref{prop:compactness_discrete_solutions_geodesic_space}, we would like to apply an Arzelà-Ascoli-type argument to obtain the existence of the limit curve $W\colon [0, T] \to \Tang[0, \varepsilon]$; the key to this end is finding some compact $K_{\mathcal{T}} \subseteq \Tang[0, \varepsilon]$ such that for all $i$ and all $t \in [0, T]$, $\mathrm{Shift}_0 \bar{\mathcal{V}}^{\tau_i, \varepsilon}_{t} \in K_{\mathcal{T}}$. As in \cref{eq:compactness_in_geodesic_space}, we need to find a compact $K_X \subseteq X$ such that for all $i$ and all $t \in [0, T]$ and for all $\sigma \in [t, t + \varepsilon]$, $\bar{\mathcal{V}}^{\tau_i, \varepsilon}_{t}(\sigma) \in K_X$. To this end, it suffices to find a uniform bound $M$ such that
	\begin{equation}
		\int_{\R^d} \bar{\mathcal{V}}^{\tau_i, \varepsilon}_{t}(\sigma)(x)^2 (1 + |x|^2) \, \diff x + \int_{\R^d} |\nabla \bar{\mathcal{V}}^{\tau_i, \varepsilon}_{t}(\sigma)(x)|^2 \, \diff x \leq M,
	\end{equation}
	since the image of such bounded sets is precompact in $L^2$. However, since $u_0 \in H^2(\R^d, \R)$ and $\dot{u}_0 \in H^1(\R^d, \R)$, we can apply \cref{lem:apriori_bounds_wave} to obtain a bound on the second summand above from \cref{eq:apriori_1_wave_a,eq:apriori_1_wave_b}; for the first summand, we may again invoke \cref{lem:apriori_3_wave} to obtain a uniform bound. The rest of the proof proceeds as in \cref{prop:compactness_discrete_solutions_geodesic_space}.
\end{proof}
\begin{proposition}[Compactness in the space of \ACt-curves, cf.~\cref{prop:compactness_discrete_solutions_ac2}]\label{prop:compactness_discrete_solutions_ac2_wave}
	Let $(\tau_i)_{i=1}^\infty$ be a sequence of step sizes such that $\tau_i \to 0$ as $i \to \infty$; let furthermore $u_0 \in H^3(\R^d,\R)$ and $\dot{u}_0\in H^2(\R^d,\R)$, both with compact support. For every $\tau_i >0$, we set $U_{\tau, 0} \coloneqq u_0$ and $U_{\tau,-1} \coloneqq u_0 -\tau \dot{u}_0$, and let $(U_{\tau,i})_{i=-1}^\infty$ be the Newtonian sequence constructed according to \cref{eq:update_wave}. Finally, let $u_\tau, |u_\tau^{(2)}|, d(Eu_\tau), \bar{\mathcal{V}}^{\tau_i, \varepsilon}, {\mathcal{V}}^{\tau_i, \varepsilon}, \tilde{\mathcal{V}}^{\tau_i, \varepsilon}$ be the functions defined as in \cref{eq:piecewise_constant_interpolation_wave,eq:acceleration_norm_interpolation_wave,eq:energy_derivative_interpolation_wave,eq:bar_velocity_interpolation_Dir,eq:velocity_interpolation_Dir,eq:fancy_velocity_interpolation_Dir}.
	Then, for any $\varepsilon > 0$, there exists---without relabeling---a subsequence $(\tau_i)_{i=1}^\infty$ and an \ACt-curve $\bar{u}$, a function $\varphi\colon [0, +\infty) \to \R$ with locally finite total variation and a function $A \in L^2_{loc}([0, +\infty), \R)$ such that
	\begin{gather}
		V^{\bar{u}}_t = \lim_{i \to \infty} \bar{\mathcal{V}}^{\tau_i, \varepsilon}_{t} = \lim_{i \to \infty} \mathcal{V}^{\tau_i, \varepsilon}_{t} = \lim_{i \to \infty} \tilde{\mathcal{V}}^{\tau_i, \varepsilon}_{t} \quad \text{ in } \mathcal{T}_{t, \varepsilon}H^1, \text{ for all } t \in [0, +\infty), \label{eq:compactness_1_wave}\\
		\bar{u}(t) = \lim_{i \to \infty} u_{\tau_i}(t) = \lim_{i \to \infty} \tilde{u}_{\tau_i}(t), \quad \text{ in } H^2, \text{ for all } t \in [0, +\infty), \label{eq:compactness_2_wave}\\
		\varphi(t) \coloneqq \lim_{i \to \infty} d(Eu_{\tau_i})(t) = \frac{\diff}{\diff T} E(V^{\bar{u}}_t(T)) \Big|_{T=t} \quad \text{ for all } t \in [0, +\infty), \label{eq:compactness_3_wave}\\
		|u^{(2)}_{\tau_i}| \rightharpoonup A, \quad A(t) \geq |\bar{u}''|(t) \quad \text{ for a.e. } t \in [0, +\infty), \label{eq:compactness_4_wave}\\
		\liminf_{i \to \infty} |G_{\tau_i}|(t) \geq |\partial^- E|(\bar{u}(t)) \quad \text{ for almost all } t \in [0, +\infty), \label{eq:compactness_5_wave}
	\end{gather}
	where $|G_{\tau_i}|$ is defined as in \cref{eq:def_G_tau_wave}.
\end{proposition}
\begin{proof}
	The proof from \cref{prop:compactness_discrete_solutions_ac2} goes through; the convergence in \cref{eq:compactness_1_wave} comes from applying \cref{prop:compactness_discrete_solutions_geodesic_space_wave} with $k=1$. For the first convergence in \cref{eq:compactness_2_wave}, we need to use the bound
	\begin{equation*}
		\| \Delta \partial_j U_{\tau,i} \|_{L^2} \leq C
	\end{equation*}
	obtained by applying \cref{lem:apriori_bounds_wave} to $\partial_j u_0$ and $\partial_j \dot{u}_0$, together with \cref{lem:apriori_3_wave} to obtain strong compactness of $u_\tau$ in $H^2$.
	For the second convergence in \cref{eq:compactness_2_wave}, we first use \cref{cor:nonincreas_total_energy} to get a global bound $D$ such that $D \geq \|\partial^{2}_{jk} \frac{U_{\tau_i, l+1} - U_{\tau_i, l}}{\tau_i}\|_{L^2}$ for all $j, k$ and $l$. For $t=l\tau_i+\delta$ with $0\leq\delta<\tau_i$, we proceed to estimate similarly to the first part of the proof of \cref{lem:apriori_bounds_2_wave}:
	\begin{align*}
		\|\partial^{2}_{jk}(u_{\tau_i}(t) - \tilde{u}_{\tau_i}(t))\|_{L^2} \leq \|\partial^2_{jk} (U_{\tau_i, l} - U_{\tau_i, l-1})\|_{L^2} + \|\partial^2_{jk}(Z_{\tau_i, l} - J_{\delta\tau_i}[Z_{\tau_i, l}] )\|_{L^2},
	\end{align*}
	The first part of the sum above is bounded by $\tau_i D$, and for the second part, we use the convexity of $E$ together with \cref{rmk:gradient_sol_MMS} to calculate:
	\begin{align*}
		\| \partial^2_{jk}(Z_{\tau_i, l} - J_{\delta\tau_i}[Z_{\tau_i, l}] )\|_{L^2} &= \| \partial^2_{jk} (Z_{\tau_i, l}) - J_{\delta\tau_i}[\partial^2_{jk} (Z_{\tau_i, l})] \|_{L^2} \leq  \| \partial^2_{jk} (Z_{\tau_i, l}) - J_{\tau_i^2}[\partial^2_{jk} (Z_{\tau_i,l})] \|_{L^2} \\
		&= \|\partial^2_{jk} (2 U_{\tau_i, l} - U_{\tau_i, l-1} - U_{\tau_i, l+1})\|_{L^2} \leq 2 \tau_i D.
	\end{align*}
	The same calculation holds for the first spatial derivatives, and thus $\tilde{u}_{\tau_i}(t) - u_{\tau_i}(t) \to 0$ in $H^2$ for all $t \in [0, +\infty)$.
	The convergence in \cref{eq:compactness_3_wave} comes from identifying, for $t\in[i\tau,(i+1)\tau)$,
	\begin{equation*}
		\begin{aligned}
			d(E u_\tau)(t) = \frac{1}{2\tau} (\|\nabla U_{\tau, i}\|_{L^2}^2 - \|\nabla U_{\tau, i-1}\|_{L^2}^2) = \langle \nabla \bar{\mathcal{V}}^{\tau, \varepsilon}_t(t), \nabla \dot \bar{\mathcal{V}}^{\tau, \varepsilon}_t \rangle_{L^2} + \frac{\tau}{2} \|\nabla \dot \bar{\mathcal{V}}^{\tau, \varepsilon}_t\|_{L^2}^2,
		\end{aligned}
	\end{equation*}
	where the first equality comes from \cref{eq:energy_derivative_interpolation_wave}, while the second comes from \cref{eq:bar_velocity_interpolation_Dir} together with the polarization formula. By \cref{eq:compactness_1_wave}, $\bar{\mathcal{V}}^{\tau, \varepsilon}_t(t)$ and $\dot \bar{\mathcal{V}}^{\tau, \varepsilon}_t$ converge strongly in $H^1$ to $\bar u(t)$ and $\dot{\bar u}(t)$, respectively, along the chosen subsequence. The first term therefore converges to $\langle\nabla\bar u(t),\nabla\dot{\bar u}(t)\rangle_{L^2}=\frac{\diff}{\diff T}E(V_t^{\bar u}(T))\big|_{T=t}$, while the second goes to zero.
\end{proof}
\begin{proposition}[Limit curves of the wave equation construction are curves of Newtonian motion, cf.~\cref{prop:limit_curves_are_newtonian_curves}]\label{prop:limit-curves-wave-equation-newtonian}
	Let $u_0 \in H^3(\R^d,\R)$ and $\dot{u}_0\in H^2(\R^d,\R)$, both with compact support. For every $\tau >0$, we set $U_{\tau, 0} \coloneqq u_0$ and $U_{\tau,-1} \coloneqq u_0 -\tau \dot{u}_0$, and let $(U_{\tau,i})_{i=-1}^\infty$ be the Newtonian sequence constructed according to \cref{eq:update_wave}.
	Then, the limit curve $\bar{u}$ from \cref{prop:compactness_discrete_solutions_ac2_wave} fulfills, for almost all $s < t$, the following equation:
	\begin{equation} \label{eq:wave_limit_is_newtonian_1}
		DE(\bar{u}_s)[\partial_t \bar{u}_s] - DE(\bar{u}_t)[\partial_t \bar{u}_t] = \int_{s}^{t} \frac{1}{2} |\partial^2_{tt} \bar{u}_r|^2_{L^2} + \frac{1}{2} |\partial E|^2(\bar{u}_r) - \bar{H}_0^E(V_r^{\bar{u}}) \, \diff r,
	\end{equation}
	where $\bar{H}_0^E$ is defined as in \cref{def:hessian_wave}.
	In particular, for all $a, b \in [0, +\infty)$ such that $a < b$, the restriction $\bar{u}|_{(a, b)}$ is a curve of Newtonian motion. Moreover, $\bar{u}$ fulfills, for a.e. $t \in [0, \infty)$, the differential equation
	\begin{align}
		\partial_{tt}^2 \bar{u}_t - \Delta \bar{u}_t &= 0. \label{eq:resulting_PDE_wave}\\
		\bar{u}_0 = u_0, \quad &\dot{\bar{u}}_0 = \dot u_0.
	\end{align}
\end{proposition}
\begin{proof}
	As in the proof of \cref{prop:ugs_wave}, we can rewrite \cref{eq:wave_limit_is_newtonian_1} into
	\begin{equation}
		\langle \nabla \bar{u}_s, \nabla \partial_t \bar{u}_s \rangle_{L^2} - \langle \nabla \bar{u}_t, \nabla \partial_t \bar{u}_t \rangle_{L^2} = \int_{s}^{t} \frac{1}{2} \|\partial^2_{tt} \bar{u}_r\|^2_{L^2} + \frac{1}{2} \| \Delta \bar{u}_r \|^2_{L^2} - \|\nabla \partial_t \bar{u}_r\|^2_{L^2} \, \diff r.
	\end{equation}
From \cref{prop:discrete_newtonian_motion_Dir}, we recall the equality
	\begin{equation}
		d(Eu_\tau)(t_\tau^i) - d(Eu_\tau)(t_\tau^j) =  \frac{1}{2} \int_{t_\tau^i}^{t_\tau^j}   |u_\tau^{(2)}|^2(r)\, \diff r + \frac{1}{2} \int_{t_\tau^i}^{t_\tau^j} |G_\tau|^2(r)\, \diff r - \int_{t_\tau^i}^{t_\tau^j} \bar{H}^E_{0}(\tilde{\mathcal{V}}^{\tau, \varepsilon}_{r}) \, \diff r,
	\end{equation}
	which, using \cref{eq:discrete_newtonian_motion_estimate_G_dir,eq:fancy_velocity_interpolation_Dir,eq:def_H0_Dir}, we can likewise rewrite as
	\begin{equation} \label{eq:discrete_newtonian_motion_Dir_rewritten}
    \begin{split}
		\langle \nabla u_\tau(t_\tau^i), \nabla \dot{u}_\tau(t_\tau^i) \rangle_{L^2} &- \langle \nabla u_\tau(t_\tau^j), \nabla \dot{u}_\tau(t_\tau^j) \rangle_{L^2}  - \frac{\| \nabla U_{\tau, i} - \nabla U_{\tau, i-1}   \|_{L^2}^2}{2\tau} + \frac{\| \nabla U_{\tau, j} - \nabla U_{\tau, j-1}   \|_{L^2}^2}{2\tau} \\   &=\int_{t_\tau^i}^{t_\tau^j} \frac{1}{2}  |u_\tau^{(2)}|^2(r) + \frac{1}{2} \|\Delta \tilde{u}_\tau(r)\|_{L^2}^2 - \| \nabla \dot{u}_{\tau}(r)\|_{L^2}^2 \, \diff r,
    \end{split}
	\end{equation}
	where we define
	\begin{equation*}
		\dot{u}_\tau(t) \coloneqq \frac{U_{\tau, i} - U_{\tau, i-1}}{\tau}, \qquad \text{for } t \in [t_i, t_{i+1}).
	\end{equation*}
	Note that the convergence in \cref{eq:compactness_1_wave} implies that $\dot{u}_\tau(t) \to \dot{\bar{u}}(t)$ in $H^1$; from \cref{prop:compactness_discrete_solutions_ac2_wave}, we furthermore have $u_\tau(t) \to \bar{u}(t)$ in $H^2$. 
    Moreover, \cref{eq:estimate_grad_laplacian} in \cref{cor:nonincreas_total_energy} implies that the last two terms on the left-hand side of \cref{eq:discrete_newtonian_motion_Dir_rewritten} vanish as $\tau\to 0$.
    We have, for any $s',s,t',t\in[0,+\infty)$ such that $s'<s<t<t'$:
	\begin{align}\label{eq:last_equation_wave}
		&\int_s^t \frac12\|\partial_{tt}^2\bar u_r\|_{L^2}^2+\frac12\|\Delta\bar u_r\|_{L^2}^2-\|\nabla\partial_t\bar u_r\|_{L^2}^2\,\diff r \nonumber\\
		&\leq\liminf_{n\to\infty}\int_s^t\frac12|u_{\tau_n}^{(2)}|^2(r)+\frac12\|\Delta\tilde u_{\tau_n}(r)\|_{L^2}^2-\|\nabla\dot u_{\tau_n}(r)\|_{L^2}^2\,\diff r \nonumber\\
		&\leq\liminf_{n\to\infty}\big(d(Eu_{\tau_n})(s')-d(Eu_{\tau_n})(t')\big)+C(t'-t)+C(s-s') \nonumber\\
		&=\langle\nabla\bar u_{s'},\nabla\partial_t\bar u_{s'}\rangle_{L^2}-\langle\nabla\bar u_{t'},\nabla\partial_t\bar u_{t'}\rangle_{L^2}+C(t'-t)+C(s-s'),
	\end{align}
	where $C$ is a uniform bound on $\|\nabla\dot u_{\tau_n}(r)\|_{L^2}^2$ stemming from the energy conservation \cref{eq:nonincreas_total_energy_deriv};
	where we used weak lower semicontinuity of $\|\cdot\|_{L^2(s,t)}^2$ and \cref{eq:compactness_4_wave} for the first term, and dominated convergence for the other two terms in the first step; the discrete identity \cref{eq:discrete_newtonian_motion_Dir} and boundedness of $\|\nabla\dot u_{\tau_n}\|_{L^2}^2$ in the second; and \cref{eq:compactness_3_wave} in the last step. 
    The bounds for dominated convergence follow from \cref{lem:apriori_bounds_wave}: To get a bound on $\|\Delta \tilde{u}_{\tau_n}(r)\|_{L^2}$, note that using convexity of $E$, we have
	\begin{equation*}
		\|\Delta \tilde{u}_{\tau_n}(r)\|_{L^2} = |\partial E|(\tilde{u}_{\tau_n}(r)) = |\partial E|(J_{\delta \cdot \tau}[Z_{\tau, i}]) \leq |\partial E|(Z_{\tau, i}) \leq 2 \|\Delta U_{\tau, i}\|_{L^2} + \|\Delta U_{\tau, i-1}\|_{L^2} \leq 3C.
	\end{equation*}
	Note that $r\mapsto DE(\bar u_r)[\partial_t\bar u_r]=\varphi(r)$ has locally finite variation by \cref{prop:compactness_discrete_solutions_ac2_wave}. Thus, for almost every $s<t$, $\varphi$ is continuous at $s$ and $t$, which allows us to let $s'\uparrow s$ and $t'\downarrow t$ in \cref{eq:last_equation_wave}. Combining the resulting inequality with \cref{prop:ugs_wave}, we see that we must have equality for almost every $s, t \in [0, +\infty)$, and that \cref{eq:resulting_PDE_wave} must hold.
	To see that the initial conditions are correct, it suffices to conclude from \cref{eq:compactness_1_wave} that $V_0(\sigma) = u_0 + \sigma \dot{u}_0$.
\end{proof}

%%%%%%%%%%%%%%%%

\newpage
\printbibliography

\bigskip
\begin{center}
  \FundingLogos
  
  \vspace{0.5em}
  \begin{tcolorbox}\centering\small
    % CHOOSE the sentence that applies to your grant/programme and replace placeholders.
    % 1) For Horizon Europe / ERC (recommended wording):
   
    Funded by the European Union. Views and opinions expressed are however those of the author(s) only and do not necessarily reflect those of the European Union or the European Research Council Executive Agency. Neither the European Union nor the granting authority can be held responsible for them. This project has received funding from the European Research Council (ERC) under the European Union’s Horizon Europe research and innovation programme (grant agreement No. 101198055, project acronym NEITALG).
    
    % 2) (Alternatively, for Horizon 2020 ERC grants use the H2020 wording:)
    % This project has received funding from the European Research Council (ERC) under the European Union's Horizon 2020 research and innovation programme (grant agreement No. <GRANT NUMBER>, project acronym <ACRONYM>).
  \end{tcolorbox}
\end{center}

\begin{table}[ht]
\centering
\caption{Disclosure of AI use in the preparation of this work.}
\label{tab:ai-disclosure}

\begingroup
\renewcommand{\arraystretch}{1.2}
\begin{tabular}{@{}lc@{}}
\toprule
\textbf{Activity} & \textbf{AI used} \\
\midrule
Asking the research question & No \\
Coming up with the approach  & No \\
Proof development            & No \\
Writing and exposition       & No \\
Checking for bugs            & Yes \\
Other supporting tasks       & Yes \\
\bottomrule
\end{tabular}
\endgroup

\captionsetup{width=0.8\textwidth}
\caption*{\textbf{Details.}
We employed Claude Opus~5.5 to find bugs in the arguments and notation inconsistencies, and for the tikz code of the figures.}
\end{table}

\appendix
\pdfbookmark[section]{Appendix}{apx:top}
 % !TEX root = main_metric.tex
\bookmarksetupnext{level=subsection}

\section{Proof of \cref{lem:cbb_interceptor_bound}} \label{sec:proof-of-CBB}
Here, we prove \cref{lem:cbb_interceptor_bound}. Let us repeat the statement:

\begin{lemma}
	Let $X$ be a complete length $CBB(\kappa)$-space for some $\kappa \in \mathbb{R}$. Then $X$ fulfills the \emph{\EpsE-interceptor bound} for some $\epsilon, \mathfrak{E} > 0$ depending only on $\kappa$.
\end{lemma}

\begin{proof}
	First, since every complete length $CBB(\kappa)$-space is also a $CBB(\kappa')$-space for every $\kappa' \leq \kappa$, we can assume without loss of generality that $\kappa < 0$. Furthermore, it suffices to show this property for the $\mathbb{M}_2(\kappa)$ model space: Assume that $\mathbb{M}_2(\kappa)$ fulfills the \EpsE-interceptor bound for some $\epsilon, \mathfrak{E} > 0$, and that two geodesics $\bar{\xi}, \bar{\zeta}$ and $s \leq t$ are given as in \cref{def:interceptor_bound}. Taking $\xi_\sigma$, $\zeta_\sigma$ the images in $\mathbb{M}_2(\kappa)$ of the boundary comparison map that goes from the triangle $[\bar{\xi}_0, \bar{\xi}_t, \bar{\zeta}_t]$ to the model triangle $[\xi_0, \xi_t, \zeta_t]$ in $\mathbb{M}_2(\kappa)$, we have that $d_X(\bar{\xi}_s, \bar{\zeta}_s) \geq d_{\mathbb{M}_2(\kappa)}(\xi_s, \zeta_s)$ ---since the boundary comparison map is a contraction---and $d_X(\bar{\xi}_t, \bar{\zeta}_t) = d_{\mathbb{M}_2(\kappa)}(\xi_t, \zeta_t)$. Thus,
	\begin{equation*}
		d_X(\bar{\xi}_t, \bar{\zeta}_t) = d_{\mathbb{M}_2(\kappa)}(\xi_t, \zeta_t) \leq \frac{d_{\mathbb{M}_2(\kappa)}(\xi_s, \zeta_s)}{s} \cdot ( t + t^3 \cdot \mathfrak{E}) \leq \frac{d_X(\bar{\xi}_s, \bar{\zeta}_s)}{s} \cdot ( t + t^3 \cdot \mathfrak{E}).
	\end{equation*}
	To prove the claim for $\mathbb{M}_2(\kappa)$, we will analyze the law of cosines in this space, given by
	\begin{equation}\label{eq:law_of_cosines}
		\cosh\left( k \cdot c \right) = \cosh\left( k \cdot a \right) \cdot \cosh\left( k \cdot b \right) - \sinh\left( k \cdot a \right) \cdot \sinh\left( k \cdot b \right) \cdot \cos(C),
	\end{equation}
	where $a, b, c$ are the sidelengths of a triangle in $\mathbb{M}_2(\kappa)$ and $C$ is the angle opposite to $c$, and where we write $k \coloneqq \sqrt{-\kappa}$ throughout the proof. We may also assume that $a \neq 0$ and $b \neq 0$, since the bound is otherwise trivial. Let us consider two geodesics $\xi, \zeta\colon [0, T]\to \mathbb{M}_2(\kappa)$ such that $\xi(0) = \zeta(0)$ and $\max\{|\xi'|, |\zeta'|\} \leq 1$. We denote by $\Phi(\xi, \zeta)$ the angle between $\xi$ and $\zeta$ at $\xi(0)=\zeta(0)$. Further, for $a \geq 0$, $b \geq 0$, and $\varphi \in [0, \pi]$, we define $\curlyvee^\kappa(a, b; \varphi)$ to be the length of the side opposite to the angle $\varphi$ in a triangle with angle $\varphi$ between two sides of length $a$ and $b$ in the model space $\mathbb{M}_2(\kappa)$. Since $\xi$ and $\zeta$ are minimising and parametrised with constant speed, we have $d_{\mathbb{M}_2}(\xi_0, \xi_\sigma) = \sigma |\xi'|$ and $d_{\mathbb{M}_2}(\xi_0, \zeta_\sigma) = \sigma |\zeta'|$, while the angle at the vertex $\xi_0$ equals $\Phi(\xi, \zeta)$ for every $\sigma$; we may therefore set
	\begin{equation*}
		d_\sigma \coloneqq d_{\mathbb{M}_2(\kappa)}(\xi_\sigma, \zeta_\sigma) = \curlyvee^\kappa\big(\sigma |\xi'|, \sigma |\zeta'|; \Phi(\xi, \zeta)\big), \qquad d^0_\sigma \coloneqq \curlyvee^0\big(\sigma |\xi'|, \sigma |\zeta'|; \Phi(\xi, \zeta)\big).
	\end{equation*}
	We divide the proof into four steps:
	\begin{enumerate}
		\item Prepare an ancillary claim.
		\item Show a cubic upper bound on $d_\sigma$.
		\item Show a compatible linear lower bound on $d_\sigma$.
		\item Combine both estimates to conclude the proof.
	\end{enumerate}
	In the rest of the proof, we define $\mathrm{sn}_\kappa(x) \coloneqq \frac{\sinh(k x)}{k}$ and $\lambda_\kappa(x) \coloneqq \frac{\mathrm{sn}_\kappa(x)}{x}$, with the convention $\lambda_\kappa(0) \coloneqq 1$. \\
	\textbf{Step 1: Ancillary claim.}
		First, we show that for $a, b \geq 0$ and $\varphi \in [0, \pi]$ and $c \coloneqq \curlyvee^\kappa(a, b; \varphi)$, we have that
		\begin{equation}\label{eq:claim_kappa}
			\sn_\kappa\left(\frac{c}{2}\right)^2 = A \sn_\kappa(P)^2 + B \sn_\kappa(Q)^2,
		\end{equation}
		where $P = \frac{a + b}{2}$, $Q = \frac{|a - b|}{2}$, $A = \sin^2\left(\frac{\varphi}{2}\right)$, and $B = 1 - A = \cos^2 \left(\frac{\varphi}{2}\right)$;
		furthermore, for $c_0 \coloneqq \curlyvee^0(a, b; \varphi)$, we have
		\begin{equation}\label{eq:claim_euclidean}
			\left(\frac{c_0}{2}\right)^2 = AP^2 + BQ^2.
		\end{equation}
		To show \cref{eq:claim_euclidean}, we use the Euclidean law of cosines:
		\begin{align*}
			AP^2 + BQ^2 &= A\frac{(a + b)^2}{4} + B\frac{(a - b)^2}{4} = \frac{(A + B)(a^2 + b^2) + 2ab(A - B)}4 \\
			&= \frac{a^2 + b^2 - 2ab \cos(\varphi)}{4} = \left(\frac{c_0}{2}\right)^2,
		\end{align*}
		where we used that $A - B = -\cos(\varphi)$ by the double-angle formula. To show \cref{eq:claim_kappa}, we rely on the law of cosines in the two half-angle forms (as stated by Lang and Schroeder, see \cite[(18)]{langToponogovsComparisonTheorem2014}), which are equivalent to \cref{eq:law_of_cosines}:
		\begin{align*}
			\sn_\kappa\left(\frac{c}{2}\right)^2 = \underbrace{\sn_\kappa\left(\frac{a + b}{2}\right)^2 - \sn_\kappa(a)\sn_\kappa(b)\cos^2\left(\frac{\varphi}{2}\right)}_{\eqqcolon Z_1} = \underbrace{\sn_\kappa\left(\frac{a-b}{2}\right)^2 + \sn_\kappa(a) \sn_\kappa(b)\sin^2\left(\frac{\varphi}{2}\right)}_{\eqqcolon Z_2}.
		\end{align*}
		We show \cref{eq:claim_kappa} by writing $\sn_\kappa\left(\frac{c}{2}\right)^2$ as a convex combination of $Z_1$ and $Z_2$ and note that the mixed terms cancel:
		\begin{equation*}
			\sn_\kappa\left(\frac{c}{2}\right)^2 = A Z_1 + B Z_2 = A \sn_\kappa\left(\frac{a + b}{2}\right)^2 + B \sn_\kappa\left(\frac{a-b}{2}\right)^2 = A \sn_\kappa(P)^2 + B \sn_\kappa(Q)^2,
		\end{equation*}
		where we used that $\sn_\kappa$ is odd, and hence its square is even, to add the absolute value around $a - b$ in the last step. \\
	\textbf{Step 2: Upper bound.}
	Here, we will show that we can find some $\epsilon, \mathfrak{E}$ independently of $\xi, \zeta$ such that for all $\sigma \leq \epsilon$, we have that
	\begin{equation}\label{eq:upper_bound_distance_goal}
		d_\sigma \leq d^0_\sigma \big(1 + \sigma^2 \cdot \mathfrak{E}\big).
	\end{equation}
	We will use \cref{eq:claim_kappa,eq:claim_euclidean} with $\varphi = \Phi(\xi, \zeta)$ and $a = \sigma |\xi'|$, $b=\sigma |\zeta'|$, so that $c = d_\sigma$ and $c_0 = d^0_\sigma$. We start by writing $\lambda_\kappa$ as its power series $\lambda_\kappa(x) = \sum_{n \geq 0} \frac{(-\kappa)^n x^{2n}}{(2n + 1)!}$ to see that $\lambda_\kappa$ is non-decreasing with $\lambda_\kappa \geq 1$, and thus we have
	\begin{equation}
		\sn_\kappa(P) = \lambda_\kappa(P) P \leq \lambda_\kappa(\sigma) P, \qquad \text{and } \sn_\kappa(Q) = \lambda_\kappa(Q) Q \leq \lambda_\kappa(\sigma) Q;
	\end{equation}
	where we use that $|\xi'| \leq 1$ and $|\zeta'| \leq 1$ imply that $P \leq \sigma$ and $Q \leq \sigma$. By substituting into \cref{eq:claim_kappa} and using \cref{eq:claim_euclidean}, we get
	\begin{equation*}
		\sn_\kappa\left(\frac{d_\sigma}{2}\right)^2 \leq \lambda_\kappa(\sigma)^2 (A P^2 + BQ^2) = \lambda_\kappa(\sigma)^2\left(\frac{d^0_\sigma}{2}\right)^2.
	\end{equation*}
	Next, we take the square root and use that $\sn_\kappa(x) \geq x$ (since $\lambda_\kappa \geq 1$) to get
	\begin{equation*}
		d_\sigma \leq \lambda_\kappa(\sigma) d^0_\sigma.
	\end{equation*}
	It remains to find a bound on $\lambda_\kappa(\sigma)$ of the form $1 + \sigma^2 \cdot \mathfrak{E}$. To this end, we define the temporary function $g(x) \coloneqq \frac{\sinh(x) - x}{x^3} = \sum_{n \geq 1} \frac{x^{2n - 2}}{(2n+1)!}$.
	Since all coefficients are non-negative, $g$ is non-decreasing, and thus we have that $g(x) \leq g(1)$ for $x \in [0, 1]$. Therefore, setting $\epsilon \coloneqq \frac{1}{k}$, we have that for $\sigma \in (0, \epsilon]$,
	\begin{equation*}
		\frac{\lambda_\kappa(\sigma) - 1}{\sigma^2} = \frac{\sinh(k \sigma) - k \sigma}{k \sigma^3} = k^2 g(k \sigma) \leq -\kappa \, g(1);
	\end{equation*}
	which yields $\lambda_\kappa(\sigma) \leq 1 - \kappa \, g(1) \sigma^2$ and thus proves \cref{eq:upper_bound_distance_goal} with $\epsilon \coloneqq \frac{1}{k}$ and $\mathfrak{E} \coloneqq -\kappa \, g(1)$.
	\textbf{Step 3: Lower bound.}
	Showing the compatible lower bound
	\begin{equation*}
		d_\sigma \geq d^0_\sigma
	\end{equation*}
	is somewhat more direct:
	It suffices to recall (\cite[1.1, (d)]{alexanderAlexandrovGeometryFoundations2023}) that the function $\kappa \mapsto \curlyvee^\kappa(a, b; \varphi)$ is non-increasing for every fixed $a, b, \varphi$. Thus, we have that
	\begin{equation}\label{eq:lower_bound_distance_model_space}
		d_\sigma = \curlyvee^\kappa(\sigma |\xi'|, \sigma |\zeta'|; \Phi(\xi, \zeta)) \geq \curlyvee^0(\sigma |\xi'|, \sigma |\zeta'|; \Phi(\xi, \zeta)) = d^0_\sigma.
	\end{equation}
	\textbf{Step 4: Conclusion.}
	Combining \cref{eq:upper_bound_distance_goal,eq:lower_bound_distance_model_space} for $s$ and $t$, where $0 < s \leq t \leq \epsilon$, and using that $\curlyvee^0$ is $1$-homogeneous in its first two arguments, we obtain
	\begin{equation*}
		d_t \leq d^0_t \cdot (1 + t^2 \cdot \mathfrak{E}) = \frac{t}{s} \cdot d^0_s \cdot (1 + t^2 \cdot \mathfrak{E}) \leq \frac{d_s}{s} \cdot ( t + t^3 \cdot \mathfrak{E}),
	\end{equation*}
	which concludes the proof.
\end{proof}

\section{Proof of \cref{lem:apriori_3_wave}} \label{sec:proof-of-apriori-3-wave}
Here, we prove \cref{lem:apriori_3_wave}. Let us repeat the statement:
\begin{lemma}
	For every $\tau \in (0, 1)$, let us consider $u_0 \in H^2(\R^d,\R)$ and $\dot{u}_0\in H^1(\R^d,\R)$, let us set $U_{\tau, 0} \coloneqq u_0$ and $U_{\tau,-1} \coloneqq u_0 -\tau \dot{u}_0$, and let $(U_{\tau,i})_{i=-1}^{N_\tau}$ be the Newtonian sequence constructed according to \cref{eq:update_wave}.
	Let $A \subset \R^d$ be such that $\supp(u_0) \cup \supp(\dot{u}_0) \subseteq A$.
	Then, for each $T > 0$, there exist $L^0$ and $v > 0$, independent of $\tau$, such that, denoting for a set $A \subseteq \R^d$ its enlargement by $r$ as $A_r \coloneqq \{x \in \R^d \mid \dist(x, A) \leq r \}$, we have that
	\begin{equation*}
		\int_{A_{v (t + \tau k)}^c} | U_{\tau, \lfloor t / \tau \rfloor }|^2(x) \, \diff x  \leq L^0 2^{-\lfloor t / \tau \rfloor  - k}, \qquad \forall t \in [0, T],\ k \in \NN.
	\end{equation*}
	In particular, if $A$ is compact, then, for each $T > 0$, there exists $M_T > 0$, independent of $\tau$, such that
	\begin{equation*}
		\int_{\R^d} |U_{\tau, i}|^2(x) (1 + |x|^2) \, \diff x \leq M_T, \qquad \int_{\R^d} |U'_{\tau, i}|^2(x) (1 + |x|^2) \, \diff x \leq M_T
	\end{equation*}
	for $\tau$ small enough and $ i \leq N_\tau$, where we set $U'_{\tau, i} \coloneqq \frac{U_{\tau, i} - U_{\tau, i-1}}{\tau}$.
\end{lemma}
\begin{proof}
	We recall from \cref{eq:update_wave,eq:Euler_Lagrange_update} that the solutions fulfill the iterative PDE
	\begin{equation*}
		- \tau^2 \Delta U_{\tau, i+1} + U_{\tau, i+1} = 2 U_{\tau, i} - U_{\tau, i-1} = Z_{\tau, i}.
	\end{equation*}
	Setting $G_\tau \coloneqq x\mapsto \tau^{-d} G(\frac{x}{\tau})$, where $G$ is the Bessel kernel of order $2$, we thus have
	\begin{equation*}
		U_{\tau, i+1} = G_\tau \star Z_{\tau, i}.
	\end{equation*}
	Our strategy to prove \cref{eq:sol_wave_Equation_decay} is to use the exponentially vanishing tail of $G_\tau$ together with an induction over $i$. \\
	First, let us note that there exist $v_0>0$, $C \geq 1$ and $c > 0$ such that
	\begin{equation*}
		\int_{|x| \geq \tau v} G_\tau(x) \, \diff x\leq C e^{-c |v|}
	\end{equation*}
	whenever $v \geq v_0$. To this end, let us recall the asymptotic formula for the Bessel kernel as $|x| \to \infty$,
	\begin{equation*}
		G(x) = \frac{e^{-|x|}}{2^\frac{d + 1}{2}\pi^\frac{d-1}{2} |x|^\frac{d-1}{2}}(1 + o(1)) \leq C' e^{-c' |x|}, \quad \text{ for some } C', c' \text{ and } |x| \text{ sufficiently large,}
	\end{equation*}
	and hence
	\begin{equation*}
		\int_{|x| \geq v} G(x) \, \diff x \leq C e^{-c v}\, \quad \text{ for some } C, c \text{ and } v \text{ big enough.}
	\end{equation*}
	Finally, we get
	\begin{equation*}
		\int_{|x| \geq \tau v} G_\tau(x) \, \diff x = \int_{|x| \geq \tau v} \tau^{-d} G\left(\frac{x}{\tau}\right)\, \diff x = \int_{|x| \geq v} G(x) \, \diff x \leq C e^{-c v}.
	\end{equation*}
	Next, we choose $v' \geq v_0$ such that $\sum_{j=0}^{\infty} 2^j C e^{-cv'j} \leq 2C$, and $C e^{-cv'j} \leq 2^{-j}$ for $j \geq 1$. Fixing $T > 0$ and $\tau > 0$, we set $L^{i, j} \coloneqq \int_{A_{\tau v' j}^c} U_{\tau, i}^2$. The next step is to show that, for $j \geq i$,
	\begin{equation}\label{eq:sol_wave_equation_decay_2}
		L^{i, j} \leq L^0 2^{-(j - i)} (10 + 20C)^i,
	\end{equation}
	where $L^0 \coloneqq \|u_0\|_{L^2}^2 + 4\|\dot{u}_0\|_{L^2}^2 + D^2$, where $D$ is chosen such that $\|U_{\tau, i}\|_{L^2} \leq D$ and $\|\frac{U_{\tau, i}- U_{\tau, i-1}}{\tau}\|_{L^2} \leq D$ for all $i \leq \lfloor T / \tau \rfloor$ by \cref{eq:apriori_1_wave_d,eq:apriori_1_wave_a}. To show \cref{eq:sol_wave_equation_decay_2}, we proceed by induction; since $G_\tau \geq 0$ and $\int_{\R^d} G_\tau = 1$, Jensen's inequality yields
	\begin{align}
		L^{i, j} &= \int_{A_{\tau v' j}^c} U_{\tau, i}^2 \leq \int_{A_{\tau v' j}^c} \int_{\R^d} G_{\tau}(x-z) Z_{\tau, i-1}^2(z) \, \diff z \, \diff x \\
		&= \int_{\R^d} Z_{\tau, i-1}^2(z) \underbrace{\int_{A_{\tau v' j}^c} G_\tau(x - z) \, \diff x}_{\eqqcolon \mathfrak{G}_{\tau, j}(z)} \, \diff z \\
		&\leq \int_{A_{\tau v'(i-1)}} Z_{\tau, i-1}^2(z) \mathfrak{G}_{\tau, j}(z) \, \diff z + \sum_{k = i}^{j} \int_{B_k} Z_{\tau, i-1}^2(z) \mathfrak{G}_{\tau, j}(z) \, \diff z + \int_{A_{\tau v'j}^c} Z_{\tau, i-1}^2(z) \mathfrak{G}_{\tau, j}(z) \, \diff z \nonumber \\
		&\leq \int_{A_{\tau v'(i-1)}} Z_{\tau, i-1}^2(z) \, \diff z \, C e^{-c v' (j - (i-1))}  + \sum_{k = i}^{j} \int_{B_k} Z_{\tau, i-1}^2(z) C e^{-cv'(j-k)} \, \diff z  \nonumber\\
		&\quad +\int_{A_{\tau v'j}^c} Z_{\tau, i-1}^2(z) \, \diff z, \label{eq:sol_wave_equation_decay_3}
	\end{align}
	where $B_k = A_{\tau v'(k-1)}^c \cap A_{\tau v'k}$. For the induction start, the cases $i=0$ and $i=1, j \leq 1$ are trivial. For $i=1, j \geq 2$, only the first summand in~\eqref{eq:sol_wave_equation_decay_3} survives, and, using the initialization $Z_{\tau,0} = u_0 + \tau \dot{u}_0$, we have, for $\tau \leq 1$,
	\begin{equation*}
		Z_{\tau,0}^2 \leq 2u_0^2 + 2\tau^2 (\dot{u}_0)^2 \leq 2u_0^2 + 2 (\dot{u}_0)^2,
	\end{equation*}
	so that
	\begin{equation*}
		L^{1, j} \leq 2(\|u_0\|_{L^2}^2 + \|\dot{u}_0\|_{L^2}^2) C e^{-cv'j} \leq L^0 2^{-(j-1)}
	\end{equation*}
	by the choice of $v'$. For the induction step, we continue from \cref{eq:sol_wave_equation_decay_3}, estimating, by the Cauchy--Schwarz inequality applied to $Z_{\tau, i-1} = 2 U_{\tau,i-1} - U_{\tau,i-2}$,
	\begin{equation*}
		Z^2_{\tau, i-1} \leq (2^2+1^2)\big(U_{\tau, i-1}^2 + U_{\tau, i-2}^2\big) = 5 U_{\tau, i-1}^2 + 5 U_{\tau, i-2}^2:
	\end{equation*}
	\begin{align}
		L^{i, j} &\leq 5 \Big( 2 D^2 C e^{-cv'(j-(i-1))} + \sum_{k=i}^{j} (L^{i-1, k-1} + L^{i-2, k-1}) C e^{-cv'(j-k)} + L^{i-1, j} + L^{i-2, j}\Big) \\
		&\leq 5 L^0 (10 + 20C)^{i-1}\Big(2^{-(j-i)} + \sum_{k=i}^j \big(2^{-(k-i)} + 2^{-(k-i+1)}\big) C e^{-cv'(j-k)} + 2^{-(j-i + 1)} + 2^{-(j-i+2)}\Big) \nonumber \\
		&= 5L^0 (10 + 20C)^{i-1} \Big(2^{-(j-i)} + \sum_{k=0}^{j-i} \big(2^{-(j-k-i)} + 2^{-(j-k-i+1)}\big) C e^{-cv'k} + 2^{-(j-i+1)} + 2^{-(j-i+2)}\Big) \nonumber\\
		&\leq 5 L^0 2^{-(j-i)} (10 + 20C)^{i-1} \Big(1 + 2 \sum_{k=0}^{j-i} 2^k C e^{-cv'k} + 1\Big) \leq L^0 2^{-(j-i)} (10 + 20C)^i,
	\end{align}
	we thus have shown \cref{eq:sol_wave_equation_decay_2}. \\
	Choosing $m\in\mathbb N$, $m\geq1$, such that $2^{m-1} \geq 2 (10 + 20 C)$, we set $v \coloneqq m v'$. We then have
	\begin{align*}
		\int_{A_{v(t + \tau k)}^c} U_{\tau, \lfloor t / \tau \rfloor }^2
		&\leq \int_{A_{v \tau (\lfloor t / \tau \rfloor + k)}^c} U_{\tau, \lfloor t / \tau \rfloor }^2 \\
		&= \int_{A_{m v' \tau (\lfloor t / \tau \rfloor + k)}^c} U_{\tau, \lfloor t / \tau \rfloor }^2
		= L^{\lfloor t / \tau \rfloor,\, m (\lfloor t / \tau \rfloor + k)} \\
		&\leq L^0 2^{- ((m-1)\lfloor t / \tau \rfloor + mk)} (10+20C)^{\lfloor t / \tau \rfloor} \\
		&\leq L^0 \cdot 2^{-k} \Big(2^{-(m-1)} (10 + 20 C)\Big)^{\lfloor t / \tau \rfloor} \\
		&\leq L^0 2^{- \lfloor t / \tau \rfloor - k},
	\end{align*}
	which finishes the proof of \cref{eq:sol_wave_Equation_decay}. Assume in addition that $0 < \tau \leq 1$. To show \cref{eq:sol_wave_equation_in_weighted_l2}, we fix $R_0$ such that $A \subseteq B_{R_0}$ and calculate
	\begin{align}
		\int_{\R^d} U_{\tau, i}^2(x) (1 + |x|^2) \diff x
		&= \|U_{\tau, i}\|_{L^2}^2 + 2 \int_0^\infty R \int_{|x| \geq R} U^2_{\tau, i}(x) \, \diff x \, \diff R \\
		&\leq D^2 (1 + 2 (R_0 + vT)^2) + 2 \int_{R_0+ vT}^\infty R \int_{A^c_{R-R_0}} U_{\tau, i}^2(x) \, \diff x \, \diff R. \label{eq:weighted_l2_u_1}
	\end{align}
	For $R \geq R_0 + vT$, define
	\[
		k_R \coloneqq \left\lfloor \frac{R-R_0}{v\tau} \right\rfloor - i.
	\]
	Since $\tau i \leq T$, we have $k_R \in \NN$, and
	\[
		v\tau (i + k_R) \leq R-R_0.
	\]
	Hence
	\[
		A_{R-R_0}^c \subseteq A_{v\tau(i+k_R)}^c,
	\]
	and by \cref{eq:sol_wave_Equation_decay},
	\[
		\int_{A^c_{R-R_0}} U_{\tau, i}^2(x) \, \diff x
		\leq \int_{A_{v\tau(i+k_R)}^c} U_{\tau, i}^2(x) \, \diff x
		\leq L^0 2^{-i-k_R}
		= L^0 2^{-\left\lfloor \frac{R-R_0}{v\tau} \right\rfloor}.
	\]
	Continuing from~\eqref{eq:weighted_l2_u_1}, we obtain
	\begin{align}
		\int_{\R^d} U_{\tau, i}^2(x) (1 + |x|^2) \diff x
		&\leq D^2 (1 + 2 (R_0 + vT)^2) + 2L^0 \int_{R_0+ vT}^\infty R\, 2^{-\left\lfloor \frac{R-R_0}{v\tau} \right\rfloor} \, \diff R \\
		&\leq D^2 (1 + 2 (R_0 + vT)^2) + 2L^0 v^2\tau^2 \int_{T/\tau}^\infty s\, 2^{-\lfloor s \rfloor} \, \diff s \\
		&\quad + 2L^0 R_0 v\tau \int_{T/\tau}^\infty 2^{-\lfloor s \rfloor} \, \diff s \leq M_T,
	\end{align}
	for some $M_T$ independent of $\tau$ and $i \leq \lfloor T/\tau \rfloor$, since both integrals in $s$ are bounded uniformly in $\tau \in (0,1]$. For the bound on $U'_{\tau,i}$, the case $i=0$ follows from the initialization $U'_{\tau,0}=\dot{u}_0$ and the compact support of $\dot{u}_0$. For $i \geq 1$, using
	\[
		(U'_{\tau,i})^2 \leq \frac{2}{\tau^2}\big(U_{\tau,i}^2 + U_{\tau,i-1}^2\big),
	\]
	we calculate
	\begin{align}
		&\int_{\R^d} (U'_{\tau, i})^2(x) (1 + |x|^2) \diff x
		= \|U'_{\tau, i}\|_{L^2}^2 + 2 \int_0^\infty R \int_{|x| \geq R} (U'_{\tau, i})^2(x) \, \diff x \, \diff R \\
		&\qquad \leq D^2 (1 + 2 (R_0 + vT)^2) + \frac{4}{\tau^2} \int_{R_0+ vT}^\infty R \int_{A^c_{R-R_0}} \big(U_{\tau,i}^2(x) + U_{\tau,i-1}^2(x)\big) \, \diff x \, \diff R. \label{eq:weighted_l2_up_1}
	\end{align}
	With the same choice
	\[
		k_R = \left\lfloor \frac{R-R_0}{v\tau} \right\rfloor - i \in \NN,
	\]
	we have
	\[
		\int_{A^c_{R-R_0}} U_{\tau, i}^2(x) \, \diff x \leq L^0 2^{-\left\lfloor \frac{R-R_0}{v\tau} \right\rfloor},
	\]
	and, since $i\ge 1$,
	\[
		\int_{A^c_{R-R_0}} U_{\tau, i-1}^2(x) \, \diff x
		\leq L^0 2^{-(i-1)-k_R}
		= 2L^0 2^{-\left\lfloor \frac{R-R_0}{v\tau} \right\rfloor}.
	\]
	Continuing from~\eqref{eq:weighted_l2_up_1}, we obtain
	\begin{align}
		\int_{\R^d} (U'_{\tau, i})^2(x) (1 + |x|^2) \diff x
		&\leq D^2 (1 + 2 (R_0 + vT)^2) + \frac{12L^0}{\tau^2} \int_{R_0+ vT}^\infty R\, 2^{-\left\lfloor \frac{R-R_0}{v\tau} \right\rfloor} \, \diff R \\
		&\leq D^2 (1 + 2 (R_0 + vT)^2) + 12L^0 v^2 \int_{T/\tau}^\infty s\, 2^{-\lfloor s \rfloor} \, \diff s \\
		&\quad + \frac{12L^0 R_0 v}{\tau} \int_{T/\tau}^\infty 2^{-\lfloor s \rfloor} \, \diff s.
	\end{align}
	Finally, since
	\[
		\int_{a}^{\infty} 2^{-\lfloor s \rfloor}\,\diff s \leq 4\cdot 2^{-a}, \qquad a \geq 1,
	\]
	we obtain
	\[
		\frac{1}{\tau}\int_{T/\tau}^{\infty} 2^{-\lfloor s \rfloor}\,\diff s \leq \frac{4}{\tau}2^{-T/\tau},
	\]
	and the right-hand side is bounded on $(0,1]$, whenever $\tau \leq T$. Hence
	\[
		\int_{\R^d} (U'_{\tau, i})^2(x) (1 + |x|^2) \diff x \leq M_T
	\]
	for some $M_T$ independent of $\tau$ and $i \leq \lfloor T/\tau \rfloor$, for all $\tau \leq T$.
\end{proof}

\end{document}